\documentclass[10pt]{amsart}
\usepackage[english,russian]{babel}
\usepackage[cp1251]{inputenc}
\usepackage{amsmath}
\usepackage{amssymb}
\usepackage{amsfonts}

\usepackage[linktocpage=true, colorlinks=true, linkcolor=blue, citecolor=blue, urlcolor=blue]{hyperref}

\begin{document}
\renewcommand{\refname}{References}
\renewcommand\contentsname{Contents}

\thispagestyle{empty}

\title[Expansion of Iterated Ito Stochastic Integrals]
{Expansion of Iterated Ito Stochastic Integrals of Arbitrary Multiplicity
Based on 
Generalized Multiple Fourier Series Converging in the Mean}
\author[D.F. Kuznetsov]{Dmitriy F. Kuznetsov}
\address{Dmitriy Feliksovich Kuznetsov
\newline\hphantom{iii} Peter the Great Saint-Petersburg Polytechnic University,
\newline\hphantom{iii} Polytechnicheskaya ul., 29,
\newline\hphantom{iii} 195251, Saint-Petersburg, Russia}%
\email{sde\_kuznetsov@inbox.ru}

\maketitle {\small
\begin{quote}
\vspace{5mm}
\noindent{\sc Abstract.} 
The article is devoted to the expansions of iterated
Ito stochastic integrals based on generalized multiple Fourier
series converging in the sense of norm in the space
$L_2([t, T]^k),$ $k\in\mathbb{N}.$ The
method of generalized multiple Fourier
series for ex\-pansion and mean-square approximation
of iterated Ito stochastic integrals of arbitrary
multiplicity $k$ ($k\in\mathbb{N}$) with respect to components
of the multidimensional Wiener process is proposed and developed.
The obtained expansions contain only one operation
of the limit transition in contrast to its existing analogues.
In the article it is also obtained the generalization 
of the proposed method for an arbitrary complete orthonormal
systems of functions in the space 
$L_2([t, T]^k),$ $k\in\mathbb{N}$ as well as for  
complete orthonormal with weight $r(t_1)\ldots r(t_k)$ 
systems of functions 
in the space 
$L_2([t, T]^k),$ $k\in\mathbb{N}$.
The comparison of the considered method with the well-known
ex\-pan\-si\-ons of iterated Ito stochastic integrals
based on the Ito formula and Hermite polynomials is given. 
The convergence in the mean of degree $2n$ $(n \in \mathbb{N})$ and with
probability 1 
of the proposed method is proved. 

\medskip
\end{quote}
}

\vspace{9mm}

\linespread{1.1}

\tableofcontents

\linespread{1.0}

\section{Introduction}

\vspace{5mm} 
  
The idea of representing of iterated Ito and Stratonovich stochastic 
integrals 
in the form of multiple stochastic integrals from 
specific discontinuous nonrandom functions of several variables and following 
expansion of these functions using multiple Fourier series in order 
to get effective mean-square approximations of the mentioned stochastic 
integrals was proposed and developed in a lot of 
publications of the author \cite{1}-\cite{Kuzh-1}
(also see related publications \cite{31b}, \cite{31bb}). Note that another
approaches to expansions of iterated stochastic integrals can be found in 
\cite{32}-\cite{35}. Specifically, 
the approach \cite{1}-\cite{Kuzh-1} appeared for the first time in 
\cite{1} (1994), \cite{2} (1996). In these
works the mentioned idea is formulated more likely at the level of 
guess (without any satisfactory grounding), and as a result the 
papers \cite{1}, \cite{2} contain rather fuzzy formulations and a number of 
incorrect conclusions. 
Note that in \cite{1}, \cite{2} we 
used the trigonometric multiple Fourier series 
converging  
in the sence of norm in the space
$L_2([t, T]^k),$ $k=1, 2, 3.$ In the final form the approach from
\cite{1}, \cite{2} has been formulated and proved for the first time
in the monograph \cite{7} (2006) (see Theorem 1 below).
It should be noted that the results 
of \cite{1}, \cite{2} are correct for a sufficiently narrow particular case 
when the numbers $i_1,...,i_k$ are pairwise different,
$i_1,\ldots,i_k = 1,\ldots,m$ (see Sect.~2 for detail).
   
Usage of Fourier series with respect to the system of Legendre polynomials 
for approximation of iterated stochastic integrals took place for 
the first time in \cite{3} (1997) \cite{4} (1998), \cite{5} (2000)
\cite{6} (2001) (also see 
\cite{7}-\cite{31bb}).

The question about what integrals (Ito or Stratonovich) are more 
suitable for expansions within the frames of the considered direction 
of researches has turned out to be rather interesting and difficult.

On the one side, Theorem 1 (see Sect. 2) conclusively 
demonstrates that the structure 
of iterated Ito stochastic integrals is rather convenient for expansions into 
multiple series with respect to the system of standard Gaussian random 
variables regardless of the 
multiplicity $k$ of iterated Ito stochastic integrals.
   
On the other side, the results of \cite{3}-\cite{6}, \cite{11}-\cite{2023xxx1},
\cite{21}-\cite{25}, \cite{27}, \cite{27aa}-\cite{28},
\cite{29a}-\cite{29aaaa}, 
\cite{31qaa}-\cite{310aaa} con\-vin\-cing\-ly demontrate
that the final formulas for 
expansions of iterated Stratonovich stochastic integrals of multiplicities 1 to 8
(the case of continuously differentiable weight functions and 
a complete orthonormal system of Legendre polynomials or 
trigonometric functions in $L_2([t, T])$) and
iterated Stratonovich stochastic integrals
of multiplicity $k,$ $k\in {\bf N}$ (the case of continuous weight functions 
and an arbitrary complete orthonormal system of functions in $L_2([t, T])$)
are more compact than their analogues for iterated
Ito stochastic integrals.

\vspace{5mm}

\section{Theorem on Expansion of Iterated Ito Stochastic Integrals of 
Arbitrary Multiplicity $k$ $(k\in\mathbb{N})$}

\vspace{5mm}

Let $(\Omega,$ ${\rm F},$ ${\sf P})$ be a complete probability space, let 
$\{{\rm F}_t, t\in[0,T]\}$ be a nondecreasing right-continous family of 
$\sigma$-algebras of ${\rm F},$
and let ${\bf f}_t$ be a standard $m$-dimensional Wiener 
stochastic process, which is
${\rm F}_t$-measurable for any $t\in[0, T].$ We assume that the components
${\bf f}_{t}^{(i)}$ $(i=1,\ldots,m)$ of this process are independent.

Let us consider 
the following iterated Ito 
stochastic integrals

\begin{equation}
\label{sodom20}
J[\psi^{(k)}]_{T,t}=\int\limits_t^T\psi_k(t_k) \ldots \int\limits_t^{t_{2}}
\psi_1(t_1) d{\bf w}_{t_1}^{(i_1)}\ldots
d{\bf w}_{t_k}^{(i_k)},
\end{equation}

\vspace{3mm}
\noindent
where every $\psi_l(\tau)$ $(l=1,\ldots,k)$ is 
a nonrandom function on $[t, T]$,
${\bf w}_{\tau}^{(i)}={\bf f}_{\tau}^{(i)}$
for $i=1,\ldots,m$ and
${\bf w}_{\tau}^{(0)}=\tau,$ 
$i_1,\ldots,i_k=0, 1,\ldots,m.$

Suppose that every $\psi_l(\tau)$ $(l=1,\ldots,k)$ is a continuous 
nonrandom function on the interval $[t, T]$ (the case $\psi_1(\tau),\ldots,\psi_k(\tau)
\in L_2([t, T])$ will be considered in Sect.~15).

Define the following function on the hypercube $[t, T]^k$

\begin{equation}
\label{ppp}
K(t_1,\ldots,t_k)=
\begin{cases}
\psi_1(t_1)\ldots \psi_k(t_k),\ &t_1<\ldots<t_k\\
~\\
~\\
0,\ &\hbox{\rm otherwise}
\end{cases}\ \ \ \ 
=\ \ \ \ 
\prod\limits_{l=1}^k
\psi_l(t_l)\ \prod\limits_{l=1}^{k-1}{\bf 1}_{\{t_l<t_{l+1}\}},\ 
\end{equation}

\vspace{4mm}
\noindent
where $t_1,\ldots,t_k\in [t, T]$ $(k\ge 2)$ and 
$K(t_1)\equiv\psi_1(t_1)$ for $t_1\in[t, T].$ Here 
${\bf 1}_A$ denotes the indicator of the set $A$.

Suppose that $\{\phi_j(x)\}_{j=0}^{\infty}$
is a complete orthonormal system of functions in 
the space $L_2([t, T])$.

The function $K(t_1,\ldots,t_k)$ is piecewise continuous in the 
hypercube $[t, T]^k.$
At this situation it is well known that the generalized 
multiple Fourier series 
of $K(t_1,\ldots,t_k)\in L_2([t, T]^k)$ is converging 
to $K(t_1,\ldots,t_k)$ in the hypercube $[t, T]^k$ in 
the mean-square sense, i.e.

\begin{equation}
\label{sos1z}
\hbox{\vtop{\offinterlineskip\halign{
\hfil#\hfil\cr
{\rm lim}\cr
$\stackrel{}{{}_{p_1,\ldots,p_k\to \infty}}$\cr
}} }\Biggl\Vert
K(t_1,\ldots,t_k)-
\sum_{j_1=0}^{p_1}\ldots \sum_{j_k=0}^{p_k}
C_{j_k\ldots j_1}\prod_{l=1}^{k} \phi_{j_l}(t_l)\Biggr\Vert_{L_2([t, T]^k)}=0,
\end{equation}

\vspace{3mm}
\noindent
where
\begin{equation}
\label{ppppa}
C_{j_k\ldots j_1}=\int\limits_{[t,T]^k}
K(t_1,\ldots,t_k)\prod_{l=1}^{k}\phi_{j_l}(t_l)dt_1\ldots dt_k
\end{equation}

\vspace{4mm}
\noindent
is the Fourier coefficient and

\vspace{-1mm}
$$
\left\Vert f\right\Vert_{L_2([t, T]^k)}=\left(\int\limits_{[t,T]^k}
f^2(t_1,\ldots,t_k)dt_1\ldots dt_k\right)^{1/2}.
$$

\vspace{4mm}

Consider the partition $\{\tau_j\}_{j=0}^N$ of the interval $[t,T]$ such that

\begin{equation}
\label{1111}
t=\tau_0<\ldots <\tau_N=T,\ \ \
\Delta_N=
\hbox{\vtop{\offinterlineskip\halign{
\hfil#\hfil\cr
{\rm max}\cr
$\stackrel{}{{}_{0\le j\le N-1}}$\cr
}} }\Delta\tau_j\to 0\ \ \hbox{if}\ \ N\to \infty,\ \ \ 
\Delta\tau_j=\tau_{j+1}-\tau_j.
\end{equation}

\vspace{3mm}

{\bf Theorem 1} \cite{7} (2006) \cite{8}-\cite{Kuzh-1}. 
{\it Suppose that
every $\psi_l(\tau)$ $(l=1,\ldots, k)$ is a continuous 
non\-ran\-dom func\-ti\-on on
the interval $[t, T]$ and
$\{\phi_j(x)\}_{j=0}^{\infty}$ is a complete orthonormal system  
of conti\-nu\-ous func\-ti\-ons in the space $L_2([t,T]).$ 
Then

$$
J[\psi^{(k)}]_{T,t}\  =\ 
\hbox{\vtop{\offinterlineskip\halign{
\hfil#\hfil\cr
{\rm l.i.m.}\cr
$\stackrel{}{{}_{p_1,\ldots,p_k\to \infty}}$\cr
}} }\sum_{j_1=0}^{p_1}\ldots\sum_{j_k=0}^{p_k}
C_{j_k\ldots j_1}\Biggl(
\prod_{l=1}^k\zeta_{j_l}^{(i_l)}\ -
\Biggr.
$$

\vspace{3mm}
\begin{equation}
\label{tyyy}
-\ \Biggl.
\hbox{\vtop{\offinterlineskip\halign{
\hfil#\hfil\cr
{\rm l.i.m.}\cr
$\stackrel{}{{}_{N\to \infty}}$\cr
}} }\sum_{(l_1,\ldots,l_k)\in {\rm G}_k}
\phi_{j_{1}}(\tau_{l_1})
\Delta{\bf w}_{\tau_{l_1}}^{(i_1)}\ldots
\phi_{j_{k}}(\tau_{l_k})
\Delta{\bf w}_{\tau_{l_k}}^{(i_k)}\Biggr),
\end{equation}

\vspace{6mm}
\noindent
where

$$
{\rm G}_k={\rm H}_k\backslash{\rm L}_k,\ \ \
{\rm H}_k=\{(l_1,\ldots,l_k):\ l_1,\ldots,l_k=0,\ 1,\ldots,N-1\},
$$

\vspace{-2mm}
$$
{\rm L}_k=\{(l_1,\ldots,l_k):\ l_1,\ldots,l_k=0,\ 1,\ldots,N-1;\
l_g\ne l_r\ (g\ne r);\ g, r=1,\ldots,k\},
$$

\vspace{5mm}
\noindent
${\rm l.i.m.}$ is a limit in the mean-square sense,
$i_1,\ldots,i_k=0,1,\ldots,m,$ 

\vspace{-1mm}
\begin{equation}
\label{rr23}
\zeta_{j}^{(i)}=
\int\limits_t^T \phi_{j}(s) d{\bf w}_s^{(i)}
\end{equation} 

\vspace{2mm}
\noindent
are independent standard Gaussian random variables
for various
$i$ or $j$ {\rm(}if $i\ne 0${\rm),}
$C_{j_k\ldots j_1}$ is the Fourier coefficient {\rm(\ref{ppppa}),}
$\Delta{\bf w}_{\tau_{j}}^{(i)}=
{\bf w}_{\tau_{j+1}}^{(i)}-{\bf w}_{\tau_{j}}^{(i)}$
$(i=0,\ 1,\ldots,m),$
$\left\{\tau_{j}\right\}_{j=0}^{N}$ is a partition of
the interval $[t,T]$ satisfying the condition {\rm (\ref{1111})}.}

\vspace{2mm}

{\bf Proof.}
At first, let us prove preparatory lemmas.

\vspace{2mm}

{\bf Lemma 1.} {\it Suppose that
every $\psi_l(\tau)$ $(l=1,\ldots, k)$ is a continuous nonrandom
function on 
the interval $[t, T]$. Then

\vspace{-2mm}
\begin{equation}
\label{30.30}
J[\psi^{(k)}]_{T,t}=
\hbox{\vtop{\offinterlineskip\halign{
\hfil#\hfil\cr
{\rm l.i.m.}\cr
$\stackrel{}{{}_{N\to \infty}}$\cr
}} }
\sum_{j_k=0}^{N-1}
\ldots \sum_{j_1=0}^{j_{2}-1}
\prod_{l=1}^k \psi_l(\tau_{j_l})\Delta{\bf w}_{\tau_
{j_l}}^{(i_l)}\ \ \ \hbox{\rm w.\ p.\ 1},
\end{equation}

\vspace{4mm}
\noindent
where $\Delta{\bf w}_{\tau_{j}}^{(i)}=
{\bf w}_{\tau_{j+1}}^{(i)}-{\bf w}_{\tau_{j}}^{(i)}$
$(i=0, 1,\ldots,m)$,
$\left\{\tau_{j}\right\}_{j=0}^{N}$ is a partition 
of the interval $[t,T]$ satisfy\-ing the condition {\rm (\ref{1111})},
hereinafter w. p. {\rm 1}  means  "with probability {\rm 1}".
}

\vspace{2mm}

{\bf Proof.}\ It is easy to notice that using the 
additive property of stochastic integrals, we can write

\begin{equation}
\label{toto}
J[\psi^{(k)}]_{T,t}=
\sum_{j_k=0}^{N-1}\ldots
\sum_{j_{1}=0}^{j_{2}-1}\prod_{l=1}^{k}
J[\psi_l]_{\tau_{j_l+1},\tau_{j_l}}+
\varepsilon_N\ \ \ \ \ \hbox{w.\ p.\ 1},
\end{equation}

\vspace{1mm}
\noindent
where

\vspace{-1mm}
$$
\varepsilon_N = 
\sum_{j_k=0}^{N-1}\int\limits_{\tau_{j_k}}^{\tau_{j_k+1}}
\psi_k(s)\int\limits_{\tau_{j_k}}^{s}
\psi_{k-1}(\tau)J[\psi^{(k-2)}]_{\tau,t}d{\bf w}_\tau^{(i_{k-1})}
d{\bf w}_s^{(i_k)} +
$$

$$
+ \sum_{r=1}^{k-3}
G[\psi_{k-r+1}^{(k)}]_{N}
\sum_{j_{k-r}=0}^{j_{k-r+1}-1}\int\limits_{\tau_{j_{k-r}}}^{\tau_{j_{k-r}+1}}
\psi_{k-r}(s)\int\limits_{\tau_{j_{k-r}}}^{s}
\psi_{k-r-1}(\tau)J[\psi^{(k-r-2)}]_{\tau,t}
d{\bf w}_\tau^{(i_{k-r-1})}
d{\bf w}_s^{(i_{k-r})}+
$$

\vspace{1mm}
$$
+ G[\psi_3^{(k)}]_{N}
\sum_{j_{2}=0}^{j_{3}-1}J[\psi^{(2)}]_{\tau_{j_{2}+1},
\tau_{j_{2}}},
$$

\vspace{5mm}
$$
G[\psi_m^{(k)}]_{N}=\sum_{j_k=0}^{N-1}\sum_{j_{k-1}=0}^{j_k-1}
\ldots \sum_{j_{m}=0}^{j_{m+1}-1}
\prod_{l=m}^{k}J[\psi_l]_{\tau_{j_l+1},\tau_{j_l}},
$$

\vspace{2mm}
$$
J[\psi_l]_{s,\theta}=\int\limits_{\theta}^s\psi_l(\tau)d{\bf w}_{\tau}
^{(i_l)},
$$

\vspace{2mm}
$$
(\psi_m,\psi_{m+1},\ldots,\psi_{k})\stackrel{\rm def}{=}\psi_m^{(k)},\ \ \
(\psi_1,\ldots,\psi_{k})\stackrel{\rm def}{=}\psi_1^{(k)}=\psi^{(k)}.
$$

\vspace{6mm}

Using the standard estimates (\ref{99.010}) for moments of stochastic 
integrals, we obtain w. p. 1

\begin{equation}
\label{999.0001}
\hbox{\vtop{\offinterlineskip\halign{
\hfil#\hfil\cr
{\rm l.i.m.}\cr
$\stackrel{}{{}_{N\to \infty}}$\cr
}} }\varepsilon_N =0.
\end{equation}

\vspace{3mm}

Comparing (\ref{toto}) and (\ref{999.0001}), we get

\begin{equation}
\label{toto1}
J[\psi^{(k)}]_{T,t}=
\hbox{\vtop{\offinterlineskip\halign{
\hfil#\hfil\cr
{\rm l.i.m.}\cr
$\stackrel{}{{}_{N\to \infty}}
$\cr
}} }
\sum_{j_k=0}^{N-1}\ldots
\sum_{j_{1}=0}^{j_{2}-1}\prod_{l=1}^{k}
J[\psi_l]_{\tau_{j_l+1},\tau_{j_l}}\ \ \ \hbox{w.\ p.\ 1}.
\end{equation}

\vspace{4mm}

Let us rewrite $J[\psi_l]_{\tau_{j_l+1},\tau_{j_l}}$ in the form 

$$
J[\psi_l]_{\tau_{j_l+1},\tau_{j_l}}=
\psi_l(\tau_{j_l})\Delta{\bf w}_{\tau_{j_l}}^{(i_l)}+
\int\limits_{\tau_{j_l}}^{\tau_{j_l+1}}
(\psi_l(\tau)-\psi_l(\tau_{j_l}))d{\bf w}_{\tau}^{(i_l)}
$$

\vspace{2mm}
\noindent
and substitute it into
(\ref{toto1}).
Then, due to the moment properties of stochastic integrals and
continuity (which means uniform continuity) of the functions 
$\psi_l(s)$ ($l=1,\ldots,k$)
it is easy to see that the 
prelimit
expression on the right-hand side of (\ref{toto1}) is a sum of 
the prelimit
expression on the right-hand side of (\ref{30.30}) and the value which 
tends to zero in the mean-square sense if 
$N\to\infty.$ Lemma 1 is proved. 

\vspace{2mm}

{\bf Remark 1.} {\it It is easy to see that if
$\Delta{\bf w}_{\tau_{j_l}}^{(i_l)}$ in {\rm (\ref{30.30})}
for some $l\in\{1,\ldots,k\}$ is replaced with 
$\left(\Delta{\bf w}_{\tau_{j_l}}^{(i_l)}\right)^p$ $(p=2,$
$i_l\ne 0),$ then
the differential $d{\bf w}_{t_{l}}^{(i_l)}$
in the integral $J[\psi^{(k)}]_{T,t}$
will be replaced with $dt_l$.
If $p=3, 4,\ldots,$ then the
right-hand side 
of the formula {\rm (\ref{30.30})}
will become zero w. p. {\rm 1}.
If we replace $\Delta{\bf w}_{\tau_{j_l}}^{(i_l)}$ in {\rm (\ref{30.30})}
for some $l\in\{1,\ldots,k\}$
with $\left(\Delta \tau_{j_l}\right)^p$ $(p=2, 3,\ldots),$
then the right-hand side of the formula
{\rm (\ref{30.30})} also 
will be equal to zero w. p. {\rm 1}.}

\vspace{2mm}

Let us define the following
multiple stochastic integral

\begin{equation}
\label{30.34}
\hbox{\vtop{\offinterlineskip\halign{
\hfil#\hfil\cr
{\rm l.i.m.}\cr
$\stackrel{}{{}_{N\to \infty}}$\cr
}} }\sum_{j_1,\ldots,j_k=0}^{N-1}
\Phi\left(\tau_{j_1},\ldots,\tau_{j_k}\right)
\prod\limits_{l=1}^k\Delta{\bf w}_{\tau_{j_l}}^{(i_l)}
\stackrel{\rm def}{=}J[\Phi]_{T,t}^{(k)},
\end{equation}

\vspace{4mm}
\noindent
where $\Phi(t_1,\ldots,t_k): [t, T]^k\to{\mathbb R},$
$\Phi(t_1,\ldots,t_k)\in C([t, T]^k)$, i.e. $\Phi(t_1,\ldots,t_k)$
is a continuous non\-ran\-dom function in $[t, T]^k.$

Denote
\begin{equation}
\label{dom1}
D_k=\{(t_1,\ldots,t_k):\ t\le t_1<\ldots <t_k\le T\}.
\end{equation}

\vspace{3mm}

We will use the same symbol $D_k$ to denote the open and closed 
domains corresponding to the domain $D_k$ defined by (\ref{dom1}).
However, we always specify what domain we consider (open or closed).

Also we will write $\Phi(t_1,\ldots,t_k)\in C(D_k)$
if 
$\Phi(t_1,\ldots,t_k)$ is a continuous nonrandom function of $k$ variables
in the closed domain $D_k$.

Let us consider the iterated Ito stochastic integral

\begin{equation}
\label{rrr29}
I[\Phi]_{T,t}^{(k)}\stackrel{\rm def}{=}
\int\limits_t^T\ldots \int\limits_t^{t_2}
\Phi(t_1,\ldots,t_k)d{\bf w}_{t_1}^{(i_1)}\ldots
d{\bf w}_{t_k}^{(i_k)},
\end{equation}

\vspace{3mm}
\noindent
where $\Phi(t_1,\ldots,t_k)\in C(D_k).$

Using the arguments which similar to the arguments used for 
the proof of Lemma 1
it is easy to demonstrate that if
$\Phi(t_1,\ldots,t_k)\in C(D_k),$ then the following equality is fulfilled

\begin{equation}
\label{30.52}
I[\Phi]_{T,t}^{(k)}=\hbox{\vtop{\offinterlineskip\halign{
\hfil#\hfil\cr
{\rm l.i.m.}\cr
$\stackrel{}{{}_{N\to \infty}}$\cr
}} }
\sum_{j_k=0}^{N-1}
\ldots \sum_{j_1=0}^{j_{2}-1}
\Phi(\tau_{j_1},\ldots,\tau_{j_k})
\prod\limits_{l=1}^k\Delta {\bf w}_{\tau_{j_l}}^{(i_l)}\ \ \ \hbox{w.\ p.\ 1}.
\end{equation}

\vspace{3mm}

In order to explain, let us check the correctness of the equality 
(\ref{30.52}) when $k=3$.
For definiteness we will suggest that 
$i_1,i_2,i_3=1,\ldots,m.$ We have

\vspace{1mm}
$$
I[\Phi]_{T,t}^{(3)}\stackrel{\rm def}{=}
\int\limits_t^T\int\limits_t^{t_3}\int\limits_t^{t_2}
\Phi(t_1,t_2,t_3)d{\bf w}_{t_1}^{(i_1)}d{\bf w}_{t_2}^{(i_2)}
d{\bf w}_{t_3}^{(i_3)}=
$$

$$
=\hbox{\vtop{\offinterlineskip\halign{
\hfil#\hfil\cr
{\rm l.i.m.}\cr
$\stackrel{}{{}_{N\to \infty}}$\cr
}} }
\sum_{j_3=0}^{N-1}
\int\limits_{t}^{\tau_{j_3}}\int\limits_t^{t_2}
\Phi(t_1,t_2,\tau_{j_3})d{\bf w}_{t_1}^{(i_1)}d{\bf w}_{t_2}^{(i_2)}
\Delta{\bf w}_{\tau_{j_3}}^{(i_3)}=
$$

\vspace{1mm}
$$
=\hbox{\vtop{\offinterlineskip\halign{
\hfil#\hfil\cr
{\rm l.i.m.}\cr
$\stackrel{}{{}_{N\to \infty}}$\cr
}} }
\sum_{j_3=0}^{N-1}\sum_{j_2=0}^{j_3-1}
\int\limits_{\tau_{j_2}}^{\tau_{j_2+1}}\int\limits_t^{t_2}
\Phi(t_1,t_2,\tau_{j_3})d{\bf w}_{t_1}^{(i_1)}d{\bf w}_{t_2}^{(i_2)}
\Delta{\bf w}_{\tau_{j_3}}^{(i_3)}=
$$

\vspace{1mm}
$$
=\hbox{\vtop{\offinterlineskip\halign{
\hfil#\hfil\cr
{\rm l.i.m.}\cr
$\stackrel{}{{}_{N\to \infty}}$\cr
}} }
\sum_{j_3=0}^{N-1}\sum_{j_2=0}^{j_3-1}
\int\limits_{\tau_{j_2}}^{\tau_{j_2+1}}
\left(\ \int\limits_t^{\tau_{j_2}}\
+\ \int\limits_{\tau_{j_2}}^{t_2}\ \right)
\Phi(t_1,t_2,\tau_{j_3})d{\bf w}_{t_1}^{(i_1)}d{\bf w}_{t_2}^{(i_2)}
\Delta{\bf w}_{\tau_{j_3}}^{(i_3)}=
$$

\vspace{1mm}
$$
=\hbox{\vtop{\offinterlineskip\halign{
\hfil#\hfil\cr
{\rm l.i.m.}\cr
$\stackrel{}{{}_{N\to \infty}}$\cr
}} }
\sum_{j_3=0}^{N-1}\sum_{j_2=0}^{j_3-1}\sum_{j_1=0}^{j_2-1}
\int\limits_{\tau_{j_2}}^{\tau_{j_2+1}}\int\limits_{\tau_{j_1}}^{\tau_{j_1+1}}
\Phi(t_1,t_2,\tau_{j_3})d{\bf w}_{t_1}^{(i_1)}d{\bf w}_{t_2}^{(i_2)}
\Delta{\bf w}_{\tau_{j_3}}^{(i_3)}+
$$

\vspace{1mm}
\begin{equation}
\label{44444.25}
+\hbox{\vtop{\offinterlineskip\halign{
\hfil#\hfil\cr
{\rm l.i.m.}\cr
$\stackrel{}{{}_{N\to \infty}}$\cr
}} }
\sum_{j_3=0}^{N-1}\sum_{j_2=0}^{j_3-1}
\int\limits_{\tau_{j_2}}^{\tau_{j_2+1}}\int\limits_{\tau_{j_2}}^{t_2}
\Phi(t_1,t_2,\tau_{j_3})d{\bf w}_{t_1}^{(i_1)}d{\bf w}_{t_2}^{(i_2)}
\Delta{\bf w}_{\tau_{j_3}}^{(i_3)}.
\end{equation}

\vspace{3mm}

Let us demonstrate that the second limit on the right-hand side 
of (\ref{44444.25}) 
equals to zero.

Actually, the second moment of its 
prelimit
expression equals to 

$$
\sum_{j_3=0}^{N-1}\sum_{j_2=0}^{j_3-1}
\int\limits_{\tau_{j_2}}^{\tau_{j_2+1}}\int\limits_{\tau_{j_2}}^{t_2}
\Phi^2(t_1,t_2,\tau_{j_3})dt_1 dt_2
\Delta\tau_{j_3}
\le M^2  \sum_{j_3=0}^{N-1}\sum_{j_2=0}^{j_3-1}
\frac{1}{2}\left(\Delta\tau_{j_2}\right)^2\Delta\tau_{j_3}\to 0
$$

\vspace{3mm}
\noindent
if $N\to\infty.$
Here $M$ is a constant, which restricts the module of
function
$\Phi(t_1,t_2,t_3)$ due to its continuity, $\Delta\tau_j=
\tau_{j+1}-\tau_j.$

Considering the obtained conclusions, we have

\vspace{1mm}
$$
I[\Phi]_{T,t}^{(3)}\stackrel{\rm def}{=}
\int\limits_t^T\int\limits_t^{t_3}\int\limits_t^{t_2}
\Phi(t_1,t_2,t_3)d{\bf w}_{t_1}^{(i_1)}d{\bf w}_{t_2}^{(i_2)}
d{\bf w}_{t_3}^{(i_3)}=
$$

\vspace{1mm}
$$
=\hbox{\vtop{\offinterlineskip\halign{
\hfil#\hfil\cr
{\rm l.i.m.}\cr
$\stackrel{}{{}_{N\to \infty}}$\cr
}} }
\sum_{j_3=0}^{N-1}\sum_{j_2=0}^{j_3-1}\sum_{j_1=0}^{j_2-1}
\int\limits_{\tau_{j_2}}^{\tau_{j_2+1}}\int\limits_{\tau_{j_1}}^{\tau_{j_1+1}}
\Phi(t_1,t_2,\tau_{j_3})d{\bf w}_{t_1}^{(i_1)}d{\bf w}_{t_2}^{(i_2)}
\Delta{\bf w}_{\tau_{j_3}}^{(i_3)}=
$$

\vspace{1mm}
$$
=\hbox{\vtop{\offinterlineskip\halign{
\hfil#\hfil\cr
{\rm l.i.m.}\cr
$\stackrel{}{{}_{N\to \infty}}$\cr
}} }
\sum_{j_3=0}^{N-1}\sum_{j_2=0}^{j_3-1}\sum_{j_1=0}^{j_2-1}
\int\limits_{\tau_{j_2}}^{\tau_{j_2+1}}\int\limits_{\tau_{j_1}}^{\tau_{j_1+1}}
\left(\Phi(t_1,t_2,\tau_{j_3})-\Phi(t_1,\tau_{j_2},\tau_{j_3})\right)
d{\bf w}_{t_1}^{(i_1)}d{\bf w}_{t_2}^{(i_2)}
\Delta{\bf w}_{\tau_{j_3}}^{(i_3)}+
$$

\vspace{1mm}
$$
+\hbox{\vtop{\offinterlineskip\halign{
\hfil#\hfil\cr
{\rm l.i.m.}\cr
$\stackrel{}{{}_{N\to \infty}}$\cr
}} }
\sum_{j_3=0}^{N-1}\sum_{j_2=0}^{j_3-1}\sum_{j_1=0}^{j_2-1}
\int\limits_{\tau_{j_2}}^{\tau_{j_2+1}}\int\limits_{\tau_{j_1}}^{\tau_{j_1+1}}
\left(\Phi(t_1,\tau_{j_2},\tau_{j_3})-
\Phi(\tau_{j_1},\tau_{j_2},\tau_{j_3})\right)
d{\bf w}_{t_1}^{(i_1)}d{\bf w}_{t_2}^{(i_2)}
\Delta{\bf w}_{\tau_{j_3}}^{(i_3)}+
$$

\vspace{2mm}
\begin{equation}
\label{4444.1}
+\hbox{\vtop{\offinterlineskip\halign{
\hfil#\hfil\cr
{\rm l.i.m.}\cr
$\stackrel{}{{}_{N\to \infty}}$\cr
}} }
\sum_{j_3=0}^{N-1}\sum_{j_2=0}^{j_3-1}\sum_{j_1=0}^{j_2-1}
\Phi(\tau_{j_1},\tau_{j_2},\tau_{j_3})
\Delta{\bf w}_{\tau_{j_1}}^{(i_1)}
\Delta{\bf w}_{\tau_{j_2}}^{(i_2)}
\Delta{\bf w}_{\tau_{j_3}}^{(i_3)}.
\end{equation}

\vspace{5mm}

In order to get the sought result, we just have to demonstrate that 
the first
two limits on the right-hand side of (\ref{4444.1}) equal to zero. 
Let us prove 
that the first one of them equals to zero (proof for the second limit 
is similar).
   
The second moment of a prelimit expression of the first limit on the 
right-hand side of (\ref{4444.1}) equals to the following expression

\begin{equation}
\label{4444.01}
\sum_{j_3=0}^{N-1}\sum_{j_2=0}^{j_3-1}\sum_{j_1=0}^{j_2-1}
\int\limits_{\tau_{j_2}}^{\tau_{j_2+1}}\int\limits_{\tau_{j_1}}^{\tau_{j_1+1}}
\left(\Phi(t_1,t_2,\tau_{j_3})-\Phi(t_1,\tau_{j_2},\tau_{j_3})\right)^2
dt_1 dt_2
\Delta\tau_{j_3}.
\end{equation}

\vspace{3mm}

Since the function $\Phi(t_1,t_2,t_3)$ is continuous in 
the closed bo\-un\-ded domain
$D_3,$ then
it is uniformly continuous in this domain. Therefore, if the 
distance between two points of the domain $D_3$ is less than 
$\delta(\varepsilon)$ ($\delta(\varepsilon)>0$ 
exists
for any $\varepsilon>0$ and it does not depend 
on mentioned points), then the cor\-res\-pond\-ing oscillation of the function 
$\Phi(t_1,t_2,t_3)$ for these two points of the domain $D_3$ is less than
$\varepsilon.$

If we assume that $\Delta\tau_j<\delta(\varepsilon)$ ($j=0, 1,\ldots,N-1$),
then the distance between points 
$(t_1,t_2,\tau_{j_3})$,\ $(t_1,\tau_{j_2},\tau_{j_3})$
is obviously less than $\delta(\varepsilon).$ In this case 

\vspace{-1mm}
$$
|\Phi(t_1,t_2,\tau_{j_3})-\Phi(t_1,\tau_{j_2},\tau_{j_3})|<\varepsilon.
$$

\vspace{3mm}

Consequently, when $\Delta\tau_j<\delta(\varepsilon)$ ($j=0,\ 1,\ldots,N-1$)
the expression (\ref{4444.01})  
is evaluated by the following value

\vspace{-1mm}
$$
\varepsilon^2
\sum_{j_3=0}^{N-1}\sum_{j_2=0}^{j_3-1}\sum_{j_1=0}^{j_2-1}
\Delta\tau_{j_1}\Delta\tau_{j_2}\Delta\tau_{j_3}<
\varepsilon^2\frac{(T-t)^3}{6}.
$$ 

\vspace{4mm}

Therefore, the first limit on the right-hand side 
of (\ref{4444.1}) equals to zero.
Similarly we can prove equality to zero of the second limit on the right-hand
side
of (\ref{4444.1}).

Consequently, the equality (\ref{30.52}) is proved for $k=3$. 
The cases $k=2$ and $k>3$ 
are analyzed absolutely similarly.

It is necessary to note that the proof of 
correctness of (\ref{30.52}) 
is similar when the nonrandom function $\Phi(t_1,\ldots,t_k)$ is 
continuous in 
the open domain $D_k$ and bounded at its boundary.

Let us consider the following multiple 
stochastic integral

\vspace{-1mm}
\begin{equation}
\label{mult11}
\hbox{\vtop{\offinterlineskip\halign{
\hfil#\hfil\cr
{\rm l.i.m.}\cr
$\stackrel{}{{}_{N\to \infty}}$\cr
}} }
\sum\limits_{\stackrel{j_1,\ldots,j_k=0}{{}_{j_q\ne j_r;\ q\ne r;\ q, r=1,\ldots, k}}}^{N-1}
\Phi\left(\tau_{j_1},\ldots,\tau_{j_k}\right)
\prod\limits_{l=1}^k
\Delta{\bf w}_{\tau_{j_l}}^{(i_l)}
\stackrel{\rm def}{=}J'[\Phi]_{T,t}^{(k)},
\end{equation}

\vspace{3mm}
\noindent
where $\Phi(t_1,\ldots,t_k): [t, T]^k\to{\mathbb R}$ is the same function as in
(\ref{30.34}).

Then, according to (\ref{30.52}) we will get the following

\begin{equation}
\label{pobeda}
J'[\Phi]_{T,t}^{(k)}=
\int\limits_t^T\ldots \int\limits_t^{t_2}
\sum\limits_{(t_1,\ldots,t_k)}\biggl(
\Phi(t_1,\ldots,t_k)
d{\bf w}_{t_1}^{(i_1)}\ldots
d{\bf w}_{t_k}^{(i_k)}\biggr),
\end{equation}

\vspace{3mm}
\noindent
where
$$
\sum\limits_{(t_1,\ldots,t_k)}
$$ 

\vspace{2mm}
\noindent
means the sum with respect to all
possible permutations
$(t_1,\ldots,t_k).$ 
At the same time permutations $(t_1,\ldots,t_k)$ 
when summing 
are performed in (\ref{pobeda}) only in the expression, which
is enclosed in pa\-ren\-the\-ses. Moreover,
the nonrandom function $\Phi(t_1,\ldots,t_k)$ is assumed 
to be continuous in the 
cor\-res\-pond\-ing closed domains of integration.
The case when the nonrandom function $\Phi(t_1,\ldots,t_k)$ is 
continuous in the open domains of integration and bounded at 
their boundaries is also possible.

It is not difficult to see that (\ref{pobeda})
can be rewritten in the form

\begin{equation}
\label{pobedaxyz}
J'[\Phi]_{T,t}^{(k)}=\sum_{(t_1,\ldots,t_k)}
\int\limits_{t}^{T}
\ldots
\int\limits_{t}^{t_2}
\Phi(t_1,\ldots,t_k)d{\bf w}_{t_1}^{(i_1)}
\ldots
d{\bf w}_{t_k}^{(i_k)},
\end{equation}

\vspace{3mm}
\noindent
where permutations $(t_1,\ldots,t_k)$ when summing are 
performed only in the values
$d{\bf w}_{t_1}^{(i_1)}
\ldots $
$d{\bf w}_{t_k}^{(i_k)}.$ At the same time the indexes near 
upper 
limits of integration in the iterated stochastic integrals are changed 
correspondently and if $t_r$ swapped with $t_q$ in the  
permutation $(t_1,\ldots,t_k)$, then $i_r$ swapped with $i_q$ in 
the permutation $(i_1,\ldots,i_k)$.

Let us consider the class ${\rm M}_2([0, T])$ of functions
$\xi: [0,T]\times\Omega\rightarrow \mathbb{R},$ which are 
measurable
in accordance with the collection of variables
$(t,\omega)$ and
${\rm F}_t$-measurable 
for all $t\in[0,T].$ Moreover, $\xi(\tau,\omega)$ is independent 
with increments ${\bf f}_{t+\Delta}-{\bf f}_{t}$
for $t\ge \tau $ $(\Delta>0),$ 

\vspace{-1mm}
$$
\int\limits_0^T{\sf M}\left\{\xi^2(t,\omega)\right\}dt
<\infty,
$$

\vspace{2mm}
\noindent
and ${\sf M}\{\xi^2(t,\omega)\}<\infty$
for all $t\in[0,T].$

It is well known \cite{36} that the Ito stochastic integral
exists in the mean-square sence for any 
$\xi\in{\rm M}_2([0, T]).$  Further, we will denote  
$\xi(\tau,\omega)$ as $\xi_{\tau}.$

\vspace{2mm}

{\bf Lemma 2.}\ {\it Suppose that $\Phi(t_1,\ldots,t_k)\in C(D_k)$ 
or $\Phi(t_1,\ldots,t_k)$ 
is a continuous nonrandom function in the open domain $D_k$ and bounded at its boundary.
Then

\vspace{1mm}
$$
{\sf M}\biggl\{\biggl|I[\Phi]_{T,t}^{(k)}\biggr|^{2}\biggr\}
\le C_{k}
\int\limits_t^T\ldots \int\limits_t^{t_2}
\Phi^{2}(t_1,\ldots,t_k)dt_1\ldots dt_k,\ \ \
C_{k}<\infty,
$$

\vspace{5mm}
\noindent
where $I[\Phi]_{T,t}^{(k)}$ is defined by the formula {\rm (\ref{rrr29})}.}

\vspace{2mm}

{\bf Proof.}\
Using standard properties and moments estimates of stochastic integrals,
we have for $\xi_{\tau}\in{\rm M}_2([t,T])$ \cite{36}

\vspace{-1mm}
\begin{equation}
\label{99.010}
{\sf M}\left\{\left|\int\limits_{t}^T \xi_\tau
df_\tau\right|^{2}\right\} =
\int\limits_{t}^T {\sf M}\{|\xi_\tau|^{2}\}d\tau,\ \ \ \
{\sf M}\left\{\left|\int\limits_{t}^T \xi_\tau
d\tau\right|^{2}\right\} \le (T-t)
\int\limits_{t}^T {\sf M}\{|\xi_\tau|^{2}\}d\tau.
\end{equation}

\vspace{5mm}

Let us denote

\vspace{-1mm}
$$
\xi[\Phi]_{t_{l+1},\ldots,t_k,t}^{(l)}=
\int\limits_t^{t_{l+1}}\ldots \int\limits_t^{t_2}
\Phi(t_1,\ldots,t_k)
d{\bf w}_{t_1}^{(i_1)}\ldots
d{\bf w}_{t_{l}}^{(i_{l})},
$$ 

\vspace{5mm}
\noindent
where $l=1,\ldots,$ $k-1$ and

\vspace{-1mm}
$$
\xi[\Phi]_{t_{1},\ldots,t_k,t}^{(0)}\stackrel{\rm def}{=}
\Phi(t_1,\ldots,t_k).
$$

\vspace{5mm}

In accordance with the induction it is easy to demonstrate that

\vspace{1mm}
$$
\xi[\Phi]_{t_{l+1},\ldots,t_k,t}^{(l)}\in{\rm M}_2([t,T])
$$

\vspace{5mm}
\noindent
with respect to the variable $t_{l+1}.$
Further, using the estimates (\ref{99.010})
repeatedly we obtain the statement of Lemma 2.

It is not difficult to see that in the case $i_1,\ldots,i_k=1,\dots,m$
from the proof of Lemma 2 we obtain

\begin{equation}
\label{dobav1}
{\sf M}\biggl\{\biggl|I[\Phi]_{T,t}^{(k)}\biggr|^{2}\biggr\}
=
\int\limits_t^T\ldots \int\limits_t^{t_2}
\Phi^{2}(t_1,\ldots,t_k)dt_1\ldots dt_k.
\end{equation}

\vspace{4mm}

{\bf Lemma 3.}\ {\it Suppose that every $\varphi_l(s)$
$(l=1,\ldots,k)$ is a continuous nonrandom function on the interval $[t, T]$.
Then
\begin{equation}
\label{30.39}
\prod_{l=1}^k 
J[\varphi_l]_{T,t}=J[\Phi]_{T,t}^{(k)}\ \ \ \hbox{\rm w. p. 1},
\end{equation}

\vspace{3mm}
\noindent
where 
$$
J[\varphi_l]_{T,t}
=\int\limits_t^T \varphi_l(s) d{\bf w}_{s}^{(i_l)},\ \ \ 
\Phi(t_1,\ldots,t_k)=\prod\limits_{l=1}^k\varphi_l(t_l)
$$

\vspace{5mm}
\noindent
and the integral $J[\Phi]_{T,t}^{(k)}$ 
is defined
by the equality
{\rm (\ref{30.34})}.
}

\vspace{3mm}

{\bf Proof.}
Let at first $i_l\ne 0,$ $l=1,\ldots,k.$
Denote

\vspace{2mm}
$$
J[\varphi_l]_{N}\stackrel{\rm def}{=}\sum\limits_{j=0}^{N-1}
\varphi_l(\tau_j)\Delta{\bf w}_{\tau_j}^{(i_l)}.
$$

\vspace{4mm}

Since
$$
\prod_{l=1}^k J[\varphi_l]_{N}-\prod_{l=1}^k J[\varphi_l]_{T,t}
=
$$

\vspace{1mm}
\begin{equation}
\label{df2}
=\sum_{l=1}^k \left(\prod_{g=1}^{l-1} J[\varphi_g]_{T,t}\right)
\biggl(J[\varphi_l]_{N}-J[\varphi_l]_{T,t}
\biggr)\left(\prod_{g=l+1}^k J[\varphi_g]_{N}\right),
\end{equation}

\vspace{4mm}
\noindent
then because of the Minkowski inequality and the inequality 
of Cauchy-Bunyakovsky we obtain

\begin{equation}
\label{30.42}
\left({\sf M}\left\{\left|\prod_{l=1}^k J[\varphi_l]_{N}
-\prod_{l=1}^k J[\varphi_l]_{T,t}\right|^2
\right\}\right)^{1/2}\le C_k
\sum_{l=1}^k
\left({\sf M}\left\{
\biggl|J[\varphi_l]_{N}-J[\varphi_l]_{T,t}\biggr|^4\right\}\right)
^{1/4},
\end{equation}

\vspace{5mm}
\noindent 
where $C_k$ is a constant.

Note that

\vspace{-1mm}
$$
J[\varphi_l]_{N}-J[\varphi_l]_{T,t}=\sum\limits_{j=0}^{N-1}
J[\Delta\varphi_l]_{\tau_{j+1},\tau_j},\ \ \ 
J[\Delta\varphi_l]_{\tau_{j+1},\tau_j}
=\int\limits_{\tau_j}^{\tau_{j+1}}\left(
\varphi_l(\tau_j)-\varphi_l(s)\right)d{\bf w}_{s}^{(i_l)}.
$$

\vspace{4mm}

Since $J[\Delta\varphi_l]_{\tau_{j+1},\tau_j}$
are independent for various $j,$ then \cite{37}

\vspace{2mm}
$$
{\sf M}\left\{\left|\sum_{j=0}^{N-1}J[\Delta\varphi_l]_{\tau_{j+1},
\tau_j}\right|^4
\right\}=
\sum_{j=0}^{N-1}{\sf M}\left\{\biggl|J[\Delta\varphi_l]_{\tau_{j+1},
\tau_j}\biggr|^4
\right\}+ 
$$

\vspace{2mm}
\begin{equation}
\label{30.43}
+6 \sum_{j=0}^{N-1}{\sf M}
\left\{\biggl|J[\Delta\varphi_l]_{\tau_{j+1},\tau_j}\biggr|^2
\right\}
\sum_{q=0}^{j-1}{\sf M}\left\{\biggl|
J[\Delta\varphi_l]_{\tau_{q+1},\tau_q}\biggr|^2
\right\}.
\end{equation}

\vspace{6mm}

It is obviously that
$J[\Delta\varphi_l]_{\tau_{j+1},\tau_j}$ are Gaussian random variables.
Then
we have

\vspace{1mm}
$$
{\sf M}\left\{\biggl|J[\Delta\varphi_l]_{\tau_{j+1},\tau_j}\biggr|^2\right\}=
\int\limits_{\tau_j}^{\tau_{j+1}}(\varphi_l(\tau_j)-\varphi_l(s))^2ds,
$$

$$
{\sf M}\left\{\biggl|J[\Delta\varphi_l]_{\tau_{j+1},\tau_j}\biggr|^4\right\}=
3\left(\int\limits_{\tau_j}^{\tau_{j+1}}(\varphi_l(\tau_j)-\varphi_l(s))^2ds
\right)^2.
$$

\vspace{4mm}

Using this relations and continuity (which means uniform continuity) 
of the functions $\varphi_l(s),$ we get

$$
{\sf M}\left\{\left|\sum_{j=0}^{N-1}J[\Delta\varphi_l]_{\tau_{j+1},
\tau_j}\right|^4
\right\}\le \varepsilon^4\left(
3 \sum_{j=0}^{N-1}(\Delta\tau_{j})^2+
6 \sum_{j=0}^{N-1}\Delta\tau_{j}
\sum_{q=0}^{j-1}\Delta\tau_{q}\right)
<
$$

\vspace{2mm}
$$
<3\varepsilon^4\left(\delta(\varepsilon) (T-t)+(T-t)^2\right),
$$

\vspace{5mm}
\noindent
where $\Delta\tau_{j}<\delta(\varepsilon),$
$j=0,1,\ldots,N-1$  ($\delta(\varepsilon)>0$ exists
for any $\varepsilon>0$ and it does not
depend on points of the interval $[t, T]$).
Then the right-hand side of the formula 
(\ref{30.43}) tends to zero when $N\to \infty.$ 

Considering this fact as well 
as (\ref{30.42}), we come to (\ref{30.39}). 

If ${\bf w}_{t_l}^{(i_l)}=t_l$ for some $l\in\{1,\ldots,k\},$
then
the proof of Lemma 3 becomes obviously simpler and  
it is performed similarly. Lemma 3 is proved.

\vspace{2mm}

{\bf Remark 2.} {\it It is easy to see that if
$\Delta{\bf w}_{\tau_{j_l}}^{(i_l)}$ in {\rm (\ref{30.39})}
for some $l\in\{1,\ldots,k\}$ is replaced with 
$\left(\Delta{\bf w}_{\tau_{j_l}}^{(i_l)}\right)^p$ $(p=2,$
$i_l\ne 0),$ then
the differential $d{\bf w}_{t_{l}}^{(i_l)}$
in the integral $J[\Phi^{(k)}]_{T,t}$
will be replaced with $dt_l$.
If $p=3, 4,\ldots,$ then the
right-hand side 
of the formula {\rm (\ref{30.39})}
will become zero w. p. {\rm 1}.}

\vspace{2mm}

Let us consider the case $p=2$ in detail.
Let $\Delta{\bf w}_{\tau_{j_l}}^{(i_l)}$ in {\rm (\ref{30.39})}
for some $l\in\{1,\ldots,k\}$ is replaced with 
$\left(\Delta{\bf w}_{\tau_{j_l}}^{(i_l)}\right)^2$ $(i_l\ne 0)$ 
and

\vspace{1mm}
$$
J[\varphi_l]_{N}\stackrel{\rm def}{=}\sum\limits_{j=0}^{N-1}
\varphi_l(\tau_j)\left(\Delta{\bf w}_{\tau_j}^{(i_l)}\right)^2,\ \ \ 
J[\varphi_l]_{T,t}\stackrel{\rm def}{=}
\int\limits_t^T\varphi_l(s)ds.
$$

\vspace{3mm}

We have

\vspace{1mm}
$$
\left({\sf M}\left\{
\biggl|J[\varphi_l]_{N}-J[\varphi_l]_{T,t}\biggr|^4\right\}\right)
^{1/4}=\left({\sf M}\left\{
\left|\sum\limits_{j=0}^{N-1}
\varphi_l(\tau_j)\left(\Delta{\bf w}_{\tau_j}^{(i_l)}\right)^2-
\int\limits_t^T\varphi_l(s)ds
\right|^4\right\}\right)
^{1/4}=
$$

\vspace{1mm}
$$
=\left({\sf M}\left\{
\left|\sum\limits_{j=0}^{N-1}\left(
\varphi_l(\tau_j)\left(\Delta{\bf w}_{\tau_j}^{(i_l)}\right)^2-
\int\limits_{\tau_j}^{\tau_{j+1}}\varphi_l(s)ds\right)
\right|^4\right\}\right)
^{1/4}\le
$$

\vspace{1mm}
\begin{equation}
\label{df1}
\le \left({\sf M}\left\{
\left|\sum\limits_{j=0}^{N-1}
\varphi_l(\tau_j)\left(\left(\Delta{\bf w}_{\tau_j}^{(i_l)}\right)^2-
\Delta\tau_j\right)
\right|^4\right\}\right)
^{1/4}+
\left|\sum\limits_{j=0}^{N-1}\int\limits_{\tau_j}^{\tau_{j+1}}\left(
\varphi_l(\tau_j)-\varphi_l(s)\right)ds
\right|.
\end{equation}

\vspace{7mm}

From the relation, which is similar to {\rm (\ref{30.43})}, we obtain

\vspace{2mm}
$$
{\sf M}\left\{
\left|\sum\limits_{j=0}^{N-1}
\varphi_l(\tau_j)\left(\left(\Delta{\bf w}_{\tau_j}^{(i_l)}\right)^2-
\Delta\tau_j\right)
\right|^4\right\}=
\sum\limits_{j=0}^{N-1}
\left(\varphi_l(\tau_j)\right)^4
{\sf M}\left\{
\left(\left(\Delta{\bf w}_{\tau_j}^{(i_l)}\right)^2-
\Delta\tau_j\right)^4\right\}+
$$

\vspace{2mm}
$$
+6\sum\limits_{j=0}^{N-1}
\left(\varphi_l(\tau_j)\right)^2
{\sf M}\left\{
\left(\left(\Delta{\bf w}_{\tau_j}^{(i_l)}\right)^2-
\Delta\tau_j\right)^2\right\}
\sum\limits_{q=0}^{j-1}
\left(\varphi_l(\tau_q)\right)^2
{\sf M}\left\{
\left(\left(\Delta{\bf w}_{\tau_q}^{(i_l)}\right)^2-
\Delta\tau_q\right)^2\right\}=
$$

\vspace{4mm}
$$
=60\sum\limits_{j=0}^{N-1}
\left(\varphi_l(\tau_j)\right)^4
\left(\Delta\tau_j\right)^4+24
\sum\limits_{j=0}^{N-1}
\left(\varphi_l(\tau_j)\right)^2
\left(\Delta\tau_j\right)^2
\sum\limits_{q=0}^{j-1}\left(\varphi_l(\tau_q)\right)^2
\left(\Delta\tau_q\right)^2\le C\left(\Delta_N\right)^2\ \to 0
$$

\vspace{5mm}
\noindent
if $N\to \infty$, where constant $C$ does not depend on $N.$

The second term on the right-hand side of 
{\rm (\ref{df1})} tends to zero if $N\to \infty$ 
due to continuity (which means uniform continuity)
of the function $\varphi_l(s)$ on the interval $[t, T].$
Then, taking into account {\rm (\ref{df2})}, {\rm (\ref{30.42})},
we come to the affirmation of Remark {\rm 2.}

Let us prove Theorem 1.
According to Lemma 1, we have

\vspace{2mm}
$$
J[\psi^{(k)}]_{T,t}=
$$

\vspace{2mm}

$$
=
\hbox{\vtop{\offinterlineskip\halign{
\hfil#\hfil\cr
{\rm l.i.m.}\cr
$\stackrel{}{{}_{N\to \infty}}$\cr
}} }\sum_{l_k=0}^{N-1}\ldots\sum_{l_1=0}^{l_2-1}
\psi_1(\tau_{l_1})\ldots\psi_k(\tau_{l_k})
\Delta{\bf w}_{\tau_{l_1}}^{(i_1)}
\ldots
\Delta{\bf w}_{\tau_{l_k}}^{(i_k)}=
$$

\vspace{2mm}
$$
=\hbox{\vtop{\offinterlineskip\halign{
\hfil#\hfil\cr
{\rm l.i.m.}\cr
$\stackrel{}{{}_{N\to \infty}}$\cr
}} }\sum_{l_k=0}^{N-1}\ldots\sum_{l_1=0}^{l_2-1}
K(\tau_{l_1},\ldots,\tau_{l_k})
\Delta{\bf w}_{\tau_{l_1}}^{(i_1)}
\ldots
\Delta{\bf w}_{\tau_{l_k}}^{(i_k)}=
$$

\vspace{2mm}
$$
=
\hbox{\vtop{\offinterlineskip\halign{
\hfil#\hfil\cr
{\rm l.i.m.}\cr
$\stackrel{}{{}_{N\to \infty}}$\cr
}} }\sum_{l_k=0}^{N-1}\ldots\sum_{l_1=0}^{N-1}
K(\tau_{l_1},\ldots,\tau_{l_k})
\Delta{\bf w}_{\tau_{l_1}}^{(i_1)}
\ldots
\Delta{\bf w}_{\tau_{l_k}}^{(i_k)}=
$$

\vspace{2mm}
$$
=\hbox{\vtop{\offinterlineskip\halign{
\hfil#\hfil\cr
{\rm l.i.m.}\cr
$\stackrel{}{{}_{N\to \infty}}$\cr
}} }
\sum\limits_{\stackrel{l_1,\ldots,l_k=0}{{}_{l_q\ne l_r;\ q\ne r;\ q, r=1,\ldots, k}}}^{N-1}
K(\tau_{l_1},\ldots,\tau_{l_k})
\Delta{\bf w}_{\tau_{l_1}}^{(i_1)}
\ldots
\Delta{\bf w}_{\tau_{l_k}}^{(i_k)}=
$$

\vspace{2mm}
\begin{equation}
\label{hehe100}
=
\int\limits_{t}^{T}
\ldots
\int\limits_{t}^{t_2}
\sum_{(t_1,\ldots,t_k)}\left(
K(t_1,\ldots,t_k)d{\bf w}_{t_1}^{(i_1)}
\ldots
d{\bf w}_{t_k}^{(i_k)}\right),
\end{equation}

\vspace{6mm}
\noindent
where permutations 
$(t_1,\ldots,t_k)$ when summing
are performed only 
in the expression, which is enclosed in pa\-ren\-the\-ses.

It is easy to see that (\ref{hehe100})
can be rewritten in the form

$$
J[\psi^{(k)}]_{T,t}=\sum_{(t_1,\ldots,t_k)}
\int\limits_{t}^{T}
\ldots
\int\limits_{t}^{t_2}
K(t_1,\ldots,t_k)d{\bf w}_{t_1}^{(i_1)}
\ldots
d{\bf w}_{t_k}^{(i_k)},
$$

\vspace{3mm}
\noindent
where permutations
$(t_1,\ldots,t_k)$ when summing
are performed only in 
the values
$d{\bf w}_{t_1}^{(i_1)}
\ldots $
$d{\bf w}_{t_k}^{(i_k)}$. At the same time the indexes near upper 
limits of integration in the iterated stochastic integrals are changed 
correspondently and if $t_r$ swapped with $t_q$ in the  
permutation $(t_1,\ldots,t_k)$, then $i_r$ swapped with $i_q$ in 
the permutation $(i_1,\ldots,i_k)$.

Since the integration of a bounded function with respect to the
set of measure zero for Riemann or Lebesgue integrals gives zero result, then the 
following formula is correct for these integrals

\vspace{1mm}
$$
\int\limits_{[t, T]^k}|G(t_1,\ldots,t_k)|dt_1\ldots dt_k=
$$

\vspace{1mm}
\begin{equation}
\label{riemann}
=
\sum_{(t_1,\ldots,t_k)}
\int\limits_{t}^{T}
\ldots
\int\limits_{t}^{t_2}
|G(t_1,\ldots,t_k)|dt_1\ldots dt_k,
\end{equation}

\vspace{5mm}
\noindent
where permutations $(t_1,\ldots,t_k)$ when summing
are performed only 
in the 
va\-lues $dt_1,\ldots, dt_k$. At the same time the indexes near upper 
limits of integration are changed correspondently
and the function $|G(t_1,\ldots,t_k)|$ is assumed to be integrable in 
the hypercube $[t, T]^k.$

According to Lemmas 1, 3 and (\ref{pobeda}), (\ref{hehe100}),
we get the following representation w. p. 1

\vspace{2mm}
$$
J[\psi^{(k)}]_{T,t}=
$$

\vspace{3mm}
$$
=
\sum_{j_1=0}^{p_1}\ldots
\sum_{j_k=0}^{p_k}
C_{j_k\ldots j_1}
\int\limits_{t}^{T}
\ldots
\int\limits_{t}^{t_2}
\sum_{(t_1,\ldots,t_k)}\left(
\phi_{j_1}(t_1)
\ldots
\phi_{j_k}(t_k)
d{\bf w}_{t_1}^{(i_1)}
\ldots
d{\bf w}_{t_k}^{(i_k)}\right)+
$$

\vspace{5mm}
$$
+R_{T,t}^{p_1,\ldots,p_k}=
$$

\vspace{7mm}
$$
=\sum_{j_1=0}^{p_1}\ldots
\sum_{j_k=0}^{p_k}
C_{j_k\ldots j_1}             
\hbox{\vtop{\offinterlineskip\halign{
\hfil#\hfil\cr
{\rm l.i.m.}\cr
$\stackrel{}{{}_{N\to \infty}}$\cr
}} }
\sum\limits_{\stackrel{l_1,\ldots,l_k=0}{{}_{l_q\ne l_r;\ 
q\ne r;\ q, r=1,\ldots, k}}}^{N-1}
\phi_{j_1}(\tau_{l_1})\ldots
\phi_{j_k}(\tau_{l_k})
\Delta{\bf w}_{\tau_{l_1}}^{(i_1)}
\ldots
\Delta{\bf w}_{\tau_{l_k}}^{(i_k)}+
$$

\vspace{5mm}
\begin{equation}
\label{ziko1500}
+R_{T,t}^{p_1,\ldots,p_k}=
\end{equation}

\vspace{7mm}
$$
=\sum_{j_1=0}^{p_1}\ldots
\sum_{j_k=0}^{p_k}
C_{j_k\ldots j_1}\left(
\hbox{\vtop{\offinterlineskip\halign{
\hfil#\hfil\cr
{\rm l.i.m.}\cr
$\stackrel{}{{}_{N\to \infty}}$\cr
}} }\sum_{l_1,\ldots,l_k=0}^{N-1}
\phi_{j_1}(\tau_{l_1})
\ldots
\phi_{j_k}(\tau_{l_k})
\Delta{\bf w}_{\tau_{l_1}}^{(i_1)}
\ldots
\Delta{\bf w}_{\tau_{l_k}}^{(i_k)}
-\right.
$$

\vspace{5mm}
$$
-\left.
\hbox{\vtop{\offinterlineskip\halign{
\hfil#\hfil\cr
{\rm l.i.m.}\cr
$\stackrel{}{{}_{N\to \infty}}$\cr
}} }\sum_{(l_1,\ldots,l_k)\in {\rm G}_k}
\phi_{j_{1}}(\tau_{l_1})
\Delta{\bf w}_{\tau_{l_1}}^{(i_1)}\ldots
\phi_{j_{k}}(\tau_{l_k})
\Delta{\bf w}_{\tau_{l_k}}^{(i_k)}\right)+
$$

\vspace{5mm}
$$
+R_{T,t}^{p_1,\ldots,p_k}=
$$

\vspace{7mm}
$$
=\sum_{j_1=0}^{p_1}\ldots\sum_{j_k=0}^{p_k}
C_{j_k\ldots j_1}\left(
\prod_{l=1}^k\zeta_{j_l}^{(i_l)}-
\hbox{\vtop{\offinterlineskip\halign{
\hfil#\hfil\cr
{\rm l.i.m.}\cr
$\stackrel{}{{}_{N\to \infty}}$\cr
}} }\sum_{(l_1,\ldots,l_k)\in {\rm G}_k}
\phi_{j_{1}}(\tau_{l_1})
\Delta{\bf w}_{\tau_{l_1}}^{(i_1)}\ldots
\phi_{j_{k}}(\tau_{l_k})
\Delta{\bf w}_{\tau_{l_k}}^{(i_k)}\right)+
$$

\vspace{5mm}
\begin{equation}
\label{zara}
+R_{T,t}^{p_1,\ldots,p_k},
\end{equation}

\vspace{5mm}
\noindent
where

$$
R_{T,t}^{p_1,\ldots,p_k}
=\sum_{(t_1,\ldots,t_k)}
\int\limits_{t}^{T}
\ldots
\int\limits_{t}^{t_2}
\left(K(t_1,\ldots,t_k)-
\sum_{j_1=0}^{p_1}\ldots
\sum_{j_k=0}^{p_k}
C_{j_k\ldots j_1}
\prod_{l=1}^k\phi_{j_l}(t_l)\right)\times
$$

\vspace{1mm}
\begin{equation}
\label{y007}
\times
d{\bf w}_{t_1}^{(i_1)}
\ldots
d{\bf w}_{t_k}^{(i_k)},
\end{equation}

\vspace{3mm}
\noindent
where permutations $(t_1,\ldots,t_k)$ when summing are performed only 
in the values $d{\bf w}_{t_1}^{(i_1)}
\ldots $
$d{\bf w}_{t_k}^{(i_k)}$. At the same time the indexes near 
upper limits of integration in the iterated stochastic integrals 
are changed correspondently and if $t_r$ swapped with $t_q$ in the  
permutation $(t_1,\ldots,t_k)$, then $i_r$ swapped with $i_q$ in the 
permutation $(i_1,\ldots,i_k)$.

Let us estimate the remainder
$R_{T,t}^{p_1,\ldots,p_k}$ of the series.

According to Lemma 2 and (\ref{riemann}), we have

\vspace{1mm}
$$
{\sf M}\left\{\left(R_{T,t}^{p_1,\ldots,p_k}\right)^2\right\}
\le 
C_k
\hspace{-2mm}\sum_{(t_1,\ldots,t_k)}
\int\limits_{t}^{T}
\ldots
\int\limits_{t}^{t_2}
\left(K(t_1,\ldots,t_k)-
\sum_{j_1=0}^{p_1}\ldots
\sum_{j_k=0}^{p_k}
C_{j_k\ldots j_1}
\prod_{l=1}^k\phi_{j_l}(t_l)\right)^2
\hspace{-2mm}dt_1
\ldots
dt_k=
$$

\vspace{2mm}
\begin{equation}
\label{obana1}
=C_k\int\limits_{[t,T]^k}
\left(K(t_1,\ldots,t_k)-
\sum_{j_1=0}^{p_1}\ldots
\sum_{j_k=0}^{p_k}
C_{j_k\ldots j_1}
\prod_{l=1}^k\phi_{j_l}(t_l)\right)^2
dt_1
\ldots
dt_k\to 0
\end{equation}

\vspace{5mm}
\noindent
if $p_1,\ldots,p_k\to\infty,$ where constant $C_k$ 
depends only
on the multiplicity $k$ of the iterated Ito stochastic integral. 
Theorem 1 is proved.

Note that from (\ref{ziko1500}) and (\ref{obana1})
it follows that

\vspace{1mm}
\begin{equation}
\label{2023abc11}
J[\psi^{(k)}]_{T,t}=
\hbox{\vtop{\offinterlineskip\halign{
\hfil#\hfil\cr
{\rm l.i.m.}\cr
$\stackrel{}{{}_{p_1,\ldots,p_k\to \infty}}$\cr
}} }\sum_{j_1=0}^{p_1}\ldots\sum_{j_k=0}^{p_k}
C_{j_k\ldots j_1}J'[\phi_{j_1}\ldots \phi_{j_k}]_{T,t}^{(i_1\ldots i_k)},
\end{equation}

\vspace{4mm}
\noindent
where $J'[\phi_{j_1}\ldots \phi_{j_k}]_{T,t}^{(i_1\ldots i_k)}$ is defined by
(\ref{mult11}).

It is not difficult to see that for the case of pairwise different numbers
$i_1,\ldots,i_k=1,\ldots,m$ from Theorem 1 we obtain

$$
J[\psi^{(k)}]_{T,t}=
\hbox{\vtop{\offinterlineskip\halign{
\hfil#\hfil\cr
{\rm l.i.m.}\cr
$\stackrel{}{{}_{p_1,\ldots,p_k\to \infty}}$\cr
}} }\sum_{j_1=0}^{p_1}\ldots\sum_{j_k=0}^{p_k}
C_{j_k\ldots j_1}\zeta_{j_1}^{(i_1)}\ldots \zeta_{j_k}^{(i_k)}.
$$

\vspace{4mm}

In order to evaluate the significance of Theorem 1 for practice we will
demonstrate its transformed particular cases (see Remark 2) for 
$k=1,\ldots,7$ \cite{7}-\cite{Kuzh-1}

\vspace{1mm}
\begin{equation}
\label{a1}
J[\psi^{(1)}]_{T,t}
=\hbox{\vtop{\offinterlineskip\halign{
\hfil#\hfil\cr
{\rm l.i.m.}\cr
$\stackrel{}{{}_{p_1\to \infty}}$\cr
}} }\sum_{j_1=0}^{p_1}
C_{j_1}\zeta_{j_1}^{(i_1)},
\end{equation}

\vspace{2mm}
\begin{equation}
\label{a2}
J[\psi^{(2)}]_{T,t}
=\hbox{\vtop{\offinterlineskip\halign{
\hfil#\hfil\cr
{\rm l.i.m.}\cr
$\stackrel{}{{}_{p_1,p_2\to \infty}}$\cr
}} }\sum_{j_1=0}^{p_1}\sum_{j_2=0}^{p_2}
C_{j_2j_1}\Biggl(\zeta_{j_1}^{(i_1)}\zeta_{j_2}^{(i_2)}
-{\bf 1}_{\{i_1=i_2\ne 0\}}
{\bf 1}_{\{j_1=j_2\}}\Biggr),
\end{equation}

\vspace{7mm}
$$
J[\psi^{(3)}]_{T,t}=
\hbox{\vtop{\offinterlineskip\halign{
\hfil#\hfil\cr
{\rm l.i.m.}\cr
$\stackrel{}{{}_{p_1,\ldots,p_3\to \infty}}$\cr
}} }\sum_{j_1=0}^{p_1}\sum_{j_2=0}^{p_2}\sum_{j_3=0}^{p_3}
C_{j_3j_2j_1}\Biggl(
\zeta_{j_1}^{(i_1)}\zeta_{j_2}^{(i_2)}\zeta_{j_3}^{(i_3)}
-\Biggr.
$$

\begin{equation}
\label{a3}
-\Biggl.
{\bf 1}_{\{i_1=i_2\ne 0\}}
{\bf 1}_{\{j_1=j_2\}}
\zeta_{j_3}^{(i_3)}
-{\bf 1}_{\{i_2=i_3\ne 0\}}
{\bf 1}_{\{j_2=j_3\}}
\zeta_{j_1}^{(i_1)}-
{\bf 1}_{\{i_1=i_3\ne 0\}}
{\bf 1}_{\{j_1=j_3\}}
\zeta_{j_2}^{(i_2)}\Biggr),
\end{equation}

\vspace{7mm}
$$
J[\psi^{(4)}]_{T,t}
=
\hbox{\vtop{\offinterlineskip\halign{
\hfil#\hfil\cr
{\rm l.i.m.}\cr
$\stackrel{}{{}_{p_1,\ldots,p_4\to \infty}}$\cr
}} }\sum_{j_1=0}^{p_1}\ldots\sum_{j_4=0}^{p_4}
C_{j_4\ldots j_1}\Biggl(
\prod_{l=1}^4\zeta_{j_l}^{(i_l)}
\Biggr.
-
$$

$$
-
{\bf 1}_{\{i_1=i_2\ne 0\}}
{\bf 1}_{\{j_1=j_2\}}
\zeta_{j_3}^{(i_3)}
\zeta_{j_4}^{(i_4)}
-
{\bf 1}_{\{i_1=i_3\ne 0\}}
{\bf 1}_{\{j_1=j_3\}}
\zeta_{j_2}^{(i_2)}
\zeta_{j_4}^{(i_4)}-
$$

$$
-
{\bf 1}_{\{i_1=i_4\ne 0\}}
{\bf 1}_{\{j_1=j_4\}}
\zeta_{j_2}^{(i_2)}
\zeta_{j_3}^{(i_3)}
-
{\bf 1}_{\{i_2=i_3\ne 0\}}
{\bf 1}_{\{j_2=j_3\}}
\zeta_{j_1}^{(i_1)}
\zeta_{j_4}^{(i_4)}-
$$

$$
-
{\bf 1}_{\{i_2=i_4\ne 0\}}
{\bf 1}_{\{j_2=j_4\}}
\zeta_{j_1}^{(i_1)}
\zeta_{j_3}^{(i_3)}
-
{\bf 1}_{\{i_3=i_4\ne 0\}}
{\bf 1}_{\{j_3=j_4\}}
\zeta_{j_1}^{(i_1)}
\zeta_{j_2}^{(i_2)}+
$$

$$
+
{\bf 1}_{\{i_1=i_2\ne 0\}}
{\bf 1}_{\{j_1=j_2\}}
{\bf 1}_{\{i_3=i_4\ne 0\}}
{\bf 1}_{\{j_3=j_4\}}
+
{\bf 1}_{\{i_1=i_3\ne 0\}}
{\bf 1}_{\{j_1=j_3\}}
{\bf 1}_{\{i_2=i_4\ne 0\}}
{\bf 1}_{\{j_2=j_4\}}+
$$
\begin{equation}
\label{a4}
+\Biggl.
{\bf 1}_{\{i_1=i_4\ne 0\}}
{\bf 1}_{\{j_1=j_4\}}
{\bf 1}_{\{i_2=i_3\ne 0\}}
{\bf 1}_{\{j_2=j_3\}}\Biggr),
\end{equation}

\vspace{9mm}
$$
J[\psi^{(5)}]_{T,t}
=\hbox{\vtop{\offinterlineskip\halign{
\hfil#\hfil\cr
{\rm l.i.m.}\cr
$\stackrel{}{{}_{p_1,\ldots,p_5\to \infty}}$\cr
}} }\sum_{j_1=0}^{p_1}\ldots\sum_{j_5=0}^{p_5}
C_{j_5\ldots j_1}\Biggl(
\prod_{l=1}^5\zeta_{j_l}^{(i_l)}
-\Biggr.
$$

\vspace{1mm}
$$
-
{\bf 1}_{\{i_1=i_2\ne 0\}}
{\bf 1}_{\{j_1=j_2\}}
\zeta_{j_3}^{(i_3)}
\zeta_{j_4}^{(i_4)}
\zeta_{j_5}^{(i_5)}-
{\bf 1}_{\{i_1=i_3\ne 0\}}
{\bf 1}_{\{j_1=j_3\}}
\zeta_{j_2}^{(i_2)}
\zeta_{j_4}^{(i_4)}
\zeta_{j_5}^{(i_5)}-
$$

$$
-
{\bf 1}_{\{i_1=i_4\ne 0\}}
{\bf 1}_{\{j_1=j_4\}}
\zeta_{j_2}^{(i_2)}
\zeta_{j_3}^{(i_3)}
\zeta_{j_5}^{(i_5)}-
{\bf 1}_{\{i_1=i_5\ne 0\}}
{\bf 1}_{\{j_1=j_5\}}
\zeta_{j_2}^{(i_2)}
\zeta_{j_3}^{(i_3)}
\zeta_{j_4}^{(i_4)}-
$$

$$
-
{\bf 1}_{\{i_2=i_3\ne 0\}}
{\bf 1}_{\{j_2=j_3\}}
\zeta_{j_1}^{(i_1)}
\zeta_{j_4}^{(i_4)}
\zeta_{j_5}^{(i_5)}-
{\bf 1}_{\{i_2=i_4\ne 0\}}
{\bf 1}_{\{j_2=j_4\}}
\zeta_{j_1}^{(i_1)}
\zeta_{j_3}^{(i_3)}
\zeta_{j_5}^{(i_5)}-
$$

$$
-
{\bf 1}_{\{i_2=i_5\ne 0\}}
{\bf 1}_{\{j_2=j_5\}}
\zeta_{j_1}^{(i_1)}
\zeta_{j_3}^{(i_3)}
\zeta_{j_4}^{(i_4)}
-{\bf 1}_{\{i_3=i_4\ne 0\}}
{\bf 1}_{\{j_3=j_4\}}
\zeta_{j_1}^{(i_1)}
\zeta_{j_2}^{(i_2)}
\zeta_{j_5}^{(i_5)}-
$$

$$
-
{\bf 1}_{\{i_3=i_5\ne 0\}}
{\bf 1}_{\{j_3=j_5\}}
\zeta_{j_1}^{(i_1)}
\zeta_{j_2}^{(i_2)}
\zeta_{j_4}^{(i_4)}
-{\bf 1}_{\{i_4=i_5\ne 0\}}
{\bf 1}_{\{j_4=j_5\}}
\zeta_{j_1}^{(i_1)}
\zeta_{j_2}^{(i_2)}
\zeta_{j_3}^{(i_3)}+
$$

\vspace{-2mm}
$$
+
{\bf 1}_{\{i_1=i_2\ne 0\}}
{\bf 1}_{\{j_1=j_2\}}
{\bf 1}_{\{i_3=i_4\ne 0\}}
{\bf 1}_{\{j_3=j_4\}}\zeta_{j_5}^{(i_5)}+
{\bf 1}_{\{i_1=i_2\ne 0\}}
{\bf 1}_{\{j_1=j_2\}}
{\bf 1}_{\{i_3=i_5\ne 0\}}
{\bf 1}_{\{j_3=j_5\}}\zeta_{j_4}^{(i_4)}+
$$

\vspace{-2mm}
$$
+
{\bf 1}_{\{i_1=i_2\ne 0\}}
{\bf 1}_{\{j_1=j_2\}}
{\bf 1}_{\{i_4=i_5\ne 0\}}
{\bf 1}_{\{j_4=j_5\}}\zeta_{j_3}^{(i_3)}+
{\bf 1}_{\{i_1=i_3\ne 0\}}
{\bf 1}_{\{j_1=j_3\}}
{\bf 1}_{\{i_2=i_4\ne 0\}}
{\bf 1}_{\{j_2=j_4\}}\zeta_{j_5}^{(i_5)}+
$$

\vspace{-2mm}
$$
+
{\bf 1}_{\{i_1=i_3\ne 0\}}
{\bf 1}_{\{j_1=j_3\}}
{\bf 1}_{\{i_2=i_5\ne 0\}}
{\bf 1}_{\{j_2=j_5\}}\zeta_{j_4}^{(i_4)}+
{\bf 1}_{\{i_1=i_3\ne 0\}}
{\bf 1}_{\{j_1=j_3\}}
{\bf 1}_{\{i_4=i_5\ne 0\}}
{\bf 1}_{\{j_4=j_5\}}\zeta_{j_2}^{(i_2)}+
$$

\vspace{-2mm}
$$
+
{\bf 1}_{\{i_1=i_4\ne 0\}}
{\bf 1}_{\{j_1=j_4\}}
{\bf 1}_{\{i_2=i_3\ne 0\}}
{\bf 1}_{\{j_2=j_3\}}\zeta_{j_5}^{(i_5)}+
{\bf 1}_{\{i_1=i_4\ne 0\}}
{\bf 1}_{\{j_1=j_4\}}
{\bf 1}_{\{i_2=i_5\ne 0\}}
{\bf 1}_{\{j_2=j_5\}}\zeta_{j_3}^{(i_3)}+
$$

\vspace{-2mm}
$$
+
{\bf 1}_{\{i_1=i_4\ne 0\}}
{\bf 1}_{\{j_1=j_4\}}
{\bf 1}_{\{i_3=i_5\ne 0\}}
{\bf 1}_{\{j_3=j_5\}}\zeta_{j_2}^{(i_2)}+
{\bf 1}_{\{i_1=i_5\ne 0\}}
{\bf 1}_{\{j_1=j_5\}}
{\bf 1}_{\{i_2=i_3\ne 0\}}
{\bf 1}_{\{j_2=j_3\}}\zeta_{j_4}^{(i_4)}+
$$

\vspace{-2mm}
$$
+
{\bf 1}_{\{i_1=i_5\ne 0\}}
{\bf 1}_{\{j_1=j_5\}}
{\bf 1}_{\{i_2=i_4\ne 0\}}
{\bf 1}_{\{j_2=j_4\}}\zeta_{j_3}^{(i_3)}+
{\bf 1}_{\{i_1=i_5\ne 0\}}
{\bf 1}_{\{j_1=j_5\}}
{\bf 1}_{\{i_3=i_4\ne 0\}}
{\bf 1}_{\{j_3=j_4\}}\zeta_{j_2}^{(i_2)}+
$$

\vspace{-2mm}
$$
+
{\bf 1}_{\{i_2=i_3\ne 0\}}
{\bf 1}_{\{j_2=j_3\}}
{\bf 1}_{\{i_4=i_5\ne 0\}}
{\bf 1}_{\{j_4=j_5\}}\zeta_{j_1}^{(i_1)}+
{\bf 1}_{\{i_2=i_4\ne 0\}}
{\bf 1}_{\{j_2=j_4\}}
{\bf 1}_{\{i_3=i_5\ne 0\}}
{\bf 1}_{\{j_3=j_5\}}\zeta_{j_1}^{(i_1)}+
$$
\begin{equation}
\label{a5}
+\Biggl.
{\bf 1}_{\{i_2=i_5\ne 0\}}
{\bf 1}_{\{j_2=j_5\}}
{\bf 1}_{\{i_3=i_4\ne 0\}}
{\bf 1}_{\{j_3=j_4\}}\zeta_{j_1}^{(i_1)}\Biggr),
\end{equation}

\vspace{10mm}
$$
J[\psi^{(6)}]_{T,t}
=\hbox{\vtop{\offinterlineskip\halign{
\hfil#\hfil\cr
{\rm l.i.m.}\cr
$\stackrel{}{{}_{p_1,\ldots,p_6\to \infty}}$\cr
}} }\sum_{j_1=0}^{p_1}\ldots\sum_{j_6=0}^{p_6}
C_{j_6\ldots j_1}\Biggl(
\prod_{l=1}^6
\zeta_{j_l}^{(i_l)}
-\Biggr.
$$

\vspace{1mm}
$$
-
{\bf 1}_{\{i_1=i_6\ne 0\}}
{\bf 1}_{\{j_1=j_6\}}
\zeta_{j_2}^{(i_2)}
\zeta_{j_3}^{(i_3)}
\zeta_{j_4}^{(i_4)}
\zeta_{j_5}^{(i_5)}-
{\bf 1}_{\{i_2=i_6\ne 0\}}
{\bf 1}_{\{j_2=j_6\}}
\zeta_{j_1}^{(i_1)}
\zeta_{j_3}^{(i_3)}
\zeta_{j_4}^{(i_4)}
\zeta_{j_5}^{(i_5)}-
$$

$$
-
{\bf 1}_{\{i_3=i_6\ne 0\}}
{\bf 1}_{\{j_3=j_6\}}
\zeta_{j_1}^{(i_1)}
\zeta_{j_2}^{(i_2)}
\zeta_{j_4}^{(i_4)}
\zeta_{j_5}^{(i_5)}-
{\bf 1}_{\{i_4=i_6\ne 0\}}
{\bf 1}_{\{j_4=j_6\}}
\zeta_{j_1}^{(i_1)}
\zeta_{j_2}^{(i_2)}
\zeta_{j_3}^{(i_3)}
\zeta_{j_5}^{(i_5)}-
$$

$$
-
{\bf 1}_{\{i_5=i_6\ne 0\}}
{\bf 1}_{\{j_5=j_6\}}
\zeta_{j_1}^{(i_1)}
\zeta_{j_2}^{(i_2)}
\zeta_{j_3}^{(i_3)}
\zeta_{j_4}^{(i_4)}-
{\bf 1}_{\{i_1=i_2\ne 0\}}
{\bf 1}_{\{j_1=j_2\}}
\zeta_{j_3}^{(i_3)}
\zeta_{j_4}^{(i_4)}
\zeta_{j_5}^{(i_5)}
\zeta_{j_6}^{(i_6)}-
$$

$$
-
{\bf 1}_{\{i_1=i_3\ne 0\}}
{\bf 1}_{\{j_1=j_3\}}
\zeta_{j_2}^{(i_2)}
\zeta_{j_4}^{(i_4)}
\zeta_{j_5}^{(i_5)}
\zeta_{j_6}^{(i_6)}-
{\bf 1}_{\{i_1=i_4\ne 0\}}
{\bf 1}_{\{j_1=j_4\}}
\zeta_{j_2}^{(i_2)}
\zeta_{j_3}^{(i_3)}
\zeta_{j_5}^{(i_5)}
\zeta_{j_6}^{(i_6)}-
$$

$$
-
{\bf 1}_{\{i_1=i_5\ne 0\}}
{\bf 1}_{\{j_1=j_5\}}
\zeta_{j_2}^{(i_2)}
\zeta_{j_3}^{(i_3)}
\zeta_{j_4}^{(i_4)}
\zeta_{j_6}^{(i_6)}-
{\bf 1}_{\{i_2=i_3\ne 0\}}
{\bf 1}_{\{j_2=j_3\}}
\zeta_{j_1}^{(i_1)}
\zeta_{j_4}^{(i_4)}
\zeta_{j_5}^{(i_5)}
\zeta_{j_6}^{(i_6)}-
$$

$$
-
{\bf 1}_{\{i_2=i_4\ne 0\}}
{\bf 1}_{\{j_2=j_4\}}
\zeta_{j_1}^{(i_1)}
\zeta_{j_3}^{(i_3)}
\zeta_{j_5}^{(i_5)}
\zeta_{j_6}^{(i_6)}-
{\bf 1}_{\{i_2=i_5\ne 0\}}
{\bf 1}_{\{j_2=j_5\}}
\zeta_{j_1}^{(i_1)}
\zeta_{j_3}^{(i_3)}
\zeta_{j_4}^{(i_4)}
\zeta_{j_6}^{(i_6)}-
$$

$$
-
{\bf 1}_{\{i_3=i_4\ne 0\}}
{\bf 1}_{\{j_3=j_4\}}
\zeta_{j_1}^{(i_1)}
\zeta_{j_2}^{(i_2)}
\zeta_{j_5}^{(i_5)}
\zeta_{j_6}^{(i_6)}-
{\bf 1}_{\{i_3=i_5\ne 0\}}
{\bf 1}_{\{j_3=j_5\}}
\zeta_{j_1}^{(i_1)}
\zeta_{j_2}^{(i_2)}
\zeta_{j_4}^{(i_4)}
\zeta_{j_6}^{(i_6)}-
$$

$$
-
{\bf 1}_{\{i_4=i_5\ne 0\}}
{\bf 1}_{\{j_4=j_5\}}
\zeta_{j_1}^{(i_1)}
\zeta_{j_2}^{(i_2)}
\zeta_{j_3}^{(i_3)}
\zeta_{j_6}^{(i_6)}+
$$

\vspace{-1.5mm}
$$
+
{\bf 1}_{\{i_1=i_2\ne 0\}}
{\bf 1}_{\{j_1=j_2\}}
{\bf 1}_{\{i_3=i_4\ne 0\}}
{\bf 1}_{\{j_3=j_4\}}
\zeta_{j_5}^{(i_5)}
\zeta_{j_6}^{(i_6)}+
{\bf 1}_{\{i_1=i_2\ne 0\}}
{\bf 1}_{\{j_1=j_2\}}
{\bf 1}_{\{i_3=i_5\ne 0\}}
{\bf 1}_{\{j_3=j_5\}}
\zeta_{j_4}^{(i_4)}
\zeta_{j_6}^{(i_6)}+
$$

\vspace{-2.7mm}
$$
+
{\bf 1}_{\{i_1=i_2\ne 0\}}
{\bf 1}_{\{j_1=j_2\}}
{\bf 1}_{\{i_4=i_5\ne 0\}}
{\bf 1}_{\{j_4=j_5\}}
\zeta_{j_3}^{(i_3)}
\zeta_{j_6}^{(i_6)}
+
{\bf 1}_{\{i_1=i_3\ne 0\}}
{\bf 1}_{\{j_1=j_3\}}
{\bf 1}_{\{i_2=i_4\ne 0\}}
{\bf 1}_{\{j_2=j_4\}}
\zeta_{j_5}^{(i_5)}
\zeta_{j_6}^{(i_6)}+
$$

\vspace{-2.7mm}
$$
+
{\bf 1}_{\{i_1=i_3\ne 0\}}
{\bf 1}_{\{j_1=j_3\}}
{\bf 1}_{\{i_2=i_5\ne 0\}}
{\bf 1}_{\{j_2=j_5\}}
\zeta_{j_4}^{(i_4)}
\zeta_{j_6}^{(i_6)}
+{\bf 1}_{\{i_1=i_3\ne 0\}}
{\bf 1}_{\{j_1=j_3\}}
{\bf 1}_{\{i_4=i_5\ne 0\}}
{\bf 1}_{\{j_4=j_5\}}
\zeta_{j_2}^{(i_2)}
\zeta_{j_6}^{(i_6)}+
$$

\vspace{-2.7mm}
$$
+
{\bf 1}_{\{i_1=i_4\ne 0\}}
{\bf 1}_{\{j_1=j_4\}}
{\bf 1}_{\{i_2=i_3\ne 0\}}
{\bf 1}_{\{j_2=j_3\}}
\zeta_{j_5}^{(i_5)}
\zeta_{j_6}^{(i_6)}
+
{\bf 1}_{\{i_1=i_4\ne 0\}}
{\bf 1}_{\{j_1=j_4\}}
{\bf 1}_{\{i_2=i_5\ne 0\}}
{\bf 1}_{\{j_2=j_5\}}
\zeta_{j_3}^{(i_3)}
\zeta_{j_6}^{(i_6)}+
$$

\vspace{-2.7mm}
$$
+
{\bf 1}_{\{i_1=i_4\ne 0\}}
{\bf 1}_{\{j_1=j_4\}}
{\bf 1}_{\{i_3=i_5\ne 0\}}
{\bf 1}_{\{j_3=j_5\}}
\zeta_{j_2}^{(i_2)}
\zeta_{j_6}^{(i_6)}
+
{\bf 1}_{\{i_1=i_5\ne 0\}}
{\bf 1}_{\{j_1=j_5\}}
{\bf 1}_{\{i_2=i_3\ne 0\}}
{\bf 1}_{\{j_2=j_3\}}
\zeta_{j_4}^{(i_4)}
\zeta_{j_6}^{(i_6)}+
$$

\vspace{-2.7mm}
$$
+
{\bf 1}_{\{i_1=i_5\ne 0\}}
{\bf 1}_{\{j_1=j_5\}}
{\bf 1}_{\{i_2=i_4\ne 0\}}
{\bf 1}_{\{j_2=j_4\}}
\zeta_{j_3}^{(i_3)}
\zeta_{j_6}^{(i_6)}
+
{\bf 1}_{\{i_1=i_5\ne 0\}}
{\bf 1}_{\{j_1=j_5\}}
{\bf 1}_{\{i_3=i_4\ne 0\}}
{\bf 1}_{\{j_3=j_4\}}
\zeta_{j_2}^{(i_2)}
\zeta_{j_6}^{(i_6)}+
$$

\vspace{-2.7mm}
$$
+
{\bf 1}_{\{i_2=i_3\ne 0\}}
{\bf 1}_{\{j_2=j_3\}}
{\bf 1}_{\{i_4=i_5\ne 0\}}
{\bf 1}_{\{j_4=j_5\}}
\zeta_{j_1}^{(i_1)}
\zeta_{j_6}^{(i_6)}
+
{\bf 1}_{\{i_2=i_4\ne 0\}}
{\bf 1}_{\{j_2=j_4\}}
{\bf 1}_{\{i_3=i_5\ne 0\}}
{\bf 1}_{\{j_3=j_5\}}
\zeta_{j_1}^{(i_1)}
\zeta_{j_6}^{(i_6)}+
$$

\vspace{-2.7mm}
$$
+
{\bf 1}_{\{i_2=i_5\ne 0\}}
{\bf 1}_{\{j_2=j_5\}}
{\bf 1}_{\{i_3=i_4\ne 0\}}
{\bf 1}_{\{j_3=j_4\}}
\zeta_{j_1}^{(i_1)}
\zeta_{j_6}^{(i_6)}
+
{\bf 1}_{\{i_6=i_1\ne 0\}}
{\bf 1}_{\{j_6=j_1\}}
{\bf 1}_{\{i_3=i_4\ne 0\}}
{\bf 1}_{\{j_3=j_4\}}
\zeta_{j_2}^{(i_2)}
\zeta_{j_5}^{(i_5)}+
$$

\vspace{-2.7mm}
$$
+
{\bf 1}_{\{i_6=i_1\ne 0\}}
{\bf 1}_{\{j_6=j_1\}}
{\bf 1}_{\{i_3=i_5\ne 0\}}
{\bf 1}_{\{j_3=j_5\}}
\zeta_{j_2}^{(i_2)}
\zeta_{j_4}^{(i_4)}
+
{\bf 1}_{\{i_6=i_1\ne 0\}}
{\bf 1}_{\{j_6=j_1\}}
{\bf 1}_{\{i_2=i_5\ne 0\}}
{\bf 1}_{\{j_2=j_5\}}
\zeta_{j_3}^{(i_3)}
\zeta_{j_4}^{(i_4)}+
$$

\vspace{-2.7mm}
$$
+
{\bf 1}_{\{i_6=i_1\ne 0\}}
{\bf 1}_{\{j_6=j_1\}}
{\bf 1}_{\{i_2=i_4\ne 0\}}
{\bf 1}_{\{j_2=j_4\}}
\zeta_{j_3}^{(i_3)}
\zeta_{j_5}^{(i_5)}
+
{\bf 1}_{\{i_6=i_1\ne 0\}}
{\bf 1}_{\{j_6=j_1\}}
{\bf 1}_{\{i_4=i_5\ne 0\}}
{\bf 1}_{\{j_4=j_5\}}
\zeta_{j_2}^{(i_2)}
\zeta_{j_3}^{(i_3)}+
$$

\vspace{-2.7mm}
$$
+
{\bf 1}_{\{i_6=i_1\ne 0\}}
{\bf 1}_{\{j_6=j_1\}}
{\bf 1}_{\{i_2=i_3\ne 0\}}
{\bf 1}_{\{j_2=j_3\}}
\zeta_{j_4}^{(i_4)}
\zeta_{j_5}^{(i_5)}
+
{\bf 1}_{\{i_6=i_2\ne 0\}}
{\bf 1}_{\{j_6=j_2\}}
{\bf 1}_{\{i_3=i_5\ne 0\}}
{\bf 1}_{\{j_3=j_5\}}
\zeta_{j_1}^{(i_1)}
\zeta_{j_4}^{(i_4)}+
$$

\vspace{-2.7mm}
$$
+
{\bf 1}_{\{i_6=i_2\ne 0\}}
{\bf 1}_{\{j_6=j_2\}}
{\bf 1}_{\{i_4=i_5\ne 0\}}
{\bf 1}_{\{j_4=j_5\}}
\zeta_{j_1}^{(i_1)}
\zeta_{j_3}^{(i_3)}
+
{\bf 1}_{\{i_6=i_2\ne 0\}}
{\bf 1}_{\{j_6=j_2\}}
{\bf 1}_{\{i_3=i_4\ne 0\}}
{\bf 1}_{\{j_3=j_4\}}
\zeta_{j_1}^{(i_1)}
\zeta_{j_5}^{(i_5)}+
$$

\vspace{-2.7mm}
$$
+
{\bf 1}_{\{i_6=i_2\ne 0\}}
{\bf 1}_{\{j_6=j_2\}}
{\bf 1}_{\{i_1=i_5\ne 0\}}
{\bf 1}_{\{j_1=j_5\}}
\zeta_{j_3}^{(i_3)}
\zeta_{j_4}^{(i_4)}
+
{\bf 1}_{\{i_6=i_2\ne 0\}}
{\bf 1}_{\{j_6=j_2\}}
{\bf 1}_{\{i_1=i_4\ne 0\}}
{\bf 1}_{\{j_1=j_4\}}
\zeta_{j_3}^{(i_3)}
\zeta_{j_5}^{(i_5)}+
$$

\vspace{-2.7mm}
$$
+
{\bf 1}_{\{i_6=i_2\ne 0\}}
{\bf 1}_{\{j_6=j_2\}}
{\bf 1}_{\{i_1=i_3\ne 0\}}
{\bf 1}_{\{j_1=j_3\}}
\zeta_{j_4}^{(i_4)}
\zeta_{j_5}^{(i_5)}
+
{\bf 1}_{\{i_6=i_3\ne 0\}}
{\bf 1}_{\{j_6=j_3\}}
{\bf 1}_{\{i_2=i_5\ne 0\}}
{\bf 1}_{\{j_2=j_5\}}
\zeta_{j_1}^{(i_1)}
\zeta_{j_4}^{(i_4)}+
$$

\vspace{-2.7mm}
$$
+
{\bf 1}_{\{i_6=i_3\ne 0\}}
{\bf 1}_{\{j_6=j_3\}}
{\bf 1}_{\{i_4=i_5\ne 0\}}
{\bf 1}_{\{j_4=j_5\}}
\zeta_{j_1}^{(i_1)}
\zeta_{j_2}^{(i_2)}
+
{\bf 1}_{\{i_6=i_3\ne 0\}}
{\bf 1}_{\{j_6=j_3\}}
{\bf 1}_{\{i_2=i_4\ne 0\}}
{\bf 1}_{\{j_2=j_4\}}
\zeta_{j_1}^{(i_1)}
\zeta_{j_5}^{(i_5)}+
$$

\vspace{-2.7mm}
$$
+
{\bf 1}_{\{i_6=i_3\ne 0\}}
{\bf 1}_{\{j_6=j_3\}}
{\bf 1}_{\{i_1=i_5\ne 0\}}
{\bf 1}_{\{j_1=j_5\}}
\zeta_{j_2}^{(i_2)}
\zeta_{j_4}^{(i_4)}
+
{\bf 1}_{\{i_6=i_3\ne 0\}}
{\bf 1}_{\{j_6=j_3\}}
{\bf 1}_{\{i_1=i_4\ne 0\}}
{\bf 1}_{\{j_1=j_4\}}
\zeta_{j_2}^{(i_2)}
\zeta_{j_5}^{(i_5)}+
$$

\vspace{-2.7mm}
$$
+
{\bf 1}_{\{i_6=i_3\ne 0\}}
{\bf 1}_{\{j_6=j_3\}}
{\bf 1}_{\{i_1=i_2\ne 0\}}
{\bf 1}_{\{j_1=j_2\}}
\zeta_{j_4}^{(i_4)}
\zeta_{j_5}^{(i_5)}
+
{\bf 1}_{\{i_6=i_4\ne 0\}}
{\bf 1}_{\{j_6=j_4\}}
{\bf 1}_{\{i_3=i_5\ne 0\}}
{\bf 1}_{\{j_3=j_5\}}
\zeta_{j_1}^{(i_1)}
\zeta_{j_2}^{(i_2)}+
$$

\vspace{-2.7mm}
$$
+
{\bf 1}_{\{i_6=i_4\ne 0\}}
{\bf 1}_{\{j_6=j_4\}}
{\bf 1}_{\{i_2=i_5\ne 0\}}
{\bf 1}_{\{j_2=j_5\}}
\zeta_{j_1}^{(i_1)}
\zeta_{j_3}^{(i_3)}
+
{\bf 1}_{\{i_6=i_4\ne 0\}}
{\bf 1}_{\{j_6=j_4\}}
{\bf 1}_{\{i_2=i_3\ne 0\}}
{\bf 1}_{\{j_2=j_3\}}
\zeta_{j_1}^{(i_1)}
\zeta_{j_5}^{(i_5)}+
$$

\vspace{-2.7mm}
$$
+
{\bf 1}_{\{i_6=i_4\ne 0\}}
{\bf 1}_{\{j_6=j_4\}}
{\bf 1}_{\{i_1=i_5\ne 0\}}
{\bf 1}_{\{j_1=j_5\}}
\zeta_{j_2}^{(i_2)}
\zeta_{j_3}^{(i_3)}
+
{\bf 1}_{\{i_6=i_4\ne 0\}}
{\bf 1}_{\{j_6=j_4\}}
{\bf 1}_{\{i_1=i_3\ne 0\}}
{\bf 1}_{\{j_1=j_3\}}
\zeta_{j_2}^{(i_2)}
\zeta_{j_5}^{(i_5)}+
$$

\vspace{-2.7mm}
$$
+
{\bf 1}_{\{i_6=i_4\ne 0\}}
{\bf 1}_{\{j_6=j_4\}}
{\bf 1}_{\{i_1=i_2\ne 0\}}
{\bf 1}_{\{j_1=j_2\}}
\zeta_{j_3}^{(i_3)}
\zeta_{j_5}^{(i_5)}
+
{\bf 1}_{\{i_6=i_5\ne 0\}}
{\bf 1}_{\{j_6=j_5\}}
{\bf 1}_{\{i_3=i_4\ne 0\}}
{\bf 1}_{\{j_3=j_4\}}
\zeta_{j_1}^{(i_1)}
\zeta_{j_2}^{(i_2)}+
$$

\vspace{-2.7mm}
$$
+
{\bf 1}_{\{i_6=i_5\ne 0\}}
{\bf 1}_{\{j_6=j_5\}}
{\bf 1}_{\{i_2=i_4\ne 0\}}
{\bf 1}_{\{j_2=j_4\}}
\zeta_{j_1}^{(i_1)}
\zeta_{j_3}^{(i_3)}
+
{\bf 1}_{\{i_6=i_5\ne 0\}}
{\bf 1}_{\{j_6=j_5\}}
{\bf 1}_{\{i_2=i_3\ne 0\}}
{\bf 1}_{\{j_2=j_3\}}
\zeta_{j_1}^{(i_1)}
\zeta_{j_4}^{(i_4)}+
$$

\vspace{-2.7mm}
$$
+
{\bf 1}_{\{i_6=i_5\ne 0\}}
{\bf 1}_{\{j_6=j_5\}}
{\bf 1}_{\{i_1=i_4\ne 0\}}
{\bf 1}_{\{j_1=j_4\}}
\zeta_{j_2}^{(i_2)}
\zeta_{j_3}^{(i_3)}
+
{\bf 1}_{\{i_6=i_5\ne 0\}}
{\bf 1}_{\{j_6=j_5\}}
{\bf 1}_{\{i_1=i_3\ne 0\}}
{\bf 1}_{\{j_1=j_3\}}
\zeta_{j_2}^{(i_2)}
\zeta_{j_4}^{(i_4)}+
$$

\vspace{-1.5mm}
$$
+
{\bf 1}_{\{i_6=i_5\ne 0\}}
{\bf 1}_{\{j_6=j_5\}}
{\bf 1}_{\{i_1=i_2\ne 0\}}
{\bf 1}_{\{j_1=j_2\}}
\zeta_{j_3}^{(i_3)}
\zeta_{j_4}^{(i_4)}-
$$

\vspace{-1.5mm}
$$
-
{\bf 1}_{\{i_6=i_1\ne 0\}}
{\bf 1}_{\{j_6=j_1\}}
{\bf 1}_{\{i_2=i_5\ne 0\}}
{\bf 1}_{\{j_2=j_5\}}
{\bf 1}_{\{i_3=i_4\ne 0\}}
{\bf 1}_{\{j_3=j_4\}}-
$$

\vspace{-2mm}
$$
-
{\bf 1}_{\{i_6=i_1\ne 0\}}
{\bf 1}_{\{j_6=j_1\}}
{\bf 1}_{\{i_2=i_4\ne 0\}}
{\bf 1}_{\{j_2=j_4\}}
{\bf 1}_{\{i_3=i_5\ne 0\}}
{\bf 1}_{\{j_3=j_5\}}-
$$

\vspace{-2mm}
$$
-
{\bf 1}_{\{i_6=i_1\ne 0\}}
{\bf 1}_{\{j_6=j_1\}}
{\bf 1}_{\{i_2=i_3\ne 0\}}
{\bf 1}_{\{j_2=j_3\}}
{\bf 1}_{\{i_4=i_5\ne 0\}}
{\bf 1}_{\{j_4=j_5\}}-
$$

\vspace{-2mm}
$$
-
{\bf 1}_{\{i_6=i_2\ne 0\}}
{\bf 1}_{\{j_6=j_2\}}
{\bf 1}_{\{i_1=i_5\ne 0\}}
{\bf 1}_{\{j_1=j_5\}}
{\bf 1}_{\{i_3=i_4\ne 0\}}
{\bf 1}_{\{j_3=j_4\}}-
$$

\vspace{-2mm}
$$
-
{\bf 1}_{\{i_6=i_2\ne 0\}}
{\bf 1}_{\{j_6=j_2\}}
{\bf 1}_{\{i_1=i_4\ne 0\}}
{\bf 1}_{\{j_1=j_4\}}
{\bf 1}_{\{i_3=i_5\ne 0\}}
{\bf 1}_{\{j_3=j_5\}}-
$$

\vspace{-2mm}
$$
-
{\bf 1}_{\{i_6=i_2\ne 0\}}
{\bf 1}_{\{j_6=j_2\}}
{\bf 1}_{\{i_1=i_3\ne 0\}}
{\bf 1}_{\{j_1=j_3\}}
{\bf 1}_{\{i_4=i_5\ne 0\}}
{\bf 1}_{\{j_4=j_5\}}-
$$

\vspace{-2mm}
$$
-
{\bf 1}_{\{i_6=i_3\ne 0\}}
{\bf 1}_{\{j_6=j_3\}}
{\bf 1}_{\{i_1=i_5\ne 0\}}
{\bf 1}_{\{j_1=j_5\}}
{\bf 1}_{\{i_2=i_4\ne 0\}}
{\bf 1}_{\{j_2=j_4\}}-
$$

\vspace{-2mm}
$$
-
{\bf 1}_{\{i_6=i_3\ne 0\}}
{\bf 1}_{\{j_6=j_3\}}
{\bf 1}_{\{i_1=i_4\ne 0\}}
{\bf 1}_{\{j_1=j_4\}}
{\bf 1}_{\{i_2=i_5\ne 0\}}
{\bf 1}_{\{j_2=j_5\}}-
$$

\vspace{-2mm}
$$
-
{\bf 1}_{\{i_3=i_6\ne 0\}}
{\bf 1}_{\{j_3=j_6\}}
{\bf 1}_{\{i_1=i_2\ne 0\}}
{\bf 1}_{\{j_1=j_2\}}
{\bf 1}_{\{i_4=i_5\ne 0\}}
{\bf 1}_{\{j_4=j_5\}}-
$$

\vspace{-2mm}
$$
-
{\bf 1}_{\{i_6=i_4\ne 0\}}
{\bf 1}_{\{j_6=j_4\}}
{\bf 1}_{\{i_1=i_5\ne 0\}}
{\bf 1}_{\{j_1=j_5\}}
{\bf 1}_{\{i_2=i_3\ne 0\}}
{\bf 1}_{\{j_2=j_3\}}-
$$

\vspace{-2mm}
$$
-
{\bf 1}_{\{i_6=i_4\ne 0\}}
{\bf 1}_{\{j_6=j_4\}}
{\bf 1}_{\{i_1=i_3\ne 0\}}
{\bf 1}_{\{j_1=j_3\}}
{\bf 1}_{\{i_2=i_5\ne 0\}}
{\bf 1}_{\{j_2=j_5\}}-
$$

\vspace{-2mm}
$$
-
{\bf 1}_{\{i_6=i_4\ne 0\}}
{\bf 1}_{\{j_6=j_4\}}
{\bf 1}_{\{i_1=i_2\ne 0\}}
{\bf 1}_{\{j_1=j_2\}}
{\bf 1}_{\{i_3=i_5\ne 0\}}
{\bf 1}_{\{j_3=j_5\}}-
$$

\vspace{-2mm}
$$
-
{\bf 1}_{\{i_6=i_5\ne 0\}}
{\bf 1}_{\{j_6=j_5\}}
{\bf 1}_{\{i_1=i_4\ne 0\}}
{\bf 1}_{\{j_1=j_4\}}
{\bf 1}_{\{i_2=i_3\ne 0\}}
{\bf 1}_{\{j_2=j_3\}}-
$$

\vspace{-2mm}
$$
-
{\bf 1}_{\{i_6=i_5\ne 0\}}
{\bf 1}_{\{j_6=j_5\}}
{\bf 1}_{\{i_1=i_2\ne 0\}}
{\bf 1}_{\{j_1=j_2\}}
{\bf 1}_{\{i_3=i_4\ne 0\}}
{\bf 1}_{\{j_3=j_4\}}-
$$
\begin{equation}
\label{a6}
\Biggl.-
{\bf 1}_{\{i_6=i_5\ne 0\}}
{\bf 1}_{\{j_6=j_5\}}
{\bf 1}_{\{i_1=i_3\ne 0\}}
{\bf 1}_{\{j_1=j_3\}}
{\bf 1}_{\{i_2=i_4\ne 0\}}
{\bf 1}_{\{j_2=j_4\}}\Biggr),
\end{equation}

\vspace{12mm}
$$
J[\psi^{(7)}]_{T,t}
=
\hbox{\vtop{\offinterlineskip\halign{
\hfil#\hfil\cr
{\rm l.i.m.}\cr
$\stackrel{}{{}_{p_1,\ldots, p_7\to \infty}}$\cr
}} }
\sum_{j_1=0}^{p_1}\ldots\sum_{j_7=0}^{p_7}
C_{j_7\ldots j_1}\Biggl(
\prod_{l=1}^7
\zeta_{j_l}^{(i_l)}
-\Biggr.
$$

\vspace{2mm}
$$
-
{\bf 1}_{\{i_1=i_6\ne 0,j_1=j_6\}}
\prod_{\stackrel{l=1}{{}_{l\ne 1, 6}}}^7\zeta_{j_l}^{(i_l)}
-
{\bf 1}_{\{i_2=i_6\ne 0,j_2=j_6\}}
\prod_{\stackrel{l=1}{{}_{l\ne 2, 6}}}^7\zeta_{j_l}^{(i_l)}
-
{\bf 1}_{\{i_3=i_6\ne 0,j_3=j_6\}}
\prod_{\stackrel{l=1}{{}_{l\ne 3, 6}}}^7\zeta_{j_l}^{(i_l)}-
$$

$$
-
{\bf 1}_{\{i_4=i_6\ne 0,j_4=j_6\}}
\prod_{\stackrel{l=1}{{}_{l\ne 4, 6}}}^7\zeta_{j_l}^{(i_l)}
-
{\bf 1}_{\{i_5=i_6\ne 0,j_5=j_6\}}
\prod_{\stackrel{l=1}{{}_{l\ne 5, 6}}}^7\zeta_{j_l}^{(i_l)}
-
{\bf 1}_{\{i_1=i_2\ne 0,j_1=j_2\}}
\prod_{\stackrel{l=1}{{}_{l\ne 1, 2}}}^7\zeta_{j_l}^{(i_l)}-
$$

$$
-
{\bf 1}_{\{i_1=i_3\ne 0,j_1=j_3\}}
\prod_{\stackrel{l=1}{{}_{l\ne 1, 3}}}^7\zeta_{j_l}^{(i_l)}
-
{\bf 1}_{\{i_1=i_4\ne 0,j_1=j_4\}}
\prod_{\stackrel{l=1}{{}_{l\ne 1, 4}}}^7\zeta_{j_l}^{(i_l)}
-
{\bf 1}_{\{i_1=i_5\ne 0,j_1=j_5\}}
\prod_{\stackrel{l=1}{{}_{l\ne 1, 5}}}^7\zeta_{j_l}^{(i_l)}-
$$

$$
-
{\bf 1}_{\{i_2=i_3\ne 0,j_2=j_3\}}
\prod_{\stackrel{l=1}{{}_{l\ne 2, 3}}}^7\zeta_{j_l}^{(i_l)}
-
{\bf 1}_{\{i_2=i_4\ne 0,j_2=j_4\}}
\prod_{\stackrel{l=1}{{}_{l\ne 2, 4}}}^7\zeta_{j_l}^{(i_l)}
-
{\bf 1}_{\{i_2=i_5\ne 0,j_2=j_5\}}
\prod_{\stackrel{l=1}{{}_{l\ne 2, 5}}}^7\zeta_{j_l}^{(i_l)}-
$$

$$
-
{\bf 1}_{\{i_3=i_4\ne 0,j_3=j_4\}}
\prod_{\stackrel{l=1}{{}_{l\ne 3, 4}}}^7\zeta_{j_l}^{(i_l)}
-
{\bf 1}_{\{i_3=i_5\ne 0,j_3=j_5\}}
\prod_{\stackrel{l=1}{{}_{l\ne 3, 5}}}^7\zeta_{j_l}^{(i_l)}
-
{\bf 1}_{\{i_4=i_5\ne 0,j_4=j_5\}}
\prod_{\stackrel{l=1}{{}_{l\ne 4, 5}}}^7\zeta_{j_l}^{(i_l)}-
$$

$$
-
{\bf 1}_{\{i_7=i_1\ne 0,j_7=j_1\}}
\prod_{\stackrel{l=1}{{}_{l\ne 1, 7}}}^7\zeta_{j_l}^{(i_l)}
-
{\bf 1}_{\{i_7=i_2\ne 0,j_7=j_2\}}
\prod_{\stackrel{l=1}{{}_{l\ne 2, 7}}}^7\zeta_{j_l}^{(i_l)}
-
{\bf 1}_{\{i_7=i_3\ne 0,j_7=j_3\}}
\prod_{\stackrel{l=1}{{}_{l\ne 3, 7}}}^7\zeta_{j_l}^{(i_l)}-
$$

$$
-
{\bf 1}_{\{i_7=i_4\ne 0,j_7=j_4\}}
\prod_{\stackrel{l=1}{{}_{l\ne 4, 7}}}^7\zeta_{j_l}^{(i_l)}
-
{\bf 1}_{\{i_7=i_5\ne 0,j_7=j_5\}}
\prod_{\stackrel{l=1}{{}_{l\ne 7, 5}}}^7\zeta_{j_l}^{(i_l)}
-
{\bf 1}_{\{i_7=i_6\ne 0,j_7=j_6\}}
\prod_{\stackrel{l=1}{{}_{l\ne 7, 6}}}^7\zeta_{j_l}^{(i_l)}+
$$

\vspace{3mm}
$$
+
{\bf 1}_{\{i_1=i_2\ne 0,j_1=j_2,i_3=i_4\ne 0,j_3=j_4\}}
\prod_{l=5,6,7}\zeta_{j_l}^{(i_l)}
+
{\bf 1}_{\{i_1=i_2\ne 0,j_1=j_2,i_3=i_5\ne 0,j_3=j_5\}}
\prod_{l=4,6,7}\zeta_{j_l}^{(i_l)}+
$$

$$
+
{\bf 1}_{\{i_1=i_2\ne 0,j_1=j_2,i_4=i_5\ne 0,j_4=j_5\}}
\prod_{l=3,6,7}\zeta_{j_l}^{(i_l)}
+
{\bf 1}_{\{i_1=i_3\ne 0,j_1=j_3,i_2=i_4\ne 0,j_2=j_4\}}
\prod_{l=5,6,7}\zeta_{j_l}^{(i_l)}+
$$

$$
+
{\bf 1}_{\{i_1=i_3\ne 0,j_1=j_3,i_2=i_5\ne 0,j_2=j_5\}}
\prod_{l=4,6,7}\zeta_{j_l}^{(i_l)}
+
{\bf 1}_{\{i_1=i_3\ne 0,j_1=j_3,i_4=i_5\ne 0,j_4=j_5\}}
\prod_{l=2,6,7}\zeta_{j_l}^{(i_l)}+
$$

$$
+
{\bf 1}_{\{i_1=i_4\ne 0,j_1=j_4,i_2=i_3\ne 0,j_2=j_3\}}
\prod_{l=5,6,7}\zeta_{j_l}^{(i_l)}
+
{\bf 1}_{\{i_1=i_4\ne 0,j_1=j_4,i_2=i_5\ne 0,j_2=j_5\}}
\prod_{l=3,6,7}\zeta_{j_l}^{(i_l)}+
$$

$$
+
{\bf 1}_{\{i_1=i_4\ne 0,j_1=j_4,i_3=i_5\ne 0,j_3=j_5\}}
\prod_{l=2,6,7}\zeta_{j_l}^{(i_l)}
+
{\bf 1}_{\{i_1=i_5\ne 0,j_1=j_5,i_2=i_3\ne 0,j_2=j_3\}}
\prod_{l=4,6,7}\zeta_{j_l}^{(i_l)}+
$$

$$
+
{\bf 1}_{\{i_1=i_5\ne 0,j_1=j_5,i_2=i_4\ne 0,j_2=j_4\}}
\prod_{l=3,6,7}\zeta_{j_l}^{(i_l)}
+
{\bf 1}_{\{i_1=i_5\ne 0,j_1=j_5,i_3=i_4\ne 0,j_3=j_4\}}
\prod_{l=2,6,7}\zeta_{j_l}^{(i_l)}+
$$

$$
+
{\bf 1}_{\{i_2=i_3\ne 0,j_2=j_3,i_4=i_5\ne 0,j_4=j_5\}}
\prod_{l=1,6,7}\zeta_{j_l}^{(i_l)}
+
{\bf 1}_{\{i_2=i_4\ne 0,j_2=j_4,i_3=i_5\ne 0,j_3=j_5\}}
\prod_{l=1,6,7}\zeta_{j_l}^{(i_l)}+
$$

$$
+
{\bf 1}_{\{i_2=i_5\ne 0,j_2=j_5,i_3=i_4\ne 0,j_3=j_4\}}
\prod_{l=1,6,7}\zeta_{j_l}^{(i_l)}
+
{\bf 1}_{\{i_6=i_1\ne 0,j_6=j_1,i_3=i_4\ne 0,j_3=j_4\}}
\prod_{l=2,5,7}\zeta_{j_l}^{(i_l)}+
$$

$$
+
{\bf 1}_{\{i_6=i_1\ne 0,j_6=j_1,i_3=i_5\ne 0,j_3=j_5\}}
\prod_{l=2,4,7}\zeta_{j_l}^{(i_l)}
+
{\bf 1}_{\{i_6=i_1\ne 0,j_6=j_1,i_2=i_5\ne 0,j_2=j_5\}}
\prod_{l=3,4,7}\zeta_{j_l}^{(i_l)}+
$$

$$
+
{\bf 1}_{\{i_6=i_1\ne 0,j_6=j_1,i_2=i_4\ne 0,j_2=j_4\}}
\prod_{l=3,5,7}\zeta_{j_l}^{(i_l)}
+
{\bf 1}_{\{i_6=i_1\ne 0,j_6=j_1,i_4=i_5\ne 0,j_4=j_5\}}
\prod_{l=2,3,7}\zeta_{j_l}^{(i_l)}+
$$

$$
+
{\bf 1}_{\{i_6=i_1\ne 0,j_6=j_1,i_2=i_3\ne 0,j_2=j_3\}}
\prod_{l=4,5,7}\zeta_{j_l}^{(i_l)}
+
{\bf 1}_{\{i_6=i_2\ne 0,j_6=j_2,i_3=i_5\ne 0,j_3=j_5\}}
\prod_{l=1,4,7}\zeta_{j_l}^{(i_l)}+
$$

$$
+
{\bf 1}_{\{i_6=i_2\ne 0,j_6=j_2,i_4=i_5\ne 0,j_4=j_5\}}
\prod_{l=1,3,7}\zeta_{j_l}^{(i_l)}
+
{\bf 1}_{\{i_6=i_2\ne 0,j_6=j_2,i_3=i_4\ne 0,j_3=j_4\}}
\prod_{l=1,5,7}\zeta_{j_l}^{(i_l)}+
$$

$$
+
{\bf 1}_{\{i_6=i_2\ne 0,j_6=j_2,i_1=i_5\ne 0,j_1=j_5\}}
\prod_{l=3,4,7}\zeta_{j_l}^{(i_l)}
+
{\bf 1}_{\{i_6=i_2\ne 0,j_6=j_2,i_1=i_4\ne 0,j_1=j_4\}}
\prod_{l=3,5,7}\zeta_{j_l}^{(i_l)}+
$$

$$
+
{\bf 1}_{\{i_6=i_2\ne 0,j_6=j_2,i_1=i_3\ne 0,j_1=j_3\}}
\prod_{l=4,5,7}\zeta_{j_l}^{(i_l)}
+
{\bf 1}_{\{i_6=i_3\ne 0,j_6=j_3,i_2=i_5\ne 0,j_2=j_5\}}
\prod_{l=1,4,7}\zeta_{j_l}^{(i_l)}+
$$

$$
+
{\bf 1}_{\{i_6=i_3\ne 0,j_6=j_3,i_4=i_5\ne 0,j_4=j_5\}}
\prod_{l=1,2,7}\zeta_{j_l}^{(i_l)}
+
{\bf 1}_{\{i_6=i_3\ne 0,j_6=j_3,i_2=i_4\ne 0,j_2=j_4\}}
\prod_{l=1,5,7}\zeta_{j_l}^{(i_l)}+
$$

$$
+
{\bf 1}_{\{i_6=i_3\ne 0,j_6=j_3,i_1=i_5\ne 0,j_1=j_5\}}
\prod_{l=2,4,7}\zeta_{j_l}^{(i_l)}
+
{\bf 1}_{\{i_6=i_3\ne 0,j_6=j_3,i_1=i_4\ne 0,j_1=j_4\}}
\prod_{l=2,5,7}\zeta_{j_l}^{(i_l)}+
$$

$$
+
{\bf 1}_{\{i_6=i_3\ne 0,j_6=j_3,i_1=i_2\ne 0,j_1=j_2\}}
\prod_{l=4,5,7}\zeta_{j_l}^{(i_l)}
+
{\bf 1}_{\{i_6=i_4\ne 0,j_6=j_4,i_3=i_5\ne 0,j_3=j_5\}}
\prod_{l=1,2,7}\zeta_{j_l}^{(i_l)}+
$$

$$
+
{\bf 1}_{\{i_6=i_4\ne 0,j_6=j_4,i_2=i_5\ne 0,j_2=j_5\}}
\prod_{l=1,3,7}\zeta_{j_l}^{(i_l)}
+
{\bf 1}_{\{i_6=i_4\ne 0,j_6=j_4,i_2=i_3\ne 0,j_2=j_3\}}
\prod_{l=1,5,7}\zeta_{j_l}^{(i_l)}+
$$

$$
+
{\bf 1}_{\{i_6=i_4\ne 0,j_6=j_4,i_1=i_5\ne 0,j_1=j_5\}}
\prod_{l=2,3,7}\zeta_{j_l}^{(i_l)}
+
{\bf 1}_{\{i_6=i_4\ne 0,j_6=j_4,i_1=i_3\ne 0,j_1=j_3\}}
\prod_{l=2,5,7}\zeta_{j_l}^{(i_l)}+
$$

$$
+
{\bf 1}_{\{i_6=i_4\ne 0,j_6=j_4,i_1=i_2\ne 0,j_1=j_2\}}
\prod_{l=3,5,7}\zeta_{j_l}^{(i_l)}
+
{\bf 1}_{\{i_6=i_5\ne 0,j_6=j_5,i_3=i_4\ne 0,j_3=j_4\}}
\prod_{l=1,2,7}\zeta_{j_l}^{(i_l)}+
$$

$$
+
{\bf 1}_{\{i_6=i_5\ne 0,j_6=j_5,i_2=i_4\ne 0,j_2=j_4\}}
\prod_{l=1,3,7}\zeta_{j_l}^{(i_l)}
+
{\bf 1}_{\{i_6=i_5\ne 0,j_6=j_5,i_2=i_3\ne 0,j_2=j_3\}}
\prod_{l=1,4,7}\zeta_{j_l}^{(i_l)}+
$$

$$
+
{\bf 1}_{\{i_6=i_5\ne 0,j_6=j_5,i_1=i_4\ne 0,j_1=j_4\}}
\prod_{l=2,3,7}\zeta_{j_l}^{(i_l)}
+
{\bf 1}_{\{i_6=i_5\ne 0,j_6=j_5,i_1=i_3\ne 0,j_1=j_3\}}
\prod_{l=2,4,7}\zeta_{j_l}^{(i_l)}+
$$

$$
+
{\bf 1}_{\{i_6=i_5\ne 0,j_6=j_5,i_1=i_2\ne 0,j_1=j_2\}}
\prod_{l=3,4,7}\zeta_{j_l}^{(i_l)}
+
{\bf 1}_{\{i_7=i_1\ne 0,j_7=j_1,i_2=i_3\ne 0,j_2=j_3\}}
\prod_{l=4,5,6}\zeta_{j_l}^{(i_l)}+
$$

$$
+
{\bf 1}_{\{i_7=i_1\ne 0,j_7=j_1,i_2=i_4\ne 0,j_2=j_4\}}
\prod_{l=3,5,6}\zeta_{j_l}^{(i_l)}
+
{\bf 1}_{\{i_7=i_1\ne 0,j_7=j_1,i_2=i_5\ne 0,j_2=j_5\}}
\prod_{l=3,4,6}\zeta_{j_l}^{(i_l)}+
$$

$$
+
{\bf 1}_{\{i_7=i_1\ne 0,j_7=j_1,i_2=i_6\ne 0,j_2=j_6\}}
\prod_{l=3,4,5}\zeta_{j_l}^{(i_l)}
+
{\bf 1}_{\{i_7=i_1\ne 0,j_7=j_1,i_3=i_4\ne 0,j_3=j_4\}}
\prod_{l=2,5,6}\zeta_{j_l}^{(i_l)}+
$$

$$
+
{\bf 1}_{\{i_7=i_1\ne 0,j_7=j_1,i_3=i_5\ne 0,j_3=j_5\}}
\prod_{l=2,4,6}\zeta_{j_l}^{(i_l)}
+
{\bf 1}_{\{i_7=i_1\ne 0,j_7=j_1,i_3=i_6\ne 0,j_3=j_6\}}
\prod_{l=2,4,5}\zeta_{j_l}^{(i_l)}+
$$

$$
+
{\bf 1}_{\{i_7=i_1\ne 0,j_7=j_1,i_4=i_5\ne 0,j_4=j_5\}}
\prod_{l=2,3,6}\zeta_{j_l}^{(i_l)}
+
{\bf 1}_{\{i_7=i_1\ne 0,j_7=j_1,i_4=i_6\ne 0,j_4=j_6\}}
\prod_{l=2,3,5}\zeta_{j_l}^{(i_l)}+
$$

$$
+
{\bf 1}_{\{i_1=i_2\ne 0,j_7=j_1,i_7=i_1\ne 0,j_5=j_6\}}
\prod_{l=2,3,4}\zeta_{j_l}^{(i_l)}
+
{\bf 1}_{\{i_7=i_2\ne 0,j_7=j_2,i_1=i_3\ne 0,j_1=j_3\}}
\prod_{l=4,5,6}\zeta_{j_l}^{(i_l)}+
$$

$$
+
{\bf 1}_{\{i_7=i_2\ne 0,j_7=j_2,i_1=i_4\ne 0,j_1=j_4\}}
\prod_{l=3,5,6}\zeta_{j_l}^{(i_l)}
+
{\bf 1}_{\{i_7=i_2\ne 0,j_7=j_2,i_1=i_5\ne 0,j_1=j_5\}}
\prod_{l=3,4,6}\zeta_{j_l}^{(i_l)}+
$$

$$
+
{\bf 1}_{\{i_7=i_2\ne 0,j_7=j_2,i_1=i_6\ne 0,j_1=j_6\}}
\prod_{l=3,4,5}\zeta_{j_l}^{(i_l)}
+
{\bf 1}_{\{i_7=i_2\ne 0,j_7=j_2,i_3=i_4\ne 0,j_3=j_4\}}
\prod_{l=1,5,6}\zeta_{j_l}^{(i_l)}+
$$

$$
+
{\bf 1}_{\{i_7=i_2\ne 0,j_7=j_2,i_3=i_5\ne 0,j_3=j_5\}}
\prod_{l=1,4,6}\zeta_{j_l}^{(i_l)}
+
{\bf 1}_{\{i_7=i_2\ne 0,j_7=j_2,i_3=i_6\ne 0,j_3=j_6\}}
\prod_{l=1,4,5}\zeta_{j_l}^{(i_l)}+
$$

$$
+
{\bf 1}_{\{i_7=i_2\ne 0,j_7=j_2,i_4=i_5\ne 0,j_4=j_5\}}
\prod_{l=1,3,6}\zeta_{j_l}^{(i_l)}
+
{\bf 1}_{\{i_7=i_2\ne 0,j_7=j_2,i_4=i_6\ne 0,j_4=j_6\}}
\prod_{l=1,3,5}\zeta_{j_l}^{(i_l)}+
$$

$$
+
{\bf 1}_{\{i_7=i_2\ne 0,j_7=j_2,i_5=i_6\ne 0,j_5=j_6\}}
\prod_{l=1,3,4}\zeta_{j_l}^{(i_l)}
+
{\bf 1}_{\{i_7=i_3\ne 0,j_7=j_3,i_1=i_2\ne 0,j_1=j_2\}}
\prod_{l=4,5,6}\zeta_{j_l}^{(i_l)}+
$$

$$
+
{\bf 1}_{\{i_7=i_3\ne 0,j_7=j_3,i_1=i_4\ne 0,j_1=j_4\}}
\prod_{l=2,3,5}\zeta_{j_l}^{(i_l)}
+
{\bf 1}_{\{i_7=i_3\ne 0,j_7=j_3,i_1=i_5\ne 0,j_1=j_5\}}
\prod_{l=2,4,6}\zeta_{j_l}^{(i_l)}+
$$

$$
+
{\bf 1}_{\{i_7=i_3\ne 0,j_7=j_3,i_1=i_6\ne 0,j_1=j_6\}}
\prod_{l=4,2,5}\zeta_{j_l}^{(i_l)}
+
{\bf 1}_{\{i_7=i_3\ne 0,j_7=j_3,i_2=i_4\ne 0,j_2=j_4\}}
\prod_{l=3,5,6}\zeta_{j_l}^{(i_l)}+
$$

$$
+
{\bf 1}_{\{i_7=i_3\ne 0,j_7=j_3,i_2=i_5\ne 0,j_2=j_5\}}
\prod_{l=1,4,6}\zeta_{j_l}^{(i_l)}
+
{\bf 1}_{\{i_7=i_3\ne 0,j_7=j_3,i_2=i_6\ne 0,j_2=j_6\}}
\prod_{l=1,4,5}\zeta_{j_l}^{(i_l)}+
$$

$$
+
{\bf 1}_{\{i_7=i_3\ne 0,j_7=j_3,i_4=i_5\ne 0,j_4=j_5\}}
\prod_{l=1,2,6}\zeta_{j_l}^{(i_l)}
+
{\bf 1}_{\{i_7=i_3\ne 0,j_7=j_3,i_4=i_6\ne 0,j_4=j_6\}}
\prod_{l=1,2,5}\zeta_{j_l}^{(i_l)}+
$$

$$
+
{\bf 1}_{\{i_7=i_3\ne 0,j_7=j_3,i_5=i_6\ne 0,j_5=j_6\}}
\prod_{l=1,2,4}\zeta_{j_l}^{(i_l)}
+
{\bf 1}_{\{i_7=i_4\ne 0,j_7=j_4,i_1=i_2\ne 0,j_1=j_2\}}
\prod_{l=3,5,6}\zeta_{j_l}^{(i_l)}+
$$

$$
+
{\bf 1}_{\{i_7=i_4\ne 0,j_7=j_4,i_1=i_3\ne 0,j_1=j_3\}}
\prod_{l=2,5,6}\zeta_{j_l}^{(i_l)}
+
{\bf 1}_{\{i_7=i_4\ne 0,j_7=j_4,i_1=i_5\ne 0,j_1=j_5\}}
\prod_{l=2,3,6}\zeta_{j_l}^{(i_l)}+
$$

$$
+
{\bf 1}_{\{i_7=i_4\ne 0,j_7=j_4,i_1=i_6\ne 0,j_1=j_6\}}
\prod_{l=2,3,5}\zeta_{j_l}^{(i_l)}
+
{\bf 1}_{\{i_7=i_4\ne 0,j_7=j_4,i_2=i_3\ne 0,j_2=j_3\}}
\prod_{l=1,5,6}\zeta_{j_l}^{(i_l)}+
$$

$$
+
{\bf 1}_{\{i_7=i_4\ne 0,j_7=j_4,i_2=i_5\ne 0,j_2=j_5\}}
\prod_{l=1,3,6}\zeta_{j_l}^{(i_l)}
+
{\bf 1}_{\{i_7=i_4\ne 0,j_7=j_4,i_2=i_6\ne 0,j_2=j_6\}}
\prod_{l=1,3,5}\zeta_{j_l}^{(i_l)}+
$$

$$
+
{\bf 1}_{\{i_7=i_4\ne 0,j_7=j_4,i_3=i_5\ne 0,j_3=j_5\}}
\prod_{l=1,2,6}\zeta_{j_l}^{(i_l)}
+
{\bf 1}_{\{i_7=i_4\ne 0,j_7=j_4,i_3=i_6\ne 0,j_3=j_6\}}
\prod_{l=1,2,5}\zeta_{j_l}^{(i_l)}+
$$

$$
+
{\bf 1}_{\{i_7=i_4\ne 0,j_7=j_4,i_5=i_6\ne 0,j_5=j_6\}}
\prod_{l=1,2,3}\zeta_{j_l}^{(i_l)}
+
{\bf 1}_{\{i_7=i_5\ne 0,j_7=j_5,i_1=i_2\ne 0,j_1=j_2\}}
\prod_{l=3,4,6}\zeta_{j_l}^{(i_l)}+
$$

$$
+
{\bf 1}_{\{i_7=i_5\ne 0,j_7=j_5,i_1=i_3\ne 0,j_1=j_3\}}
\prod_{l=2,4,6}\zeta_{j_l}^{(i_l)}
+
{\bf 1}_{\{i_7=i_5\ne 0,j_7=j_5,i_1=i_4\ne 0,j_1=j_4\}}
\prod_{l=2,3,6}\zeta_{j_l}^{(i_l)}+
$$

$$
+
{\bf 1}_{\{i_7=i_5\ne 0,j_7=j_5,i_1=i_6\ne 0,j_1=j_6\}}
\prod_{l=2,3,4}\zeta_{j_l}^{(i_l)}
+
{\bf 1}_{\{i_7=i_5\ne 0,j_7=j_5,i_2=i_3\ne 0,j_2=j_3\}}
\prod_{l=1,4,6}\zeta_{j_l}^{(i_l)}+
$$

$$
+
{\bf 1}_{\{i_7=i_5\ne 0,j_7=j_5,i_2=i_4\ne 0,j_2=j_4\}}
\prod_{l=1,3,6}\zeta_{j_l}^{(i_l)}
+
{\bf 1}_{\{i_7=i_5\ne 0,j_7=j_5,i_2=i_6\ne 0,j_2=j_6\}}
\prod_{l=1,3,5}\zeta_{j_l}^{(i_l)}+
$$

$$
+
{\bf 1}_{\{i_7=i_5\ne 0,j_7=j_5,i_3=i_4\ne 0,j_3=j_4\}}
\prod_{l=1,2,6}\zeta_{j_l}^{(i_l)}
+
{\bf 1}_{\{i_7=i_5\ne 0,j_7=j_5,i_3=i_6\ne 0,j_3=j_6\}}
\prod_{l=1,2,4}\zeta_{j_l}^{(i_l)}+
$$

$$
+
{\bf 1}_{\{i_7=i_5\ne 0,j_7=j_5,i_4=i_6\ne 0,j_4=j_6\}}
\prod_{l=1,2,3}\zeta_{j_l}^{(i_l)}
+
{\bf 1}_{\{i_7=i_6\ne 0,j_7=j_6,i_1=i_2\ne 0,j_1=j_2\}}
\prod_{l=3,4,5}\zeta_{j_l}^{(i_l)}+
$$

$$
+
{\bf 1}_{\{i_7=i_6\ne 0,j_7=j_6,i_1=i_3\ne 0,j_1=j_3\}}
\prod_{l=2,4,5}\zeta_{j_l}^{(i_l)}
+
{\bf 1}_{\{i_7=i_6\ne 0,j_7=j_6,i_1=i_4\ne 0,j_1=j_4\}}
\prod_{l=2,3,5}\zeta_{j_l}^{(i_l)}+
$$

$$
+
{\bf 1}_{\{i_7=i_6\ne 0,j_7=j_6,i_1=i_5\ne 0,j_1=j_5\}}
\prod_{l=2,3,4}\zeta_{j_l}^{(i_l)}
+
{\bf 1}_{\{i_7=i_6\ne 0,j_7=j_6,i_2=i_3\ne 0,j_2=j_3\}}
\prod_{l=1,4,5}\zeta_{j_l}^{(i_l)}+
$$

$$
+
{\bf 1}_{\{i_7=i_6\ne 0,j_7=j_6,i_2=i_4\ne 0,j_2=j_4\}}
\prod_{l=1,3,5}\zeta_{j_l}^{(i_l)}
+
{\bf 1}_{\{i_7=i_6\ne 0,j_7=j_6,i_2=i_5\ne 0,j_2=j_5\}}
\prod_{l=1,3,4}\zeta_{j_l}^{(i_l)}+
$$

$$
+
{\bf 1}_{\{i_7=i_6\ne 0,j_7=j_6,i_3=i_5\ne 0,j_3=j_5\}}
\prod_{l=1,2,4}\zeta_{j_l}^{(i_l)}
+
{\bf 1}_{\{i_7=i_6\ne 0,j_7=j_6,i_4=i_5\ne 0,j_4=j_5\}}
\prod_{l=1,2,3}\zeta_{j_l}^{(i_l)}+
$$

$$
+
{\bf 1}_{\{i_7=i_6\ne 0,j_7=j_6,i_3=i_4\ne 0,j_3=j_4\}}
\prod_{l=1,2,5}\zeta_{j_l}^{(i_l)} -
$$
$$
-\Biggl(  
{\bf 1}_{\{i_2=i_3\ne 0,j_2=j_3,i_4=i_5\ne 0,j_4=j_5,i_6=i_7\ne 0,j_6=j_7\}}
+\Biggr.
{\bf 1}_{\{i_2=i_3\ne 0,j_2=j_3,i_4=i_6\ne 0,j_4=j_6,i_5=i_7\ne 0,j_5=j_7\}}+
$$
$$
+
{\bf 1}_{\{i_2=i_3\ne 0,j_2=j_3,i_4=i_7\ne 0,j_4=j_7,i_5=i_6\ne 0,j_5=j_6\}}
+
{\bf 1}_{\{i_2=i_4\ne 0,j_2=j_4,i_3=i_5\ne 0,j_3=j_5,i_6=i_7\ne 0,j_6=j_7\}}+
$$

\vspace{-4mm}
$$
+
{\bf 1}_{\{i_2=i_4\ne 0,j_2=j_4,i_3=i_6\ne 0,j_3=j_6,i_5=i_7\ne 0,j_5=j_7\}}
+
{\bf 1}_{\{i_2=i_4\ne 0,j_2=j_4,i_3=i_7\ne 0,j_3=j_7,i_5=i_6\ne 0,j_5=j_6\}}+
$$

\vspace{-4mm}
$$
+
{\bf 1}_{\{i_2=i_5\ne 0,j_2=j_5,i_3=i_4\ne 0,j_3=j_4,i_6=i_7\ne 0,j_6=j_7\}}
+
{\bf 1}_{\{i_2=i_5\ne 0,j_2=j_5,i_3=i_6\ne 0,j_3=j_6,i_4=i_7\ne 0,j_4=j_7\}}+
$$

\vspace{-4mm}
$$
+
{\bf 1}_{\{i_2=i_5\ne 0,j_2=j_5,i_3=i_7\ne 0,j_3=j_7,i_4=i_6\ne 0,j_4=j_6\}}
+
{\bf 1}_{\{i_2=i_6\ne 0,j_2=j_6,i_3=i_4\ne 0,j_3=j_4,i_5=i_7\ne 0,j_5=j_7\}}+
$$

\vspace{-4mm}
$$
+
{\bf 1}_{\{i_2=i_6\ne 0,j_2=j_6,i_3=i_5\ne 0,j_3=j_5,i_4=i_7\ne 0,j_4=j_7\}}
+
{\bf 1}_{\{i_2=i_6\ne 0,j_2=j_6,i_3=i_7\ne 0,j_3=j_7,i_4=i_5\ne 0,j_4=j_5\}}+
$$

\vspace{-4mm}
$$
+
{\bf 1}_{\{i_2=i_7\ne 0,j_2=j_7,i_3=i_4\ne 0,j_3=j_4,i_5=i_6\ne 0,j_5=j_6\}}
+
{\bf 1}_{\{i_2=i_7\ne 0,j_2=j_7,i_3=i_5\ne 0,j_3=j_5,i_4=i_6\ne 0,j_4=j_6\}}+
$$
$$
\Biggl.
+{\bf 1}_{\{i_2=i_7\ne 0,j_2=j_7,i_3=i_6\ne 0,j_3=j_6,i_4=i_5\ne 0,j_4=j_5\}}
\Biggr)\zeta_{j_1}^{(i_1)} -
$$
$$
-\Biggl(
{\bf 1}_{\{i_1=i_3\ne 0,j_1=j_3,i_4=i_7\ne 0,j_4=j_7,i_5=i_6\ne 0,j_5=j_6\}}
+
{\bf 1}_{\{i_1=i_3\ne 0,j_1=j_3,i_4=i_5\ne 0,j_4=j_5,i_6=i_7\ne 0,j_6=j_7\}}+
\Biggr.
$$
$$
+
{\bf 1}_{\{i_1=i_3\ne 0,j_1=j_3,i_4=i_6\ne 0,j_4=j_6,i_5=i_7\ne 0,j_5=j_7\}}
+
{\bf 1}_{\{i_1=i_4\ne 0,j_1=j_4,i_3=i_5\ne 0,j_3=j_5,i_6=i_7\ne 0,j_6=j_7\}}+
$$

\vspace{-4mm}
$$
+
{\bf 1}_{\{i_1=i_4\ne 0,j_1=j_4,i_3=i_6\ne 0,j_3=j_6,i_5=i_7\ne 0,j_5=j_7\}}
+
{\bf 1}_{\{i_1=i_4\ne 0,j_1=j_4,i_3=i_7\ne 0,j_3=j_7,i_5=i_6\ne 0,j_5=j_6\}}+
$$

\vspace{-4mm}
$$
+
{\bf 1}_{\{i_1=i_5\ne 0,j_1=j_5,i_3=i_4\ne 0,j_3=j_4,i_6=i_7\ne 0,j_6=j_7\}}
+
{\bf 1}_{\{i_1=i_5\ne 0,j_1=j_5,i_3=i_6\ne 0,j_3=j_6,i_4=i_7\ne 0,j_4=j_7\}}+
$$

\vspace{-4mm}
$$
+
{\bf 1}_{\{i_1=i_5\ne 0,j_1=j_5,i_3=i_7\ne 0,j_3=j_7,i_4=i_6\ne 0,j_4=j_6\}}
+
{\bf 1}_{\{i_1=i_6\ne 0,j_1=j_6,i_3=i_4\ne 0,j_3=j_4,i_5=i_7\ne 0,j_5=j_7\}}+
$$

\vspace{-4mm}
$$
+
{\bf 1}_{\{i_6=i_1\ne 0,j_6=j_1,i_3=i_5\ne 0,j_3=j_5,i_4=i_7\ne 0,j_4=j_7\}}
+
{\bf 1}_{\{i_6=i_1\ne 0,j_6=j_1,i_3=i_7\ne 0,j_3=j_7,i_4=i_5\ne 0,j_4=j_5\}}+
$$

\vspace{-4mm}
$$
+
{\bf 1}_{\{i_1=i_7\ne 0,j_1=j_7,i_3=i_4\ne 0,j_3=j_4,i_5=i_6\ne 0,j_5=j_6\}}
+
{\bf 1}_{\{i_1=i_7\ne 0,j_1=j_7,i_3=i_5\ne 0,j_3=j_5,i_4=i_6\ne 0,j_4=j_6\}}+
$$
$$
\Biggl.
+{\bf 1}_{\{i_1=i_7\ne 0,j_1=j_7,i_3=i_6\ne 0,j_3=j_6,i_4=i_5\ne 0,j_4=j_5\}}
\Biggr)\zeta_{j_2}^{(i_2)} -
$$
$$
-\Biggl(
{\bf 1}_{\{i_1=i_2\ne 0,j_1=j_2,i_4=i_5\ne 0,j_4=j_5,i_6=i_7\ne 0,j_6=j_7\}}
+
{\bf 1}_{\{i_1=i_2\ne 0,j_1=j_2,i_4=i_6\ne 0,j_4=j_6,i_5=i_7\ne 0,j_5=j_7\}}+
\Biggr.
$$
$$
+
{\bf 1}_{\{i_1=i_2\ne 0,j_1=j_2,i_4=i_7\ne 0,j_4=j_7,i_5=i_6\ne 0,j_5=j_6\}}
+
{\bf 1}_{\{i_1=i_4\ne 0,j_1=j_4,i_2=i_5\ne 0,j_2=j_5,i_6=i_7\ne 0,j_6=j_7\}}+
$$

\vspace{-4mm}
$$
+
{\bf 1}_{\{i_1=i_4\ne 0,j_1=j_4,i_2=i_6\ne 0,j_2=j_6,i_5=i_7\ne 0,j_5=j_7\}}
+
{\bf 1}_{\{i_1=i_4\ne 0,j_1=j_4,i_2=i_7\ne 0,j_2=j_7,i_5=i_6\ne 0,j_5=j_6\}}+
$$

\vspace{-4mm}
$$
+
{\bf 1}_{\{i_1=i_5\ne 0,j_1=j_5,i_2=i_4\ne 0,j_2=j_4,i_6=i_7\ne 0,j_6=j_7\}}
+
{\bf 1}_{\{i_1=i_5\ne 0,j_1=j_5,i_2=i_6\ne 0,j_2=j_6,i_4=i_7\ne 0,j_4=j_7\}}+
$$

\vspace{-4mm}
$$
+
{\bf 1}_{\{i_1=i_5\ne 0,j_1=j_5,i_2=i_7\ne 0,j_2=j_7,i_4=i_6\ne 0,j_4=j_6\}}
+
{\bf 1}_{\{i_6=i_1\ne 0,j_6=j_1,i_2=i_4\ne 0,j_2=j_4,i_5=i_7\ne 0,j_5=j_7\}}+
$$

\vspace{-4mm}
$$
+
{\bf 1}_{\{i_6=i_1\ne 0,j_6=j_1,i_2=i_5\ne 0,j_2=j_5,i_4=i_7\ne 0,j_4=j_7\}}
+
{\bf 1}_{\{i_6=i_1\ne 0,j_6=j_1,i_2=i_7\ne 0,j_2=j_7,i_4=i_5\ne 0,j_4=j_5\}}+
$$

\vspace{-4mm}
$$
+
{\bf 1}_{\{i_1=i_7\ne 0,j_1=j_7,i_2=i_4\ne 0,j_2=j_4,i_5=i_6\ne 0,j_5=j_6\}}
+
{\bf 1}_{\{i_1=i_7\ne 0,j_1=j_7,i_2=i_5\ne 0,j_2=j_5,i_4=i_6\ne 0,j_4=j_6\}}+
$$
$$
\Biggl.+
{\bf 1}_{\{i_1=i_7\ne 0,j_1=j_7,i_2=i_6\ne 0,j_2=j_6,i_4=i_5\ne 0,j_4=j_5\}}
\Biggr)\zeta_{j_3}^{(i_3)} -
$$
$$
-\Biggl(
{\bf 1}_{\{i_1=i_2\ne 0,j_1=j_2,i_3=i_5\ne 0,j_3=j_5,i_6=i_7\ne 0,j_6=j_7\}}
+\Biggr.
{\bf 1}_{\{i_1=i_2\ne 0,j_1=j_2,i_3=i_6\ne 0,j_3=j_6,i_5=i_7\ne 0,j_5=j_7\}}+
$$
$$
+
{\bf 1}_{\{i_1=i_2\ne 0,j_1=j_2,i_3=i_7\ne 0,j_3=j_7,i_5=i_6\ne 0,j_5=j_6\}}
+
{\bf 1}_{\{i_1=i_3\ne 0,j_1=j_3,i_2=i_5\ne 0,j_2=j_5,i_6=i_7\ne 0,j_6=j_7\}}+
$$

\vspace{-4mm}
$$
+
{\bf 1}_{\{i_1=i_3\ne 0,j_1=j_3,i_2=i_6\ne 0,j_2=j_6,i_5=i_7\ne 0,j_5=j_7\}}
+
{\bf 1}_{\{i_1=i_3\ne 0,j_1=j_3,i_2=i_7\ne 0,j_2=j_7,i_5=i_6\ne 0,j_5=j_6\}}+
$$

\vspace{-4mm}
$$
+
{\bf 1}_{\{i_1=i_5\ne 0,j_1=j_5,i_2=i_3\ne 0,j_2=j_3,i_6=i_7\ne 0,j_6=j_7\}}
+
{\bf 1}_{\{i_1=i_5\ne 0,j_1=j_5,i_2=i_6\ne 0,j_2=j_6,i_3=i_7\ne 0,j_3=j_7\}}+
$$

\vspace{-4mm}
$$
+
{\bf 1}_{\{i_1=i_5\ne 0,j_1=j_5,i_2=i_7\ne 0,j_2=j_7,i_3=i_6\ne 0,j_3=j_6\}}
+
{\bf 1}_{\{i_6=i_1\ne 0,j_6=j_1,i_2=i_3\ne 0,j_2=j_3,i_5=i_7\ne 0,j_5=j_7\}}+
$$

\vspace{-4mm}
$$
+
{\bf 1}_{\{i_6=i_1\ne 0,j_6=j_1,i_2=i_5\ne 0,j_2=j_5,i_3=i_7\ne 0,j_3=j_7\}}
+
{\bf 1}_{\{i_6=i_1\ne 0,j_6=j_1,i_2=i_7\ne 0,j_2=j_7,i_3=i_5\ne 0,j_3=j_5\}}+
$$

\vspace{-4mm}
$$
+
{\bf 1}_{\{i_7=i_1\ne 0,j_7=j_1,i_2=i_3\ne 0,j_2=j_3,i_5=i_6\ne 0,j_5=j_6\}}
+
{\bf 1}_{\{i_7=i_1\ne 0,j_7=j_1,i_2=i_5\ne 0,j_2=j_5,i_3=i_6\ne 0,j_3=j_6\}}+
$$
$$
\Biggl.
+{\bf 1}_{\{i_7=i_1\ne 0,j_7=j_1,i_2=i_6\ne 0,j_2=j_6,i_3=i_5\ne 0,j_3=j_5\}}
\Biggr)\zeta_{j_4}^{(i_4)} -
$$
$$
-\Biggl(
{\bf 1}_{\{i_1=i_2\ne 0,j_1=j_2,i_3=i_4\ne 0,j_3=j_4,i_6=i_7\ne 0,j_6=j_7\}}
+\Biggr.
{\bf 1}_{\{i_1=i_2\ne 0,j_1=j_2,i_3=i_6\ne 0,j_3=j_6,i_4=i_7\ne 0,j_4=j_7\}}+
$$
$$
+
{\bf 1}_{\{i_1=i_2\ne 0,j_1=j_2,i_3=i_7\ne 0,j_3=j_7,i_4=i_6\ne 0,j_4=j_6\}}
+
{\bf 1}_{\{i_1=i_3\ne 0,j_1=j_3,i_2=i_4\ne 0,j_2=j_4,i_6=i_7\ne 0,j_6=j_7\}}+
$$

\vspace{-4mm}
$$
+
{\bf 1}_{\{i_1=i_3\ne 0,j_1=j_3,i_2=i_6\ne 0,j_2=j_6,i_4=i_7\ne 0,j_4=j_7\}}
+
{\bf 1}_{\{i_1=i_3\ne 0,j_1=j_3,i_2=i_7\ne 0,j_2=j_7,i_4=i_6\ne 0,j_4=j_6\}}+
$$

\vspace{-4mm}
$$
+
{\bf 1}_{\{i_1=i_4\ne 0,j_1=j_4,i_2=i_3\ne 0,j_2=j_3,i_6=i_7\ne 0,j_6=j_7\}}
+
{\bf 1}_{\{i_1=i_4\ne 0,j_1=j_4,i_2=i_6\ne 0,j_2=j_6,i_3=i_7\ne 0,j_3=j_7\}}+
$$

\vspace{-4mm}
$$
+
{\bf 1}_{\{i_1=i_4\ne 0,j_1=j_4,i_2=i_7\ne 0,j_2=j_7,i_3=i_6\ne 0,j_3=j_6\}}
+
{\bf 1}_{\{i_6=i_1\ne 0,j_6=j_1,i_2=i_3\ne 0,j_2=j_3,i_4=i_7\ne 0,j_4=j_7\}}+
$$

\vspace{-4mm}
$$
+
{\bf 1}_{\{i_6=i_1\ne 0,j_6=j_1,i_2=i_4\ne 0,j_2=j_4,i_3=i_7\ne 0,j_3=j_7\}}
+
{\bf 1}_{\{i_6=i_1\ne 0,j_6=j_1,i_2=i_7\ne 0,j_2=j_7,i_3=i_4\ne 0,j_3=j_4\}}+
$$

\vspace{-4mm}
$$
+
{\bf 1}_{\{i_1=i_7\ne 0,j_1=j_7,i_2=i_3\ne 0,j_2=j_3,i_4=i_6\ne 0,j_4=j_6\}}
+
{\bf 1}_{\{i_1=i_7\ne 0,j_1=j_7,i_2=i_4\ne 0,j_2=j_4,i_3=i_6\ne 0,j_3=j_6\}}+
$$
$$
\Biggl.+
{\bf 1}_{\{i_7=i_1\ne 0,j_7=j_1,i_2=i_6\ne 0,j_2=j_6,i_3=i_4\ne 0,j_3=j_4\}}
\Biggr)\zeta_{j_5}^{(i_5)} -
$$
$$
-\Biggl(
{\bf 1}_{\{i_1=i_2\ne 0,j_1=j_2,i_3=i_4\ne 0,j_3=j_4,i_5=i_7\ne 0,j_5=j_7\}}
+\Biggr.
{\bf 1}_{\{i_1=i_2\ne 0,j_1=j_2,i_3=i_5\ne 0,j_3=j_5,i_4=i_7\ne 0,j_4=j_7\}}+
$$
$$
+
{\bf 1}_{\{i_1=i_2\ne 0,j_1=j_2,i_3=i_7\ne 0,j_3=j_7,i_4=i_5\ne 0,j_4=j_5\}}
+
{\bf 1}_{\{i_1=i_3\ne 0,j_1=j_3,i_2=i_4\ne 0,j_2=j_4,i_5=i_7\ne 0,j_5=j_7\}}+
$$

\vspace{-4mm}
$$
+
{\bf 1}_{\{i_1=i_3\ne 0,j_1=j_3,i_2=i_5\ne 0,j_2=j_5,i_4=i_7\ne 0,j_4=j_7\}}
+
{\bf 1}_{\{i_1=i_3\ne 0,j_1=j_3,i_2=i_7\ne 0,j_2=j_7,i_4=i_5\ne 0,j_4=j_5\}}+
$$

\vspace{-4mm}
$$
+
{\bf 1}_{\{i_1=i_4\ne 0,j_1=j_4,i_2=i_3\ne 0,j_2=j_3,i_5=i_7\ne 0,j_5=j_7\}}
+
{\bf 1}_{\{i_1=i_4\ne 0,j_1=j_4,i_2=i_5\ne 0,j_2=j_5,i_3=i_7\ne 0,j_3=j_7\}}+
$$

\vspace{-4mm}
$$
+
{\bf 1}_{\{i_1=i_4\ne 0,j_1=j_4,i_2=i_7\ne 0,j_2=j_7,i_3=i_5\ne 0,j_3=j_5\}}
+
{\bf 1}_{\{i_1=i_5\ne 0,j_1=j_5,i_2=i_3\ne 0,j_2=j_3,i_4=i_7\ne 0,j_4=j_7\}}+
$$

\vspace{-4mm}
$$
+
{\bf 1}_{\{i_1=i_5\ne 0,j_1=j_5,i_2=i_4\ne 0,j_2=j_4,i_3=i_7\ne 0,j_3=j_7\}}
+
{\bf 1}_{\{i_1=i_5\ne 0,j_1=j_5,i_2=i_7\ne 0,j_2=j_7,i_3=i_4\ne 0,j_3=j_4\}}+
$$

\vspace{-4mm}
$$
+
{\bf 1}_{\{i_7=i_1\ne 0,j_7=j_1,i_2=i_3\ne 0,j_2=j_3,i_4=i_5\ne 0,j_4=j_5\}}
+
{\bf 1}_{\{i_7=i_1\ne 0,j_7=j_1,i_2=i_4\ne 0,j_2=j_4,i_3=i_5\ne 0,j_3=j_5\}}+
$$
$$
\Biggl.+
{\bf 1}_{\{i_7=i_1\ne 0,j_7=j_1,i_2=i_5\ne 0,j_2=j_5,i_3=i_4\ne 0,j_3=j_4\}}
\Biggr)\zeta_{j_6}^{(i_6)} - 
$$
$$
-\Biggl(
{\bf 1}_{\{i_1=i_2\ne 0,j_1=j_2,i_3=i_4\ne 0,j_3=j_4,i_5=i_6\ne 0,j_5=j_6\}}
+\Biggr.
{\bf 1}_{\{i_1=i_2\ne 0,j_1=j_2,i_3=i_5\ne 0,j_3=j_5,i_4=i_6\ne 0,j_4=j_6\}}+
$$
$$
+
{\bf 1}_{\{i_1=i_2\ne 0,j_1=j_2,i_3=i_6\ne 0,j_3=j_6,i_4=i_5\ne 0,j_4=j_5\}}
+
{\bf 1}_{\{i_1=i_3\ne 0,j_1=j_3,i_2=i_4\ne 0,j_2=j_4,i_5=i_6\ne 0,j_5=j_6\}}+
$$

\vspace{-4mm}
$$
+
{\bf 1}_{\{i_1=i_3\ne 0,j_1=j_3,i_2=i_5\ne 0,j_2=j_5,i_4=i_6\ne 0,j_4=j_6\}}
+
{\bf 1}_{\{i_1=i_3\ne 0,j_1=j_3,i_2=i_6\ne 0,j_2=j_6,i_4=i_5\ne 0,j_4=j_5\}}+
$$

\vspace{-4mm}
$$
+
{\bf 1}_{\{i_4=i_1\ne 0,j_4=j_1,i_2=i_3\ne 0,j_2=j_3,i_5=i_6\ne 0,j_5=j_6\}}
+
{\bf 1}_{\{i_4=i_1\ne 0,j_4=j_1,i_2=i_5\ne 0,j_2=j_5,i_3=i_6\ne 0,j_3=j_6\}}+
$$

\vspace{-4mm}
$$
+
{\bf 1}_{\{i_4=i_1\ne 0,j_4=j_1,i_2=i_6\ne 0,j_2=j_6,i_3=i_5\ne 0,j_3=j_5\}}
+
{\bf 1}_{\{i_5=i_1\ne 0,j_5=j_1,i_2=i_3\ne 0,j_2=j_3,i_4=i_6\ne 0,j_4=j_6\}}+
$$

\vspace{-4mm}
$$
+
{\bf 1}_{\{i_5=i_1\ne 0,j_5=j_1,i_2=i_4\ne 0,j_2=j_4,i_3=i_6\ne 0,j_3=j_6\}}
+
{\bf 1}_{\{i_5=i_1\ne 0,j_5=j_1,i_2=i_6\ne 0,j_2=j_6,i_3=i_4\ne 0,j_3=j_4\}}+
$$

\vspace{-4mm}
$$
+
{\bf 1}_{\{i_6=i_1\ne 0,j_6=j_1,i_2=i_3\ne 0,j_2=j_3,i_4=i_5\ne 0,j_4=j_5\}}
+
{\bf 1}_{\{i_6=i_1\ne 0,j_6=j_1,i_2=i_4\ne 0,j_2=j_4,i_3=i_5\ne 0,j_3=j_5\}}+
$$
\begin{equation}
\label{a7}
\Biggl.\Biggl.
+{\bf 1}_{\{i_6=i_1\ne 0,j_6=j_1,i_2=i_5\ne 0,j_2=j_5,i_3=i_4\ne 0,j_3=j_4\}}
\Biggr)\zeta_{j_7}^{(i_7)}\Biggr),
\end{equation}

\vspace{7mm}
\noindent
where ${\bf 1}_A$ is the indicator of the set $A$.

Consider the generalization of the formulas (\ref{a1})--(\ref{a7}) 
for the case of arbitrary multiplicity $k$ of 
$J[\psi^{(k)}]_{T,t}$.
In order to do this, let us
consider the unordered
set $\{1, 2, \ldots, k\}$ 
and separate it into two parts:
the first part consists of $r$ unordered 
pairs (sequence order of these pairs is also unimportant) and the 
second one consists of the 
remaining $k-2r$ numbers.
So, we have

\vspace{1mm}
\begin{equation}
\label{leto5007}
(\{
\underbrace{\{g_1, g_2\}, \ldots, 
\{g_{2r-1}, g_{2r}\}}_{\small{\hbox{part 1}}}
\},
\{\underbrace{q_1, \ldots, q_{k-2r}}_{\small{\hbox{part 2}}}
\}),
\end{equation}

\vspace{4mm}
\noindent
where 

\vspace{-3mm}
$$
\{g_1, g_2, \ldots, 
g_{2r-1}, g_{2r}, q_1, \ldots, q_{k-2r}\}=\{1, 2, \ldots, k\},
$$

\vspace{3mm}
\noindent
braces   
mean an unordered 
set, and pa\-ren\-the\-ses mean an ordered set.

We will say that (\ref{leto5007}) is the partition 
and consider the sum with respect to all possible
partitions

\begin{equation}
\label{leto5008}
\sum_{\stackrel{(\{\{g_1, g_2\}, \ldots, 
\{g_{2r-1}, g_{2r}\}\}, \{q_1, \ldots, q_{k-2r}\})}
{{}_{\{g_1, g_2, \ldots, 
g_{2r-1}, g_{2r}, q_1, \ldots, q_{k-2r}\}=\{1, 2, \ldots, k\}}}}
a_{g_1 g_2, \ldots, 
g_{2r-1} g_{2r}, q_1 \ldots q_{k-2r}},
\end{equation}

\vspace{4mm}
\noindent
where $a_{g_1 g_2, \ldots, 
g_{2r-1} g_{2r}, q_1 \ldots q_{k-2r}}\in\mathbb{R}.$

Below there are several examples of sums in the form (\ref{leto5008})

\vspace{3mm}
$$
\sum_{\stackrel{(\{g_1, g_2\})}{{}_{\{g_1, g_2\}=\{1, 2\}}}}
a_{g_1 g_2}=a_{12},
$$

\vspace{5mm}
$$
\sum_{\stackrel{(\{\{g_1, g_2\}, \{g_3, g_4\}\})}
{{}_{\{g_1, g_2, g_3, g_4\}=\{1, 2, 3, 4\}}}}
a_{g_1 g_2, g_3 g_4}=a_{12,34} + a_{13,24} + a_{23,14},
$$

\vspace{5mm}
$$
\sum_{\stackrel{(\{g_1, g_2\}, \{q_1, q_{2}\})}
{{}_{\{g_1, g_2, q_1, q_{2}\}=\{1, 2, 3, 4\}}}}
a_{g_1 g_2, q_1 q_{2}}=
$$

$$
=a_{12,34}+a_{13,24}+a_{14,23}
+a_{23,14}+a_{24,13}+a_{34,12},
$$

\vspace{5mm}
$$
\sum_{\stackrel{(\{g_1, g_2\}, \{q_1, q_{2}, q_3\})}
{{}_{\{g_1, g_2, q_1, q_{2}, q_3\}=\{1, 2, 3, 4, 5\}}}}
a_{g_1 g_2, q_1 q_{2}q_3}
=
$$

$$
=a_{12,345}+a_{13,245}+a_{14,235}
+a_{15,234}+a_{23,145}+a_{24,135}+
$$

$$
+a_{25,134}+a_{34,125}+a_{35,124}+a_{45,123},
$$

\vspace{5mm}
$$
\sum_{\stackrel{(\{\{g_1, g_2\}, \{g_3, g_{4}\}\}, \{q_1\})}
{{}_{\{g_1, g_2, g_3, g_{4}, q_1\}=\{1, 2, 3, 4, 5\}}}}
a_{g_1 g_2, g_3 g_{4},q_1}
=
$$

$$
=
a_{12,34,5}+a_{13,24,5}+a_{14,23,5}+
a_{12,35,4}+a_{13,25,4}+a_{15,23,4}+
$$

$$
+a_{12,54,3}+a_{15,24,3}+a_{14,25,3}+a_{15,34,2}+a_{13,54,2}+a_{14,53,2}+
$$

$$
+
a_{52,34,1}+a_{53,24,1}+a_{54,23,1}.
$$

\vspace{9mm}

Now, we can formulate Theorem 1 
(see (\ref{tyyy})) 
using the alternative form.

\vspace{2mm}

{\bf Theorem 2} \cite{10} (2009) (also see \cite{11}-\cite{16}, 
\cite{19}-\cite{20aa}, \cite{26}, \cite{31a}-\cite{31aaa}). 
{\it Under the conditions of Theorem {\rm 1} 
the following expansion

\vspace{1mm}

$$
J[\psi^{(k)}]_{T,t}=
\hbox{\vtop{\offinterlineskip\halign{
\hfil#\hfil\cr
{\rm l.i.m.}\cr
$\stackrel{}{{}_{p_1,\ldots,p_k\to \infty}}$\cr
}} }
\sum\limits_{j_1=0}^{p_1}\ldots
\sum\limits_{j_k=0}^{p_k}
C_{j_k\ldots j_1}\Biggl(
\prod_{l=1}^k\zeta_{j_l}^{(i_l)}+\sum\limits_{r=1}^{[k/2]}
(-1)^r \times
\Biggr.
$$

\vspace{2mm}
\begin{equation}
\label{leto6000}
\times
\sum_{\stackrel{(\{\{g_1, g_2\}, \ldots, 
\{g_{2r-1}, g_{2r}\}\}, \{q_1, \ldots, q_{k-2r}\})}
{{}_{\{g_1, g_2, \ldots, 
g_{2r-1}, g_{2r}, q_1, \ldots, q_{k-2r}\}=\{1, 2, \ldots, k\}}}}
\prod\limits_{s=1}^r
{\bf 1}_{\{i_{g_{{}_{2s-1}}}=~i_{g_{{}_{2s}}}\ne 0\}}
\Biggl.{\bf 1}_{\{j_{g_{{}_{2s-1}}}=~j_{g_{{}_{2s}}}\}}
\prod_{l=1}^{k-2r}\zeta_{j_{q_l}}^{(i_{q_l})}\Biggr)
\end{equation}

\vspace{5mm}
\noindent
con\-verg\-ing in the mean-square sense is valid$,$ where $i_1,\ldots,i_k=0,1,\ldots,m,$
$[x]$ is an integer part of a real number $x,$
$\prod\limits_{\emptyset}
\stackrel{\sf def}{=}1,$ 
$\sum\limits_{\emptyset}
\stackrel{\sf def}{=}0;$ 
another notations are the same as in Theorem~{\rm 1.}}

\vspace{2mm}

{\bf Proof.}\ The equality (\ref{leto6000}) will be proved by induction in Sect.~18
(see the proof of Theorem~21).

\vspace{2mm}
\noindent
In particular, from (\ref{leto6000}) for $k=5$ we obtain

\vspace{3mm}

$$
J[\psi^{(5)}]_{T,t}=
\hbox{\vtop{\offinterlineskip\halign{
\hfil#\hfil\cr
{\rm l.i.m.}\cr
$\stackrel{}{{}_{p_1,\ldots,p_5\to \infty}}$\cr
}} }\sum_{j_1=0}^{p_1}\ldots\sum_{j_5=0}^{p_5}
C_{j_5\ldots j_1}\Biggl(
\prod_{l=1}^5\zeta_{j_l}^{(i_l)}-\Biggr.
$$

\vspace{2mm}
$$
-
\sum\limits_{\stackrel{(\{g_1, g_2\}, \{q_1, q_{2}, q_3\})}
{{}_{\{g_1, g_2, q_{1}, q_{2}, q_3\}=\{1, 2, 3, 4, 5\}}}}
{\bf 1}_{\{i_{g_{{}_{1}}}=~i_{g_{{}_{2}}}\ne 0\}}
{\bf 1}_{\{j_{g_{{}_{1}}}=~j_{g_{{}_{2}}}\}}
\prod_{l=1}^{3}\zeta_{j_{q_l}}^{(i_{q_l})}+
$$

\vspace{2mm}
$$
+
\sum_{\stackrel{(\{\{g_1, g_2\}, 
\{g_{3}, g_{4}\}\}, \{q_1\})}
{{}_{\{g_1, g_2, g_{3}, g_{4}, q_1\}=\{1, 2, 3, 4, 5\}}}}
{\bf 1}_{\{i_{g_{{}_{1}}}=~i_{g_{{}_{2}}}\ne 0\}}
{\bf 1}_{\{j_{g_{{}_{1}}}=~j_{g_{{}_{2}}}\}}
\Biggl.{\bf 1}_{\{i_{g_{{}_{3}}}=~i_{g_{{}_{4}}}\ne 0\}}
{\bf 1}_{\{j_{g_{{}_{3}}}=~j_{g_{{}_{4}}}\}}
\zeta_{j_{q_1}}^{(i_{q_1})}\Biggr).
$$

\vspace{6mm}

The last equality obviously agrees with
(\ref{a5}).

It is now appropriate
to make a remark about the structure 
of the formulas (\ref{a1})--(\ref{a7}) and           
(\ref{leto6000}). Using (\ref{ziko1500}), (\ref{2023abc11}),
(\ref{a1})--(\ref{a7}), (\ref{leto6000}), we obtain

\vspace{1mm}
$$
J'[\phi_{j_1}\ldots \phi_{j_k}]_{T,t}^{(i_1\ldots i_k)}=
\prod_{l=1}^k\zeta_{j_l}^{(i_l)}+\sum\limits_{r=1}^{[k/2]}
(-1)^r \times
\Biggr.
$$

\vspace{3mm}
\begin{equation}
\label{leto60001a1b}
\times
\sum_{\stackrel{(\{\{g_1, g_2\}, \ldots, 
\{g_{2r-1}, g_{2r}\}\}, \{q_1, \ldots, q_{k-2r}\})}
{{}_{\{g_1, g_2, \ldots, 
g_{2r-1}, g_{2r}, q_1, \ldots, q_{k-2r}\}=\{1, 2, \ldots, k\}}}}
\prod\limits_{s=1}^r
{\bf 1}_{\{i_{g_{{}_{2s-1}}}=~i_{g_{{}_{2s}}}\ne 0\}}
\Biggl.{\bf 1}_{\{j_{g_{{}_{2s-1}}}=~j_{g_{{}_{2s}}}\}}
\prod_{l=1}^{k-2r}\zeta_{j_{q_l}}^{(i_{q_l})}\Biggr)
\end{equation}

\vspace{5mm}
\noindent
w.~p.~1, where the multiple stochastic integral
$J'[\phi_{j_1}\ldots \phi_{j_k}]_{T,t}^{(i_1\ldots i_k)}$
is defined by (\ref{mult11}); another notations in 
(\ref{leto60001a1b}) are the same as in Theorem~2.

The stochastic integral with respect to the scalar standard Wiener process
($i_1=\ldots=i_k\ne 0$)
and similar to (\ref{mult11}) was considered in \cite{ito1951} (1951)
and is called the multiple Wiener stochastic integral \cite{ito1951}.
Note that $\Phi(t_1,\ldots,t_k)\in L_2([t, T]^k)$ in \cite{ito1951}
(this case will be considered in Sect.~15--18).

As we will see in Sect.~14, 15, 18, the expression on the right-hand side
of (\ref{leto60001a1b}) is the Wick polynomial with arguments
$\zeta_{j_1}^{(i_1)},\ldots,\zeta_{j_k}^{(i_k)}.$
Moreover, the given expression is an explicit representation
of the Wick polynomial, in contrast to its 
representation in the form of a product of Hermite
polynomials (see Sect.~14, 15, 18) or its another representation
(or definition) using a recurrence relation (see (\ref{recur1})).

To best of our knowledge, the representation of the multiple Wiener
stochastic integral in the form of a Wick polynomial 
(see (\ref{leto60001a1b})) for the case of a multidimensional
Wiener process ($i_1,\ldots,i_k=0,1,\ldots,m$)
and the case $j_1,\ldots,j_k=0,1,2,\ldots $ was first obtained
in our monographs \cite{7} (2006), \cite{9} (2007), and
\cite{10} (2009). More precisely, 
the formula (\ref{leto60001a1b}) is obtained
in our monograph \cite{10} (2009) as part of the formula
(5.30) (see \cite{10}, p.~220).
Moreover, partiular cases $k=1,\ldots,5$ (see (\ref{a1})--(\ref{a5})) 
of the formula (\ref{leto60001a1b})
were obtained in \cite{4} (2006) as parts of the formulas 
on the pages 243-244 and partiular cases $k=1,\ldots,7$ 
(see (\ref{a1})--(\ref{a7})) 
of the formula (\ref{leto60001a1b})
were obtained in \cite{9} (2007) as parts of the formulas 
on the pages 208-218.

The indicated formulas are obtained for the case
when $\psi_1(\tau),\ldots,\psi_k(\tau)$ 
are conti\-nu\-ous nonrandom functions on the interval $[t, T]$
and 
$\{\phi_j(x)\}_{j=0}^{\infty}$ is a complete orthonormal system  
of piecewise continuous functions in the space $L_2([t,T])$ 
(see Sect.~2, 4 in this article and \cite{7} (2006), \cite{9} (2007), and
\cite{10} (2009)). Note that the generality
of the above results is even too great when
applied to the numerical integration
of Ito stochastic differential equations.

It should be noted that in \cite{fox} (1987)
an $L_2$--version of the formula (\ref{leto60001a1b})
was obtained, but only for the special case
$j_1=\ldots=j_k$. The above result in \cite{fox} (Proposition~5.1)
is obtained using diagrams, i.e. (unlike our results)
in an implicit form
(see Sect.~18 (below Remark~15) for details).

Let us turn to the comparison
of the formula (\ref{leto60001a1b}) with another interesting work \cite{major2} (2019).
An $L_2$-version of (\ref{leto60001a1b}) was obtained in \cite{major2} in terms 
of Wick polynomials and for the case of vector valued random measures 
(see \cite{major2}, Theorem~7.2, p.~69). In earlier works of this author
(see for example \cite{major1}) only the case of scalar valued random measures 
was considered (see Sect.~18 (below Remark~15) for details).

In Sect.~18 (Theorems~20, 21) 
we consider $L_2$--versions of the formula 
(\ref{leto60001a1b}). At that, to prove 
Theorems~20 and 21 we use only
the Ito formula, in contrast to the diagram method from 
\cite{major2}.

\vspace{5mm}

\section{Comparison of Theorem 2 With Representations of Iterated
Ito Stochastic Integrals Based on Hermite Polynomials}

\vspace{5mm}

Note that the correctness of the formulas (\ref{a1})--(\ref{a7}) 
can be 
verified 
by the fact that if 
$i_1=\ldots=i_7=i=1,\ldots,m$
and $\psi_1(s),\ldots,\psi_7(s)\equiv \psi(s)$,
then we can derive from (\ref{a1})--(\ref{a7}) 
\cite{9} (2007) (also see \cite{10}-\cite{16}, \cite{19}-\cite{20aa})
the well-known
equalities

\vspace{0.5mm}
$$
J[\psi^{(1)}]_{T,t}
=\frac{1}{1!}\delta_{T,t},
$$

\vspace{1.5mm}
$$
J[\psi^{(2)}]_{T,t}
=\frac{1}{2!}\left(\delta^2_{T,t}-\Delta_{T,t}\right),\
$$

\vspace{1.5mm}
$$
J[\psi^{(3)}]_{T,t}
=\frac{1}{3!}\left(\delta_{T,t}^3-3\delta_{T,t}\Delta_{T,t}\right),
$$

\vspace{1.5mm}
$$
J[\psi^{(4)}]_{T,t}
=\frac{1}{4!}\left(\delta^4_{T,t}-6\delta_{T,t}^2\Delta_{T,t}
+3\Delta^2_{T,t}\right),\
$$

\vspace{1.5mm}
$$
J[\psi^{(5)}]_{T,t}
=\frac{1}{5!}\left(\delta^5_{T,t}-10\delta_{T,t}^3\Delta_{T,t}
+15\delta_{T,t}\Delta^2_{T,t}\right),
$$

\vspace{1.5mm}
$$
J[\psi^{(6)}]_{T,t}
=\frac{1}{6!}\left(\delta^6_{T,t}-15\delta_{T,t}^4\Delta_{T,t}
+45\delta_{T,t}^2\Delta^2_{T,t}-15\Delta_{T,t}^3\right),
$$

\vspace{1.5mm}
$$
J[\psi^{(7)}]_{T,t}
=\frac{1}{7!}\left(\delta^7_{T,t}-21\delta_{T,t}^5\Delta_{T,t}
+105\delta_{T,t}^3\Delta^2_{T,t}-105\delta_{T,t}\Delta_{T,t}^3\right),
$$

\vspace{4mm}
\noindent
which fulfilled w. p. 1, where

\vspace{-1mm} 
$$
\delta_{T,t}=\int\limits_t^T\psi(s)d{\bf f}_s^{(i)},\ \ \ \
\Delta_{T,t}=\int\limits_t^T\psi^2(s)ds.
$$

\vspace{3mm}

The above equalities can be independently  
obtained using the Ito formula and Hermite polynomials.

When $k=1$ everything is evident. Let us consider the cases  
$k=2,$ $3.$ When $k=2$ for the case $p_1=p_2=p$ we obtain

\vspace{-1mm}
$$
J[\psi^{(2)}]_{T,t}=
\hbox{\vtop{\offinterlineskip\halign{
\hfil#\hfil\cr
{\rm l.i.m.}\cr
$\stackrel{}{{}_{p\to \infty}}$\cr
}} }\left(
\sum_{j_1,j_2=0}^{p}
C_{j_2j_1}\zeta_{j_1}^{(i)}\zeta_{j_2}^{(i)}-
\sum_{j_1=0}^{p}
C_{j_1j_1}\right)=
$$

\vspace{2mm}
$$
=
\hbox{\vtop{\offinterlineskip\halign{
\hfil#\hfil\cr
{\rm l.i.m.}\cr
$\stackrel{}{{}_{p\to \infty}}$\cr
}} }\left(
\sum_{j_1=0}^{p}\sum_{j_2=0}^{j_1-1}\biggl(
C_{j_2j_1}+C_{j_1j_2}\biggr)
\zeta_{j_1}^{(i)}\zeta_{j_2}^{(i)}+
\sum_{j_1=0}^{p}
C_{j_1j_1}\left(\left(\zeta_{j_1}^{(i)}\right)^2-1\right)\right)
=
$$

\vspace{2mm}
$$
=\hbox{\vtop{\offinterlineskip\halign{
\hfil#\hfil\cr
{\rm l.i.m.}\cr
$\stackrel{}{{}_{p\to \infty}}$\cr
}} }\left(
\sum_{j_1=0}^{p}\sum_{j_2=0}^{j_1-1}
C_{j_1}C_{j_2}
\zeta_{j_1}^{(i)}\zeta_{j_2}^{(i)}+
\frac{1}{2}\sum_{j_1=0}^{p}
C_{j_1}^2\left(\left(\zeta_{j_1}^{(i)}\right)^2-1\right)\right)=
$$

\vspace{2mm}
$$
=\hbox{\vtop{\offinterlineskip\halign{
\hfil#\hfil\cr
{\rm l.i.m.}\cr
$\stackrel{}{{}_{p\to \infty}}$\cr
}} }\left(
\frac{1}{2}
\sum_{\stackrel{j_1,j_2=0}{{}_{j_1\ne j_2}}}^{p}
C_{j_1}C_{j_2}
\zeta_{j_1}^{(i)}\zeta_{j_2}^{(i)}+
\frac{1}{2}\sum_{j_1=0}^{p}
C_{j_1}^2\left(\left(\zeta_{j_1}^{(i)}\right)^2-1\right)\right)=
$$

\vspace{2mm}
$$
=
\hbox{\vtop{\offinterlineskip\halign{
\hfil#\hfil\cr
{\rm l.i.m.}\cr
$\stackrel{}{{}_{p\to \infty}}$\cr
}} }\left(
\frac{1}{2}
\left(\sum_{j_1=0}^{p}
C_{j_1}\zeta_{j_1}^{(i)}\right)^2-
\frac{1}{2}\sum_{j_1=0}^{p}
C_{j_1}^2\right)
=
$$

\vspace{2mm}
\begin{equation}
\label{pipi20}
=\frac{1}{2!}\left(\delta^2_{T,t}-\Delta_{T,t}\right).
\end{equation}

\vspace{5mm}

Let us explain the last 
step in (\ref{pipi20}). For the Ito stochastic 
integrals the following estimate is valid \cite{36}

\vspace{-1mm}
\begin{equation}
\label{pipi}
{\sf M}\left\{\left|\int\limits_t^T \xi_\tau df_\tau\right|^q\right\}
\le K_q {\sf M}\left\{\left(\int\limits_t^T|\xi_\tau|^2 d\tau
\right)^{q/2}\right\},
\end{equation}

\vspace{4mm}
\noindent
where $q>0$ is a fixed number, $f_\tau$ is a scalar
standard Wiener process, $\xi_\tau\in{\rm M}_2([t, T])$,
$K_q$ is a constant depending only on $q$,

\vspace{-1mm}
$$
\int\limits_t^T|\xi_\tau|^2 d\tau<\infty\ \ \ {\rm w.\ p.\ 1},
$$
$$
{\sf M}\left\{\left(\int\limits_t^T|\xi_\tau|^2 d\tau
\right)^{q/2}\right\}<\infty.
$$

\vspace{2mm}

Since

\vspace{-1mm}
$$
\delta_{T,t}-\sum\limits_{j_1=0}^{p}
C_{j_1}\zeta_{j_1}^{(i)}=
\int\limits_t^T\biggl(\psi(s)-
\sum\limits_{j_1=0}^{p}
C_{j_1}\phi_{j_1}(s)\biggr)d{\bf f}_s^{(i)},
$$

\vspace{4mm}
\noindent
then 
using the 
estimate
(\ref{pipi}) to the right-hand side of this expression and considering that

$$
\int\limits_t^T\biggl(\psi(s)-
\sum\limits_{j_1=0}^{p}
C_{j_1}\phi_{j_1}(s)\biggr)^2ds\ \to 0
$$

\vspace{3mm}
\noindent
if $p \to \infty$, we obtain

\vspace{-1mm}
\begin{equation}
\label{pipi21}
\int\limits_t^T\psi(s)d{\bf f}_s^{(i)}=
\hbox{\vtop{\offinterlineskip\halign{
\hfil#\hfil\cr
$q$~-~{\rm l.i.m.}\cr
$\stackrel{}{{}_{p\to \infty}}$\cr
}} }\sum_{j_1=0}^{p}
C_{j_1}\zeta_{j_1}^{(i)},\ \ \ q>0,
\end{equation}

\vspace{3mm}
\noindent
where $\hbox{\vtop{\offinterlineskip\halign{
\hfil#\hfil\cr
$q$~-~{\rm l.i.m.}\cr
$\stackrel{}{{}_{p\to \infty}}$\cr
}} }$ is a limit in the mean of degree $q.$
Hence, if $q=4$, then it is easy to conclude that 

$$
\hbox{\vtop{\offinterlineskip\halign{
\hfil#\hfil\cr
{\rm l.i.m.}\cr
$\stackrel{}{{}_{p\to \infty}}$\cr
}} }
\left(\sum\limits_{j_1=0}^{p}
C_{j_1}\zeta_{j_1}^{(i)}\right)^2=\delta^2_{T,t}.
$$

\vspace{4mm}

This equality as well as Parseval's
equality were used in the last transition of the formula (\ref{pipi20}).

If $k=3$ for the case $p_1=p_2=p_3=p$ we have

\vspace{2mm}
$$
J[\psi^{(3)}]_{T,t}=
$$

\vspace{1mm}
$$
=
\hbox{\vtop{\offinterlineskip\halign{
\hfil#\hfil\cr
{\rm l.i.m.}\cr
$\stackrel{}{{}_{p\to \infty}}$\cr
}} }\left(
\sum_{j_1,j_2,j_3=0}^{p}
C_{j_3j_2j_1}\zeta_{j_1}^{(i)}
\zeta_{j_2}^{(i)}
\zeta_{j_3}^{(i)}
-
\sum_{j_1,j_3=0}^{p}
C_{j_3j_1j_1}\zeta_{j_3}^{(i)}-
\sum_{j_1,j_2=0}^{p}
C_{j_2j_2j_1}\zeta_{j_1}^{(i)}-
\sum_{j_1,j_2=0}^{p}
C_{j_1j_2j_1}\zeta_{j_2}^{(i)}\right)=
$$

\vspace{4mm}
$$
=
\hbox{\vtop{\offinterlineskip\halign{
\hfil#\hfil\cr
{\rm l.i.m.}\cr
$\stackrel{}{{}_{p\to \infty}}$\cr
}} }\left(
\sum_{j_1,j_2,j_3=0}^{p}
C_{j_3j_2j_1}\zeta_{j_1}^{(i)}
\zeta_{j_2}^{(i)}
\zeta_{j_3}^{(i)}
-\sum_{j_1,j_3=0}^{p}\biggl(
C_{j_3j_1j_1}+
C_{j_1j_1j_3}+
C_{j_1j_3j_1}\biggr)\zeta_{j_3}^{(i)}\right)=
$$

\vspace{6mm}
$$
=
\hbox{\vtop{\offinterlineskip\halign{
\hfil#\hfil\cr
{\rm l.i.m.}\cr
$\stackrel{}{{}_{p\to \infty}}$\cr
}} }\left(
\sum_{j_1=0}^{p}
\sum_{j_2=0}^{j_1-1}
\sum_{j_3=0}^{j_2-1}
\biggl(C_{j_3j_2j_1}+
C_{j_3j_1j_2}+C_{j_2j_1j_3}+
C_{j_2j_3j_1}+
C_{j_1j_2j_3}+
C_{j_1j_3j_2}\biggr)
\zeta_{j_1}^{(i)}
\zeta_{j_2}^{(i)}
\zeta_{j_3}^{(i)}+\right.
$$

\vspace{4mm}
$$
+\sum_{j_1=0}^{p}\sum_{j_3=0}^{j_1-1}\biggl(
C_{j_3j_1j_3}+
C_{j_1j_3j_3}+
C_{j_3j_3j_1}\biggr)\left(\zeta_{j_3}^{(i)}\right)^2
\zeta_{j_1}^{(i)}+
$$

\vspace{4mm}
$$
+
\sum_{j_1=0}^{p}\sum_{j_3=0}^{j_1-1}\biggl(
C_{j_3j_1j_1}+
C_{j_1j_1j_3}+
C_{j_1j_3j_1}\biggr)\left(\zeta_{j_1}^{(i)}\right)^2
\zeta_{j_3}^{(i)}+
$$

\vspace{4mm}
$$
\left.+
\sum_{j_1=0}^p C_{j_1j_1j_1}
\left(\zeta_{j_1}^{(i)}\right)^3-
\sum_{j_1,j_3=0}^{p}\left(
C_{j_3j_1j_1}+
C_{j_1j_1j_3}+
C_{j_1j_3j_1}\right)\zeta_{j_3}^{(i)}\right)=
$$

\vspace{8mm}
$$
=
\hbox{\vtop{\offinterlineskip\halign{
\hfil#\hfil\cr
{\rm l.i.m.}\cr
$\stackrel{}{{}_{p\to \infty}}$\cr
}} }\left(
\sum_{j_1=0}^{p}
\sum_{j_2=0}^{j_1-1}
\sum_{j_3=0}^{j_2-1}
C_{j_1}
C_{j_2}C_{j_3}
\zeta_{j_1}^{(i)}
\zeta_{j_2}^{(i)}
\zeta_{j_3}^{(i)}+\right.
$$

\vspace{4mm}
$$
+
\frac{1}{2}\sum_{j_1=0}^{p}\sum_{j_3=0}^{j_1-1}
C_{j_3}^2C_{j_1}\left(\zeta_{j_3}^{(i)}\right)^2
\zeta_{j_1}^{(i)}
+\frac{1}{2}\sum_{j_1=0}^{p}\sum_{j_3=0}^{j_1-1}
C_{j_1}^2C_{j_3}\left(\zeta_{j_1}^{(i)}\right)^2
\zeta_{j_3}^{(i)}+
$$

\vspace{4mm}
$$
\left.+\frac{1}{6}\sum_{j_1=0}^p C_{j_1}^3
\left(\zeta_{j_1}^{(i)}\right)^3-
\frac{1}{2}\sum_{j_1,j_3=0}^{p}
C_{j_1}^2C_{j_3}\zeta_{j_3}^{(i)}\right)=
$$

\vspace{8mm}
$$
=
\hbox{\vtop{\offinterlineskip\halign{
\hfil#\hfil\cr
{\rm l.i.m.}\cr
$\stackrel{}{{}_{p\to \infty}}$\cr
}} }\left(
\frac{1}{6}\sum_{\stackrel{j_1,j_2,j_3=0}{{}_{j_1\ne j_2, j_2\ne j_3,
j_1\ne j_3}}}^{p}
C_{j_1}
C_{j_2}C_{j_3}
\zeta_{j_1}^{(i)}
\zeta_{j_2}^{(i)}
\zeta_{j_3}^{(i)}+\right.
$$

\vspace{4mm}
$$
+\frac{1}{2}\sum_{j_1=0}^{p}\sum_{j_3=0}^{j_1-1}
C_{j_3}^2C_{j_1}\left(\zeta_{j_3}^{(i)}\right)^2
\zeta_{j_1}^{(i)}
+\frac{1}{2}\sum_{j_1=0}^{p}\sum_{j_3=0}^{j_1-1}
C_{j_1}^2C_{j_3}\left(\zeta_{j_1}^{(i)}\right)^2
\zeta_{j_3}^{(i)}+
$$

\vspace{4mm}
$$
\left.+\frac{1}{6}\sum_{j_1=0}^p C_{j_1}^3
\left(\zeta_{j_1}^{(i)}\right)^3-
\frac{1}{2}\sum_{j_1,j_3=0}^{p}
C_{j_1}^2C_{j_3}\zeta_{j_3}^{(i)}\right)
=
$$

\vspace{8mm}
$$
=
\hbox{\vtop{\offinterlineskip\halign{
\hfil#\hfil\cr
{\rm l.i.m.}\cr
$\stackrel{}{{}_{p\to \infty}}$\cr
}} }\left(
\frac{1}{6}\sum_{j_1,j_2,j_3=0}^{p}
C_{j_1}C_{j_2}C_{j_3}\zeta_{j_1}^{(i)}
\zeta_{j_2}^{(i)}\zeta_{j_3}^{(i)}-\right.
$$

\vspace{4mm}
$$
-\frac{1}{6}\left(
3\sum_{j_1=0}^{p}\sum_{j_3=0}^{j_1-1}
C_{j_3}^2C_{j_1}\left(\zeta_{j_3}^{(i)}\right)^2
\zeta_{j_1}^{(i)}
+3\sum_{j_1=0}^{p}\sum_{j_3=0}^{j_1-1}
C_{j_1}^2C_{j_3}\left(\zeta_{j_1}^{(i)}\right)^2
\zeta_{j_3}^{(i)}+
\sum_{j_1=0}^p C_{j_1}^3
\left(\zeta_{j_1}^{(i)}\right)^3\right)+
$$

\vspace{4mm}
$$
+\frac{1}{2}\sum_{j_1=0}^{p}\sum_{j_3=0}^{j_1-1}
C_{j_3}^2C_{j_1}\left(\zeta_{j_3}^{(i)}\right)^2
\zeta_{j_1}^{(i)}
+\frac{1}{2}\sum_{j_1=0}^{p}\sum_{j_3=0}^{j_1-1}
C_{j_1}^2C_{j_3}\left(\zeta_{j_1}^{(i)}\right)^2
\zeta_{j_3}^{(i)}+
$$

\vspace{4mm}
$$
\left.+\frac{1}{6}\sum_{j_1=0}^p C_{j_1}^3
\left(\zeta_{j_1}^{(i)}\right)^3-
\frac{1}{2}\sum_{j_1,j_3=0}^{p}
C_{j_1}^2C_{j_3}\zeta_{j_3}^{(i)}\right)=
$$

\vspace{8mm}
$$
=
\hbox{\vtop{\offinterlineskip\halign{
\hfil#\hfil\cr
{\rm l.i.m.}\cr
$\stackrel{}{{}_{p\to \infty}}$\cr
}} }\left(
\frac{1}{6}\left(\sum_{j_1=0}^{p}
C_{j_1}
\zeta_{j_1}^{(i)}\right)^3-
\frac{1}{2}\sum_{j_1=0}^{p}
C_{j_1}^2 \sum_{j_3=0}^p C_{j_3}\zeta_{j_3}^{(i)}\right)
=
$$

\vspace{4mm}
\begin{equation}
\label{pipi22}
=\frac{1}{3!}\left(\delta_{T,t}^3-3\delta_{T,t}\Delta_{T,t}\right).
\end{equation}

\vspace{6mm}

The last step in (\ref{pipi22}) follows from the Parseval equality,
Theorem 1 for $k=1$, and the equality

$$
\hbox{\vtop{\offinterlineskip\halign{
\hfil#\hfil\cr
{\rm l.i.m.}\cr
$\stackrel{}{{}_{p\to \infty}}$\cr
}} }
\left(\sum\limits_{j_1=0}^{p}
C_{j_1}\zeta_{j_1}^{(i)}\right)^3=\delta^3_{T,t},
$$

\vspace{4mm}
\noindent
which can be obtained easily when
$q=8$ (see (\ref{pipi21})).

In addition, we used the following relations between Fourier 
coefficients for the considered case

\vspace{-1mm}
$$
C_{j_1j_2}+C_{j_2j_1}=C_{j_1}C_{j_2},\ \ \
2C_{j_1j_1}=C_{j_1}^2,
$$

$$
C_{j_1j_2j_3}+
C_{j_1j_3j_2}+
C_{j_2j_3j_1}+
C_{j_2j_1j_3}+
C_{j_3j_2j_1}+
C_{j_3j_1j_2}=C_{j_1}C_{j_2}C_{j_3},
$$

$$
2\left(C_{j_1j_1j_3}+
C_{j_1j_3j_1}+
C_{j_3j_1j_1}\right)=C_{j_1}^2C_{j_3},
$$

$$
6C_{j_1j_1j_1}=C_{j_1}^3.
$$

\vspace{5mm}

\section{On Usage of Discontinuous Complete Orthonormal Systems 
of Functions in Theorem 1}

\vspace{5mm}

Analyzing the proof of Theorem 1, we can ask a natural question: can we 
weaken the condition of continuity of the functions 
$\phi_j(x),$ $j=1, 2,\ldots$?

\vspace{2mm}
   
{\it We will say that the function $f(x):$ $[t, T]\to \mathbb{R}$ 
satisfies the condition {\rm (}$\star ${\rm )} if it is 
continuous on the interval $[t, T]$ except may be
for the finite number of points 
of the finite discontinuity as well as it is 
right-continuous 
on the interval  $[t, T].$}

\vspace{2mm}

Furthermore, let us suppose that $\{\phi_j(x)\}_{j=0}^{\infty}$
is a complete orthonormal system of functions in the 
space $L_2([t, T])$, each function $\phi_j(x)$ of which
for $j<\infty$ satisfies the 
condition $(\star)$.

It is easy to see that continuity of the functions $\phi_j(x)$ was used 
substantially in the proof of Theorem 1 in two places: Lemma 3 and 
the formula (\ref{30.52}).
It is clear that without the loss of generality the partition 
$\{\tau_j\}_{j=0}^N$ of the interval $[t, T]$ in Lemma 3 and 
in the formula (\ref{30.52}) can be taken so "dense"\  that 
among the 
points $\tau_j$ of this partition will be all points 
of jumps of the functions
$\varphi_1(\tau)=\phi_{j_1}(\tau),$ $\ldots,$ 
$\varphi_k(\tau)=\phi_{j_k}(\tau)$
($j_1,\ldots,j_k<\infty$)
and among the points $(\tau_{j_1},\ldots,\tau_{j_k})$ for which
$0\le j_1<\ldots<j_k\le N-1$
there will be all points of jumps
of the function
$\Phi(t_1,\ldots,t_k)$.

Let us demonstrate how to modify the proofs of Lemma 3 and the formula 
(\ref{30.52}) 
in the case when $\{\phi_j(x)\}_{j=0}^{\infty}$ is a complete
orthonormal 
system of functions in the space $L_2([t, T])$,  
each function $\phi_j(x)$ 
of which for $j<\infty$ satisfies the condition $(\star)$.

At first, consider Lemma 3. In the proof of this lemma we
obtained the following 
relations

$$
{\sf M}\left\{\left|\sum_{j=0}^{N-1}J[\Delta\varphi_l]_{\tau_{j+1},
\tau_j}\right|^4
\right\}=
\sum_{j=0}^{N-1}{\sf M}\left\{\biggl|J[\Delta\varphi_l]_{\tau_{j+1},
\tau_j}\biggr|^4
\right\}+
$$

\vspace{2mm}
\begin{equation}
\label{4444.02}
+ 6 \sum_{j=0}^{N-1}{\sf M}
\left\{\biggl|J[\Delta\varphi_l]_{\tau_{j+1},\tau_j}\biggr|^2
\right\}
\sum_{q=0}^{j-1}{\sf M}\left\{\biggl|
J[\Delta\varphi_l]_{\tau_{q+1},\tau_q}\biggr|^2
\right\},
\end{equation}

\vspace{3mm}

$$
{\sf M}\left\{\left|J[\Delta\varphi_l]_{\tau_{j+1},\tau_j}\right|^2\right\}=
\int\limits_{\tau_j}^{\tau_{j+1}}(\varphi_l(\tau_j)-\varphi_l(s))^2ds,\
$$

\vspace{1mm}
$$
{\sf M}\left\{\left|J[\Delta\varphi_l]_{\tau_{j+1},\tau_j}\right|^4\right\}=
3\left(\int\limits_{\tau_j}^{\tau_{j+1}}(\varphi_l(\tau_j)-\varphi_l(s))^2ds
\right)^2.
$$

\vspace{6mm}

Suppose that the functions $\varphi_l(s)$ ($l=1,\ldots,k$)  satisfy 
the condition $(\star)$ and the partition $\{\tau_j\}_{j=0}^{N}$ 
includes all points of 
jumps of the functions
$\varphi_l(s)$ ($l=1,\ldots,k$).  It means that for the integral

$$
\int\limits_{\tau_j}^{\tau_{j+1}}(\varphi_l(\tau_j)-\varphi_l(s))^2ds
$$

\vspace{3mm}
\noindent
the integrand function is 
continuous 
at the interval
$[\tau_j, \tau_{j+1}],$ except possibly the point $\tau_{j+1}$ of
finite 
discontinuity.

Let $\mu\in (0, \Delta\tau_j)$ be fixed. Then, due to
continuity (which means uniform continuity) of the functions
$\varphi_l(s)$ ($l=1,\ldots,k$) 
on the interval $[\tau_j, \tau_{j+1}-\mu]$ we have

\vspace{1mm}
$$
\int\limits_{\tau_j}^{\tau_{j+1}}(\varphi_l(\tau_j)-\varphi_l(s))^2ds=
$$
\begin{equation}
\label{4444.90}
=
\int\limits_{\tau_j}^{\tau_{j+1}-\mu}(\varphi_l(\tau_j)-\varphi_l(s))^2ds 
+ \int\limits_{\tau_{j+1}-\mu}^{\tau_{j+1}}(\varphi_l(\tau_j)-\varphi_l(s))^2ds
<\varepsilon^2 (\Delta\tau_j-\mu)+M^2\mu.
\end{equation}

\vspace{5mm}

Obtaining the inequality (\ref{4444.90}), we proposed that
$\Delta\tau_j<\delta(\varepsilon)$ for $j=0,\ 1,\ldots,N-1$ 
($\delta(\varepsilon)>0$ exists for any $\varepsilon>0$ 
and it does not depend 
on $s$), 

$$
|\varphi_l(\tau_j)-\varphi_l(s)|<\varepsilon
$$

\vspace{3mm}
\noindent
if $s\in [\tau_{j}, \tau_{j+1}-\mu]$
(due to uniform continuity of the 
functions $\varphi_l(s)$ ($l=1,\ldots,k$)),

$$
|\varphi_l(\tau_j)-\varphi_l(s)|<M
$$

\vspace{3mm}
\noindent
if
$s\in [\tau_{j+1}-\mu, \tau_{j+1}]$, $M$ is a constant
(potential point of discontinuity 
of the function $\varphi_l(s)$ 
is supposed in the 
point $\tau_{j+1}$).

Performing the 
passage to the limit
in the inequality (\ref{4444.90}) when $\mu\to +0$, we get

$$
\int\limits_{\tau_j}^{\tau_{j+1}}(\varphi_l(\tau_j)-\varphi_l(s))^2ds\le
\varepsilon^2 \Delta\tau_j.
$$

\vspace{3mm}

Using this estimate for the right-hand side of 
(\ref{4444.02}),
we obtain

\vspace{1mm}
$$
{\sf M}\left\{\left|\sum_{j=0}^{N-1}J[\Delta\varphi_l]_{\tau_{j+1},
\tau_j}\right|^4
\right\}\le 
\varepsilon^4\left(
3 \sum_{j=0}^{N-1}(\Delta\tau_{j})^2+
6 \sum_{j=0}^{N-1}\Delta\tau_{j}
\sum_{q=0}^{j-1}\Delta\tau_{q}\right)<
$$

\vspace{2mm}
\begin{equation}
\label{4444.92}
<
3\varepsilon^4\left(\delta(\varepsilon) (T-t)+(T-t)^2\right).
\end{equation}

\vspace{5mm}

This implies that
$$
{\sf M}\left\{\left|\sum\limits_{j=0}^{N-1}J[\Delta\varphi_l]_{\tau_{j+1}
\tau_j}\right|^4\right\} \to 0
$$

\vspace{3mm}
\noindent
if $N\to\infty$. So, Lemma 3 remains valid.

Now, let us present explanations concerning the correctness of the formula 
(\ref{30.52}) when $\{\phi_j(x)\}_{j=0}^{\infty}$ is a 
complete
orthonormal system of functions in the space $L_2([t, T])$, 
each function
$\phi_j(x)$ of which for $j<\infty$ satisfies the condition $(\star)$.

Let us consider the case $k=3$ and the representation (\ref{4444.1}). 
We can demonstrate that in the studied case the first limit 
on the right-hand side of (\ref{4444.1}) equals to zero 
(similarly we demonstrate that the second limit on the right-hand side 
of (\ref{4444.1}) equals to zero; 
proof of the second limit 
equality to zero on the right-hand side
of the formula
(\ref{44444.25}) is the same as for the case of continuous functions
$\phi_j(x),$ $j=0, 1,\ldots).$

The second moment of the prelimit expression of first limit on the 
right-hand side of (\ref{4444.1}) looks as follows

\vspace{-1mm}
$$
\sum_{j_3=0}^{N-1}\sum_{j_2=0}^{j_3-1}\sum_{j_1=0}^{j_2-1}
\int\limits_{\tau_{j_2}}^{\tau_{j_2+1}}\int\limits_{\tau_{j_1}}^{\tau_{j_1+1}}
\left(\Phi(t_1,t_2,\tau_{j_3})-\Phi(t_1,\tau_{j_2},\tau_{j_3})\right)^2
dt_1 dt_2
\Delta\tau_{j_3}.
$$

\vspace{5mm}
       
Further, for the fixed $\mu\in(0, \Delta\tau_{j_2})$ 
and $\rho\in(0, \Delta\tau_{j_1})$ we have

\vspace{2mm}
$$
\int\limits_{\tau_{j_2}}^{\tau_{j_2+1}}\int\limits_{\tau_{j_1}}^{\tau_{j_1+1}}
\left(\Phi(t_1,t_2,\tau_{j_3})-\Phi(t_1,\tau_{j_2},\tau_{j_3})\right)^2
dt_1 dt_2=
$$

\vspace{1mm}
$$
=\left(\int\limits_{\tau_{j_2}}^{\tau_{j_2+1}-\mu}
+\int\limits_{\tau_{j_2+1}-\mu}^{\tau_{j_2+1}}\right)
\left(
\int\limits_{\tau_{j_1}}^{\tau_{j_1+1}-\rho}+
\int\limits_{\tau_{j_1+1}-\rho}^{\tau_{j_1+1}}\right)
\left(\Phi(t_1,t_2,\tau_{j_3})-\Phi(t_1,\tau_{j_2},\tau_{j_3})\right)^2
dt_1 dt_2=
$$

\vspace{1mm}
$$
=\left(\int\limits_{\tau_{j_2}}^{\tau_{j_2+1}-\mu}
\int\limits_{\tau_{j_1}}^{\tau_{j_1+1}-\rho}+
\int\limits_{\tau_{j_2}}^{\tau_{j_2+1}-\mu}
\int\limits_{\tau_{j_1+1}-\rho}^{\tau_{j_1+1}}+
\int\limits_{\tau_{j_2+1}-\mu}^{\tau_{j_2+1}}
\int\limits_{\tau_{j_1}}^{\tau_{j_1+1}-\rho}+
\int\limits_{\tau_{j_2+1}-\mu}^{\tau_{j_2+1}}
\int\limits_{\tau_{j_1+1}-\rho}^{\tau_{j_1+1}}\right)\times
$$

\vspace{3mm}
$$
\times
\left(\Phi(t_1,t_2,\tau_{j_3})-\Phi(t_1,\tau_{j_2},\tau_{j_3})\right)^2
dt_1 dt_2<
$$

\begin{equation}
\label{4444.54}
<\varepsilon^2\left(\Delta\tau_{j_2}-\mu\right)
\left(\Delta\tau_{j_1}-\rho\right)
+M^2\rho\left(\Delta\tau_{j_2}-\mu\right)
+M^2\mu\left(\Delta\tau_{j_1}-\rho\right)
+M^2\mu\rho,
\end{equation}

\vspace{5mm}
\noindent
where $M$ is a constant,
$\Delta\tau_j<\delta(\varepsilon)$ for $j=0, 1,\ldots,N-1$
($\delta(\varepsilon)>0$ exists for any $\varepsilon>0$ and it does not  
depend on points
$(t_1,t_2,\tau_{j_3}),$ $(t_1,\tau_{j_2},\tau_{j_3})$).
We suppose here that the partition 
$\{\tau_j\}_{j=0}^{N}$ contains all discontinuity points
of the function
$\Phi(t_1,t_2,t_3)$ 
as points $\tau_j$ (for every variable).
When obtaining
(\ref{4444.54}),  
we also supposed that potential discontinuity points
of this function (for every variable) are contained among the points
$\tau_{j_1+1}, \tau_{j_2+1}, \tau_{j_3+1}$.

Let us explain in detail how we obtained the inequality 
(\ref{4444.54}). 
Since the function $\Phi(t_1,t_2,t_3)$ is continuous on the closed 
bounded set  

$$
Q_3=\biggl\{(t_1, t_2, t_3): t_1\in[\tau_{j_1}, \tau_{j_1+1}-\rho],
t_2\in[\tau_{j_2}, \tau_{j_2+1}-\mu],
t_3\in [\tau_{j_3}, \tau_{j_3+1}-\nu]\biggr\},
$$

\vspace{3mm}
\noindent 
where $\rho, \mu, \nu$
are fixed small positive numbers such that

$$
\nu\in(0, \Delta\tau_{j_3}),\ \ \ 
\mu\in(0, \Delta\tau_{j_2}),\ \ \ 
\rho\in(0, \Delta\tau_{j_1}),
$$

\vspace{3mm}
\noindent
then this function is also uniformly continous 
on this set and bounded on the closed set 
$D_3$.

Since the distance between the points
$(t_1,t_2,\tau_{j_3})$,
$(t_1,\tau_{j_2},\tau_{j_3})$ $\in $ $Q_3$ 
is obviously less than
$\delta(\varepsilon)$ ($\Delta\tau_j<\delta(\varepsilon)$ for
$j=0, 1,\ldots,N-1$),
then 

$$
|\Phi(t_1,t_2,\tau_{j_3})-\Phi(t_1,\tau_{j_2},\tau_{j_3})|
<\varepsilon.
$$

\vspace{3mm}
 
This inequality was used to 
estimate the first double integral in (\ref{4444.54}). 
Estimating 
the three remaining double integrals, we used the property of 
boundedness of the function 
$\Phi(t_1,t_2,t_3)$ in form of inequality  

\vspace{-2mm}
$$
|\Phi(t_1,t_2,\tau_{j_3})-\Phi(t_1,\tau_{j_2},\tau_{j_3})|
<M.
$$

\vspace{4mm}

Performing the 
passage to the limit 
in the inequality (\ref{4444.54})
if $\mu, \rho\to +0$, we obtain the estimate

$$
\int\limits_{\tau_{j_2}}^{\tau_{j_2+1}}\int\limits_{\tau_{j_1}}^{\tau_{j_1+1}}
\left(\Phi(t_1,t_2,\tau_{j_3})-\Phi(t_1,\tau_{j_2},\tau_{j_3})\right)^2
dt_1 dt_2\le
\varepsilon^2\Delta\tau_{j_2}
\Delta\tau_{j_1}.
$$

\vspace{3mm}

Usage of this estimate provides

\vspace{1mm}
$$
\sum_{j_3=0}^{N-1}\sum_{j_2=0}^{j_3-1}\sum_{j_1=0}^{j_2-1}
\int\limits_{\tau_{j_2}}^{\tau_{j_2+1}}\int\limits_{\tau_{j_1}}^{\tau_{j_1+1}}
\left(\Phi(t_1,t_2,\tau_{j_3})-\Phi(t_1,\tau_{j_2},\tau_{j_3})\right)^2
dt_1 dt_2
\Delta\tau_{j_3}\le
$$

\vspace{2mm}
$$
\le\varepsilon^2
\sum_{j_3=0}^{N-1}\sum_{j_2=0}^{j_3-1}\sum_{j_1=0}^{j_2-1}
\Delta\tau_{j_1}\Delta\tau_{j_2}\Delta\tau_{j_3}<
\varepsilon^2 \frac{(T-t)^3}{6}.
$$

\vspace{5mm}

The last estimate means that in the considered case the first 
limit on the right-hand side of (\ref{4444.1}) equals to zero 
(similarly we can
demonstrate that the second limit on the right-hand side of (\ref{4444.1})
equals to zero).

Consequently, the formula (\ref{30.52}) is correct when 
$k=3$ in the considered case. 
Similarly, we perform argumentation for the cases
$k=2$ and $k>3.$

Consequently, in Theorem 1 we can use complete orthonormal systems of 
functions $\{\phi_j(x)\}_{j=0}^{\infty}$ in the space 
$L_2([t, T]),$ each function 
$\phi_j(x)$ of which for $j<\infty$ satisfies the condition
$(\star)$.    

One of the examples of such systems of functions is
a complete orthonormal 
system of Haar functions in the space 
$L_2([t, T])$ 

\vspace{1mm}
$$
\phi_0(x)=\frac{1}{\sqrt{T-t}},\ \ \ \
\phi_{nj}(x)=\frac{1}{\sqrt{T-t}}\varphi_{nj}\biggl(\frac{x-t}{T-t}\biggr),
$$

\vspace{5mm}
\noindent
where
$n=0, 1,\ldots,$ $j=1, 2,\ldots, 2^n,$
and the functions $\varphi_{nj}(x)$ have the following form

\vspace{3mm}
$$
\varphi_{nj}(x)=
\begin{cases}
2^{n/2},\ &x\in[(j-1)/2^n,\ (j-1)/2^n+
1/2^{n+1})\cr\cr\cr
-2^{n/2},\ &x\in[(j-1)/2^n+1/2^{n+1},\
j/2^n)\cr\cr\cr
0,\  &\hbox{otherwise}
\end{cases},
$$

\vspace{6mm}
\noindent
where $n=0, 1,\ldots,$ $j=1, 2,\ldots, 2^n$ 
(we choose the values of Haar functions 
in the points of discontinuity in such a way that these functions
will be 
right-continuous).

The other example of similar system of functions is a complete orthonormal 
system of Rademacher--Walsh functions in the space  $L_2([t, T])$

\vspace{1mm}
$$
\phi_0(x)=\frac{1}{\sqrt{T-t}},
$$

\vspace{2mm}
$$
\phi_{m_1\ldots m_k}(x)=
\frac{1}{\sqrt{T-t}}\varphi_{m_1}\biggl(\frac{x-t}{T-t}\biggr)
\ldots \varphi_{m_k}
\biggl(\frac{x-t}{T-t}\biggr),
$$

\vspace{7mm}
\noindent
where $0<m_1<\ldots<m_k,$ $m_1,\ldots,m_k=1, 2,\ldots,$ $k=1, 2,\ldots,$

\vspace{2mm}
$$
\varphi_m(x)=(-1)^{[2^m x]},
$$

\vspace{5mm}
\noindent
$x\in [0, 1]$, $ m=1, 2,\ldots,$ $[y]$ is an integer part of a real number
$y.$

\vspace{5mm}

\section{Remark on Usage of Complete Orthonormal Systems 
of Functions in Theorem 1}

\vspace{5mm}

Note that actually the 
func\-ti\-ons $\phi_j(s)$ of complete orthonormal system 
of func\-ti\-ons $\{\phi_j(s)\}_{j=0}^{\infty}$ in the space 
$L_2([t, T])$ depend not only on $s$, but 
on $t$ and $T.$
   
For example, the complete orthonormal systems of Legendre polynomials 
and trigonometric func\-ti\-ons in the space $L_2([t, T])$ have the 
following form

\vspace{1mm}
$$
\phi_j(s,t,T)=\sqrt{\frac{2j+1}{T-t}}P_j\left(\left(
s-\frac{T+t}{2}\right)\frac{2}{T-t}\right),
$$

\vspace{5mm}
where
$P_j(s)$ $(j=0, 1, 2,\ldots)$ is the Legendre polynomial,

\vspace{2mm}

$$
\phi_j(s,t,T)=\frac{1}{\sqrt{T-t}}
\left\{
\begin{matrix}
1,\ & j=0\cr\cr\cr
\sqrt{2}{\rm sin} \left(2\pi r(s-t)/(T-t)\right),\ & j=2r-1\cr\cr\cr
\sqrt{2}{\rm cos} \left(2\pi r(s-t)/(T-t)\right),\ & j=2r
\end{matrix},\right.
$$

\vspace{5mm}
\noindent
where $r=1, 2,\ldots $

Note that the specified systems of functions are assumed to be used in
the context of implementing of numerical methods 
for Ito stochastic differential 
equations for the sequences of time intervals 

$$
[T_0, T_1],\
[T_1, T_2],\ [T_2, T_3],\ \ldots\ ,
$$

\vspace{2mm}
\noindent
and spaces 

\vspace{-1mm}
$$
L_2([T_0, T_1]),\ L_2([T_1, T_2]),\
L_2([T_2, T_3]),\ \ldots
$$

\vspace{2mm}

We can explain that the dependence of functions 
$\phi_j(s,t,T)$ from $t$ and $T$ (hereinafter these constants will 
mean fixed moments of time) will not affect the main properties 
of independence of the random variables

\vspace{-1mm}
\begin{equation}
\label{a102}
\zeta_{(j)T,t}^{(i)}=
\int\limits_t^T \phi_{j}(s,t,T) d{\bf w}_s^{(i)},
\end{equation}

\vspace{3mm}
\noindent
where $i=1,\ldots,m$ and $j=0, 1, 2,\ldots $

Indeed, for fixed  $t$ and $T$ due to orthonormality of the mentioned systems 
of functions we have

\vspace{1mm}
\begin{equation}
\label{a101}
{\sf M}\left\{\zeta_{(j)T,t}^{(i)}\zeta_{(g)T,t}^{(r)}\right\}=
{\bf 1}_{\{i=r\}}
{\bf 1}_{\{j=g\}},
\end{equation}

\vspace{3mm}
\noindent
where 
$$
\zeta_{(j)T,t}^{(i)}=
\int\limits_t^T \phi_{j}(s,t,T) d{\bf w}_s^{(i)},\ \ \
i, r=1,\ldots,m,\ \ \ j, g=0, 1, 2,\ldots
$$

\vspace{3mm}

Note that (\ref{a101}) means the independence of random variables
(\ref{a102}) for various $i$ or $j$.

On the other side,
the random
variables

$$
\zeta_{(j)T_1,t_1}^{(i)}=
\int\limits_{t_1}^{T_1}\phi_{j}(s,t_1,T_1) d{\bf w}_s^{(i)},\ \ \ 
\zeta_{(j)T_2,t_2}^{(i)}=
\int\limits_{t_2}^{T_2}\phi_{j}(s,t_2,T_2) d{\bf w}_s^{(i)}
$$

\vspace{3mm}
\noindent
are independent if
$[t_1, T_1]\cap [t_2, T_2]=\emptyset$
(the case $T_1=t_2$ is possible)
according to 
the property
of the Ito stochastic integral.

Therefore, two important characteristics of random variables
$\zeta_{(j)T,t}^{(i)}$, which are the basic motive of their 
usage, are saved.

\vspace{5mm}

\section{Convergence in the 
Mean of Degree $2n$ ($n\in\mathbb{N}$) of Expansion of Iterated
Ito Stochastic Integrals From Theorem 1}

\vspace{5mm}

Constructing the expansions of iterated
Ito stochastic integrals from Theorem 1 we 
saved all 
information about these integrals. That is why it is 
natural to expect that the mentioned expansions will converge
not only in the mean-square sense but in the stronger probabilistic senses.

We will obtain the general estimate which prove convergence in 
the mean
of degree $2n$ ($n\in {\bf N}$) of expansion from Theorem 1.

According to the notations of Theorem 1 (see (\ref{pobedaxyz}), (\ref{y007})), we have

\vspace{2mm}
$$
R_{T,t}^{p_1,\ldots, p_k}
=J[\psi^{(k)}]_{T,t}-
J[\psi^{(k)}]_{T,t}^{p_1,\ldots,p_k}=J'[R_{p_1\ldots p_k}]_{T,t}^{(i_1\ldots i_k)}=
$$

\vspace{1mm}
\begin{equation}
\label{jjjye}
=\sum_{(t_1,\ldots,t_k)}
\int\limits_{t}^{T}
\ldots
\int\limits_{t}^{t_2}
R_{p_1\ldots p_k}(t_1,\ldots,t_k)
d{\bf f}_{t_1}^{(i_1)}
\ldots
d{\bf f}_{t_k}^{(i_k)},
\end{equation}

\vspace{3mm}
\noindent
where
\begin{equation}
\label{ch12026ll12}
R_{p_1\ldots p_k}(t_1,\ldots,t_k)\stackrel{{\rm def}}{=}
K(t_1,\ldots,t_k)-
\sum_{j_1=0}^{p_1}\ldots
\sum_{j_k=0}^{p_k}
C_{j_k\ldots j_1}
\prod_{l=1}^k\phi_{j_l}(t_l),
\end{equation}

\vspace{3mm}
\noindent
$J[\psi^{(k)}]_{T,t}$ is the stochastic integral {\rm (\ref{sodom20}),}
$J[\psi^{(k)}]_{T,t}^{p_1,\ldots,p_k}$ is the 
expression on the right-hand side of {\rm (\ref{tyyy})} before
passing to the limit 
$$
\hbox{\vtop{\offinterlineskip\halign{
\hfil#\hfil\cr
{\rm l.i.m.}\cr
$\stackrel{}{{}_{p_1,\ldots,p_k\to \infty}}$\cr
}} }.
$$ 

\vspace{3mm}

Note that for definiteness we consider 
the case $i_1,\ldots,i_k=1,\ldots,m$ in this section.
Another 
notations from this section are the same as in the
formulation and proof of Theorem 1.

When proving Theorem 1 we 
obtained the following estimate (see (\ref{obana1})) 

\vspace{1mm}
$$
{\sf M}\left\{\left(J'[R_{p_1\ldots p_k}]_{T,t}^{(i_1\ldots i_k)}\right)^2\right\}\le C_k 
\int\limits_{[t, T]^k}
R_{p_1\ldots p_k}^2(t_1,\ldots,t_k)
dt_1
\ldots
dt_k,
$$

\vspace{4mm}
\noindent
where $C_k$ is a constant.
Obviously, $C_k=k!$ for the case $i_1,\ldots,i_k=1,\ldots,m$.

First, we note that the iterated Ito stochastic integral (\ref{sodom20})
can be considered as a multiple Wiener stochastic integral 
with respect to the components of a multidimensional
Wiener process. The multiple Wiener stochastic integral
with respect to a scalar Wiener process was first considered
in \cite{ito1951} (1951). 
Multiple Wiener stochastic integrals, including integrals
with respect to the components of a 
multidimensional Wiener process, are discussed in detail
in Sect.~14, 15, 18 (also see Sect.~13).
In fact, we have already considered the 
multiple Wiener stochastic integral 
(see (\ref{mult11})).

Let $H(s): [t,T]\to {\mathbb R}$.
Let $(\Omega,{\rm F},{\sf P})$ is a probability space, where
${\rm F}$ is the smallest $\sigma$-algebra such that
the random variables
$$
\int\limits_t^T H(s)dw_s
$$

\vspace{2mm}

\noindent
are ${\rm F}$-measurable for every $H(s)\in L_2([t,T]),$
where 
$w_s$ is a standard Wiener process,

\vspace{-1mm}
$$
\int\limits_t^T H(s)dw_s
$$

\vspace{3mm}
\noindent
is a usual Wiener--Ito stochastic integral.

It is well known (see Theorem~9.7.1, Theorem~9.7.3 and Theorem~9.6.7 in \cite{Kuo}) that

$$
L_2(\Omega,{\rm F},{\sf P})=\bigoplus\limits_{k=0}^{\infty}{\mathbb H}_k,
$$

\vspace{4mm}
\noindent
where ${\mathbb H}_0$ contains only constants,
the space ${\mathbb H}_k$ is the so-called $k$th $(k\ge 1)$ homogeneous Wiener chaos
which consists of all random variables of the form

\vspace{-1mm}
\begin{equation}
\label{ch12026oooo1}
n!\int\limits_t^{T}\ldots \int\limits_t^{t_2}
\hat H(t_1,\ldots,t_k)dw_{t_1}\ldots dw_{t_k},
\end{equation}

\vspace{3mm}
\noindent
where (\ref{ch12026oooo1}) is a representation for the 
multiple Wiener stochastic integral with respect to 
the scalar Wiener process (see (\ref{mult11}) or (\ref{WiI}) for the case 
$i_1=\ldots=i_k=i\in\{1,\ldots,m\}$ and with $w_s$ instead of ${\bf w}^{(i)}_s$),
$\hat H(t_1,\ldots,t_k)$ is a symmetrization
of $H(t_1,\ldots,t_k),$ $H(t_1,\ldots,t_k)\in L_2(D_k),$
$D_k=\{(t_1,\ldots,t_k):  t<t_1<\ldots <t_k<T\}.$

Note that (see \cite{Nurd2}, Corollary~2.8.4)

\begin{equation}
\label{ch12026777}
{\sf M}\left\{\xi^{2n}\right\}\le (2n-1)^{kn}\left({\sf M}\left\{\xi^2\right\}\right)^n,
\end{equation}

\vspace{3mm}
\noindent
where $\xi\in {\mathbb H}_k,$ $n, k\in {\mathbb N}.$

The following estimate 
for the multiple Wiener stochastic integral
with respect to a scalar Wiener process 
is a consequence of the inequality (\ref{ch12026777})
\cite{Nual1}, \cite{Nurd1}

\begin{equation}
\label{2026ch1001}
{\sf M}\left\{\left(
J'[\Phi]_{T,t}^{(i_1\ldots i_k)}\right)^{2n}\right\}
\le (2n-1)^{nk} \left({\sf M}\left\{\left(
J'[\Phi]_{T,t}^{(i_1\ldots i_k)}\right)^{2}\right\}\right)^{n},
\end{equation}

\vspace{3mm}
\noindent
where $n\in{\bf N},$ 
$J'[\Phi]_{T,t}^{(i_1\ldots i_k)}$ is a multiple
Wiener stochastic integral defined as in Sect.~2 (see (\ref{mult11})) or as in Sect.~15
(see (\ref{WiI}))
but for the case of a scalar Wiener process $(i_1=\ldots=i_k=i\in\{1,\ldots,m\})$,
$\Phi(t_1,\ldots,t_k)\in L_2([t, T]^k)$ in (\ref{WiI}) and  
$\Phi(t_1,\ldots,t_k)\in C([t, T]^k)\subset L_2([t, T]^k)$ in (\ref{mult11}),
$J'[\Phi]_{T,t}^{(i_1\ldots i_k)}\in {\mathbb H}_k$ 
(see (\ref{pobeda}), (\ref{pobedaxyz}), (\ref{Wi110})). We also note that

$$
J'[\Phi]_{T,t}^{(i_1\ldots i_k)}=J'[\hat\Phi]_{T,t}^{(i_1\ldots i_k)}\ \ \ \hbox{w.~p.~1,}
$$

\vspace{3mm}
\noindent
where $i_1=\ldots=i_k=i\in\{1,\ldots,m\}$ and $\hat\Phi(t_1,\ldots,t_k)$ is a symmetrization
of the function $\Phi(t_1,\ldots,t_k)$.

Consider the elementary inequality

\vspace{-2mm}
\begin{equation}
\label{ch12026m300}
\left(a_1+a_2+\ldots+a_p\right)^2 \le
p\left(a_1^2+a_2^2+\ldots+a_p^2\right),\ \ \ p\in {\bf N}.
\end{equation}

\vspace{3mm}

Using the inequality (\ref{ch12026m300}) and (\ref{jjjye}), we obtain 

\vspace{1mm}
$$
{\sf M}\left\{\left(
J'[R_{p_1\ldots p_k}]_{T,t}^{(i_1\ldots i_k)}\right)^{2}\right\}=
$$

\vspace{2mm}
$$
=
{\sf M}\left\{\left(
\sum_{(t_1,\ldots,t_k)}
\int\limits_{t}^{T}
\ldots
\int\limits_{t}^{t_2}
R_{p_1\ldots p_k}(t_1,\ldots,t_k)
d{\bf f}_{t_1}^{(i_1)}
\ldots
d{\bf f}_{t_k}^{(i_k)}\right)^2\right\}\le
$$

\vspace{2mm}
$$
\le k!
\sum_{(t_1,\ldots,t_k)}
{\sf M}\left\{\left(
\int\limits_{t}^{T}
\ldots
\int\limits_{t}^{t_2}
R_{p_1\ldots p_k}(t_1,\ldots,t_k)
d{\bf f}_{t_1}^{(i_1)}
\ldots
d{\bf f}_{t_k}^{(i_k)}\right)^2\right\}=
$$

\vspace{2mm}
$$
= k!
\sum_{(t_1,\ldots,t_k)}
\int\limits_{t}^{T}
\ldots
\int\limits_{t}^{t_2}
R^2_{p_1\ldots p_k}(t_1,\ldots,t_k)
dt_1
\ldots
dt_k=
$$

\vspace{3mm}
\begin{equation}
\label{ch12026102}
=k!
\int\limits_{[t, T]^k}
R^2_{p_1\ldots p_k}(t_1,\ldots,t_k)
dt_1
\ldots
dt_k,
\end{equation}

\vspace{3mm}
\noindent
where $i_1,\ldots,i_k=1,\ldots,m.$

Suppose that $\{\phi_j(x)\}_{j=0}^{\infty}$
is a complete orthonormal system of functions in 
the space $L_2([t, T])$. 
Using the orthonormality of the functions $\phi_j(x)$ $(j=0, 1, 2,\ldots),$ 
we obtain

\vspace{2mm}
$$
\int\limits_{[t, T]^k}
R^2_{p_1\ldots p_k}(t_1,\dots,t_k)dt_1\ldots dt_k
=
$$

\vspace{3mm}
$$
=\int\limits_{[t,T]^k}
\Biggl(K(t_1,\ldots,t_k)-
\sum_{j_1=0}^{p_1}\ldots
\sum_{j_k=0}^{p_k}
C_{j_k\ldots j_1}
\prod_{l=1}^k\phi_{j_l}(t_l)\Biggr)^2
dt_1
\ldots
dt_k=
$$

\vspace{3mm}
$$
=\int\limits_{[t,T]^k}
K^2(t_1,\ldots,t_k)
dt_1\ldots dt_k -
$$

\vspace{3mm}
$$
- 2\int\limits_{[t,T]^k}
K(t_1,\ldots,t_k)\sum_{j_1=0}^{p_1}\ldots
\sum_{j_k=0}^{p_k}
C_{j_k\ldots j_1}
\prod_{l=1}^k\phi_{j_l}(t_l)
dt_1\ldots dt_k+
$$

\vspace{3mm}
$$
+\int\limits_{[t,T]^k}
\Biggl(\sum_{j_1=0}^{p_1}\ldots
\sum_{j_k=0}^{p_k}
C_{j_k\ldots j_1}
\prod_{l=1}^k\phi_{j_l}(t_l)\Biggr)^2
dt_1\ldots dt_k=
$$

\vspace{3mm}
$$
=\int\limits_{[t,T]^k}
K^2(t_1,\ldots,t_k)
dt_1\ldots dt_k - 
$$

\vspace{3mm}
$$
- 2 \sum_{j_1=0}^{p_1}\ldots
\sum_{j_k=0}^{p_k}C_{j_k\ldots j_1}\int\limits_{[t,T]^k}
K(t_1,\ldots,t_k)
\prod_{l=1}^k\phi_{j_l}(t_l)
dt_1\ldots dt_k+
$$

\vspace{3mm}
$$
+
\sum_{j_1=0}^{p_1}\sum_{j_1'=0}^{p_1}\ldots
\sum_{j_k=0}^{p_k}\sum_{j_k'=0}^{p_k}
C_{j_k\ldots j_1}C_{j_k'\ldots j_1'}
\prod_{l=1}^k\int\limits_{\stackrel{~}{t}}^T\phi_{j_l}(t_l)
\phi_{j_l'}(t_l)dt_l=
$$

\vspace{3mm}
$$
=\int\limits_{[t,T]^k}
K^2(t_1,\ldots,t_k)
dt_1\ldots dt_k 
-2 \sum_{j_1=0}^{p_1}\ldots
\sum_{j_k=0}^{p_k}C^2_{j_k\ldots j_1}
+\sum_{j_1=0}^{p_1}\ldots
\sum\limits_{j_k=0}^{p_k}
C^2_{j_k\ldots j_1}=
$$

\vspace{3mm}
\begin{equation}
\label{dobav100}
=\int\limits_{[t,T]^k}
K^2(t_1,\ldots,t_k)
dt_1\ldots dt_k -\sum_{j_1=0}^{p_1}\ldots
\sum_{j_k=0}^{p_k}C^2_{j_k\ldots j_1}.
\end{equation}

\vspace{4mm}

Let us substitute (\ref{dobav100}) into (\ref{ch12026102})

\vspace{1mm}
$$
{\sf M}\left\{\left(
J'[R_{p_1\ldots p_k}]_{T,t}^{(i_1\ldots i_k)}\right)^{2}\right\}\le
$$

\vspace{1mm}
\begin{equation}
\label{ch12026103}
\le  k!
\left(~\int\limits_{[t,T]^k}
K^2(t_1,\ldots,t_k)
dt_1\ldots dt_k -\sum_{j_1=0}^{p_1}\ldots
\sum_{j_k=0}^{p_k}C^2_{j_k\ldots j_1}\right),
\end{equation}

\vspace{5mm}
\noindent
where $i_1,\ldots,i_k=1,\ldots,m.$

Due to Parseval's equality

$$
\int\limits_{[t, T]^k}
R^2_{p_1\ldots p_k}(t_1,\dots,t_k)dt_1\ldots dt_k=
$$

\vspace{1mm}
\begin{equation}
\label{ziko11000}
=\int\limits_{[t,T]^k}
K^2(t_1,\ldots,t_k)
dt_1\ldots dt_k -\sum_{j_1=0}^{p_1}\ldots
\sum_{j_k=0}^{p_k}C^2_{j_k\ldots j_1}\  \to \  0
\end{equation}

\vspace{4mm}
\noindent
if $p_1,\ldots,p_k\to\infty.$

Combining (\ref{2026ch1001}) and (\ref{ch12026103}), we get

\vspace{1mm}
$$
{\sf M}\left\{\left(J'[R_{p_1\ldots p_k}]_{T,t}^{(i_1\ldots i_k)}\right)^{2n}\right\}\le
$$

\vspace{4mm}
$$
\le
(k!)^{n}(2n-1)^{nk}\ \times
$$

\begin{equation}
\label{2026ch1001s}
\times\ 
\left(~
\int\limits_{[t,T]^k}
K^2(t_1,\ldots,t_k)
dt_1\ldots dt_k -\sum_{j_1=0}^{p_1}\ldots
\sum_{j_k=0}^{p_k}C^2_{j_k\ldots j_1}
\right)^n,
\end{equation}

\vspace{4mm}
\noindent
or
$$
{\sf M}\left\{\left(J[\psi^{(k)}]_{T,t}-
J[\psi^{(k)}]_{T,t}^{p_1,\ldots,p_k}\right)^{2n}\right\}\le
$$

\vspace{4mm}
$$
\le
(k!)^{n} (2n-1)^{nk}\ \times
$$

\begin{equation}
\label{dima2ye100}
\times\ 
\left(~
\int\limits_{[t,T]^k}
K^2(t_1,\ldots,t_k)
dt_1\ldots dt_k -\sum_{j_1=0}^{p_1}\ldots
\sum_{j_k=0}^{p_k}C^2_{j_k\ldots j_1}
\right)^n,
\end{equation}

\vspace{4mm}
\noindent
where $n\in{\mathbb N}$ and $i_1=\ldots=i_k=i\in\{1,\ldots,m\}.$

The inequality (\ref{2026ch1001s}) (or (\ref{dima2ye100}))
means that the expansion of 
iterated Ito stochastic integral obtained using Theorem 1
(the case $k\in{\mathbb N},$ $i_1=\ldots=i_k=i\in\{1,\ldots,m\}$)
converges in the 
mean
of degree $2n$ ($n\in {\mathbb N}$) to the appropriate 
iterated Ito stochastic integral.

Now we consider the case of a multidimensional Wiener process
and obtain an estimate of type (\ref{dima2ye100}) for the case $k=2,$ $i_1,i_2=1,\ldots,m.$

Suppose that $\{\phi_j(x)\}_{j=0}^{\infty}$
is a complete orthonormal system of continuous functions in 
the space $L_2([t, T])$ and $\psi_1(\tau),\psi_2(\tau)$
are continuous functions on $[t, T].$

Applying the Minkowski inequality and (\ref{jjjye}), we obtain 

$$
{\sf M}\left\{\left(J[\psi^{(2)}]_{T,t}-
J[\psi^{(2)}]_{T,t}^{p_1,p_2}\right)^{2n}\right\}=
{\sf M}\left\{\left(
J'[R_{p_1 p_2}]_{T,t}^{(i_1 i_2)}\right)^{2n}\right\}=
$$

\vspace{2mm}
$$
=
{\sf M}\left\{\left(
\sum_{(t_1,t_2)}
\int\limits_{t}^{T}
\int\limits_{t}^{t_2}
R_{p_1 p_2}(t_1,t_2)
d{\bf f}_{t_1}^{(i_1)}
d{\bf f}_{t_2}^{(i_2)}\right)^{2n}\right\}\le
$$

\vspace{2mm}
\begin{equation}
\label{ch12026e201}
\le
\left(\sum_{(t_1,t_2)}
\left({\sf M}\left\{\left(
\int\limits_{t}^{T}
\int\limits_{t}^{t_2}
R_{p_1 p_2}(t_1,t_2)
d{\bf f}_{t_1}^{(i_1)}
d{\bf f}_{t_2}^{(i_2)}\right)^{2n}\right\}\right)^{1/2n}\right)^{2n},
\end{equation}

\vspace{4mm}
\noindent
where $n\in{\mathbb N}.$

Let us evaluate 
$$
{\sf M}\left\{\left(
\int\limits_{t}^{T}
\int\limits_{t}^{t_2}
R_{p_1 p_2}(t_1,t_2)
d{\bf f}_{t_1}^{(i_1)}
d{\bf f}_{t_2}^{(i_2)}\right)^{2n}\right\}.
$$

\vspace{3mm}

Denote
\begin{equation}
\label{ch1w1004}
\eta_s=\int\limits_{t}^{s}
\int\limits_{t}^{t_2}
R_{p_1 p_2}(t_1,t_2)
d{\bf f}_{t_1}^{(i_1)}
d{\bf f}_{t_2}^{(i_2)},\ \ \ s\in[t, T].
\end{equation}

\vspace{3mm}

We have
$$
d\eta_s=\xi_s d{\bf f}_{s}^{(i_2)},
$$

\vspace{2mm}
\noindent
where
\begin{equation}
\label{ch12026ner11}
\xi_s=\int\limits_{t}^s
R_{p_1 p_2}(t_1,s)
d{\bf f}_{t_1}^{(i_1)}.
\end{equation}

\vspace{3mm}

Using the Ito formula it is easy to demonstrate that \cite{36} 

\vspace{1mm}
$$
{\sf M}\left\{(\eta_{\tau})^{2n}
\right\}=n(2n-1)
{\sf M}\left\{\int\limits_{t}^{\tau}
\left(\eta_s\right)^{2n-2}\xi_s^2 ds\right\}=
$$

\vspace{2mm}
\begin{equation}
\label{ch12026qqq1}
=
n(2n-1)
\int\limits_{t}^{\tau}
{\sf M}\left\{\left(\eta_s\right)^{2n-2}\xi_s^2\right\}ds.
\end{equation}

\vspace{4mm}

The last step in (\ref{ch12026qqq1}) is carried out
on the basis of a consequence from Fubini's Theorem \cite{Shir999},
since (as we will see later) 
\begin{equation}
\label{ch1ww2001}
{\sf M}\left\{\int\limits_{t}^{\tau}\left(\eta_s\right)^{2n-2}\xi_s^2 ds\right\}<\infty\ \ \
\hbox{for}\ \ \ p_1,p_2<\infty.
\end{equation}

\vspace{3mm}

Using the H\"{o}lder inequality (under the integral sign 
on the right-hand side of (\ref{ch12026qqq1})) for 
$p=n/(n-1)$, $q=n$ $(n>1)$ and using the non-decreasing property of 
the value ${\sf M}\left\{(\eta_{\tau})^{2n}\right\}$
with the 
growth of $\tau$ (see (\ref{ch12026qqq1})), we get

\vspace{-1mm}
$$
{\sf M}\left\{\left(\eta_{\tau}\right)^{2n}
\right\}\le
n(2n-1)
\left({\sf M}\left\{\left(\eta_{\tau}\right)^{2n}
\right\}\right)^{(n-1)/n}
\int\limits_{t}^{\tau}\left({\sf M}\left\{(\xi_s)^{2n}\right\}\right)^{1/n}ds.
$$

\vspace{2mm}

After raising to power $n$ the obtained inequality and 
dividing the result by 

$$
\left({\sf M}\left\{\left(\eta_{\tau}\right)^{2n}
\right\}\right)^{n-1},
$$

\vspace{1mm}
\noindent
we get the following estimate

\vspace{-1mm}
\begin{equation}
\label{neogidal}
{\sf M}\left\{\left(\eta_{\tau}\right)^{2n}
\right\}\le
(n(2n-1))^n
\left(
\int\limits_{t}^{\tau}\left({\sf M}\left\{(\xi_s)^{2n}\right\}\right)^{1/n}ds
\right)^n.
\end{equation}

\vspace{3mm}

Note that
\begin{equation}
\label{ch12026r67}
{\sf M}\left\{(\xi_s)^{2n}\right\}=
(2n-1)!!\left(
\int\limits_t^{s}
R^2_{p_1 p_2}(t_1,s)dt_1
\right)^n,
\end{equation} 

\vspace{3mm}
\noindent
since the randon variable $\xi_s$ has a Gaussian distribution
and

\vspace{-2mm}
\begin{equation}
\label{ch12026art34}
{\sf M}\left\{\xi_s^{2}\right\}={\sf M}\left\{\left(
\int\limits_{t}^s
R_{p_1 p_2}(t_1,s)
d{\bf f}_{t_1}^{(i_1)}\right)^2
\right\}=
\int\limits_{t}^s
R^2_{p_1 p_2}(t_1,s)
dt_1.
\end{equation}

\vspace{3mm}

Combining (\ref{neogidal}) and (\ref{ch12026r67}), we obtain

\vspace{-1mm}
\begin{equation}
\label{ch12026tttt1}
{\sf M}\left\{\left(\eta_{\tau}\right)^{2n}
\right\}\le
(n(2n-1))^n (2n-1)!!
\left(\int\limits_{t}^{\tau}
\int\limits_t^{s}
R^2_{p_1 p_2}(t_1,s)dt_1 ds
\right)^n.
\end{equation}

\vspace{3mm}

Then
\begin{equation}
\label{ch12026r67d}
{\sf M}\left\{\left(\eta_{T}\right)^{2n}
\right\}\le
(n(2n-1))^n (2n-1)!!
\left(~\int\limits_{[t, T]^2}
R^2_{p_1 p_2}(t_1,t_2)dt_1 dt_2
\right)^n.
\end{equation}

\vspace{3mm}

Finally, using (\ref{ch12026e201}) and (\ref{ch12026r67d}), we have

\vspace{1mm}
$$
{\sf M}\left\{\left(J[\psi^{(2)}]_{T,t}-
J[\psi^{(2)}]_{T,t}^{p_1 p_2}\right)^{2n}\right\}\le
$$

\vspace{2mm}
$$
\le
\left(\sum_{(t_1,t_2)}
\left({\sf M}\left\{\left(\eta_T\right)^{2n}
\right\}\right)^{1/2n}\right)^{2n}\le
$$

\vspace{2mm}
$$
\le
2^{2n}(n(2n-1))^n (2n-1)!! 
\left(~\int\limits_{[t, T]^2}
R^2_{p_1 p_2}(t_1,t_2)dt_1 dt_2
\right)^{n}=
$$

\vspace{4mm}
$$
=
2^{2n}(n(2n-1))^n (2n-1)!!\times
$$

\begin{equation}
\label{ch12026ri12}
\times
\left(~
\int\limits_{[t,T]^2}
K^2(t_1,t_2)
dt_1dt_2 -\sum_{j_1=0}^{p_1}
\sum_{j_2=0}^{p_2}C^2_{j_2 j_1}
\right)^{n},
\end{equation}

\vspace{4mm}
\noindent
where $n\in{\mathbb N}.$

Let us show that for $p_1,p_2<\infty$ the following inequality

\vspace{-3mm}
\begin{equation}
\label{ch120204567}
{\sf M}\left\{(\eta_{\tau})^{2n}\right\}\le C<\infty
\end{equation}

\vspace{2.5mm}
\noindent
is satisfied, where $\tau\in [t, T]$ and $C$ is a constant
that depends on $p_1,p_2,T,n.$

Consider the following well known estimate
for the moments of the Ito stochastic integral \cite{36}

\vspace{-2mm}
\begin{equation}
\label{ch12026pupol1}
{\sf M}\left\{\left|\int\limits_{t}^T \phi_\tau
dw_\tau\right|^{2n}\right\} \le (T-t)^{n-1}\left(n(2n-1)\right)^n
\int\limits_{t}^T {\sf M}\left\{(\phi_\tau)^{2n}\right\}d\tau,
\end{equation}

\vspace{2mm}
\noindent
where the process $\phi_{\tau}$ is such that
$\left(\phi_{\tau}\right)^n\in{\rm M}_2
([t,T])$ and $w_{\tau}$ is a scalar standard Wiener 
process,
$n=1, 2,\ldots$ (definition of the class 
${\rm M}_2([t,T])$ see in Sect.~2).

Applying (\ref{ch12026pupol1}), we obtain

\vspace{-3mm}
\begin{equation}
\label{ch12026pupol12}
{\sf M}\left\{(\eta_{\tau})^{2n}\right\} \le (T-t)^{n-1}\left(n(2n-1)\right)^n
\int\limits_{t}^{\tau} {\sf M}\left\{(\xi_s)^{2n}\right\}ds,
\end{equation}

\vspace{1mm}
\noindent
where $\xi_\tau$ is defined by (\ref{ch12026ner11}).

Combining (\ref{ch12026pupol12}) and (\ref{ch12026r67}), we get

\vspace{-1mm}
\begin{equation}
\label{ch12026pupol122}
{\sf M}\left\{(\eta_{\tau})^{2n}\right\} \le (T-t)^{n-1}\left(n(2n-1)\right)^n
(2n-1)!!
\int\limits_{t}^{\tau} \left(
\int\limits_t^{s}
R^2_{p_1 p_2}(t_1,s)dt_1
\right)^n ds.
\end{equation}

\vspace{3.5mm}

Under the conditions of Theorem~1, 
the integrals

$$
\int\limits_t^{s}
R^2_{p_1 p_2}(t_1,s)dt_1,\ \ \ \int\limits_{t}^{\tau} \left(
\int\limits_t^{s}
R^2_{p_1 p_2}(t_1,s)dt_1
\right)^n ds
$$

\vspace{3mm}
\noindent
are continuous functions with respect to $s$ and $\tau$, respectively.

Thus, the estimate (\ref{ch120204567}) is proved (see (\ref{ch12026pupol122})).
Then, the inequality (\ref{ch1ww2001}) holds (see (\ref{ch12026qqq1})).
This means that for $p_1,p_2<\infty$
we can apply the consequence from Fubini's Theorem \cite{Shir999}
in (\ref{ch12026qqq1}) and therefore for $p_1,p_2<\infty$ the estimate (\ref{ch12026ri12})
will be true. The proof of the estimate (\ref{ch12026ri12}) is completed.

Let us explain why this approach cannot be generalized
to the case $k\ge 3.$
Let $k=3.$ Now

$$
\eta_s=\int\limits_{t}^{s}
\int\limits_{t}^{t_3}
\int\limits_{t}^{t_2}
R_{p_1 p_2 p_3}(t_1,t_2,t_3)
d{\bf f}_{t_1}^{(i_1)}
d{\bf f}_{t_2}^{(i_2)}d{\bf f}_{t_3}^{(i_3)},\ \ \ s\in[t, T],
$$

\vspace{2mm}
\noindent
and
$$
d\eta_s=\xi_s
d{\bf f}_{s}^{(i_3)},
$$

\vspace{3mm}
\noindent
where
$$
\xi_s=\int\limits_{t}^{s}
\int\limits_{t}^{t_2}
R_{p_1 p_2 p_3}(t_1,t_2,s)
d{\bf f}_{t_1}^{(i_1)}
d{\bf f}_{t_2}^{(i_2)}.
$$

\vspace{3mm}

In the next step, we need that the stochastic differential
$d\xi_s$ to have the form

\vspace{-1mm}
\begin{equation}
\label{ch12020yyy1}
d\xi_s=\mu_s d{\bf f}_{s}^{(i_2)}.
\end{equation}

\vspace{4mm}

According to (\ref{ch12026ll12}), $\xi_s$ is a finite
linear combination of integrals of the form

$$
\rho_s=h(s)\int\limits_{t}^{s}g(t_2)
\int\limits_{t}^{t_2}
q(t_1)d{\bf f}_{t_1}^{(i_1)}
d{\bf f}_{t_2}^{(i_2)},
$$

\vspace{3mm}
\noindent
where $h(s), g(s), q(s)$ are some continuous functions on $[t, T].$

If we assume that the function $h(s)$ is continuously differentiable,
then according to the Ito formula the stochastic differential
$d\rho_s$ will have a non-zero drift coefficient.
This means that the stochastic differential $d\xi_s$
will have a more compex form than (\ref{ch12020yyy1})
(with a non-zero drift coefficient), which makes
it impossible to generalize this approach to the case $k=3.$

Let us generalize the estimate (\ref{dima2ye100}) to the case $k\in{\mathbb N},$
$i_1,\ldots,i_k=1,\ldots,m.$
Let $H_1(s),\ldots,H_m(s): [0,\infty)\to {\mathbb R}$.
Let $(\Omega,{\rm F},{\sf P})$ is a probability space, where
${\rm F}$ is the smallest $\sigma$-algebra such that
the random variables

\vspace{-2mm}
$$
\sum\limits_{i=1}^m\int\limits_0^{\infty} H_i(s)d{\bf f}_s^{(i)}
$$

\vspace{3mm}
\noindent
are ${\rm F}$-measurable for every $H_1(s),\ldots,H_m(s)\in L_2([0,\infty)),$
where 
${\bf f}_s^{(i)}$ are independent standard Wiener processes, $i=1,\ldots,m,$

\vspace{-1mm}
$$
\int\limits_0^{\infty} H_i(s)d{\bf w}_s^{(i)}
$$

\vspace{3mm}
\noindent
is a usual Wiener--Ito stochastic integral.

It is well known (see Theorem~2.4 and Lemma~2.6 in \cite{hairer1}) that

\vspace{1mm}
$$
L_2(\Omega,{\rm F},{\sf P})=\bigoplus\limits_{k=0}^{\infty}{\mathbb H}_k,
$$

\vspace{5mm}
\noindent
where ${\mathbb H}_0$ contains only constants,
the space ${\mathbb H}_k$ is the so-called $k$th $(k\ge 1)$ homogeneous Wiener chaos
which consists of all random variables of the form

\vspace{1mm}
$$
\sum\limits_{i_1,\ldots,i_k=1}^m 
\int\limits_0^{\infty}\int\limits_0^{t_k}\ldots \int\limits_0^{t_2}
H_{i_1\ldots i_k}(t_1,\ldots,t_k)d{\bf f}_{t_1}^{(i_1)}\ldots d{\bf f}_{t_k}^{(i_k)},
$$

\vspace{5.5mm}
\noindent
where $H_{i_1\ldots i_k}(t_1,\ldots,t_k)\in L_2(D_k),$
$D_k=\{(t_1,\ldots,t_k):\  0<t_1<\ldots <t_k\}.$

Let 

$$
H_{i_1\ldots i_k}(t_1,\ldots,t_k)=
\Phi(t_1,\ldots,t_k){\bf 1}_{\{0\le t<t_1<\ldots <t_k\le T\}}
{\bf 1}_{\{i_1=j_1,\ldots,i_k=j_k\}},
$$

\vspace{6.5mm}
\noindent
where ${\bf 1}_A$ denotes the indicator of the set $A,$
$\Phi(t_1,\ldots,t_k)\in L_2([t,T]^k),$
$j_1,\ldots,j_k$ are some fixed numbers from the set
$\{1,\ldots,m\}.$

Then, we have w.~p.~1

\vspace{1mm}
$$
\sum\limits_{i_1,\ldots,i_k=1}^m 
\int\limits_0^{\infty}\int\limits_0^{t_k}\ldots \int\limits_0^{t_2}
H_{i_1\ldots i_k}(t_1,\ldots,t_k)d{\bf f}_{t_1}^{(i_1)}\ldots {\bf f}_{t_k}^{(i_k)}=
$$

\vspace{3mm}
$$
=\int\limits_t^{T}\ldots \int\limits_t^{t_2}
\Phi(t_1,\ldots,t_k)d{\bf f}_{t_1}^{(j_1)}\ldots 
d{\bf f}_{t_k}^{(j_k)}\in {\mathbb H}_k.
$$

\vspace{5mm}

Obviously (see (\ref{pobeda}), (\ref{pobedaxyz}), (\ref{Wi110})),

\vspace{1mm}
$$
J'[\Phi]_{T,t}^{(i_1\ldots i_k)}=\sum_{(t_1,\ldots,t_k)}
\int\limits_{t}^{T}
\ldots
\int\limits_{t}^{t_2}
\Phi(t_1,\ldots,t_k)
d{\bf f}_{t_1}^{(i_1)}
\ldots
d{\bf f}_{t_k}^{(i_k)}=
$$

\vspace{3mm}
$$
=
\int\limits_t^T\ldots \int\limits_t^{t_2}
\sum\limits_{(t_1,\ldots,t_k)}\biggl(
\Phi(t_1,\ldots,t_k)
d{\bf f}_{t_1}^{(i_1)}\ldots
d{\bf f}_{t_k}^{(i_k)}\biggr)
\in {\mathbb H}_k,
$$

\vspace{5mm}
\noindent
where 
$J'[\Phi]_{T,t}^{(i_1\ldots i_k)}$ is a multiple
Wiener stochastic integral 
with respect to components of a multi\-di\-men\-si\-o\-nal 
Wiener process
defined as in Sect.~2 (see (\ref{mult11})) or as in Sect.~15
(see (\ref{WiI})),
$\Phi(t_1,\ldots,t_k)\in L_2([t, T]^k)$ in (\ref{WiI}) and  
$\Phi(t_1,\ldots,t_k)\in C([t, T]^k)\subset L_2([t, T]^k)$ in (\ref{mult11}),
$i_1,\ldots,i_k=1,\ldots,m.$

It is well known that (see \cite{hairer1})

$$
{\sf M}\left\{\xi^{2n}\right\}\le (2n-1)^{kn}\left({\sf M}\left\{\xi^2\right\}\right)^n,
$$

\vspace{3mm}
\noindent
where $\xi\in {\mathbb H}_k,$ $n, k\in {\mathbb N}.$

Then, we have the following estimate for 
$J'[\Phi]_{T,t}^{(i_1\ldots i_k)}\in {\mathbb H}_k$

\begin{equation}
\label{2026ch1001dd}
{\sf M}\left\{\left(
J'[\Phi]_{T,t}^{(i_1\ldots i_k)}\right)^{2n}\right\}
\le (2n-1)^{nk} \left({\sf M}\left\{\left(
J'[\Phi]_{T,t}^{(i_1\ldots i_k)}\right)^{2}\right\}\right)^{n},
\end{equation}

\vspace{3.5mm}
\noindent
where $n\in{\mathbb N}$ and $k\in{\mathbb N},$ $i_1,\ldots,i_k=1,\ldots,m.$

Combining (\ref{jjjye}), (\ref{ch12026103}), and (\ref{2026ch1001dd}), we get

\vspace{2mm}
$$
{\sf M}\left\{\left(J'[R_{p_1\ldots p_k}]_{T,t}^{(i_1\ldots i_k)}\right)^{2n}\right\}\le
$$

\vspace{4.5mm}
$$
\le
(k!)^{n}(2n-1)^{nk}\ \times
$$

\vspace{1mm}
$$
\times\ 
\left(~
\int\limits_{[t,T]^k}
K^2(t_1,\ldots,t_k)
dt_1\ldots dt_k -\sum_{j_1=0}^{p_1}\ldots
\sum_{j_k=0}^{p_k}C^2_{j_k\ldots j_1}
\right)^n,
$$

\vspace{6mm}
\noindent
or

\vspace{-1mm}
$$
{\sf M}\left\{\left(J[\psi^{(k)}]_{T,t}-
J[\psi^{(k)}]_{T,t}^{p_1,\ldots,p_k}\right)^{2n}\right\}\le
$$

\vspace{4.5mm}
$$
\le
(k!)^{n} (2n-1)^{nk}\ \times
$$

\vspace{1mm}
\begin{equation}
\label{2026ch1001s11}
\times\ 
\left(~
\int\limits_{[t,T]^k}
K^2(t_1,\ldots,t_k)
dt_1\ldots dt_k -\sum_{j_1=0}^{p_1}\ldots
\sum_{j_k=0}^{p_k}C^2_{j_k\ldots j_1}
\right)^n,
\end{equation}

\vspace{6mm}
\noindent
where $n\in{\mathbb N}$ and $k\in{\mathbb N},$ $i_1,\ldots,i_k=1,\ldots,m.$

The inequality (\ref{2026ch1001s11})
means that the expansion of 
iterated Ito stochastic integral obtained using Theorem 1
(the case $k\in{\mathbb N},$ $i_1,\ldots,i_k=1,\ldots,m$)
converges in the 
mean
of degree $2n$ ($n\in {\mathbb N}$) to the appropriate 
iterated Ito stochastic integral.

\vspace{5mm}

\section{Estimate for the Mean-Square Error of Approximation
of Iterated Ito Stochastic Integrals Based on Theorem 1}

\vspace{5mm}

In this section, we prove the useful estimate for the
mean-square error of approximation in Theorem 1.

{\bf Theorem 3} \cite{19}-\cite{2023xxx1}, \cite{26}. 
{\it Suppose that
every $\psi_l(\tau)$ $(l=1,\ldots, k)$ is a continuous nonrandom function on 
the interval $[t, T]$ and
$\{\phi_j(x)\}_{j=0}^{\infty}$ is a complete orthonormal system  
of functions in the space $L_2([t,T]),$ 
each function of which 
for finite $j$ satisfies the condition 
$(\star)$ {\rm (}see Sect. {\rm 4)}. 
Then
the estimate

$$
{\sf M}\left\{\left(
J[\psi^{(k)}]_{T,t}-J[\psi^{(k)}]_{T,t}^{p_1,\ldots,p_k}
\right)^2\right\}
\le 
$$

\vspace{2mm}
\begin{equation}
\label{z1}
\le k!\left(\int\limits_{[t,T]^k}
K^2(t_1,\ldots,t_k)
dt_1\ldots dt_k -\sum_{j_1=0}^{p_1}\ldots
\sum_{j_k=0}^{p_k}C^2_{j_k\ldots j_1}\right)
\end{equation}

\vspace{5mm}
\noindent
is valid for the following cases{\rm :}

\vspace{2mm}

{\rm 1.}\ $i_1,\ldots,i_k=1,\ldots,m$\ \ and\ \ $0<T-t<\infty,$

\vspace{1mm}

{\rm 2.}\ $i_1,\ldots,i_k=0, 1,\ldots,m,$\ \ $i_1^2+\ldots+i_k^2>0,$\ \
and\ \ $0<T-t<1,$

\vspace{3mm}
\noindent
where $J[\psi^{(k)}]_{T,t}$ is the stochastic integral {\rm (\ref{sodom20}),}
$J[\psi^{(k)}]_{T,t}^{p_1,\ldots,p_k}$ is the 
expression on the right-hand side of {\rm (\ref{tyyy})} before
passing to the limit 
$\hbox{\vtop{\offinterlineskip\halign{
\hfil#\hfil\cr
{\rm l.i.m.}\cr
$\stackrel{}{{}_{p_1,\ldots,p_k\to \infty}}$\cr
}} };$ another 
notations are the same as in Theorem {\rm 1}.
}

\vspace{5mm}

{\bf Proof.}\ Prooving Theorem 1,
we obtained w. p. 1 the 
following re\-pre\-sen\-ta\-ti\-on (see (\ref{zara}), (\ref{y007}))

\vspace{2mm}
$$
J[\psi^{(k)}]_{T,t}=J[\psi^{(k)}]_{T,t}^{p_1,\ldots,p_k}+
R_{T,t}^{p_1,\ldots,p_k},
$$

\vspace{5mm}
\noindent
where $J[\psi^{(k)}]_{T,t}^{p_1,\ldots,p_k}$
is the 
expression on the right-hand side of {\rm (\ref{tyyy})} before
passing to the limit 
$\hbox{\vtop{\offinterlineskip\halign{
\hfil#\hfil\cr
{\rm l.i.m.}\cr
$\stackrel{}{{}_{p_1,\ldots,p_k\to \infty}}$\cr
}} }$ and

\vspace{1mm}
$$
R_{T,t}^{p_1,\ldots,p_k}=
\sum_{(t_1,\ldots,t_k)}
\int\limits_{t}^{T}
\ldots
\int\limits_{t}^{t_2}
\Biggl(K(t_1,\ldots,t_k)-
\sum_{j_1=0}^{p_1}\ldots
\sum_{j_k=0}^{p_k}
C_{j_k\ldots j_1}
\prod_{l=1}^k\phi_{j_l}(t_l)\Biggr)\times
$$

\vspace{2mm}
\begin{equation}
\label{y007a}
\times d{\bf w}_{t_1}^{(i_1)}
\ldots
d{\bf w}_{t_k}^{(i_k)},
\end{equation}

\vspace{5mm}
\noindent
where
$$
\sum_{(t_1,\ldots,t_k)}
$$

\vspace{2mm}
\noindent
means the sum with respect to all
possible permutations
$(t_1,\ldots,t_k),$ which are
performed only 
in the values $d{\bf w}_{t_1}^{(i_1)}
\ldots $
$d{\bf w}_{t_k}^{(i_k)}$. At the same time the indexes near 
upper limits of integration in the iterated stochastic integrals 
are changed correspondently and if $t_r$ swapped with $t_q$ in the  
permutation $(t_1,\ldots,t_k)$, then $i_r$ swapped with $i_q$ in the 
permutation $(i_1,\ldots,i_k)$.

In the case of 
any fixed $k$ and numbers $i_1,\ldots,i_k=1,\ldots,m$
the integrals on the right-hand side of (\ref{y007a})
will be dependent in a stochastic sense. 
Let us estimate the second moment of $R_{T,t}^{p_1,\ldots,p_k}.$ 
From (\ref{99.010}), (\ref{y007a})
and elementary inequality

\begin{equation}
\label{y5}
\left(a_1+a_2+\ldots+a_p\right)^2 \le
p\left(a_1^2+a_2^2+\ldots+a_p^2\right),\ \ \ p\in \mathbb{N},
\end{equation}

\vspace{4mm}
\noindent
we obtain
the following estimate for the case $i_1,\ldots,i_k=1,\dots,m$
$(0<T-t<\infty)$

\vspace{3mm}
$$
{\sf M}\left\{\left(R_{T,t}^{p_1,\ldots,p_k}\right)^2\right\}
\le 
$$

\vspace{4mm}
$$
\le k!\left(
\sum_{(t_1,\ldots,t_k)}
\int\limits_{t}^{T}
\ldots
\int\limits_{t}^{t_2}
\Biggl(K(t_1,\ldots,t_k)-
\sum_{j_1=0}^{p_1}\ldots
\sum_{j_k=0}^{p_k}
C_{j_k\ldots j_1}
\prod_{l=1}^k\phi_{j_l}(t_l)\Biggr)^2
dt_1\ldots dt_k\right)=
$$

\vspace{4mm}
$$
=k!\int\limits_{[t,T]^k}
\Biggl(K(t_1,\ldots,t_k)-
\sum_{j_1=0}^{p_1}\ldots
\sum_{j_k=0}^{p_k}
C_{j_k\ldots j_1}
\prod_{l=1}^k\phi_{j_l}(t_l)\Biggr)^2
dt_1
\ldots
dt_k=
$$

\vspace{3mm}
\begin{equation}
\label{star00011}
= k!\left(\int\limits_{[t,T]^k}
K^2(t_1,\ldots,t_k)
dt_1\ldots dt_k -\sum_{j_1=0}^{p_1}\ldots
\sum_{j_k=0}^{p_k}C^2_{j_k\ldots j_1}\right).
\end{equation}

\vspace{6mm}

For the case of 
any fixed $k$ and numbers $i_1,\ldots,i_k=0, 1,\ldots,m$\ 
$(i_1^2+\ldots+i_k^2>0)$
from (\ref{99.010}),  (\ref{y007a}), (\ref{y5}) we obtain

\vspace{1.5mm}
$$
{\sf M}\left\{\left(R_{T,t}^{p_1,\ldots,p_k}\right)^2\right\}
\le 
$$

\vspace{4mm}
$$
\le 
C_k
\sum_{(t_1,\ldots,t_k)}
\int\limits_{t}^{T}
\ldots
\int\limits_{t}^{t_2}
\Biggl(K(t_1,\ldots,t_k)-
\sum_{j_1=0}^{p_1}\ldots
\sum_{j_k=0}^{p_k}
C_{j_k\ldots j_1}
\prod_{l=1}^k\phi_{j_l}(t_l)\Biggr)^2
dt_1
\ldots
dt_k=
$$

\vspace{4mm}
$$
=C_k\int\limits_{[t,T]^k}
\Biggl(K(t_1,\ldots,t_k)-
\sum_{j_1=0}^{p_1}\ldots
\sum_{j_k=0}^{p_k}
C_{j_k\ldots j_1}
\prod_{l=1}^k\phi_{j_l}(t_l)\Biggr)^2
dt_1
\ldots
dt_k=
$$

\vspace{3mm}
$$
=C_k\left(\int\limits_{[t,T]^k}
K^2(t_1,\ldots,t_k)
dt_1\ldots dt_k -\sum_{j_1=0}^{p_1}\ldots
\sum_{j_k=0}^{p_k}C^2_{j_k\ldots j_1}\right),
$$

\vspace{6.5mm}
\noindent
where $C_k$ is a constant.

It is not difficult to see that the constant $C_k$ depends on 
$k$ ($k$ is the multiplicity 
of the iterated Ito stochastic integral) and $T-t$ ($T-t$ is the length
of integration interval of the iterated Ito stochastic integral).
Moreover, $C_k$ has the following form

$$
C_k=k!\cdot{\rm max}\Bigl\{
(T-t)^{\alpha_1},\ (T-t)^{\alpha_2},\ \ldots,\ (T-t)^{\alpha_{k!}}
\Bigr\},
$$

\vspace{3mm}
\noindent
where $\alpha_1, \alpha_2, \ldots, \alpha_{k!}=0, 1,\ldots, k-1.$

However, $T-t$ is supposed as an integration step of
numerical procedures 
for Ito stochastic differential equations, which is 
a rather small value. For example $0<T-t<1.$ Then $C_k\le k!$

It means, that for the case
of any fixed $k$ and 
$i_1,\ldots,i_k=0, 1,\ldots,m,$  $i_1^2+\ldots+i_k^2>0$
$(0<T-t<1)$
we can
write (\ref{z1}). Theorem 3 is proved.

\vspace{5mm}

\section{Expansion of Iterated Ito Stochastic Integrals 
Based on Generalized Multiple Fourier Series.
The Case of Complete Orthonormal With We\-ight $r(t_1)\ldots r(t_k)\ge 0$  
Systems of Functions 
in the Space $L_2([t, T]^k)$}

\vspace{5mm}

In this section, we consider the modification of Theorem 1 for 
the case of complete orthonormal with weight $r(t_1)\ldots r(t_k)\ge 0$ 
systems of functions 
in the space $L_2([t, T]^k)$ ($k\in\mathbb{N}$).

Let $\{\Psi_j(x)\}_{j=0}^{\infty}$ be a complete orthonormal 
with weight $r(x)\ge 0$ 
system of functions in the space $L_2([t, T]).$ It is well known that the
Fourier
series with respect to the system $\{\Psi_j(x)\}_{j=0}^{\infty}$
of function 

$$
f(x)\ \ \ \left(f(x)\sqrt{r(x)}\in L_2([t, T])\right)
$$

\vspace{3mm}
\noindent
converges 
to the function $f(x)$ in the
mean-square sense with weight $r(x),$ i.e.

\begin{equation}
\label{g1}
\lim\limits_{p\to\infty}
\int\limits_t^T\biggl(f(x)-\sum\limits_{j=0}^p 
{\tilde C}_j \Psi_j(x)\biggr)^2 r(x)dx = 0,
\end{equation}

\vspace{3mm}
\noindent
where
\begin{equation}
\label{h1}
{\tilde C}_j=\int\limits_t^T f(x)\Psi_j(x)r(x)dx
\end{equation}

\vspace{2mm}
\noindent
is the Fourier coefficient.

Obviously, the relation (\ref{g1}) can be obtained if we will 
expand the function
$f(x)\sqrt{r(x)}\in L_2([t, T])$ into a usual Fourier series with respect
to the complete orthonormal with weight $1$ system of functions

$$
\left\{\Psi_j(x)\sqrt{r(x)}\right\}_{j=0}^{\infty}
$$ 

\vspace{2mm}
\noindent
in 
the space $L_2([t, T]).$ Then

\vspace{1mm}
$$
\lim\limits_{p\to\infty}
\int\limits_t^T\biggl(f(x)\sqrt{r(x)}-\sum\limits_{j=0}^p {\tilde C}_j 
\Psi_j(x)\sqrt{r(x)}\biggr)^2dx = 
$$

\vspace{1mm}
\begin{equation}
\label{g2}
=\lim\limits_{p\to\infty}
\int\limits_t^T\biggl(f(x)-
\sum\limits_{j=0}^p {\tilde C}_j \Psi_j(x)\biggr)^2 r(x)dx = 0,
\end{equation}

\vspace{4mm}
\noindent
where ${\tilde C}_j$ has the form (\ref{h1}).

Let us consider an obvious generalization of this approach 
to the case of several
variables.
Let us expand the function $K(t_1,\ldots,t_k)$ such that

$$
K(t_1,\ldots,t_k)\prod\limits_{l=1}^k \sqrt{r(t_l)}\in L_2([t, T]^k)
$$

\vspace{3mm}
\noindent
using the complete orthonormal system of functions 
$$
\prod\limits_{l=1}^k \Psi_{j_l}(t_l)\sqrt{r(t_l)},\ \ \ 
j_l=0, 1, 2, \ldots,\ \ \  l=1,\ldots,k
$$

\vspace{2mm}
\noindent
in the space $L_2([t, T]^k)$ into the generalized multiple Fourier 
series. 

It is well known that the mentioned
generalized multiple Fourier series converges in the mean-square sense,
i.e.

$$
\lim\limits_{p_1,\ldots,p_k\to\infty}
\int\limits_{[t,T]^k}
\left(K(t_1,\ldots,t_k)\prod\limits_{l=1}^k \sqrt{r(t_l)}-
\sum\limits_{j_1=0}^{p_1}\ldots\sum\limits_{j_k=0}^{p_k}
{\tilde C}_{j_k\ldots j_1}
\prod\limits_{l=1}^k \Psi_{j_l}(t_l)\sqrt{r(t_l)}\right)^2
dt_1\ldots dt_k=
$$

\vspace{2mm}
\begin{equation}
\label{z1aaa}
=\lim\limits_{p_1,\ldots,p_k\to\infty}
\int\limits_{[t,T]^k}
\left(K(t_1,\ldots,t_k)-
\sum\limits_{j_1=0}^{p_1}\ldots\sum\limits_{j_k=0}^{p_k}
{\tilde C}_{j_k\ldots j_1}\prod\limits_{l=1}^k \Psi_{j_l}(t_l)\right)^2 
\left(\prod\limits_{l=1}^k r(t_l)\right)
dt_1\ldots dt_k=0,
\end{equation}

\vspace{4mm}
\noindent
where

\vspace{-2mm}
$$
{\tilde C}_{j_k\ldots j_1}=\int\limits_{[t,T]^k}
K(t_1,\ldots,t_k)\prod\limits_{l=1}^k 
\biggl(\Psi_{j_l}(t_l)r(t_l)\biggr)dt_1\ldots dt_k.
$$

\vspace{5mm}

Let us consider 
the following iterated Ito 
stochastic integrals

\begin{equation}
\label{ito-rr}
{\tilde J}[\psi^{(k)}]_{T,t}=\int\limits_t^T\psi_k(t_k)\sqrt{r(t_k)} 
\ldots \int\limits_t^{t_{2}}
\psi_1(t_1)\sqrt{r(t_1)} d{\bf w}_{t_1}^{(i_1)}\ldots
d{\bf w}_{t_k}^{(i_k)},
\end{equation}

\vspace{3mm}
\noindent
where every $\psi_l(\tau)$ $(l=1,\ldots,k)$ is 
a nonrandom function on $[t, T]$,
${\bf w}_{\tau}^{(i)}={\bf f}_{\tau}^{(i)}$
for $i=1,\ldots,m$ and
${\bf w}_{\tau}^{(0)}=\tau,$ 
$i_1,\ldots,i_k=0, 1,\ldots,m.$

So, we obtain the following version of Theorem 1.

\vspace{3mm}

{\bf Theorem 4}\ \cite{20a}-\cite{2023xxx1} (also see \cite{20}, \cite{27a}).
{\it Suppose that
every $\psi_l(\tau)$ $(l=$ $1,\ldots, k)$ is a continuous 
nonrandom function on 
the interval $[t, T]$. Moreover, let
$\{\Psi_j(x)\sqrt{r(x)}\}_{j=0}^{\infty}$ $(r(x)\ge 0)$ 
is a complete orthonormal 
system of functions in the space $L_2([t,T])$, each function 
$\Psi_j(x)\sqrt{r(x)}$ of which 
for finite $j$ satisfies the condition 
$(\star)$ {\rm (}see Sect. {\rm 4)}.
Then

\vspace{2mm}
$$
{\tilde J}[\psi^{(k)}]_{T,t} =
\hbox{\vtop{\offinterlineskip\halign{
\hfil#\hfil\cr
{\rm l.i.m.}\cr
$\stackrel{}{{}_{p_1,\ldots,p_k\to \infty}}$\cr
}} }\sum_{j_1=0}^{p_1}\ldots\sum_{j_k=0}^{p_k}
{\tilde C}_{j_k\ldots j_1}\Biggl(
\prod_{l=1}^k{\tilde \zeta}_{j_l}^{(i_l)} -
\Biggr.
$$

\vspace{3mm}
\begin{equation}
\label{tyyy-rr}
-\Biggl.
\hbox{\vtop{\offinterlineskip\halign{
\hfil#\hfil\cr
{\rm l.i.m.}\cr
$\stackrel{}{{}_{N\to \infty}}$\cr
}} }\sum_{(l_1,\ldots,l_k)\in {\rm G}_k}
\Psi_{j_{1}}(\tau_{l_1})\sqrt{r(\tau_{l_1})}
\Delta{\bf w}_{\tau_{l_1}}^{(i_1)}\ldots
\Psi_{j_{k}}(\tau_{l_k})\sqrt{r(\tau_{l_k})}
\Delta{\bf w}_{\tau_{l_k}}^{(i_k)}\Biggr),
\end{equation}

\vspace{6mm}
\noindent
where

$$
{\rm G}_k={\rm H}_k\backslash{\rm L}_k,\ \ \
{\rm H}_k=\left\{(l_1,\ldots,l_k):\ l_1,\ldots,l_k=0,\ 1,\ldots,N-1\right\},
$$

$$
{\rm L}_k=\left\{(l_1,\ldots,l_k):\ l_1,\ldots,l_k=0,\ 1,\ldots,N-1;\
l_g\ne l_r\ (g\ne r);\ g, r=1,\ldots,k\right\},
$$

\vspace{4mm}
\noindent
${\rm l.i.m.}$ is a limit in the mean-square sense,
$i_1,\ldots,i_k=0,1,\ldots,m,$ 

$$
{\tilde \zeta}_{j}^{(i)}=
\int\limits_t^T \Psi_{j}(s)\sqrt{r(s)}d{\bf w}_s^{(i)}
$$

\vspace{3mm}
\noindent
are independent standard Gaussian random variables
for various
$i$ or $j$ {\rm(}in the case when $i\ne 0${\rm),}
$\Delta{\bf w}_{\tau_{j}}^{(i)}=
{\bf w}_{\tau_{j+1}}^{(i)}-{\bf w}_{\tau_{j}}^{(i)}$
$(i=0, 1,\ldots,m),$
$\left\{\tau_{j}\right\}_{j=0}^{N}$ is a partition of
the interval $[t,T],$ which satisfies the condition {\rm (\ref{1111})},

$$
{\tilde C}_{j_k\ldots j_1}=\int\limits_{[t,T]^k}
K(t_1,\ldots,t_k)
\prod_{l=1}^{k}\biggl(\Psi_{j_l}(t_l)r(t_l)\biggr)dt_1\ldots dt_k
$$

\vspace{3mm}
is the Fourier coefficient,

$$
K(t_1,\ldots,t_k)=
\left\{
\begin{matrix}
\psi_1(t_1)\ldots \psi_k(t_k),\ &t_1<\ldots<t_k\cr\cr\cr
0,\ &\hbox{otherwise}
\end{matrix}\right.,\ \  t_1,\ldots,t_k\in[t, T],\ \  k\ge 2,
$$

\vspace{5mm}
\noindent
and 
$K(t_1)\equiv\psi_1(t_1)$ for $t_1\in[t, T].$  
}

\vspace{2mm}

{\bf Proof.}
According to Lemmas 1, 3 and (\ref{pobeda}), (\ref{hehe100}),
we get the following representation w. p. 1

\vspace{2mm}
$$
{\tilde J}[\psi^{(k)}]_{T,t}=
\sum_{(t_1,\ldots,t_k)}
\int\limits_{t}^{T}
\ldots
\int\limits_{t}^{t_2}
K(t_1,\ldots,t_k)\prod\limits_{l=1}^k \sqrt{r(t_l)}d{\bf w}_{t_1}^{(i_1)}
\ldots
d{\bf w}_{t_k}^{(i_k)}=
$$

\vspace{7mm}
$$
=
\sum_{j_1=0}^{p_1}\ldots
\sum_{j_k=0}^{p_k}
{\tilde C}_{j_k\ldots j_1}
\int\limits_{t}^{T}
\ldots
\int\limits_{t}^{t_2}
\sum_{(t_1,\ldots,t_k)}\left(
\prod\limits_{l=1}^k\left(\Psi_{j_l}(t_l)\sqrt{r(t_l)}\right)
d{\bf w}_{t_1}^{(i_1)}
\ldots
d{\bf w}_{t_k}^{(i_k)}\right)
{\tilde R}_{T,t}^{p_1,\ldots,p_k}=
$$

\vspace{7mm}
$$
=\sum_{j_1=0}^{p_1}\ldots
\sum_{j_k=0}^{p_k}
{\tilde C}_{j_k\ldots j_1}
\hbox{\vtop{\offinterlineskip\halign{
\hfil#\hfil\cr
{\rm l.i.m.}\cr
$\stackrel{}{{}_{N\to \infty}}$\cr
}} }
\sum\limits_{\stackrel{l_1,\ldots,l_k=0}{{}_{l_q\ne l_r;\ 
q\ne r;\ q, r=1,\ldots, k}}}^{N-1}
\Psi_{j_1}(\tau_{l_1})\sqrt{r(\tau_{l_1})}\Delta{\bf w}_{\tau_{l_1}}^{(i_1)}
\ldots
\Psi_{j_k}(\tau_{l_k})\sqrt{(\tau_{l_k})}
\Delta{\bf w}_{\tau_{l_k}}^{(i_k)}+
$$

\vspace{2mm}
$$
+{\tilde R}_{T,t}^{p_1,\ldots,p_k}=
$$

\vspace{7mm}
$$
=\sum_{j_1=0}^{p_1}\ldots
\sum_{j_k=0}^{p_k}
{\tilde C}_{j_k\ldots j_1}
\left(
\hbox{\vtop{\offinterlineskip\halign{
\hfil#\hfil\cr
{\rm l.i.m.}\cr
$\stackrel{}{{}_{N\to \infty}}$\cr
}} }\sum_{l_1,\ldots,l_k=0}^{N-1}
\Psi_{j_1}(\tau_{l_1})\sqrt{r(\tau_{l_1})}\Delta{\bf w}_{\tau_{l_1}}^{(i_1)}
\ldots
\Psi_{j_k}(\tau_{l_k})\sqrt{(\tau_{l_k})}
\Delta{\bf w}_{\tau_{l_k}}^{(i_k)}
-\right.
$$

\vspace{2mm}
$$
-\left.
\hbox{\vtop{\offinterlineskip\halign{
\hfil#\hfil\cr
{\rm l.i.m.}\cr
$\stackrel{}{{}_{N\to \infty}}$\cr
}} }\sum_{(l_1,\ldots,l_k)\in {\rm G}_k}
\Psi_{j_1}(\tau_{l_1})\sqrt{r(\tau_{l_1})}\Delta{\bf w}_{\tau_{l_1}}^{(i_1)}
\ldots
\Psi_{j_k}(\tau_{l_k})\sqrt{(\tau_{l_k})}
\Delta{\bf w}_{\tau_{l_k}}^{(i_k)}
\right)
+{\tilde R}_{T,t}^{p_1,\ldots,p_k}=
$$

\vspace{7mm}
$$
=\sum_{j_1=0}^{p_1}\ldots\sum_{j_k=0}^{p_k}
{\tilde C}_{j_k\ldots j_1}\times
$$

\vspace{2mm}

\begin{equation}
\label{novoe2}
\times
\left(
\prod_{l=1}^k {\tilde \zeta}_{j_l}^{(i_l)}-
\hbox{\vtop{\offinterlineskip\halign{
\hfil#\hfil\cr
{\rm l.i.m.}\cr
$\stackrel{}{{}_{N\to \infty}}$\cr
}} }\sum_{(l_1,\ldots,l_k)\in {\rm G}_k}
\Psi_{j_1}(\tau_{l_1})\sqrt{r(\tau_{l_1})}\Delta{\bf w}_{\tau_{l_1}}^{(i_1)}
\ldots
\Psi_{j_k}(\tau_{l_k})\sqrt{(\tau_{l_k})}
\Delta{\bf w}_{\tau_{l_k}}^{(i_k)}
\right)+
{\tilde R}_{T,t}^{p_1,\ldots,p_k},
\end{equation}

\vspace{7mm}
\noindent
where

\vspace{-4mm}
$$
{\tilde R}_{T,t}^{p_1,\ldots,p_k}
=\sum_{(t_1,\ldots,t_k)}
\int\limits_{t}^{T}
\ldots
\int\limits_{t}^{t_2}
\left(K(t_1,\ldots,t_k)\prod_{l=1}^k\sqrt{r(t_l)}-\right.
$$

\vspace{2mm}
\begin{equation}
\label{y007ggg}
\left.
-\sum_{j_1=0}^{p_1}\ldots
\sum_{j_k=0}^{p_k}
{\tilde C}_{j_k\ldots j_1}
\prod_{l=1}^k\left(\Psi_{j_l}(t_l)\sqrt{r(t_l)}\right)\right)
d{\bf w}_{t_1}^{(i_1)}
\ldots
d{\bf w}_{t_k}^{(i_k)},
\end{equation}

\vspace{5mm}
\noindent
where permutations $(t_1,\ldots,t_k)$ when summing are performed only 
in the values $d{\bf w}_{t_1}^{(i_1)}
\ldots $
$d{\bf w}_{t_k}^{(i_k)}$. At the same time the indexes near 
upper limits of integration in the iterated stochastic integrals 
are changed correspondently and if $t_r$ swapped with $t_q$ in the  
permutation $(t_1,\ldots,t_k)$, then $i_r$ swapped with $i_q$ in the 
permutation $(i_1,\ldots,i_k)$.

Let us evaluate the remainder
${\tilde R}_{T,t}^{p_1,\ldots,p_k}$ of the series.

According to Lemma 2 and (\ref{riemann}), we have

\vspace{1mm}
$$
{\sf M}\left\{\left({\tilde R}_{T,t}^{p_1,\ldots,p_k}\right)^2\right\}
\le C_k
\sum_{(t_1,\ldots,t_k)}
\int\limits_{t}^{T}
\ldots
\int\limits_{t}^{t_2}
\left(K(t_1,\ldots,t_k)\prod_{l=1}^k\sqrt{r(t_l)}\right.-
$$

\vspace{2mm}
$$
\left.-
\sum_{j_1=0}^{p_1}\ldots
\sum_{j_k=0}^{p_k}
{\tilde C}_{j_k\ldots j_1}
\prod_{l=1}^k\left(\Psi_{j_l}(t_l)\sqrt{r(t_l)}\right)\right)^2
dt_1
\ldots
dt_k=
$$

\vspace{4mm}
\begin{equation}
\label{obana1ggg}
=C_k\int\limits_{[t,T]^k}
\left(K(t_1,\ldots,t_k)-
\sum_{j_1=0}^{p_1}\ldots
\sum_{j_k=0}^{p_k}
{\tilde C}_{j_k\ldots j_1}
\prod_{l=1}^k\Psi_{j_l}(t_l)\right)^2
\left(\prod_{l=1}^k r(t_l)\right)
dt_1 \ldots
dt_k\to 0
\end{equation}

\vspace{6mm}
\noindent
if $p_1,\ldots,p_k\to\infty,$ where constant $C_k$ 
depends only
on the multiplicity $k$ of the iterated Ito stochastic integral
(\ref{ito-rr}). 
Theorem 4 is proved.

Let us formulate the version of Theorem 3.

{\bf Theorem 5}\ \cite{20a}-\cite{2023xxx1}, \cite{27a}. 
{\it Suppose that every $\psi_l(\tau)$ $(l=$ $1,\ldots, k)$ 
is a continuous 
nonrandom function on 
$[t, T].$ Moreover, let 
$\{\Psi_j(x)\sqrt{r(x)}\}_{j=0}^{\infty}$ $(r(x)\ge 0)$
is a complete orthonormal 
system of functions in the space $L_2([t,T])$, each function 
$\Psi_j(x)\sqrt{r(x)}$
of which 
for finite $j$ satisfies the condition 
$(\star)$ {\rm (}see Sect. {\rm 4)}.
Then the estimate

\vspace{1mm}
$$
{\sf M}\left\{\left(
{\tilde J}[\psi^{(k)}]_{T,t}-{\tilde J}[\psi^{(k)}]_{T,t}^{p_1,\ldots,p_k}
\right)^2\right\}
\le 
$$

\vspace{2mm}
\begin{equation}
\label{z1-ura}
\le k!\left(\int\limits_{[t,T]^k}
K^2(t_1,\ldots,t_k)\left(\prod_{l=1}^k r(t_l)\right)
dt_1\ldots dt_k -\sum_{j_1=0}^{p_1}\ldots
\sum_{j_k=0}^{p_k}{\tilde C}^2_{j_k\ldots j_1}\right)
\end{equation}

\vspace{5mm}
\noindent
is valid for the following cases{\rm :}

\vspace{2mm}

{\rm 1.}\ $i_1,\ldots,i_k=1,\ldots,m$\ \ and\ \ $0<T-t<\infty,$

\vspace{1mm}

{\rm 2.}\ $i_1,\ldots,i_k=0, 1,\ldots,m,$\ \ $i_1^2+\ldots+i_k^2>0,$\ \
and\ \ $0<T-t<1,$

\vspace{3mm}
\noindent
where ${\tilde J}[\psi^{(k)}]_{T,t}$ is the 
stochastic integral {\rm (\ref{ito-rr}),}
${\tilde J}[\psi^{(k)}]_{T,t}^{p_1,\ldots,p_k}$ is the 
expression on the right-hand side of {\rm (\ref{tyyy-rr})} before
passing to the limit 
$\hbox{\vtop{\offinterlineskip\halign{
\hfil#\hfil\cr
{\rm l.i.m.}\cr
$\stackrel{}{{}_{p_1,\ldots,p_k\to \infty}}$\cr
}} };$ another 
notations are the same as in Theorem {\rm 4}.
}

\vspace{5mm}

\section{Convergence With Probability 1
of Expansion of Iterated
Ito Stochastic Integrals in Theorem 1 for the Case of Multiplicity
$k$ $(k\in\mathbb{N})$}

\vspace{5mm}

In this section, we formulate and prove the theorem on 
convergence with probability 1 (w. p. 1) of expansions 
of iterated Ito stochastic integrals in Theorem 1
for the case of multiplicity $k$ $(k\in{\mathbb N})$.
This section is written on the base of Sect.~1.7.2 
from \cite{20a}.

Let us remind the well-known fact from the mathematical analysis,
which is connected to existence
of iterated limits.

\vspace{2mm}

{\bf Proposition 1.} {\it Let $\bigl\{x_{n,m}\bigr\}_{n,m=1}^{\infty}$
be a double sequence and let there exists the limit

\vspace{-1mm}
$$
\lim\limits_{n,m\to\infty}x_{n,m}=a<\infty.
$$

\vspace{2mm}

Moreover, let there exist the limits

\vspace{-1mm}
$$
\lim\limits_{n\to\infty}x_{n,m}<\infty\ \ \hbox{for any}\ \
m,\ \ \ 
\lim\limits_{m\to\infty}x_{n,m}<\infty\ \ \hbox{for any}\ \ n.
$$

\vspace{2mm}

Then, there exist the iterated limits

\vspace{-1mm}
$$
\lim\limits_{n\to\infty}\lim\limits_{m\to\infty}x_{n,m},\ \ \ 
\lim\limits_{m\to\infty}\lim\limits_{n\to\infty}x_{n,m}
$$

\vspace{2mm}
and moreover,

\vspace{-1mm}
$$
\lim\limits_{n\to\infty}\lim\limits_{m\to\infty}x_{n,m}=
\lim\limits_{m\to\infty}\lim\limits_{n\to\infty}x_{n,m}=a.
$$
}

\vspace{3mm}

{\bf Theorem 6} \cite{20a}-\cite{2023xxx1}, \cite{22}.
{\it Let 
$\psi_l(\tau)$ $(l=1,\ldots, k)$ are 
continuously differentiable nonrandom functions on the interval
$[t, T]$ and $\{\phi_j(x)\}_{j=0}^{\infty}$ is a complete
orthonormal system of Legendre poly\-no\-mi\-als or 
trigonometric functions in the space $L_2([t, T]).$
Then 

\vspace{2mm}
$$
J[\psi^{(k)}]_{T,t}^{p,\ldots,p}\ \to \ J[\psi^{(k)}]_{T,t}\ \ \ 
\hbox{if}\ \ \ p\to \infty
$$

\vspace{5mm}
\noindent
w.\ p.\ {\rm 1,} where $J[\psi^{(k)}]_{T,t}^{p,\ldots,p}$
is the expression on the right-hand side of {\rm (\ref{tyyy})}
before passing to the limit 
$\hbox{\vtop{\offinterlineskip\halign{
\hfil#\hfil\cr
{\rm l.i.m.}\cr
$\stackrel{}{{}_{p_1,\ldots,p_k\to \infty}}$\cr
}} }$ for the case $p_1=\ldots=p_k=p,$ i.e. {\rm (}see Theorem {\rm 1)}

\vspace{2mm}
$$
J[\psi^{(k)}]_{T,t}^{p,\ldots,p}=
\sum_{j_1=0}^{p}\ldots\sum_{j_k=0}^{p}
C_{j_k\ldots j_1}\Biggl(
\prod_{l=1}^k\zeta_{j_l}^{(i_l)}\ -
\Biggr.
$$

\vspace{3mm}
$$
-\ \Biggl.
\hbox{\vtop{\offinterlineskip\halign{
\hfil#\hfil\cr
{\rm l.i.m.}\cr
$\stackrel{}{{}_{N\to \infty}}$\cr
}} }\sum_{(l_1,\ldots,l_k)\in {\rm G}_k}
\phi_{j_{1}}(\tau_{l_1})
\Delta{\bf w}_{\tau_{l_1}}^{(i_1)}\ldots
\phi_{j_{k}}(\tau_{l_k})
\Delta{\bf w}_{\tau_{l_k}}^{(i_k)}\Biggr),
$$

\vspace{6mm}
\noindent
where $i_1,\ldots,i_k=1,\ldots,m$.}

\vspace{2mm}

{\bf Proof.} Let us consider the Parseval equality

\begin{equation}
\label{par1}
\int\limits_{[t,T]^k}K^2(t_1,\ldots,t_k)dt_1\ldots dt_k=
\lim\limits_{p_1,\ldots,p_k\to\infty}
\sum_{j_1=0}^{p_1}\ldots \sum_{j_k=0}^{p_k}
C_{j_k\ldots j_1}^2,
\end{equation}

\vspace{3mm}
\noindent
where
\begin{equation}
\label{pppx}
K(t_1,\ldots,t_k)=
\begin{cases}
\psi_1(t_1)\ldots \psi_k(t_k),\ &t_1<\ldots<t_k\\
~\\
~\\
0,\ &\hbox{\rm otherwise}
\end{cases}\ \ \ \ 
=\ \ \ \ 
\prod\limits_{l=1}^k
\psi_l(t_l)\ \prod\limits_{l=1}^{k-1}{\bf 1}_{\{t_l<t_{l+1}\}},\ 
\end{equation}

\vspace{3mm}
\noindent
where $t_1,\ldots,t_k\in [t, T]$ for $k\ge 2$ and 
$K(t_1)\equiv\psi_1(t_1)$ for $t_1\in[t, T],$ 
${\bf 1}_A$ denotes the indicator of the set $A$,
\begin{equation}
\label{ppppax}
C_{j_k\ldots j_1}=\int\limits_{[t,T]^k}
K(t_1,\ldots,t_k)\prod_{l=1}^{k}\phi_{j_l}(t_l)dt_1\ldots dt_k
\end{equation}

\vspace{2mm}
\noindent
is the Fourier coefficient.

Using (\ref{pppx}), we obtain
$$
C_{j_k\ldots j_1}=
\int\limits_t^T
\phi_{j_k}(t_k)\psi_k(t_k)\ldots \int\limits_t^{t_2}
\phi_{j_1}(t_1)\psi_1(t_1)dt_1\ldots dt_k.
$$

\vspace{2mm}

Further, we denote

\vspace{-1mm}
$$
\lim\limits_{p_1,\ldots,p_k\to\infty}
\sum_{j_1=0}^{p_1}\ldots \sum_{j_k=0}^{p_k}
C_{j_k\ldots j_1}^2\stackrel{\sf def}{=}
\sum_{j_1,\ldots,j_k=0}^{\infty}
C_{j_k\ldots j_1}^2.
$$

\vspace{4mm}

If $p_1=\ldots=p_k=p,$ then we also write

\vspace{1mm}
$$
\lim\limits_{p\to\infty}
\sum_{j_1=0}^{p}\ldots \sum_{j_k=0}^{p}
C_{j_k\ldots j_1}^2\stackrel{\sf def}{=}
\sum_{j_1,\ldots,j_k=0}^{\infty}
C_{j_k\ldots j_1}^2.
$$

\vspace{4mm}

From the other hand, for iterated limits we write

\vspace{1mm}
$$
\lim\limits_{p_1\to\infty}\ldots \lim\limits_{p_k\to\infty}
\sum_{j_1=0}^{p_1}\ldots \sum_{j_k=0}^{p_k}
C_{j_k\ldots j_1}^2\stackrel{\sf def}{=}
\sum_{j_1=0}^{\infty}\ldots
\sum_{j_k=0}^{\infty}
C_{j_k\ldots j_1}^2,
$$

\vspace{3mm}
$$
\lim\limits_{p_1\to\infty}\lim\limits_{p_2,\ldots,p_k\to\infty}
\sum_{j_1=0}^{p_1}\ldots \sum_{j_k=0}^{p_k}
C_{j_k\ldots j_1}^2\stackrel{\sf def}{=}
\sum_{j_1=0}^{\infty}
\sum_{j_2,\ldots,j_k=0}^{\infty}
C_{j_k\ldots j_1}^2
$$

\vspace{2mm}
\noindent
and so on.

Let us consider the following lemma.

\vspace{2mm}

{\bf Lemma 4.}\ {\it The following equalities are fulfilled

\vspace{1mm}
$$
\sum_{j_1,\ldots,j_k=0}^{\infty}
C_{j_k\ldots j_1}^2=
\sum_{j_1=0}^{\infty}\ldots
\sum_{j_k=0}^{\infty}
C_{j_k\ldots j_1}^2=
$$

\vspace{2mm}
\begin{equation}
\label{lem1}
=\sum_{j_k=0}^{\infty}\ldots
\sum_{j_1=0}^{\infty}
C_{j_k\ldots j_1}^2=
\sum_{j_{q_1}=0}^{\infty}\ldots
\sum_{j_{q_k}=0}^{\infty}
C_{j_k\ldots j_1}^2
\end{equation}

\vspace{5mm}
\noindent
for any permutation $(q_1,\ldots,q_k)$ such that
$\{q_1,\ldots,q_k\}=\{1,\ldots,k\}.$}

\vspace{2mm}

{\bf Proof.} Let us consider the value

\vspace{-3mm}
\begin{equation}
\label{21}
\sum_{j_{q_l}=0}^{p}\ldots
\sum_{j_{q_k}=0}^{p}
C_{j_k\ldots j_1}^2
\end{equation}

\vspace{3mm}
\noindent
for any permutation $(q_l,\ldots,q_k)$, where $l=1,2,\ldots,k$,
$\{q_1,\ldots,q_k\}=\{1,\ldots,k\}.$

Obviously, (\ref{21}) 
is the non-decreasing sequence with respect to $p$.
Moreover,

\vspace{1mm}
$$
\sum_{j_{q_l}=0}^{p}\ldots
\sum_{j_{q_k}=0}^{p}
C_{j_k\ldots j_1}^2\le 
\sum_{j_{q_1}=0}^{p}\sum_{j_{q_2}=0}^{p}\ldots
\sum_{j_{q_k}=0}^{p}
C_{j_k\ldots j_1}^2\le 
$$

\vspace{4mm}
$$
\le
\sum_{j_1,\ldots,j_k=0}^{\infty}
C_{j_k\ldots j_1}^2<\infty.
$$

\vspace{4mm}

Then, the following limit

\vspace{1mm}
$$
\lim\limits_{p\to\infty}\sum\limits_{j_{q_l}=0}^p \ldots 
\sum\limits_{j_{q_k}=0}^{p}
C_{j_k\ldots j_1}^2=
\sum_{j_{q_l},\ldots,j_{q_k}=0}^{\infty}
C_{j_k\ldots j_1}^2
$$

\vspace{3mm}
\noindent
exists.

Let $p_l,\ldots,p_k$ simultaneously tend to infinity.
Then $g, r\to \infty$, where $g=\min\{p_l,\ldots,p_k\}$ and
$r=\max\{p_l,\ldots,p_k\}$. Moreover,

\vspace{2mm}
$$
\sum_{j_{q_l}=0}^{g}\ldots
\sum_{j_{q_k}=0}^{g}
C_{j_k\ldots j_1}^2\le 
\sum_{j_{q_l}=0}^{p_l}\ldots
\sum_{j_{q_k}=0}^{p_k}
C_{j_k\ldots j_1}^2\le
\sum_{j_{q_l}=0}^{r}\ldots
\sum_{j_{q_k}=0}^{r}
C_{j_k\ldots j_1}^2.
$$

\vspace{5mm}

This means that the existence of the limit 

\begin{equation}
\label{1c1c}
\lim\limits_{p\to\infty}\sum_{j_{q_l}=0}^{p}\ldots
\sum_{j_{q_k}=0}^{p}
C_{j_k\ldots j_1}^2
\end{equation}

\vspace{3mm}
\noindent
implies the existence of the limit 
\begin{equation}
\label{1d1d}
\lim\limits_{p_l,\ldots,p_k\to\infty}\sum_{j_{q_l}=0}^{p_l}\ldots
\sum_{j_{q_k}=0}^{p_k}
C_{j_k\ldots j_1}^2
\end{equation}

\vspace{3mm}
\noindent
and equality of the limits (\ref{1c1c}) and (\ref{1d1d}).
 
Taking into account the above reasoning, we have 

$$
\lim\limits_{p,q\to\infty}\sum_{j_{q_l}=0}^{q}\sum_{j_{q_{l+1}}=0}^{p}\ldots
\sum_{j_{q_k}=0}^{p}
C_{j_k\ldots j_1}^2=
\lim\limits_{p\to\infty}\sum_{j_{q_l}=0}^{p}\ldots
\sum_{j_{q_k}=0}^{p}
C_{j_k\ldots j_1}^2=
$$

\vspace{1mm}
\begin{equation}
\label{1h1h}
=\lim\limits_{p_l,\ldots,p_k\to\infty}\sum_{j_{q_l}=0}^{p_l}\ldots
\sum_{j_{q_k}=0}^{p_k}
C_{j_k\ldots j_1}^2.
\end{equation}

\vspace{4mm}

Since the limit
$$
\sum_{j_1,\ldots,j_k=0}^{\infty}
C_{j_k\ldots j_1}^2
$$

\vspace{3mm}
\noindent
exists (see the Parseval equality (\ref{par1})), then from Proposition 1
we have

\vspace{1mm}
$$
\sum_{j_{q_1}=0}^{\infty}\sum_{j_{q_2},\ldots,j_{q_k}=0}^{\infty}
C_{j_k\ldots j_1}^2=
\lim\limits_{q\to\infty}
\lim\limits_{p\to\infty}
\sum_{j_{q_1}=0}^{q}\sum_{j_{q_2}=0}^p \ldots \sum_{j_{q_k}=0}^{p}
C_{j_k\ldots j_1}^2=
$$

\vspace{3mm}
\begin{equation}
\label{1b1b}
=\lim\limits_{q,p\to\infty}
\sum_{j_{q_1}=0}^{q}\sum_{j_{q_2}=0}^p \ldots \sum_{j_{q_k}=0}^{p}
C_{j_k\ldots j_1}^2=
\sum_{j_1,\ldots,j_k=0}^{\infty}
C_{j_k\ldots j_1}^2.
\end{equation}

\vspace{5mm}

Using (\ref{1h1h}) and Proposition 1, we obtain

$$
\sum_{j_{q_2}=0}^{\infty}\sum_{j_{q_3},\ldots,j_{q_k}=0}^{\infty}
C_{j_k\ldots j_1}^2=
\lim\limits_{q\to\infty}
\lim\limits_{p\to\infty}
\sum_{j_{q_2}=0}^{q}\sum_{j_{q_3}=0}^p \ldots \sum_{j_{q_k}=0}^{p}
C_{j_k\ldots j_1}^2=
$$

\vspace{3mm}
\begin{equation}
\label{1a1a}
=\lim\limits_{q,p\to\infty}
\sum_{j_{q_2}=0}^{q}\sum_{j_{q_3}=0}^p \ldots \sum_{j_{q_k}=0}^{p}
C_{j_k\ldots j_1}^2=
\sum_{j_{q_2},\ldots,j_{q_k}=0}^{\infty}
C_{j_k\ldots j_1}^2.
\end{equation}

\vspace{5mm}

Combining (\ref{1a1a}) and (\ref{1b1b}), we get

$$
\sum_{j_{q_1}=0}^{\infty}\sum_{j_{q_2}=0}^{\infty}
\sum_{j_{q_3},\ldots,j_{q_k}=0}^{\infty}
C_{j_k\ldots j_1}^2=
\sum_{j_{1},\ldots,j_{k}=0}^{\infty}
C_{j_k\ldots j_1}^2.
$$

\vspace{3mm}

Repeating the above steps, we complete the proof of Lemma
4.

Further, let us show that for $s=1,\ldots,k$

\vspace{2mm}
$$
\sum_{j_1=0}^{\infty}\ldots
\sum_{j_{s-1}=0}^{\infty}
\sum_{j_s=p+1}^{\infty}\sum_{j_{s+1}=0}^{\infty}\ldots \sum_{j_k=0}^{\infty}
C_{j_k\ldots j_1}^2=
$$

\vspace{3mm}
\begin{equation}
\label{d11}
=
\sum_{j_s=p+1}^{\infty}\sum_{j_{s-1}=0}^{\infty}\ldots
\sum_{j_{1}=0}^{\infty}
\sum_{j_{s+1}=0}^{\infty}\ldots \sum_{j_k=0}^{\infty}
C_{j_k\ldots j_1}^2.
\end{equation}

\vspace{5mm}

Using the arguments which we used when proving Lemma 4, we have

\vspace{2mm}
$$
\lim\limits_{n\to\infty}
\sum_{j_1=0}^{n}\ldots
\sum_{j_{s-1}=0}^{n}
\sum_{j_s=0}^{p}\sum_{j_{s+1}=0}^{n}\ldots \sum_{j_k=0}^{n}
C_{j_k\ldots j_1}^2=
$$

\vspace{3mm}
\begin{equation}
\label{ura0}
=\sum_{j_s=0}^{p}\ \sum_{j_{1},\ldots, j_{s-1}, 
j_{s+1},\ldots,j_k=0}^{\infty}
C_{j_k\ldots j_1}^2
=\sum_{j_s=0}^{p}\sum_{j_{q_1}=0}^{\infty}\ldots
\sum_{j_{q_{k-1}}=0}^{\infty}
C_{j_k\ldots j_1}^2
\end{equation}

\vspace{5mm}
\noindent
for any permutation $(q_1,\ldots,q_{k-1})$ such that
$\{q_1,\ldots,q_{k-1}\}=\{1,\ldots,s-1,s+1,\ldots,k\}$,
where $p$ is a fixed natural number.

Obviously, we have

\vspace{2mm}
$$
\sum_{j_s=0}^{p}\sum_{j_{q_1}=0}^{\infty}\ldots
\sum_{j_{q_{k-1}}=0}^{\infty}
C_{j_k\ldots j_1}^2=
\sum_{j_{q_1}=0}^{\infty}\ldots \sum_{j_s=0}^{p} \ldots
\sum_{j_{q_{k-1}}=0}^{\infty}C_{j_k\ldots j_1}^2 = \ldots =
$$

\vspace{3mm}
\begin{equation}
\label{ura1}
=
\sum_{j_{q_1}=0}^{\infty}\ldots 
\sum_{j_{q_{k-1}}=0}^{\infty}
\sum_{j_s=0}^{p}
C_{j_k\ldots j_1}^2.
\end{equation}

\vspace{5mm}

Using (\ref{ura0}), (\ref{ura1}), and Lemma 4, we obtain

\vspace{2mm}
$$
\sum_{j_1=0}^{\infty}\ldots
\sum_{j_{s-1}=0}^{\infty}
\sum_{j_s=p+1}^{\infty}\sum_{j_{s+1}=0}^{\infty}\ldots \sum_{j_k=0}^{\infty}
C_{j_k\ldots j_1}^2=
\sum_{j_1=0}^{\infty}\ldots
\sum_{j_{s-1}=0}^{\infty}
\sum_{j_s=0}^{\infty}\sum_{j_{s+1}=0}^{\infty}\ldots \sum_{j_k=0}^{\infty}
C_{j_k\ldots j_1}^2-
$$

\vspace{3mm}
$$
-\sum_{j_1=0}^{\infty}\ldots
\sum_{j_{s-1}=0}^{\infty}
\sum_{j_s=0}^{p}\sum_{j_{s+1}=0}^{\infty}\ldots \sum_{j_k=0}^{\infty}
C_{j_k\ldots j_1}^2=
$$

\vspace{3mm}
$$
=
\sum_{j_s=0}^{\infty}
\sum_{j_{s-1}=0}^{\infty}\ldots
\sum_{j_1=0}^{\infty}\sum_{j_{s+1}=0}^{\infty}\ldots \sum_{j_k=0}^{\infty}
C_{j_k\ldots j_1}^2-
\sum_{j_s=0}^{p}
\sum_{j_{s-1}=0}^{\infty}\ldots
\sum_{j_1=0}^{\infty}\sum_{j_{s+1}=0}^{\infty}\ldots \sum_{j_k=0}^{\infty}
C_{j_k\ldots j_1}^2=
$$

\vspace{3mm}
$$
=
\sum_{j_s=p+1}^{\infty}
\sum_{j_{s-1}=0}^{\infty}\ldots
\sum_{j_1=0}^{\infty}\sum_{j_{s+1}=0}^{\infty}\ldots \sum_{j_k=0}^{\infty}
C_{j_k\ldots j_1}^2.
$$

\vspace{5mm}

The equality (\ref{d11}) is proved.

Using the Parseval equality and Lemma 4, we obtain

\vspace{2mm}
$$
\int\limits_{[t,T]^k}K^2(t_1,\ldots,t_k)dt_1\ldots dt_k-
\sum_{j_1=0}^{p}\ldots \sum_{j_k=0}^{p}
C_{j_k\ldots j_1}^2=
$$

\vspace{3mm}

$$
=\sum_{j_1,\ldots,j_k=0}^{\infty}
C_{j_k\ldots j_1}^2-
\sum_{j_1=0}^{p}\ldots \sum_{j_k=0}^{p}
C_{j_k\ldots j_1}^2=
$$

\vspace{5mm}

$$
=
\sum_{j_1=0}^{\infty}\ldots \sum_{j_k=0}^{\infty}
C_{j_k\ldots j_1}^2-
\sum_{j_1=0}^{p}\ldots \sum_{j_k=0}^{p}
C_{j_k\ldots j_1}^2=
$$

\vspace{5mm}

$$
=\sum_{j_1=0}^{p}\sum_{j_2=0}^{\infty}\ldots \sum_{j_k=0}^{\infty}
C_{j_k\ldots j_1}^2+
\sum_{j_1=p+1}^{\infty}\sum_{j_2=0}^{\infty}\ldots \sum_{j_k=0}^{\infty}
C_{j_k\ldots j_1}^2-
\sum_{j_1=0}^{p}\ldots \sum_{j_k=0}^{p}
C_{j_k\ldots j_1}^2=
$$

\vspace{6mm}

$$
=\sum_{j_1=0}^{p}\sum_{j_2=0}^{p}\sum_{j_3=0}^{\infty}
\ldots \sum_{j_k=0}^{\infty}
C_{j_k\ldots j_1}^2+
\sum_{j_1=0}^{p}\sum_{j_2=p+1}^{\infty}
\sum_{j_3=0}^{\infty}
\ldots \sum_{j_k=0}^{\infty}+
$$

\vspace{4mm}

$$
+\sum_{j_1=p+1}^{\infty}\sum_{j_2=0}^{\infty}\ldots \sum_{j_k=0}^{\infty}
C_{j_k\ldots j_1}^2-
\sum_{j_1=0}^{p}\ldots \sum_{j_k=0}^{p}
C_{j_k\ldots j_1}^2=\ldots =
$$

\vspace{6mm}

$$
=\sum_{j_1=p+1}^{\infty}\sum_{j_2=0}^{\infty}\ldots \sum_{j_k=0}^{\infty}
C_{j_k\ldots j_1}^2+
\sum_{j_1=0}^p
\sum_{j_2=p+1}^{\infty}\sum_{j_2=0}^{\infty}\ldots \sum_{j_k=0}^{\infty}
C_{j_k\ldots j_1}^2+
$$

\vspace{6mm}

$$
+\sum_{j_1=0}^p\sum_{j_2=0}^p
\sum_{j_3=p+1}^{\infty}\sum_{j_4=0}^{\infty}\ldots \sum_{j_k=0}^{\infty}
C_{j_k\ldots j_1}^2+ \ldots +
\sum_{j_1=0}^p\ldots \sum_{j_{k-1}=0}^p
\sum_{j_k=p+1}^{\infty}C_{j_k\ldots j_1}^2\le
$$

\vspace{6mm}

$$
\le\sum_{j_1=p+1}^{\infty}\sum_{j_2=0}^{\infty}\ldots \sum_{j_k=0}^{\infty}
C_{j_k\ldots j_1}^2+
\sum_{j_1=0}^{\infty}
\sum_{j_2=p+1}^{\infty}\sum_{j_2=0}^{\infty}\ldots \sum_{j_k=0}^{\infty}
C_{j_k\ldots j_1}^2+
$$

\vspace{4mm}

$$
+\sum_{j_1=0}^{\infty}\sum_{j_2=0}^{\infty}
\sum_{j_3=p+1}^{\infty}\sum_{j_4=0}^{\infty}\ldots \sum_{j_k=0}^{\infty}
C_{j_k\ldots j_1}^2+ \ldots +
\sum_{j_1=0}^{\infty}\ldots \sum_{j_{k-1}=0}^{\infty}
\sum_{j_k=p+1}^{\infty}C_{j_k\ldots j_1}^2=
$$

\vspace{5mm}

\begin{equation}
\label{aaap}
=\sum\limits_{s=1}^k \left(\sum_{j_1=0}^{\infty}\ldots
\sum_{j_{s-1}=0}^{\infty}
\sum_{j_s=p+1}^{\infty}\sum_{j_{s+1}=0}^{\infty}\ldots \sum_{j_k=0}^{\infty}
C_{j_k\ldots j_1}^2\right).
\end{equation}

\vspace{5mm}

Note that deriving (\ref{aaap}) we use the following

\vspace{1mm}
$$
\sum_{j_1=0}^{p}\ldots
\sum_{j_{s-1}=0}^{p}
\sum_{j_s=p+1}^{\infty}\sum_{j_{s+1}=0}^{\infty}\ldots \sum_{j_k=0}^{\infty}
C_{j_k\ldots j_1}^2\le
$$

\vspace{3mm}
$$
\le
\sum_{j_1=0}^{m_1}\ldots
\sum_{j_{s-1}=0}^{m_{s-1}}
\sum_{j_s=p+1}^{\infty}\sum_{j_{s+1}=0}^{\infty}\ldots \sum_{j_k=0}^{\infty}
C_{j_k\ldots j_1}^2\le
$$

\vspace{3mm}
$$
\le
\lim\limits_{m_{s-1}\to\infty}
\sum_{j_1=0}^{m_1}\ldots
\sum_{j_{s-1}=0}^{m_{s-1}}
\sum_{j_s=p+1}^{\infty}\sum_{j_{s+1}=0}^{\infty}\ldots \sum_{j_k=0}^{\infty}
C_{j_k\ldots j_1}^2=
$$

\vspace{3mm}
$$
=
\sum_{j_1=0}^{m_1}\ldots
\sum_{j_{s-2}=0}^{m_{s-2}}\sum_{j_{s-1}=0}^{\infty}
\sum_{j_s=p+1}^{\infty}\sum_{j_{s+1}=0}^{\infty}\ldots \sum_{j_k=0}^{\infty}
C_{j_k\ldots j_1}^2\le
$$

$$
\le\ldots\le
$$

$$
\le\sum_{j_1=0}^{\infty}\ldots
\sum_{j_{s-1}=0}^{\infty}
\sum_{j_s=p+1}^{\infty}\sum_{j_{s+1}=0}^{\infty}\ldots \sum_{j_k=0}^{\infty}
C_{j_k\ldots j_1}^2,
$$

\vspace{5mm}
\noindent
where $m_1,\ldots,m_{s-1}>p.$

Denote
$$
C_{j_s\ldots j_1}(\tau)=
\int\limits_t^{\tau}
\phi_{j_s}(t_s)\psi_s(t_s)\ldots \int\limits_t^{t_2}
\phi_{j_1}(t_1)\psi_1(t_1)dt_1\ldots dt_s,
$$

\vspace{2mm}
\noindent
where
$s=1,\ldots,k-1.$

Let us remind the Dini Theorem, which we will use further.

\vspace{2mm}

{\bf Theorem (Dini).} {\it 
Let the functional sequence $u_n(x)$ 
be non-decreasing at each point of the interval $[a, b]$.
In addition, all the functions $u_n(x)$
of this sequence and the limit function $u(x)$ are continuous on the interval
$[a, b].$ Then the convergence $u_n(x)$ to 
$u(x)$ is uniform on the interval $[a,b].$}

\vspace{2mm}

For $s<k$ due to the Parseval equality and Dini Theorem as well as
(\ref{d11}) we obtain

\vspace{2mm}
$$
\sum_{j_1=0}^{\infty}\ldots
\sum_{j_{s-1}=0}^{\infty}
\sum_{j_s=p+1}^{\infty}\sum_{j_{s+1}=0}^{\infty}\ldots \sum_{j_k=0}^{\infty}
C_{j_k\ldots j_1}^2=
$$

\vspace{5mm}

$$
\stackrel{\hbox{(\ref{d11})}}{=}\ \ \ 
\sum_{j_s=p+1}^{\infty}
\sum_{j_{s-1}=0}^{\infty}\ldots
\sum_{j_{1}=0}^{\infty}
\sum_{j_{s+1}=0}^{\infty}\ldots \sum_{j_k=0}^{\infty}
C_{j_k\ldots j_1}^2=
$$

\vspace{5mm}

$$
\stackrel{\hbox{(Parseval Eq.)}}{=}\ \ \
\sum_{j_s=p+1}^{\infty}
\sum_{j_{s-1}=0}^{\infty}\ldots
\sum_{j_{1}=0}^{\infty}
\sum_{j_{s+1}=0}^{\infty}\ldots 
\sum_{j_{k-1}=0}^{\infty}
\int\limits_t^T \psi_k^2(t_k) \left(C_{j_{k-1}\ldots j_1}(t_k)\right)^2 dt_k=
$$

\vspace{5mm}

$$
\stackrel{\hbox{(Dini Th.)}}{=}\ \ \
\sum_{j_s=p+1}^{\infty}
\sum_{j_{s-1}=0}^{\infty}\ldots
\sum_{j_{1}=0}^{\infty}
\sum_{j_{s+1}=0}^{\infty}\ldots 
\sum_{j_{k-2}=0}^{\infty}
\int\limits_t^T \psi_k^2(t_k) 
\sum_{j_{k-1}=0}^{\infty}\left(C_{j_{k-1}\ldots j_1}(t_k)\right)^2 dt_k=
$$

\vspace{5mm}

$$
\stackrel{\hbox{(Parseval Eq.)}}{=}\ \ \
\sum_{j_s=p+1}^{\infty}
\sum_{j_{s-1}=0}^{\infty}\ldots
\sum_{j_{1}=0}^{\infty}
\sum_{j_{s+1}=0}^{\infty}\ldots 
\sum_{j_{k-2}=0}^{\infty}
\int\limits_t^T \psi_k^2(t_k) \int\limits_t^{t_k} \psi_{k-1}^2(t_{k-1}) 
\left(C_{j_{k-2}\ldots j_1}(t_{k-1})\right)^2\times
$$

\vspace{2mm}
$$
\times dt_{k-1}dt_k\le
$$

\vspace{5mm}
$$
\le C\sum_{j_s=p+1}^{\infty}
\sum_{j_{s-1}=0}^{\infty}\ldots
\sum_{j_{1}=0}^{\infty}
\sum_{j_{s+1}=0}^{\infty}\ldots 
\sum_{j_{k-2}=0}^{\infty}
\int\limits_t^T 
\left(C_{j_{k-2}\ldots j_1}(\tau)\right)^2 d\tau =
$$

\vspace{5mm}

$$
\stackrel{\hbox{(Dini Th.)}}{=}\ \ \
C\sum_{j_s=p+1}^{\infty}
\sum_{j_{s-1}=0}^{\infty}\ldots
\sum_{j_{1}=0}^{\infty}
\sum_{j_{s+1}=0}^{\infty}\ldots 
\sum_{j_{k-3}=0}^{\infty}
\int\limits_t^T 
\sum_{j_{k-2}=0}^{\infty}
\left(C_{j_{k-2}\ldots j_1}(\tau)\right)^2 d\tau =
$$

\vspace{5mm}

$$
\stackrel{\hbox{(Parseval Eq.)}}{=}\ \ \ C\sum_{j_s=p+1}^{\infty}
\sum_{j_{s-1}=0}^{\infty}\ldots
\sum_{j_{1}=0}^{\infty}
\sum_{j_{s+1}=0}^{\infty}\ldots 
\sum_{j_{k-3}=0}^{\infty}
\int\limits_t^T \int\limits_t^{\tau}
\psi_{k-2}^2(\theta)
\left(C_{j_{k-3}\ldots j_1}(\theta)\right)^2
d\theta d\tau\le 
$$

\vspace{5mm}

$$
\le K
\sum_{j_s=p+1}^{\infty}
\sum_{j_{s-1}=0}^{\infty}\ldots
\sum_{j_{1}=0}^{\infty}
\sum_{j_{s+1}=0}^{\infty}\ldots 
\sum_{j_{k-3}=0}^{\infty}
\int\limits_t^T
\left(C_{j_{k-3}\ldots j_1}(\tau)\right)^2
d\tau\le 
$$

\vspace{4mm}
$$
\le \ldots \le
$$

$$
\le C_k
\sum_{j_s=p+1}^{\infty}
\sum_{j_{s-1}=0}^{\infty}\ldots
\sum_{j_{1}=0}^{\infty}
\int\limits_t^T 
\left(C_{j_{s}\ldots j_1}(\tau)\right)^2 d\tau=
$$

\vspace{5mm}

\begin{equation}
\label{d14}
\stackrel{\hbox{(Dini Th.)}}{=}\ \ \ C_k
\sum_{j_s=p+1}^{\infty}
\sum_{j_{s-1}=0}^{\infty}\ldots
\sum_{j_{2}=0}^{\infty}
\int\limits_t^T  \sum_{j_{1}=0}^{\infty}
\left(C_{j_{s}\ldots j_1}(\tau)\right)^2 d\tau,
\end{equation}

\vspace{6mm}
\noindent
where constants $C,$ $K$ depend on $T-t$ and
constant $C_k$ depends on $k$ and $T-t.$

Let us explane more precisely how we obtain (\ref{d14}).
For any function $g(s)\in L_2([t,T])$ we have the following
Parseval equality

$$
\sum\limits_{j=0}^{\infty}\left(\int\limits_t^{\tau}
\phi_j(s)g(s)ds\right)^2=
\sum\limits_{j=0}^{\infty}\left(\int\limits_t^T
{\bf 1}_{\{s<\tau\}}\phi_j(s)g(s)ds\right)^2=
$$

\begin{equation}
\label{d15}
=\int\limits_t^T
\left({\bf 1}_{\{s<\tau\}}\right)^2 g^2(s)ds=
\int\limits_t^{\tau}
g^2(s)ds.
\end{equation}

\vspace{3mm}

The equality (\ref{d15}) has been applied repeatedly when we obtaining
(\ref{d14}).

Using the replacement of integrating order in Riemann integrals, we have

\vspace{1mm}
$$
C_{j_s\ldots j_1}(\tau)=
\int\limits_t^{\tau}
\phi_{j_s}(t_s)\psi_s(t_s)\ldots \int\limits_t^{t_2}
\phi_{j_1}(t_1)\psi_1(t_1)dt_1\ldots dt_s=
$$

\vspace{2mm}
$$
=\int\limits_t^{\tau}
\phi_{j_1}(t_1)\psi_1(t_1)\int\limits_{t_1}^{\tau}
\phi_{j_2}(t_2)\psi_2(t_2)
\ldots
\int\limits_{t_{s-1}}^{\tau}
\phi_{j_s}(t_s)\psi_s(t_s)dt_s\ldots dt_2dt_1
\stackrel{\sf def}{=}
$$

\vspace{4mm}
$$
\stackrel{\sf def}{=}
{\tilde C}_{j_s\ldots j_1}(\tau).
$$

\vspace{7mm}

For $l=1,\ldots,s$ we will use the following notation

\vspace{2mm}
$$
{\tilde C}_{j_s\ldots j_l}(\tau,\theta)=
\int\limits_{\theta}^{\tau}
\phi_{j_l}(t_l)\psi_l(t_l)\int\limits_{t_l}^{\tau}
\phi_{j_{l+1}}(t_{l+1})\psi_{l+1}(t_{l+1})
\ldots
\int\limits_{t_{s-1}}^{\tau}
\phi_{j_s}(t_s)\psi_s(t_s)dt_s\ldots dt_{l+1}dt_l.
$$

\vspace{5mm}

Using the Parseval equality and Dini Theorem, from (\ref{d14}) we obtain

\vspace{2mm}
$$
\sum_{j_1=0}^{\infty}\ldots
\sum_{j_{s-1}=0}^{\infty}
\sum_{j_s=p+1}^{\infty}\sum_{j_{s+1}=0}^{\infty}\ldots \sum_{j_k=0}^{\infty}
C_{j_k\ldots j_1}^2\le
$$

\vspace{5mm}

$$
\le
C_k
\sum_{j_s=p+1}^{\infty}
\sum_{j_{s-1}=0}^{\infty}\ldots
\sum_{j_{2}=0}^{\infty}
\int\limits_t^T  \sum_{j_{1}=0}^{\infty}
\left(C_{j_{s}\ldots j_1}(\tau)\right)^2 d\tau=
$$

\vspace{5mm}

$$
=C_k
\sum_{j_s=p+1}^{\infty}
\sum_{j_{s-1}=0}^{\infty}\ldots
\sum_{j_{2}=0}^{\infty}
\int\limits_t^T  \sum_{j_{1}=0}^{\infty}
\left({\tilde C}_{j_{s}\ldots j_1}(\tau)\right)^2 d\tau=
$$

\vspace{5mm}

\begin{equation}
\label{molod1}
\stackrel{\hbox{(Parseval Eq.)}}{=}\ \ \ C_k
\sum_{j_s=p+1}^{\infty}
\sum_{j_{s-1}=0}^{\infty}\ldots
\sum_{j_{2}=0}^{\infty}
\int\limits_t^T\int\limits_t^{\tau}\psi_1^2(t_1)  
\left({\tilde C}_{j_{s}\ldots j_2}(\tau,t_1)\right)^2 dt_1d\tau=
\end{equation}

\vspace{5mm}

\begin{equation}
\label{molod2}
\stackrel{\hbox{(Dini Th.)}}{=}\ \ \ C_k
\sum_{j_s=p+1}^{\infty}
\sum_{j_{s-1}=0}^{\infty}\ldots
\sum_{j_{3}=0}^{\infty}
\int\limits_t^T\int\limits_t^{\tau}\psi_1^2(t_1)  
\sum_{j_{2}=0}^{\infty}
\left({\tilde C}_{j_{s}\ldots j_2}(\tau,t_1)\right)^2 dt_1d\tau=
\end{equation}

\vspace{5mm}

$$
\stackrel{\hbox{(Parseval Eq.)}}{=}\ \ \ C_k
\sum_{j_s=p+1}^{\infty}
\sum_{j_{s-1}=0}^{\infty}\ldots
\sum_{j_{3}=0}^{\infty}
\int\limits_t^T\int\limits_t^{\tau}\psi_1^2(t_1)  
\int\limits_{t_1}^{\tau}\psi_2^2(t_2)  
\left({\tilde C}_{j_{s}\ldots j_3}(\tau,t_2)\right)^2 dt_2dt_1d\tau \le
$$

\vspace{5mm}

$$
\le C_k
\sum_{j_s=p+1}^{\infty}
\sum_{j_{s-1}=0}^{\infty}\ldots
\sum_{j_{3}=0}^{\infty}
\int\limits_t^T\int\limits_t^{\tau}\psi_1^2(t_1)  
\int\limits_{t}^{\tau}\psi_2^2(t_2)  
\left({\tilde C}_{j_{s}\ldots j_3}(\tau,t_2)\right)^2 dt_2dt_1d\tau\le
$$

\vspace{5mm}

$$
\le C^{'}_k
\sum_{j_s=p+1}^{\infty}
\sum_{j_{s-1}=0}^{\infty}\ldots
\sum_{j_{3}=0}^{\infty}
\int\limits_t^T
\int\limits_{t}^{\tau}\psi_2^2(t_2)  
\left({\tilde C}_{j_{s}\ldots j_3}(\tau,t_2)\right)^2 dt_2d\tau
\le 
$$

\vspace{4mm}

$$
\le  \ldots \le
$$

\vspace{1mm}

$$
\le C^{''}_k
\sum_{j_s=p+1}^{\infty}
\int\limits_t^T\int\limits_t^{\tau}
\psi_{s-1}^2(t_{s-1})
\left({\tilde C}_{j_{s}}(\tau,t_{s-1})\right)^2 dt_{s-1} d\tau\le
$$

\vspace{4mm}

\begin{equation}
\label{la}
\le {\tilde C}_k
\sum_{j_s=p+1}^{\infty}
\int\limits_t^T\int\limits_t^{\tau}
\left(~\int\limits_{u}^{\tau}\phi_{j_s}(\theta)
\psi_s(\theta)d\theta\right)^2 du d\tau,
\end{equation}

\vspace{6mm}
\noindent
where constants $C^{'}_k,$ $C^{''}_k,$ $\tilde C_k$
depend on $k$ and $T-t.$

Let us explane more precisely how we obtain (\ref{la}).
For any function $g(s)\in L_2([t,T])$ we have the following
Parseval equality

$$
\sum\limits_{j=0}^{\infty}\left(\int\limits_{\theta}^{\tau}
\phi_j(s)g(s)ds\right)^2=
\sum\limits_{j=0}^{\infty}\left(\int\limits_t^T
{\bf 1}_{\{\theta<s<\tau\}}\phi_j(s)g(s)ds\right)^2=
$$

\begin{equation}
\label{d22}
=\int\limits_t^T
\left({\bf 1}_{\{\theta<s<\tau\}}\right)^2 g^2(s)ds=
\int\limits_{\theta}^{\tau}
g^2(s)ds.
\end{equation}

\vspace{3mm}

The equality (\ref{d22}) has been applied repeatedly when we obtaining
(\ref{la}).

Let us explane more precisely the passing from (\ref{molod1})
to (\ref{molod2}) (the same steps have been used when we 
deriving (\ref{la})).

We have

$$
\int\limits_t^T\int\limits_t^{\tau}\psi_1^2(t_1)  
\sum_{j_{2}=0}^{\infty}
\left({\tilde C}_{j_{s}\ldots j_2}(\tau,t_1)\right)^2 dt_1d\tau -
\sum_{j_{2}=0}^{n}\int\limits_t^T\int\limits_t^{\tau}\psi_1^2(t_1)  
\left({\tilde C}_{j_{s}\ldots j_2}(\tau,t_1)\right)^2 dt_1d\tau =
$$

\vspace{4mm}

$$
=\int\limits_t^T\int\limits_t^{\tau}\psi_1^2(t_1)  
\sum_{j_{2}=n+1}^{\infty}
\left({\tilde C}_{j_{s}\ldots j_2}(\tau,t_1)\right)^2 dt_1d\tau =
$$

\vspace{4mm}

\begin{equation}
\label{molod3}
=\lim\limits_{N\to\infty}
\sum\limits_{j=0}^{N-1}\int\limits_t^{\tau_j}\psi_1^2(t_1)  
\sum_{j_{2}=n+1}^{\infty}
\left({\tilde C}_{j_{s}\ldots j_2}(\tau_j,t_1)\right)^2 dt_1 \Delta\tau_j,
\end{equation}

\vspace{4mm}
\noindent
where $\{\tau_j\}_{j=0}^{N}$ is the partition of the 
interval $[t, T],$ which satisfies the condition (\ref{1111}).

Since the non-decreasing functional sequence $u_n(\tau_j,t_1)$ and its
limit function $u(\tau_j,t_1)$ are continuous on the
interval $[t,\tau_j]\subseteq [t, T]$ with respect to $t_1$,
where

$$
u_n(\tau_j,t_1)=
\sum_{j_{2}=0}^{n}
\left({\tilde C}_{j_{s}\ldots j_2}(\tau_j,t_1)\right)^2,
$$

$$
u(\tau_j,t_1)=
\sum_{j_{2}=0}^{\infty}
\left({\tilde C}_{j_{s}\ldots j_2}(\tau_j,t_1)\right)^2=
\int\limits_{t_1}^{\tau_j}
\psi_2^2(t_2)
\left({\tilde C}_{j_{s}\ldots j_3}(\tau_j,t_2)\right)^2 dt_2,
$$

\vspace{4mm}

\noindent 
then by Dini Theorem we have the uniform convergence
of $u_n(\tau_j,t_1)$ to $u(\tau_j,t_1)$ at the interval $[t,\tau_j]\subseteq
[t, T]$
with respect to $t_1.$ As a result, we obtain

\begin{equation}
\label{molod4}
\sum_{j_{2}=n+1}^{\infty}
\left({\tilde C}_{j_{s}\ldots j_2}(\tau_j,t_1)\right)^2<\varepsilon,\ \ \ 
t_1\in [t,\tau_j]
\end{equation}

\vspace{4mm}
\noindent
for $n>N(\varepsilon)$ ($N(\varepsilon)$ exists
for any $\varepsilon>0$ and it does not depend on $t_1$).

From (\ref{molod3}) and (\ref{molod4}) we obtain

$$
\lim\limits_{N\to\infty}
\sum\limits_{j=0}^{N-1}\int\limits_t^{\tau_j}\psi_1^2(t_1)  
\sum_{j_{2}=n+1}^{\infty}
\left({\tilde C}_{j_{s}\ldots j_2}(\tau_j,t_1)\right)^2 dt_1 \Delta\tau_j
\le
\varepsilon 
\lim\limits_{N\to\infty}
\sum\limits_{j=0}^{N-1}\int\limits_t^{\tau_j}\psi_1^2(t_1)  
dt_1 \Delta\tau_j= 
$$

\vspace{3mm}
\begin{equation}
\label{molod6}
=\varepsilon \int\limits_t^T
\int\limits_t^{\tau}\psi_1^2(t_1)  
dt_1 d\tau.
\end{equation}

\vspace{4mm}

From (\ref{molod6}) we get

$$
\lim\limits_{n\to\infty}\int\limits_t^T\int\limits_t^{\tau}\psi_1^2(t_1)  
\sum_{j_{2}=n+1}^{\infty}
\left({\tilde C}_{j_{s}\ldots j_2}(\tau,t_1)\right)^2 dt_1d\tau = 0.
$$

\vspace{4mm}

This fact completes the proof of passing 
from (\ref{molod1})
to (\ref{molod2}).

Let us estimate the integral 
\begin{equation}
\label{st1}
\int\limits_{u}^{\tau}\phi_{j_s}(\theta)
\psi_s(\theta)d\theta
\end{equation}

\vspace{2mm}
\noindent
from (\ref{la}) for the cases when $\{\phi_j(s)\}_{j=0}^{\infty}$
is a complete orthonormal system of Legendre polynomials or
trigonometric functions in the space $L_2([t,T])$.

Note that the estimates for the integral

\vspace{-2mm}
\begin{equation}
\label{st2}
\int\limits_{t}^{\tau}\phi_{j}(\theta)\psi(\theta)d\theta,\ \ \ j\ge p+1,
\end{equation}

\vspace{2mm}
\noindent
where $\psi(\theta)$ is a continuously
differentiable function on the interval $[t, T]$,
have been obtained in \cite{27}
(see the formulas 
(54) (55), (60)) or in \cite{21}
(see the formulas (57), (58), (63)). The same estimates 
also can be found in early publications \cite{15}, \cite{16}, 
\cite{19}, \cite{20} and in \cite{20a}-\cite{2023xxx1} (2020, 2021, 2023).

Let us estimate the integral (\ref{st1}) using the approach from
\cite{21}, \cite{27}.

First, consider the case of Legendre polynomials.
Then $\phi_j(s)$ looks as follows

\vspace{1mm}
\begin{equation}
\label{ogo7}
\phi_j(\theta)=\sqrt{\frac{2j+1}{T-t}}P_j\left(\left(
\theta-\frac{T+t}{2}\right)\frac{2}{T-t}\right),\ \ \ j\ge 0,
\end{equation}

\vspace{5mm}
\noindent
where $P_j(x)$ $(j=0, 1, 2\ldots)$ is a complete
orthonormal system of Legendre polynomials in the space $L_2([-1,1]).$

Further, we have 

$$
\int\limits_v^x\phi_{j}(\theta)\psi(\theta)d\theta=
\frac{\sqrt{T-t}\sqrt{2j+1}}{2}
\int\limits_{z(v)}^{z(x)}P_{j}(y)
\psi(u(y))dy=
$$

$$
=\frac{\sqrt{T-t}}{2\sqrt{2j+1}}\Biggl((P_{j+1}(z(x))-
P_{j-1}(z(x)))\psi(x)-
(P_{j+1}(z(v))-
P_{j-1}(z(v)))\psi(v)-
\Biggr.
$$

\begin{equation}
\label{6000}
\Biggl.-
\frac{T-t}{2}
\int\limits_{z(v)}^{z(x)}((P_{j+1}(y)-P_{j-1}(y))
{\psi}'(u(y))dy\Biggr),
\end{equation}

\vspace{5mm}
\noindent
where $x, v\in (t, T),$ $j\ge p+1,$ 
$u(y)$ and $z(x)$ are defined by the following relations

\vspace{1mm}
$$
u(y)=\frac{T-t}{2}y+\frac{T+t}{2},\ \ \
z(x)=\left(x-\frac{T+t}{2}\right)\frac{2}{T-t},
$$

\vspace{4mm}
\noindent
${\psi}'$ is a derivative of the function $\psi(\theta)$
with respect to the variable $u(y).$

Note that in (\ref{6000}) we used the following well-known property
of the Legendre polynomials

\vspace{2mm}
$$
\frac{dP_{j+1}}{dx}(x)-\frac{dP_{j-1}}{dx}(x)=(2j+1)P_j(x),\ \ \ 
j=1, 2,\ldots
$$

\vspace{5mm}

From (\ref{6000}) and the well-known estimate for the Legendre
polynomials

\vspace{1mm}
\begin{equation}
\label{200}
|P_j(y)| <\frac{K}{\sqrt{j+1}(1-y^2)^{1/4}},\ \ \ 
y\in (-1, 1),\ \ \ j\in \mathbb{N},
\end{equation}

\vspace{5mm}
\noindent
where constant $K$ does not depend on $y$ and $j$, it follows that

\vspace{2mm}
\begin{equation}
\label{101}
\left|
\int\limits_v^x\phi_{j}(\theta)\psi(\theta)d\theta
\right| <
\frac{C}{j}\Biggl(\frac{1}{(1-(z(x))^2)^{1/4}}+
\frac{1}{(1-(z(v))^2)^{1/4}}+C_1\Biggr),
\end{equation}

\vspace{5mm}
\noindent
where $j\in \mathbb{N},$ 
$z(x), z(v)\in (-1, 1),$ $x, v\in (t, T),$
constants $C,$ $C_1$ do not depend on $j$.

From (\ref{101}) we obtain

\vspace{1mm}
\begin{equation}
\label{102}
\left(
\int\limits_v^x\phi_{j}(\theta)\psi(\theta)d\theta
\right)^2 <
\frac{C_2}{j^2}\Biggl(\frac{1}{(1-(z(x))^2)^{1/2}}+
\frac{1}{(1-(z(v))^2)^{1/2}}+C_3\Biggr),
\end{equation}

\vspace{4mm}
\noindent
where $j\in \mathbb{N},$ constants $C_2,$ $C_3$ do not depend on $j$.

Let us apply (\ref{102}) for the estimate of the right-hand side
of (\ref{la}). We have

\vspace{1mm}
$$
\int\limits_t^T\int\limits_t^{\tau}
\left(~\int\limits_{u}^{\tau}\phi_{j_s}(\theta)
\psi_s(\theta)d\theta\right)^2 du d\tau\le
$$

$$
\le \frac{K_1}{j_s^2}
\left(
\int\limits_{-1}^1
\frac{dy}{\left(1-y^2\right)^{1/2}}+
\int\limits_{-1}^1\int\limits_{-1}^x
\frac{dy}{\left(1-y^2\right)^{1/2}}dx + K_2\right)\le
$$

\begin{equation}
\label{103}
\le \frac{K_3}{j_s^2},
\end{equation}

\vspace{4mm}
\noindent
where $j_s\in \mathbb{N},$ constants $K_1, K_2, K_3$ are independent of $j_s.$

Now consider the trigonometric case.
The complete orthonormal system of trigonometric functions
in the space $L_2([t, T])$ has the following form

\begin{equation}
\label{trig11}
\phi_j(\theta)=\frac{1}{\sqrt{T-t}}
\left\{
\begin{matrix}
1,\ & j=0\cr\cr
\sqrt{2}{\rm sin} \left(2\pi r(\theta-t)/(T-t)\right),\ & j=2r-1\cr\cr
\sqrt{2}{\rm cos} \left(2\pi r(\theta-t)/(T-t)\right),\ & j=2r
\end{matrix}
,\right.
\end{equation}

\vspace{3mm}
\noindent
where $r=1, 2,\ldots $

Using the system of functions 
(\ref{trig11}), we have

\vspace{1mm}
$$
\int\limits_v^x\phi_{2r-1}(\theta)\psi(\theta)d\theta=
\sqrt{\frac{2}{T-t}}\int\limits_v^x
{\rm sin} \frac{2\pi r(\theta-t)}{T-t}\psi(\theta)d\theta=
$$

\vspace{1mm}
$$
=-\sqrt{\frac{T-t}{2}}\frac{1}{\pi r}\Biggl(
\psi(x){\rm cos}\frac{2\pi r(x-t)}{T-t}-
\psi(v){\rm cos}\frac{2\pi r(v-t)}{T-t}-\Biggr.
$$

\vspace{1mm}
\begin{equation}
\label{201}
\Biggl.-
\int\limits_v^x
{\rm cos} \frac{2\pi r(\theta-t)}{T-t}\psi'(\theta)d\theta\Biggr),
\end{equation}

\vspace{4mm}

$$
\int\limits_v^x\phi_{2r}(\theta)\psi(\theta)d\theta=
\sqrt{\frac{2}{T-t}}\int\limits_v^x
{\rm cos} \frac{2\pi r(\theta-t)}{T-t}\psi(\theta)d\theta=
$$

\vspace{1mm}
$$
=\sqrt{\frac{T-t}{2}}\frac{1}{\pi r}\Biggl(
\psi(x){\rm sin}\frac{2\pi r(x-t)}{T-t}-
\psi(v){\rm sin}\frac{2\pi r(v-t)}{T-t}-\Biggr.
$$

\vspace{1mm}
\begin{equation}
\label{202}
\Biggl.-
\int\limits_v^x
{\rm sin} \frac{2\pi r(\theta-t)}{T-t}\psi'(\theta)d\theta\Biggr),
\end{equation}

\vspace{5mm}
\noindent
where $\psi'(\theta)$ is a derivative of the function $\psi(\theta)$
with respect to the variable $\theta.$

Combining (\ref{201}) and (\ref{202}), we obtain for the
trigonometric case

\begin{equation}
\label{203}
\left(
\int\limits_v^x\phi_{j}(\theta)\psi(\theta)d\theta
\right)^2 \le 
\frac{C_4}{j^2},
\end{equation}

\vspace{4mm}
\noindent
where $j\in \mathbb{N},$ constant $C_4$ is independent of $j.$

From (\ref{203}) we finally have

\begin{equation}
\label{103x}
\int\limits_t^T\int\limits_t^{\tau}
\left(~\int\limits_{u}^{\tau}\phi_{j_s}(\theta)
\psi_s(\theta)d\theta\right)^2 du d\tau
\le \frac{K_4}{j_s^2},
\end{equation}

\vspace{4mm}
\noindent
where $j_s\in \mathbb{N},$ constant $K_4$ does not depend on $j_s.$

Combibing (\ref{la}), (\ref{103}) and (\ref{103x}), we obtain

\vspace{2mm}
$$
\sum_{j_1=0}^{\infty}\ldots
\sum_{j_{s-1}=0}^{\infty}
\sum_{j_s=p+1}^{\infty}\sum_{j_{s+1}=0}^{\infty}\ldots \sum_{j_k=0}^{\infty}
C_{j_k\ldots j_1}^2\le
$$

\vspace{2mm}
\begin{equation}
\label{fff}
\le L_k
\sum_{j_s=p+1}^{\infty}\frac{1}{j_s^2} \le 
\frac{L_k}{p},
\end{equation}

\vspace{5mm}
\noindent
where constant $L_k$ depends on $k$ and $T-t.$

Obviously, the case $s=k$ can be considered absolutely analogously to the
case $s<k$. Then from (\ref{aaap}) and (\ref{fff})
we obtain

\begin{equation}
\label{ddd1}
\int\limits_{[t,T]^k}K^2(t_1,\ldots,t_k)dt_1\ldots dt_k-
\sum_{j_1=0}^{p}\ldots \sum_{j_k=0}^{p}
C_{j_k\ldots j_1}^2\le \frac{G_k}{p},
\end{equation}

\vspace{3mm}
\noindent
where constant $G_k$ depends on $k$ and $T-t.$

For the further consideration we will use the estimate (\ref{2026ch1001s11}).
Using (\ref{ddd1}) and the estimate (\ref{2026ch1001s11})
for the case $p_1=\ldots=p_k=p$ and $n=2$,
we obtain

$$
{\sf M}\left\{\biggl(J[\psi^{(k)}]_{T,t}-
J[\psi^{(k)}]_{T,t}^{p,\ldots,p}\biggr)^{4}\right\}\le
$$

\vspace{1mm}
\begin{equation}
\label{fff5}
\le C_{2,k}
\left(
\int\limits_{[t,T]^k}
K^2(t_1,\ldots,t_k)
dt_1\ldots dt_k -\sum_{j_1=0}^{p}\ldots
\sum_{j_k=0}^{p}C^2_{j_k\ldots j_1}
\right)^2\le 
\frac{H_{2,k}}{p^2},
\end{equation}

\vspace{3mm}
\noindent
where 

\vspace{-2mm}
$$
C_{n,k}=(k!)^{n} (2n-1)^{nk} 
$$

\vspace{4mm}
\noindent
and $H_{2,k}=G_k^2{C}_{2,k}.$

Note the well known fact.

\vspace{2mm}

{\bf Lemma 5.}\ {\it If for the sequence of random variables
$\xi_p$ and for some
$\alpha>0$ the number series 

$$
\sum\limits_{p=1}^{\infty}{\sf M}\left\{\left|\xi_p\right|^{\alpha}\right\}
$$

\vspace{3mm}
\noindent
converges, then the sequence $\xi_p$ converges to zero w.~p.~{\rm 1}.}

Let $\alpha$ and $\xi_p$ in
Lemma 5 be chosen as follows

\vspace{1mm}
$$
\alpha=4,\ \ \ \xi_p=\biggl|J[\psi^{(k)}]_{T,t}-
J[\psi^{(k)}]_{T,t}^{p,\ldots,p}\biggr|.
$$

\vspace{3mm}

Then from (\ref{fff5}) we obtain

\vspace{1mm}
\begin{equation}
\label{qqq1}
\sum\limits_{p=1}^{\infty}
{\sf M}\left\{\biggl(J[\psi^{(k)}]_{T,t}-
J[\psi^{(k)}]_{T,t}^{p,\ldots,p}\biggr)^{4}\right\}\le
H_{2,k}\sum\limits_{p=1}^{\infty}\frac{1}{p^2}<\infty.
\end{equation}

\vspace{4mm}

Using Lemma 5, from (\ref{qqq1}) we obtain

\vspace{2mm}
$$
J[\psi^{(k)}]_{T,t}^{p,\ldots,p}\ \to \ J[\psi^{(k)}]_{T,t}\ \ \ 
\hbox{if}\ \ \ p\to \infty
$$

\vspace{5mm}
\noindent
w.\ p.\ 1, where (see Theorem 1)

\vspace{2mm}
$$
J[\psi^{(k)}]_{T,t}^{p,\ldots,p}=
\sum_{j_1=0}^{p}\ldots\sum_{j_k=0}^{p}
C_{j_k\ldots j_1}\Biggl(
\prod_{l=1}^k\zeta_{j_l}^{(i_l)}\ -
\Biggr.
$$

\vspace{3mm}
\begin{equation}
\label{kk0}
-\ \Biggl.
\hbox{\vtop{\offinterlineskip\halign{
\hfil#\hfil\cr
{\rm l.i.m.}\cr
$\stackrel{}{{}_{N\to \infty}}$\cr
}} }\sum_{(l_1,\ldots,l_k)\in {\rm G}_k}
\phi_{j_{1}}(\tau_{l_1})
\Delta{\bf w}_{\tau_{l_1}}^{(i_1)}\ldots
\phi_{j_{k}}(\tau_{l_k})
\Delta{\bf w}_{\tau_{l_k}}^{(i_k)}\Biggr)
\end{equation}

\vspace{5mm}
\noindent
or (see Theorem 2)

\vspace{1mm}
$$
J[\psi^{(k)}]_{T,t}^{p,\ldots,p}=
\sum\limits_{j_1=0}^{p}\ldots
\sum\limits_{j_k=0}^{p}
C_{j_k\ldots j_1}\Biggl(
\prod_{l=1}^k\zeta_{j_l}^{(i_l)}+\sum\limits_{r=1}^{[k/2]}
(-1)^r \times
\Biggr.
$$

\vspace{3mm}
\begin{equation}
\label{kk1}
\times
\sum_{\stackrel{(\{\{g_1, g_2\}, \ldots, 
\{g_{2r-1}, g_{2r}\}\}, \{q_1, \ldots, q_{k-2r}\})}
{{}_{\{g_1, g_2, \ldots, 
g_{2r-1}, g_{2r}, q_1, \ldots, q_{k-2r}\}=\{1, 2, \ldots, k\}}}}
\prod\limits_{s=1}^r
{\bf 1}_{\{i_{g_{{}_{2s-1}}}=~i_{g_{{}_{2s}}}\ne 0\}}
\Biggl.{\bf 1}_{\{j_{g_{{}_{2s-1}}}=~j_{g_{{}_{2s}}}\}}
\prod_{l=1}^{k-2r}\zeta_{j_{q_l}}^{(i_{q_l})}\Biggr),
\end{equation}

\vspace{7mm}
\noindent
where $i_1,\ldots,i_k=1,\ldots,m$ in (\ref{kk0}) and (\ref{kk1}).
Theorem 6 is proved.

\vspace{4mm}

{\bf Remark 3.} {\it From Theorem {\rm 3} and Lemma {\rm 4}
we obtain

\vspace{2mm}
$$
\lim\limits_{p_{q_1}\to \infty}
\varlimsup\limits_{p_{q_2}\to \infty}
\ldots\varlimsup\limits_{p_{q_k}\to\infty}
{\sf M}\left\{\left(
J[\psi^{(k)}]_{T,t}-J[\psi^{(k)}]_{T,t}^{p_1,\ldots,p_k}
\right)^2\right\}
\le 
$$

\vspace{3mm}
$$
\le k! \cdot \lim\limits_{p_{q_1}\to 0}\ldots\lim\limits_{p_{q_k}\to\infty}
\left(\int\limits_{[t,T]^k}
K^2(t_1,\ldots,t_k)
dt_1\ldots dt_k -\sum_{j_1=0}^{p_{1}}\ldots
\sum_{j_{k}=0}^{p_{k}}C^2_{j_k\ldots j_1}\right)=
$$

\vspace{3mm}
$$
= k! 
\left(\int\limits_{[t,T]^k}
K^2(t_1,\ldots,t_k)
dt_1\ldots dt_k -\sum_{j_{q_1}=0}^{\infty}\ldots
\sum_{j_{q_k}=0}^{\infty}C^2_{j_k\ldots j_1}\right)=0
$$

\vspace{5mm}
\noindent
for the following cases{\rm :}

\vspace{2mm}

{\rm 1.}\ $i_1,\ldots,i_k=1,\ldots,m$\ \ and\ \ $0<T-t<\infty,$

\vspace{2mm}

{\rm 2.}\ $i_1,\ldots,i_k=0, 1,\ldots,m,$\ \ $i_1^2+\ldots+i_k^2>0,$\ \
and\ \ $0<T-t<1.$

\vspace{4mm}
\noindent
At that, 
$(q_1,\ldots,q_k)$
is any permutation such that
$\{q_1,\ldots,q_k\}=\{1,\ldots,k\},$
$J[\psi^{(k)}]_{T,t}$ is the stochastic integral {\rm (\ref{sodom20}),}
$J[\psi^{(k)}]_{T,t}^{p_1,\ldots,p_k}$ is the 
expression on the right-hand side of {\rm (\ref{tyyy})} before
passing to the limit 

\vspace{-1mm}
$$
\hbox{\vtop{\offinterlineskip\halign{
\hfil#\hfil\cr
{\rm l.i.m.}\cr
$\stackrel{}{{}_{p_1,\ldots,p_k\to \infty}}$\cr
}} },
$$ 

\vspace{3mm}
\noindent
$\varlimsup$ means ${\rm lim\ sup};$ another 
notations are the same as in Theorem {\rm 1}.
}

\vspace{4mm}

{\bf Remark 4.} {\it Taking into account Theorem {\rm 3} and
the estimate {\rm (\ref{ddd1})}, we obtain the following
inequality

\vspace{-1mm}
\begin{equation}
\label{zsel1}
{\sf M}\left\{\left(
J[\psi^{(k)}]_{T,t}-J[\psi^{(k)}]_{T,t}^{p,\ldots,p}
\right)^2\right\}\le \frac{k! P_k (T-t)^k}{p},
\end{equation}

\vspace{6mm}
\noindent
where $i_1,\ldots,i_k=1,\ldots,m$ and constant $P_k$ depends only on $k$.
}

\vspace{4mm}

{\bf Remark 5.} {\it The estimates
{\rm (\ref{2026ch1001s11})} and 
{\rm (\ref{ddd1})} imply the following 
inequality

\vspace{1mm}
$$
{\sf M}\left\{\left(
J[\psi^{(k)}]_{T,t}-J[\psi^{(k)}]_{T,t}^{p,\ldots,p}
\right)^{2n}\right\}\le 
$$

\vspace{3mm}
\begin{equation}
\label{xyzyx1}
\le (k!)^{n} (2n-1)^{nk}\
\frac{\left(P_k\right)^n (T-t)^{nk}}{p^n},
\end{equation}

\vspace{6mm}
\noindent
where $i_1,\ldots,i_k=1,\ldots,m,$\ \ $n\in\mathbb{N},$ and 
constant $P_k$ depends only on $k$.}

\vspace{4mm}

{\bf Remark 6.} {\it Consider the question
on the rate of convergence w. p. {\rm 1} in Theorem {\rm 6.}
Using the inequality {\rm (\ref{xyzyx1}),} we obtain

\begin{equation}
\label{xyzyx11}
\left({\sf M}\left\{\left(
J[\psi^{(k)}]_{T,t}-J[\psi^{(k)}]_{T,t}^{p,\ldots,p}
\right)^{2n}\right\}\right)^{1/2n}\le \frac{Q_{n,k}}{\sqrt{p}},
\end{equation}

\vspace{6mm}
\noindent
where $i_1,\ldots,i_k=1,\ldots,m,$\ \ $n\in \mathbb{N},$ and 

\vspace{1mm}
$$
Q_{n,k}=(2n-1)^{k/2}\ \sqrt{k!}\
\sqrt{P_k}\ (T-t)^{k/2}.
$$

\vspace{5mm}

According to the Lyapunov inequality, we have

\begin{equation}
\label{xyzyx12}
\left({\sf M}\biggl\{\left(
J[\psi^{(k)}]_{T,t}-J[\psi^{(k)}]_{T,t}^{p,\ldots,p}
\right)^{n}\biggr\}\right)^{1/n}\le \frac{Q_{n,k}}{\sqrt{p}}
\end{equation}

\vspace{6mm}
\noindent
for all $n\in \mathbb{N}$. Following \cite{xyz1001} {\rm (}Lemma {\rm 2.1)}, we get

\vspace{2mm}
$$
\biggl|J[\psi^{(k)}]_{T,t}-
J[\psi^{(k)}]_{T,t}^{p,\ldots,p}\biggr|=
$$

\vspace{3mm}
$$
=
\frac{p^{1/2 - \varepsilon}}{p^{1/2 - \varepsilon}}\biggl|J[\psi^{(k)}]_{T,t}-
J[\psi^{(k)}]_{T,t}^{p,\ldots,p}\biggr|\le
$$

\vspace{3mm}
$$
\le 
\frac{1}{p^{1/2 - \varepsilon}}
\sup\limits_{p\in \mathbb{N}}\left(p^{1/2 - \varepsilon}
\biggl|J[\psi^{(k)}]_{T,t}-
J[\psi^{(k)}]_{T,t}^{p,\ldots,p}\biggr|\right)=
$$

\vspace{3mm}
\begin{equation}
\label{xyzyx13}
=\frac{\eta_{\varepsilon}}{p^{1/2 - \varepsilon}}
\end{equation}

\vspace{5mm}
\noindent
w. p. {\rm 1}, where

\vspace{-2mm}
$$
\eta_{\varepsilon}=
\sup\limits_{p\in \mathbb{N}}\left(p^{1/2 - \varepsilon}
\biggl|J[\psi^{(k)}]_{T,t}-
J[\psi^{(k)}]_{T,t}^{p,\ldots,p}\biggr|\right)
$$

\vspace{5mm}
\noindent
and $\varepsilon>0$ is fixed.

For $q>1/\varepsilon,$ $q\in \mathbb{N}$ we obtain \cite{xyz1001} {\rm (}see {\rm (\ref{xyzyx12}))}

\vspace{2mm}
$$
{\sf M}\left\{\left|\eta_{\varepsilon}\right|^q\right\}=
$$

\vspace{3mm}
$$
=
{\sf M}\left\{\left(\sup\limits_{p\in \mathbb{N}}\left(p^{1/2 - \varepsilon}
\biggl|J[\psi^{(k)}]_{T,t}-
J[\psi^{(k)}]_{T,t}^{p,\ldots,p}\biggr|\right)\right)^q\right\}=
$$

\vspace{3mm}
$$
=
{\sf M}\left\{\sup\limits_{p\in \mathbb{N}}\left(p^{(1/2 - \varepsilon)q}
\biggl|J[\psi^{(k)}]_{T,t}-
J[\psi^{(k)}]_{T,t}^{p,\ldots,p}\biggr|^q\right)\right\}\le
$$

\vspace{3mm}
$$
\le {\sf M}\left\{\sum\limits_{p=1}^{\infty}p^{(1/2 - \varepsilon)q}
\biggl|J[\psi^{(k)}]_{T,t}-
J[\psi^{(k)}]_{T,t}^{p,\ldots,p}\biggr|^q\right\}=
$$

\vspace{3mm}
$$
= \sum\limits_{p=1}^{\infty}p^{(1/2 - \varepsilon)q}
{\sf M}\left\{\biggl|J[\psi^{(k)}]_{T,t}-
J[\psi^{(k)}]_{T,t}^{p,\ldots,p}\biggr|^q\right\}\le
$$

\vspace{3mm}
$$
\le
\sum\limits_{p=1}^{\infty}p^{(1/2 - \varepsilon)q}
\frac{\left(Q_{q,k}\right)^q}{p^{q/2}}=
$$

\vspace{2mm}

\begin{equation}
\label{xyzyx14}
=
\left(Q_{q,k}\right)^q\sum\limits_{p=1}^{\infty}\frac{1}{p^{\varepsilon q}}<\infty.
\end{equation}

\vspace{5mm}

From {\rm (\ref{xyzyx13})} we obtain that for all $\varepsilon>0$
there exists a random variable $\eta_{\varepsilon}$ such that 
the inequality {\rm (\ref{xyzyx13})} is fulfilled w.~p.~{\rm 1} for all $p\in \mathbb{N}.$
Moreover, from the Lyapunov inequality and {\rm (\ref{xyzyx14})}, we obtain
${\sf M}\left\{\left|\eta_{\varepsilon}\right|^q\right\}<\infty$
for all $q\ge 1.$
}

\vspace{5mm}

\section{Conclusions}

\vspace{5mm}

Thus, we obtain the following useful possibilities and modifications
of the method based on Theorem 1 (in Sect.~15, we will consider the 
generalization of Theorems 1, 2 to the case of an arbitrary 
complete orthonormal system of functions in the space $L_2([t, T])$
and $\psi_1(\tau),$ $\ldots,\psi_k(\tau)\in L_2([t, T])).$

\vspace{3.5mm}

1. There is an explicit formula (see (\ref{ppppa})) for calculation 
of expansion coefficients 
of the iterated Ito stochastic integral (\ref{sodom20}) with any
fixed multiplicity $k$ $(k\in\mathbb{N})$. 

\vspace{3.5mm}

2. We have possibilities for exact calculation of the mean-square 
error of approximation 
of the iterated Ito stochastic integral (\ref{sodom20})
\cite{26} (also see \cite{19}-\cite{2023xxx1}).

\vspace{3.5mm}

3. Since the used
multiple Fourier series is a generalized in the sense
that it is built using various complete orthonormal
systems of functions in the space $L_2([t, T])$, then we 
have new possibilities 
for approximation --- we can 
use not only trigonometric functions as in \cite{32}-\cite{Zapad-4}
but Legendre polynomials.

\vspace{3.5mm}

4. As it turned out \cite{7}-\cite{310aaa} (also see early publications
\cite{3}-\cite{6}) 
it is more convenient to work 
with Legendre polynomials for building of approximations 
of the iterated Ito stochastic integrals (\ref{sodom20}). 
Approximations based on the Legendre polynomials essentially simpler 
than their analogues based on the trigonometric functions
\cite{7}-\cite{310aaa} (also see early publications
\cite{3}-\cite{6}).
Another advantages of the application of Legendre polynomials 
in the framework of the mentioned problem are considered
in \cite{20a}-\cite{2023xxx1}, \cite{29}, \cite{30}.

\vspace{3.5mm}

5. The approach to expansion of iterated 
stochastic integrals based on the Karhunen--Loeve expansion
of the Brownian bridge process \cite{32}-\cite{Zapad-4},
\cite{Zapad-7b}, \cite{Zapad-9}, \cite{Zapad-11}, \cite{Zapad-12a}
leads to 
iterated application of the operation of limit
transition (the operation of limit transition 
is implemented only once in Theorem 1)
starting from the 
second multiplicity (in the general case) 
and third multiplicity (for the case
$\psi_1(s), \psi_2(s), \psi_3(s)\equiv 1;$ 
$i_1, i_2, i_3=0,1,\ldots,m$)
of the iterated Ito stochastic integrals (\ref{sodom20}).
Multiple series (the operation of limit transition 
is implemented only once) are more convenient 
for approximation than the iterated ones
(iterated application of the operation of limit
transition), 
since the partial sums of multiple series converge for any possible case of  
convergence to infinity of their upper limits of summation 
(let us denote them as $p_1,\ldots, p_k$). 
For example, when
$p_1=\ldots=p_k=p\to\infty$.
For iterated series, the condition $p_1=\ldots=p_k=p\to\infty$ obviously 
does not guarantee the convergence of this series.
However, in 
\cite{32}
(Sect.~5.8, pp.~202--204), \cite{Zapad-2} (pp.~438-439),  
\cite{Zapad-4} (pp.~82-84),
\cite{Zapad-9} (pp.~263-264) the authors use 
(without rigorous proof)
the condition $p_1=p_2=p_3=p\to\infty$
within the frames of the mentioned approach
based on the Karhunen--Loeve expansion of the Brownian bridge
process \cite{33} together with the Wong--Zakai approximation
\cite{W-Z-1}-\cite{Watanabe} (see discussion
in Sect.~11 of this paper for detail).

\vspace{3.5mm}

6. As we mentioned above,
constructing the expansions of iterated
Ito stochastic integrals from Theorem 1 we 
saved all 
information about these integrals. That is why it is 
natural to expect that the mentioned expansions will converge
w. p. 1 and in the mean of degree $2n$ $(n\in\mathbb{N})$
(see Sect. 6 and 9 from this article).

\vspace{3.5mm}

7. The modification of Theorem 1 for complete 
orthonormal with weight  
$r(t_1)\ldots r(t_k)\ge 0$
systems of functions  in the space $L_2([t,T]^k)$ ($k\in\mathbb{N}$) 
(Theorems 4, 18)
as well 
as for some other types of iterated stochastic 
integrals (iterated stochastic integrals 
with respect to martingale Poisson measures and 
iterated stochastic integrals with respect 
to martingales) were obtained in \cite{20}, \cite{20a}-\cite{2023xxx1}, \cite{27a}.

\vspace{3.5mm}

8. The adaptation of Theorem 1 for iterated Stratonovich
stochastic integrals of multiplicities 1 to 8 was realized in 
\cite{12}-\cite{16}, \cite{19}-\cite{2023xxx1}, \cite{21}-\cite{23a}, 
\cite{24}, \cite{25}, \cite{27}, \cite{27aa}, \cite{28}-\cite{29aaaa}.

\vspace{3.5mm}

9. Application of Theorem 1 and Theorem 12 (see below) for the mean-square
approximation of iterated stochastic integrals 
with respect to the 
infinite-dimensional $Q$-Wiener process can be found
in the monographs \cite{20a}-\cite{2023xxx1} 
(Chapter 7) and in \cite{31a}-\cite{31qa}.

\vspace{5mm}

\section{Theorem 1 from Point
of View of the Wong--Zakai Approximation}

\vspace{5mm}

The iterated Ito stochastic integrals and solutions
of Ito stochastic differential equations
are complex and important func\-ti\-o\-nals
from the independent components ${\bf f}_{s}^{(i)},$
$i=1,\ldots,m$ of the multidimensional
Wiener process ${\bf f}_{s},$ $s\in[0, T].$
Let ${\bf f}_{s}^{(i)p},$ $p\in\mathbb{N}$ 
be some approximation of
${\bf f}_{s}^{(i)},$
$i=1,\ldots,m$.
Suppose that 
${\bf f}_{s}^{(i)p}$
converges to
${\bf f}_{s}^{(i)},$
$i=1,\ldots,m$ if $p\to\infty$ in some sense and has
differentiable sample trajectories.

A natural question arises: if we replace 
${\bf f}_{s}^{(i)}$
by ${\bf f}_{s}^{(i)p},$
$i=1,\ldots,m$ in the functionals
mentioned above, will the resulting
functionals converge to the original
functionals from the components 
${\bf f}_{s}^{(i)},$
$i=1,\ldots,m$ of the multidimentional
Wiener process ${\bf f}_{s}$?
The answere to this question is negative 
in the general case. However, 
in the pioneering works of Wong E. and Zakai M. \cite{W-Z-1},
\cite{W-Z-2},
it was shown that under the special conditions and 
for some types of approximations 
of the Wiener process the answere is affirmative
with one peculiarity: the convergence takes place 
to the iterated Stratonovich stochastic integrals
and solutions of Stratonovich stochastic differential equations
and not to iterated 
Ito stochastic integrals and solutions
of Ito stochastic differential equations.
The piecewise 
linear approximation 
as well as the regularization by convolution 
\cite{W-Z-1}-\cite{Watanabe} relate the 
mentioned types of approximations
of the Wiener process. The above approximation 
of stochastic integrals and solutions of stochastic differential equations
is often called the Wong--Zakai approximation.

Let ${\bf w}_{\tau},$ $\tau\in[0, T]$ is a random vector with 
an $m+1$ components: ${\bf w}_{\tau}^{(i)}={\bf f}_{\tau}^{(i)}$ 
for $i=1,\ldots,m$ and 
${\bf w}_{\tau}^{(0)}=\tau,$\ 
${\bf f}_{\tau}^{(i)}$ $(i=1,\ldots,m)$
are independent standard Wiener processes.

It is well known that the following representation 
takes place \cite{Lipt}, \cite{7e}

\begin{equation}
\label{um1x}
{\bf w}_{\tau}^{(i)}-{\bf w}_{t}^{(i)}=
\sum_{j=0}^{\infty}\int\limits_t^{\tau}
\phi_j(s)ds\ \zeta_j^{(i)},
\end{equation}

\vspace{3mm}
\noindent
where 
$$
\zeta_j^{(i)}=
\int\limits_t^T \phi_j(s)d{\bf w}_s^{(i)},
$$

\vspace{3mm}
\noindent
$\tau\in[t, T],$ $t\ge 0,$
$\{\phi_j(x)\}_{j=0}^{\infty}$ is an arbitrary complete 
orthonormal system of functions in the space $L_2([t, T]),$ and
$\zeta_j^{(i)}$ are independent standard Gaussian 
random variables for various $i$ or $j.$
Moreover, the series (\ref{um1x}) converges for any $\tau\in [t, T]$
in the mean-square sense.

Let ${\bf w}_{\tau}^{(i)p}-{\bf w}_{t}^{(i)p}$ be 
the mean-square approximation of the process
${\bf w}_{\tau}^{(i)}-{\bf w}_{t}^{(i)},$
which has the following form

\vspace{-3mm}
\begin{equation}
\label{um1xx}
{\bf w}_{\tau}^{(i)p}-{\bf w}_{t}^{(i)p}=
\sum_{j=0}^{p}\int\limits_t^{\tau}
\phi_j(s)ds\ \zeta_j^{(i)}.
\end{equation}

\vspace{3mm}

From (\ref{um1xx}) we obtain

\vspace{-4mm}
\begin{equation}
\label{um1xxx}
d{\bf w}_{\tau}^{(i)p}=
\sum_{j=0}^{p}
\phi_j(\tau)\zeta_j^{(i)} d\tau.
\end{equation}

\vspace{4mm}

Consider the following iterated Riemann--Stieltjes
integral

\begin{equation}
\label{um1xxxx}
\int\limits_t^T
\psi_k(t_k)\ldots \int\limits_t^{t_2}\psi_1(t_1)
d{\bf w}_{t_1}^{(i_1)p_1}\ldots d{\bf w}_{t_k}^{(i_k)p_k},
\end{equation}

\vspace{4mm}
\noindent
where $i_1,\ldots,i_k=0,1,\ldots,m,$\ \ $p_1,\ldots,p_k\in\mathbb{N},$

\begin{equation}
\label{um1xxx1}
d{\bf w}_{\tau}^{(i)p}=
\left\{\begin{matrix}
d{\bf f}_{\tau}^{(i)p}\ &\hbox{\rm for}\ \ \ i=1,\ldots,m\cr\cr\cr
d\tau^p\ &\hbox{\rm for}\ \ \ i=0
\end{matrix}
,\right.
\end{equation}

\vspace{4mm}
\noindent
and $d{\bf f}_{\tau}^{(i)p},$ $d\tau^p$ are defined by the relation (\ref{um1xxx}).

Let us substitute (\ref{um1xxx}) into (\ref{um1xxxx})

\begin{equation}
\label{um1xxxx1}
\int\limits_t^T
\psi_k(t_k)\ldots \int\limits_t^{t_2}\psi_1(t_1)
d{\bf w}_{t_1}^{(i_1)p_1}\ldots d{\bf w}_{t_k}^{(i_k)p_k}=
\sum\limits_{j_1=0}^{p_1}\ldots \sum\limits_{j_k=0}^{p_k}
C_{j_k \ldots j_1}\prod\limits_{l=1}^k \zeta_{j_l}^{(i_l)},
\end{equation}

\vspace{4mm}
\noindent
where 
$$
\zeta_j^{(i)}=\int\limits_t^T \phi_j(s)d{\bf w}_s^{(i)}
$$ 

\vspace{2mm}
\noindent
are independent standard Gaussian random variables for various 
$i$ or $j$ (in the case when $i\ne 0$),
${\bf w}_{s}^{(i)}={\bf f}_{s}^{(i)}$ for
$i=1,\ldots,m$ and 
${\bf w}_{s}^{(0)}=s,$

$$
C_{j_k \ldots j_1}=\int\limits_t^T\psi_k(t_k)\phi_{j_k}(t_k)\ldots
\int\limits_t^{t_2}
\psi_1(t_1)\phi_{j_1}(t_1)
dt_1\ldots dt_k
$$

\vspace{4mm}
\noindent
is the Fourier coefficient.

Consider the following iterated Stratonovich
stochastic integral

\vspace{1mm}
\begin{equation}
\label{str}
\int\limits_t^{*T}\psi_k(t_k) \ldots \int\limits_t^{*t_{2}}
\psi_1(t_1) d{\bf w}_{t_1}^{(i_1)}\ldots
d{\bf w}_{t_k}^{(i_k)},
\end{equation}

\vspace{3mm}
\noindent
where every $\psi_l(\tau)$ $(l=1,\ldots,k)$ is a continuously
differentialble 
nonrandom function 
on $[t,T],$ ${\bf w}_{s}^{(i)}={\bf f}_{s}^{(i)}$
for $i=1,\ldots,m$ and
${\bf w}_{s}^{(0)}=\tau$; $i_1,\ldots,i_k = 0, 1,\ldots,m.$

To best of our knowledge \cite{W-Z-1}-\cite{Watanabe}
the approximations of the Wiener process
in the Wong--Zakai approximation must satisfy fairly strong
restrictions
\cite{Watanabe}
(see Definition 7.1, pp.~480--481).
Moreover, approximations of the Wiener process that are
similar to (\ref{um1xx})
were not considered in \cite{W-Z-1}, \cite{W-Z-2}
(also see \cite{Watanabe}, Theorems 7.1, 7.2).
Therefore, the proof of analogs of Theorems 7.1 and 7.2 \cite{Watanabe}
for approximations of the Wiener 
process based on its series expansion (\ref{um1x})
should be carried out separately.
Thus, the mean-square convergence of the right-hand side
of (\ref{um1xxxx1}) to the iterated Stratonovich stochastic integral 
(\ref{str})
does not follow from the results of the papers
\cite{W-Z-1}, \cite{W-Z-2} (also see \cite{Watanabe},
Theorems 7.1, 7.2).

From the other hand, Theorems 1 from this paper 
and the theory built in Chapters
1 and 2 of the monographs \cite{20a}-\cite{2023xxx1}
can be considered as the proof of the
Wong--Zakai approximation for the iterated 
Stratonovich stochastic integrals (\ref{str}) of multiplicities 1 to 6
based on the 
iterated Riemann--Stieltjes
integrals
(\ref{um1xxxx})
and 
approximation (\ref{um1xx}) of the Wiener process.
At that, the Riemann--Stieltjes integrals 
(\ref{um1xxxx}) of multiplicities 1 to 6 converge
(according to Theorems 2.1--2.9 from \cite{20a}-\cite{2023xxx1})
to the appropriate Stratonovich 
stochastic integrals (\ref{str}). Recall that
$\{\phi_j(x)\}_{j=0}^{\infty}$ (see (\ref{um1x}), (\ref{um1xx}))
is a complete 
orthonormal system of Legendre polynomials or 
trigonometric functions 
in the space $L_2([t, T])$.

To illustrate the above reasoning, 
consider two examples for the case $k=2,$
$\psi_1(s),$ $\psi_2(s)\equiv 1;$ $i_1, i_2=1,\ldots,m.$

The first example relates to the piecewise linear approximation
of the multidimensional Wiener process (these approximations 
were considered in \cite{W-Z-1}-\cite{Watanabe}).

Let ${\bf b}_{\Delta}^{(i)}(t),$ $t\in[0, T]$ be the piecewise
linear approximation of the $i$th component ${\bf f}_t^{(i)}$
of the multidimensional standard Wiener process ${\bf f}_t,$
$t\in [0, T]$ with independent components
${\bf f}_t^{(i)},$ $i=1,\ldots,m,$ i.e.

$$
{\bf b}_{\Delta}^{(i)}(t)={\bf f}_{k\Delta}^{(i)}+
\frac{t-k\Delta}{\Delta}\Delta{\bf f}_{k\Delta}^{(i)},
$$

\vspace{3mm}
\noindent
where 

\vspace{-2mm}
$$
\Delta{\bf f}_{k\Delta}^{(i)}={\bf f}_{(k+1)\Delta}^{(i)}-
{\bf f}_{k\Delta}^{(i)},\ \ \
t\in[k\Delta, (k+1)\Delta),\ \ \ k=0, 1,\ldots, N-1.
$$

\vspace{4mm}

Note that w.~p.~1

\vspace{-1mm}
\begin{equation}
\label{pridum}
\frac{d{\bf b}_{\Delta}^{(i)}}{dt}(t)=
\frac{\Delta{\bf f}_{k\Delta}^{(i)}}{\Delta},\ \ \
t\in[k\Delta, (k+1)\Delta),\ \ \ k=0, 1,\ldots, N-1.
\end{equation}

\vspace{4mm}

Consider the following iterated Riemann--Stieltjes
integral

\vspace{1mm}
$$
\int\limits_0^T
\int\limits_0^{s}
d{\bf b}_{\Delta}^{(i_1)}(\tau)d{\bf b}_{\Delta}^{(i_2)}(s),\ \ \ 
i_1,i_2=1,\ldots,m.
$$

\vspace{4mm}

Using (\ref{pridum}) and additive property of Riemann--Stieltjes integrals, 
we can write w.~p.~1

\vspace{2mm}
$$
\int\limits_0^T
\int\limits_0^{s}
d{\bf b}_{\Delta}^{(i_1)}(\tau)d{\bf b}_{\Delta}^{(i_2)}(s)=
\int\limits_0^T
\int\limits_0^{s}
\frac{d{\bf b}_{\Delta}^{(i_1)}}{d\tau}(\tau)d\tau
\frac{d {\bf b}_{\Delta}^{(i_2)}}{d s}(s)
ds =
$$

\vspace{3mm}
$$
=
\sum\limits_{l=0}^{N-1}\int\limits_{l\Delta}^{(l+1)\Delta}
\left(
\sum\limits_{q=0}^{l-1}\int\limits_{q\Delta}^{(q+1)\Delta}
\frac{\Delta{\bf f}_{q\Delta}^{(i_1)}}{\Delta}d\tau+
\int\limits_{l\Delta}^{s}
\frac{\Delta{\bf f}_{l\Delta}^{(i_1)}}{\Delta}d\tau\right)
\frac{\Delta{\bf f}_{l\Delta}^{(i_2)}}{\Delta}ds=
$$

\vspace{3mm}
$$
=\sum\limits_{l=0}^{N-1}\sum\limits_{q=0}^{l-1}
\Delta{\bf f}_{q\Delta}^{(i_1)}
\Delta{\bf f}_{l\Delta}^{(i_2)}+
\frac{1}{\Delta^2}\sum\limits_{l=0}^{N-1}
\Delta{\bf f}_{l\Delta}^{(i_1)}
\Delta{\bf f}_{l\Delta}^{(i_2)}
\int\limits_{l\Delta}^{(l+1)\Delta}
\int\limits_{l\Delta}^{s}d\tau ds=
$$

\vspace{3mm}
\begin{equation}
\label{oh-ty}
=\sum\limits_{l=0}^{N-1}\sum\limits_{q=0}^{l-1}
\Delta{\bf f}_{q\Delta}^{(i_1)}
\Delta{\bf f}_{l\Delta}^{(i_2)}+
\frac{1}{2}\sum\limits_{l=0}^{N-1}
\Delta{\bf f}_{l\Delta}^{(i_1)}
\Delta{\bf f}_{l\Delta}^{(i_2)}.
\end{equation}

\vspace{6mm}

Using (\ref{oh-ty}) it 
is not difficult to show 
that

\vspace{1mm}
$$
\hbox{\vtop{\offinterlineskip\halign{
\hfil#\hfil\cr
{\rm l.i.m.}\cr
$\stackrel{}{{}_{N\to \infty}}$\cr
}} }
\int\limits_0^T
\int\limits_0^{s}
d{\bf b}_{\Delta}^{(i_1)}(\tau)d{\bf b}_{\Delta}^{(i_2)}(s)=
\int\limits_0^T
\int\limits_0^{s}
d{\bf f}_{\tau}^{(i_1)}d{\bf f}_{s}^{(i_2)}+
\frac{1}{2}{\bf 1}_{\{i_1=i_2\}}\int\limits_0^T ds=
$$

\vspace{3mm}
\begin{equation}
\label{uh-111}
=
\int\limits_0^{*T}
\int\limits_0^{*s}
d{\bf f}_{\tau}^{(i_1)}d{\bf f}_{s}^{(i_2)},
\end{equation}

\vspace{5mm}
\noindent
where $\Delta\to 0$ if $N\to\infty$ ($N\Delta=T$).

Obviously, (\ref{uh-111}) agrees with Theorem 7.1 (see \cite{Watanabe},
p.~486).

The next example relates to the approximation
of the Wiener process based on its series expansion
(\ref{um1x}) for $t=0$, where
$\{\phi_j(x)\}_{j=0}^{\infty}$ 
is a complete 
orthonormal system of Legendre polynomials or 
trigonometric functions 
in the space $L_2([0, T])$.

Consider the following iterated Riemann--Stieltjes
integral

\vspace{-1mm}
\begin{equation}
\label{abcd1}
\int\limits_0^T
\int\limits_0^{s}
d{\bf f}_{\tau}^{(i_1)p}d{\bf f}_{s}^{(i_2)p},\ \ \ 
i_1,i_2=1,\ldots,m,
\end{equation}

\vspace{3mm}
\noindent
where $d{\bf f}_{\tau}^{(i)p}$ is defined by the
relation
(\ref{um1xxx}).

Let us substitute (\ref{um1xxx}) into (\ref{abcd1}) 

\vspace{-1mm}
\begin{equation}
\label{set18}
\int\limits_0^T
\int\limits_0^{s}
d{\bf f}_{\tau}^{(i_1)p}d{\bf f}_{s}^{(i_2)p}=
\sum\limits_{j_1,j_2=0}^p
C_{j_2 j_1} \zeta_{j_1}^{(i_1)}\zeta_{j_2}^{(i_2)},
\end{equation}

\vspace{3mm}
\noindent
where 
$$
C_{j_2 j_1}=
\int\limits_0^T \phi_{j_2}(s)\int\limits_0^s
\phi_{j_1}(\tau)d\tau ds
$$

\vspace{3mm}
\noindent
is the Fourier coefficient; another notations 
are the same as in (\ref{um1xxxx1}).

As we noted above, approximations of the Wiener process that are
similar to (\ref{um1xx})
were not considered in \cite{W-Z-1}, \cite{W-Z-2}
(also see Theorems 7.1, 7.2 in \cite{Watanabe}).
Furthermore, the extension of the results of Theorems 7.1 and 7.2
\cite{Watanabe} to the case under consideration is
not obvious.

Nevertheless, in 
\cite{32}
(Sect.~5.8, pp.~202--204), \cite{Zapad-2} (pp.~438-439),  
\cite{Zapad-4} (pp.~82-84),
\cite{Zapad-9} (pp.~263-264) the authors use 
(without rigorous proof) the Wong--Zakai approximation
\cite{W-Z-1}-\cite{Watanabe} together with
the approximation of the Wiener process based on
its series expansion. 

On the other hand, we can apply the theory built in Chapters 1 and 2
of the monographs \cite{20a}-\cite{2023xxx1}. More precisely, 
using 
Theorems 2.1, 2.2 \cite{20a}-\cite{2023xxx1} 
we obtain from (\ref{set18}) the desired result

\vspace{-1mm}
$$
\hbox{\vtop{\offinterlineskip\halign{
\hfil#\hfil\cr
{\rm l.i.m.}\cr
$\stackrel{}{{}_{p\to \infty}}$\cr
}} }
\int\limits_0^T
\int\limits_0^{s}
d{\bf f}_{\tau}^{(i_1)p}d{\bf f}_{s}^{(i_2)p}=
\hbox{\vtop{\offinterlineskip\halign{
\hfil#\hfil\cr
{\rm l.i.m.}\cr
$\stackrel{}{{}_{p\to \infty}}$\cr
}} }
\sum\limits_{j_1,j_2=0}^p
C_{j_2 j_1} \zeta_{j_1}^{(i_1)}\zeta_{j_2}^{(i_2)}=
$$

\vspace{2mm}
\begin{equation}
\label{umen-bl}
=
\int\limits_0^{*T}
\int\limits_0^{*s}
d{\bf f}_{\tau}^{(i_1)}d{\bf f}_{s}^{(i_2)}.
\end{equation}

\vspace{5mm}

From the other hand, by Theorem 1 from this paper
(see (\ref{a2})) for the case
$k=2$ we obtain from (\ref{set18}) the following relation

$$
\hbox{\vtop{\offinterlineskip\halign{
\hfil#\hfil\cr
{\rm l.i.m.}\cr
$\stackrel{}{{}_{p\to \infty}}$\cr
}} }
\int\limits_0^T
\int\limits_0^{s}
d{\bf f}_{\tau}^{(i_1)p}d{\bf f}_{s}^{(i_2)p}=
\hbox{\vtop{\offinterlineskip\halign{
\hfil#\hfil\cr
{\rm l.i.m.}\cr
$\stackrel{}{{}_{p\to \infty}}$\cr
}} }
\sum\limits_{j_1,j_2=0}^p
C_{j_2 j_1} \zeta_{j_1}^{(i_1)}\zeta_{j_2}^{(i_2)}=
$$

\vspace{2mm}
$$
=
\hbox{\vtop{\offinterlineskip\halign{
\hfil#\hfil\cr
{\rm l.i.m.}\cr
$\stackrel{}{{}_{p\to \infty}}$\cr
}} }
\sum\limits_{j_1,j_2=0}^p
C_{j_2 j_1} \biggl(\zeta_{j_1}^{(i_1)}\zeta_{j_2}^{(i_2)}-
{\bf 1}_{\{i_1=i_2\}}{\bf 1}_{\{j_1=j_2\}}\biggr)+
{\bf 1}_{\{i_1=i_2\}}\sum\limits_{j_1=0}^{\infty}
C_{j_1 j_1}=
$$

\vspace{2mm}
\begin{equation}
\label{umen-blx}
=
\int\limits_0^T
\int\limits_0^{s}
d{\bf f}_{\tau}^{(i_1)}d{\bf f}_{s}^{(i_2)}+
{\bf 1}_{\{i_1=i_2\}}\sum\limits_{j_1=0}^{\infty}
C_{j_1 j_1}.
\end{equation}

\vspace{6mm}

Since

\vspace{-4mm}
$$
\sum\limits_{j_1=0}^{\infty}
C_{j_1 j_1}=\frac{1}{2}\sum\limits_{j_1=0}^{\infty}
\left(\int\limits_0^T \phi_j(\tau)d\tau\right)^2
=
$$

\vspace{3mm}
$$
=\frac{1}{2}
\left(\int\limits_0^T \phi_0(\tau)d\tau\right)^2=\frac{1}{2}
\int\limits_0^T ds,
$$

\vspace{6mm}
\noindent
then from the standard relation between Stratonovich and
Ito stochastic integrals and (\ref{umen-blx}) we obtain (\ref{umen-bl}).

\vspace{5mm}

\section{Modification of Theorem 1 for the Case
of the In\-teg\-ra\-tion Interval $[t, s]$ $(s\in (t, T])$ 
of Iterated Ito Sto\-chas\-tic Integrals}

\vspace{5mm}

Suppose that every $\psi_l(\tau)$ $(l=1,\ldots,k)$ is a continuous 
nonrandom
function on $[t, T]$. 
Define the following function on the hypercube $[t, T]^k$

\vspace{-1mm}
$$
\bar K(t_1,\ldots,t_k,s)={\bf 1}_{\{t_k<s\}}K(t_1,\ldots,t_k),
$$

\vspace{3mm}
\noindent
where the function $K(t_1,\ldots,t_k)$ is defined by 
(\ref{ppp}), $s\in (t, T]$ ($s$ is fixed), 
and ${\bf 1}_A$ is the indicator of the set $A.$
So we have

$$
\bar K(t_1,\ldots,t_k,s)=
{\bf 1}_{\{t_1<\ldots <t_k<s\}}\psi_1(t_1)\ldots \psi_k(t_k)=
$$

\vspace{-1mm}
\begin{equation}
\label{pppxyz}
=
\begin{cases}
\psi_1(t_1)\ldots \psi_k(t_k),\ &t_1<\ldots<t_k<s\\
~\\
~\\
0,\ &\hbox{\rm otherwise}
\end{cases},
\end{equation}

\vspace{4mm}
\noindent
where $k\ge 1, $ $t_1,\ldots,t_k\in [t, T],$ and 
$s\in (t, T]$.

Suppose that $\{\phi_j(x)\}_{j=0}^{\infty}$
is a complete orthonormal system of functions in 
the space $L_2([t, T])$.

The function $\bar K(t_1,\ldots,t_k,s)$ defined by
(\ref{pppxyz})
is piecewise continuous in the 
hypercube $[t, T]^k.$
At this situation it is well known that the generalized 
multiple Fourier series 
of $\bar K(t_1,\ldots,t_k,s)$ $\in L_2([t, T]^k)$ is converging 
to $\bar K(t_1,\ldots,t_k,s)$ in the hypercube $[t, T]^k$ in 
the mean-square sense, i.e.

\vspace{1mm}
\begin{equation}
\label{sos1zxyz}
\hbox{\vtop{\offinterlineskip\halign{
\hfil#\hfil\cr
{\rm lim}\cr
$\stackrel{}{{}_{p_1,\ldots,p_k\to \infty}}$\cr
}} }\Biggl\Vert
\bar K(t_1,\ldots,t_k,s)-
\sum_{j_1=0}^{p_1}\ldots \sum_{j_k=0}^{p_k}
C_{j_k\ldots j_1}(s)
\prod_{l=1}^{k} \phi_{j_l}(t_l)\Biggr\Vert_{L_2([t, T]^k)}=0,
\end{equation}

\vspace{4mm}
\noindent
where

\vspace{-3mm}
$$
C_{j_k\ldots j_1}(s)=\int\limits_{[t,T]^k}
\bar K(t_1,\ldots,t_k,s)\prod_{l=1}^{k}\phi_{j_l}(t_l)dt_1\ldots dt_k=
$$

\begin{equation}
\label{ppppaxyz}
=\int\limits_t^s\psi_k(t_k)\phi_{j_k}(t_k)\ldots
\int\limits_t^{t_2}
\psi_1(t_1)\phi_{j_1}(t_1)
dt_1\ldots dt_k
\end{equation}

\vspace{5mm}
\noindent
is the Fourier coefficient, and

$$
\left\Vert f\right\Vert_{L_2([t, T]^k)}=\left(\int\limits_{[t,T]^k}
f^2(t_1,\ldots,t_k)dt_1\ldots dt_k\right)^{1/2}.
$$

\vspace{5mm}

Note that

\vspace{-2mm}
\begin{equation}
\label{opr22}
J[\psi^{(k)}]_{s,t}=\int\limits_t^s\psi_k(t_k) \ldots \int\limits_t^{t_{2}}
\psi_1(t_1) d{\bf w}_{t_1}^{(i_1)}\ldots
d{\bf w}_{t_k}^{(i_k)}=
\end{equation}

\vspace{1mm}
$$
=
\int\limits_t^T {\bf 1}_{\{t_k<s\}}\psi_k(t_k) \ldots \int\limits_t^{t_{2}}
\psi_1(t_1) d{\bf w}_{t_1}^{(i_1)}\ldots
d{\bf w}_{t_k}^{(i_k)}\ \ \ \hbox{w.~p.~1},
$$

\vspace{5mm}
\noindent
where $s\in (t, T]$ ($s$ is fixed), $i_1,\ldots,i_k=0,1,\ldots,m.$

Consider the partition $\{\tau_j\}_{j=0}^N$ of $[t,T]$ such that

\vspace{-2mm}
\begin{equation}
\label{1111xxx}
t=\tau_0<\ldots <\tau_N=T,\ \ \
\Delta_N=
\hbox{\vtop{\offinterlineskip\halign{
\hfil#\hfil\cr
{\rm max}\cr
$\stackrel{}{{}_{0\le j\le N-1}}$\cr
}} }\Delta\tau_j\to 0\ \ \hbox{if}\ \ N\to \infty,\ \ \ 
\Delta\tau_j=\tau_{j+1}-\tau_j.
\end{equation}

\vspace{4mm}

{\bf Theorem 7} \cite{20a}, \cite{zzzzzs1}, \cite{2023xxx1}.\
{\it Suppose that
every $\psi_l(\tau)$ $(l=$ $1,\ldots, k)$ is a continuous 
non\-ran\-dom function on 
$[t, T]$ and
$\{\phi_j(x)\}_{j=0}^{\infty}$ is a complete orthonormal system  
of functions in the space $L_2([t,T]),$ 
each function $\phi_j(x)$ of which 
for finite $j$ satisfies the condition 
$(\star)$ {\rm (}see Sect.~{\rm 4)}.
Then

\vspace{1mm}
$$
J[\psi^{(k)}]_{s,t} =
\hbox{\vtop{\offinterlineskip\halign{
\hfil#\hfil\cr
{\rm l.i.m.}\cr
$\stackrel{}{{}_{p_1,\ldots,p_k\to \infty}}$\cr
}} }\sum_{j_1=0}^{p_1}\ldots\sum_{j_k=0}^{p_k}
C_{j_k\ldots j_1}(s)\Biggl(
\prod_{l=1}^k\zeta_{j_l}^{(i_l)} -
\Biggr.
$$

\vspace{3mm}
\begin{equation}
\label{form1}
-\Biggl.
\hbox{\vtop{\offinterlineskip\halign{
\hfil#\hfil\cr
{\rm l.i.m.}\cr
$\stackrel{}{{}_{N\to \infty}}$\cr
}} }\sum_{(l_1,\ldots,l_k)\in {\rm G}_k}
\phi_{j_{1}}(\tau_{l_1})
\Delta{\bf w}_{\tau_{l_1}}^{(i_1)}\ldots
\phi_{j_{k}}(\tau_{l_k})
\Delta{\bf w}_{\tau_{l_k}}^{(i_k)}\Biggr),
\end{equation}

\vspace{6mm}
\noindent
where $J[\psi^{(k)}]_{s,t}$ is the iterated Ito
stochastic integral {\rm (\ref{opr22}),}
$s\in (t, T]$ {\rm ($s$ is fixed),}

$$
{\rm G}_k={\rm H}_k\backslash{\rm L}_k,\ \ \
{\rm H}_k=\bigl\{(l_1,\ldots,l_k):\ l_1,\ldots,l_k=0,\ 1,\ldots,N-1\bigr\},
$$

\vspace{-1mm}
$$
{\rm L}_k=\bigl\{(l_1,\ldots,l_k):\ l_1,\ldots,l_k=0,\ 1,\ldots,N-1;\
l_g\ne l_r\ (g\ne r);\ g, r=1,\ldots,k\bigr\},
$$

\vspace{5mm}
\noindent
${\rm l.i.m.}$ is a limit in the mean-square sense,
$i_1,\ldots,i_k=0,1,\ldots,m,$

$$
\zeta_{j}^{(i)}=
\int\limits_t^T \phi_{j}(s) d{\bf w}_s^{(i)}
$$

\vspace{3mm}
\noindent
are independent standard Gaussian random variables
for various
$i$ or $j$ {\rm(}in the case when $i\ne 0${\rm),}
$C_{j_k\ldots j_1}(s)$ is the Fourier coefficient {\rm(\ref{ppppaxyz}),}
$\Delta{\bf w}_{\tau_{j}}^{(i)}=
{\bf w}_{\tau_{j+1}}^{(i)}-{\bf w}_{\tau_{j}}^{(i)}$
$(i=0,\ 1,\ldots,m),$\
$\left\{\tau_{j}\right\}_{j=0}^{N}$ is a partition of
$[t,T],$ which satisfies the condition {\rm (\ref{1111xxx})}.}

\vspace{2mm}

{\bf Proof.}\ Let us consider the multiple 
stochastic integrals (\ref{30.34}), (\ref{mult11}).
We will write $J[\Phi]_{s,t}^{(k)}$ and $J'[\Phi]_{s,t}^{(k)}$ $(s\in (t, T],$ $s$ is fixed)
if the function
$\Phi(t_1,\ldots,t_k)$ in (\ref{30.34}) and (\ref{mult11})
is replaced by the function ${\bf 1}_{\{t_1,\ldots,t_k<s\}}\Phi(t_1,\ldots,t_k).$

By analogy with (\ref{pobeda}), we have

\begin{equation}
\label{pobedaxxx}
J'[\Phi]_{s,t}^{(k)}=
\int\limits_t^T\ldots \int\limits_t^{t_2}
{\bf 1}_{\{t_k<s\}} \sum\limits_{(t_1,\ldots,t_k)}\biggl(
\Phi(t_1,\ldots,t_k)
d{\bf w}_{t_1}^{(i_1)}\ldots
d{\bf w}_{t_k}^{(i_k)}\biggr)\ \ \ \hbox{w.\ p.\ 1},
\end{equation}

\vspace{3mm}
\noindent
where $J'[\Phi]_{s,t}^{(k)}$ is defined by (\ref{mult11}) and

\vspace{-1mm}
$$
\sum\limits_{(t_1,\ldots,t_k)}
$$ 

\vspace{3mm}
\noindent
means the sum with respect to all
possible permutations
$(t_1,\ldots,t_k).$ 
At the same time permutations $(t_1,\ldots,t_k)$ 
when summing
are performed in (\ref{pobedaxxx}) only in the expression, which
is enclosed in pa\-ren\-the\-ses.
Moreover, 
the nonrandom function $\Phi(t_1,\ldots,t_k)$ is assumed 
to be continuous in the 
cor\-res\-pond\-ing closed domains of integration. The case
when the nonrandom function $\Phi(t_1,\ldots,t_k)$ is 
continuous in 
the open domains of integration and bounded at their boundaries is also possible.

Let us write (\ref{pobedaxxx}) as

\begin{equation}
\label{pobedaxyz1}
J'[\Phi]_{s,t}^{(k)}=
\int\limits_t^T\ldots \int\limits_t^{t_2}
\sum\limits_{(t_1,\ldots,t_k)}\biggl(
{\bf 1}_{\{t_k<s\}}
\Phi(t_1,\ldots,t_k)
d{\bf w}_{t_1}^{(i_1)}\ldots
d{\bf w}_{t_k}^{(i_k)}\biggr)\ \ \ \hbox{w.\ p.\ 1},
\end{equation}

\vspace{3mm}
\noindent
where permutations $(t_1,\ldots,t_k)$
when summing 
are performed in (\ref{pobedaxyz1}) only in the expression

$$
\Phi(t_1,\ldots,t_k)
d{\bf w}_{t_1}^{(i_1)}\ldots
d{\bf w}_{t_k}^{(i_k)}.
$$

\vspace{4mm}

It is not difficult to notice that (\ref{pobedaxxx}),
(\ref{pobedaxyz1}) can be rewritten in the form (see (\ref{pobedaxyz}))

\vspace{1mm}
\begin{equation}
\label{s2sxxx}
J'[\Phi]_{s,t}^{(k)}=\sum_{(t_1,\ldots,t_k)}
\int\limits_{t}^{T}
\ldots
\int\limits_{t}^{t_2}
\Phi(t_1,\ldots,t_k){\bf 1}_{\{t_k<s\}}d{\bf w}_{t_1}^{(i_1)}
\ldots
d{\bf w}_{t_k}^{(i_k)}\ \ \ \hbox{w.\ p.\ 1},
\end{equation}

\vspace{5mm}
\noindent
where permutations $(t_1,\ldots,t_k)$ when summing are 
performed only in the values

\vspace{1mm}
$$
{\bf 1}_{\{t_k<s\}}d{\bf w}_{t_1}^{(i_1)}
\ldots d{\bf w}_{t_k}^{(i_k)}.
$$ 

\vspace{5mm}
\noindent
At the same time the indices near 
upper 
limits of integration in the iterated stochastic integrals are changed 
correspondently and if $t_r$ swapped with $t_q$ in the  
permutation $(t_1,\ldots,t_k)$, then $i_r$ swapped with $i_q$ in 
the permutation $(i_1,\ldots,i_k)$.

According to Lemma 1, we have

\vspace{1mm}
$$
J[\psi^{(k)}]_{s,t}=
\hbox{\vtop{\offinterlineskip\halign{
\hfil#\hfil\cr
{\rm l.i.m.}\cr
$\stackrel{}{{}_{N\to \infty}}$\cr
}} }\sum_{l_k=0}^{N-1}\ldots\sum_{l_1=0}^{l_2-1}
{\bf 1}_{\{\tau_{l_k}<s\}}\psi_1(\tau_{l_1})\ldots\psi_k(\tau_{l_k})
\Delta{\bf w}_{\tau_{l_1}}^{(i_1)}
\ldots
\Delta{\bf w}_{\tau_{l_k}}^{(i_k)}=
$$

\vspace{3mm}
$$
=
\hbox{\vtop{\offinterlineskip\halign{
\hfil#\hfil\cr
{\rm l.i.m.}\cr
$\stackrel{}{{}_{N\to \infty}}$\cr
}} }\sum_{l_k=0}^{N-1}\ldots\sum_{l_1=0}^{N-1}
{\bf 1}_{\{\tau_{l_k}<s\}}K(\tau_{l_1},\ldots,\tau_{l_k})
\Delta{\bf w}_{\tau_{l_1}}^{(i_1)}
\ldots
\Delta{\bf w}_{\tau_{l_k}}^{(i_k)}=
$$

\vspace{3mm}
$$
=\hbox{\vtop{\offinterlineskip\halign{
\hfil#\hfil\cr
{\rm l.i.m.}\cr
$\stackrel{}{{}_{N\to \infty}}$\cr
}} }
\sum\limits_{\stackrel{l_1,\ldots,l_k=0}{{}_{l_q\ne l_r;\ 
q\ne r;\ q, r=1,\ldots, k}}}^{N-1}
{\bf 1}_{\{\tau_{l_k}<s\}}K(\tau_{l_1},\ldots,\tau_{l_k})
\Delta{\bf w}_{\tau_{l_1}}^{(i_1)}
\ldots
\Delta{\bf w}_{\tau_{l_k}}^{(i_k)}=
$$

\vspace{3mm}
\begin{equation}
\label{hehe100xyz}
=
\int\limits_{t}^{T}
\ldots
\int\limits_{t}^{t_2}
\sum_{(t_1,\ldots,t_k)}\left(
{\bf 1}_{\{t_k<s\}}K(t_1,\ldots,t_k)d{\bf w}_{t_1}^{(i_1)}
\ldots
d{\bf w}_{t_k}^{(i_k)}\right)\ \ \ \hbox{w.\ p.\ 1},
\end{equation}

\vspace{5mm}
\noindent
where $K(t_1,\ldots,t_k)$ is defined by 
(\ref{ppp}) and permutations 
$(t_1,\ldots,t_k)$ when summing
are performed only 
in the expression 

\vspace{1mm}
$$
K(t_1,\ldots,t_k)d{\bf w}_{t_1}^{(i_1)}
\ldots
d{\bf w}_{t_k}^{(i_k)}.
$$

\vspace{7mm}

According to Lemmas 1, 3 and (\ref{pobedaxyz}), (\ref{s2sxxx}),
(\ref{hehe100xyz}),
we get the following representation

\newpage
\noindent
$$
J[\psi^{(k)}]_{s,t}=
$$

\vspace{1mm}
$$
=
\sum_{j_1=0}^{p_1}\ldots
\sum_{j_k=0}^{p_k}
C_{j_k\ldots j_1}(s)
\int\limits_{t}^{T}
\ldots
\int\limits_{t}^{t_2}
\sum_{(t_1,\ldots,t_k)}\left(
\phi_{j_1}(t_1)
\ldots
\phi_{j_k}(t_k)
d{\bf w}_{t_1}^{(i_1)}
\ldots
d{\bf w}_{t_k}^{(i_k)}\right)
+
$$

\vspace{3mm}
$$
+R_{T,t,s}^{p_1,\ldots,p_k}=
$$

\vspace{7mm}
$$
=\sum_{j_1=0}^{p_1}\ldots
\sum_{j_k=0}^{p_k}
C_{j_k\ldots j_1}(s)\times
$$

\vspace{3mm}
$$
\times\ 
\hbox{\vtop{\offinterlineskip\halign{
\hfil#\hfil\cr
{\rm l.i.m.}\cr
$\stackrel{}{{}_{N\to \infty}}$\cr
}} }
\sum\limits_{\stackrel{l_1,\ldots,l_k=0}{{}_{l_q\ne l_r;\ 
q\ne r;\ q, r=1,\ldots, k}}}^{N-1}
\phi_{j_1}(\tau_{l_1})\ldots
\phi_{j_k}(\tau_{l_k})
\Delta{\bf w}_{\tau_{l_1}}^{(i_1)}
\ldots
\Delta{\bf w}_{\tau_{l_k}}^{(i_k)}+
$$

\vspace{3mm}
$$
+R_{T,t,s}^{p_1,\ldots,p_k}=
$$

\vspace{9mm}
$$
=\sum_{j_1=0}^{p_1}\ldots
\sum_{j_k=0}^{p_k}
C_{j_k\ldots j_1}(s)\left(
\hbox{\vtop{\offinterlineskip\halign{
\hfil#\hfil\cr
{\rm l.i.m.}\cr
$\stackrel{}{{}_{N\to \infty}}$\cr
}} }\sum_{l_1,\ldots,l_k=0}^{N-1}
\phi_{j_1}(\tau_{l_1})
\ldots
\phi_{j_k}(\tau_{l_k})
\Delta{\bf w}_{\tau_{l_1}}^{(i_1)}
\ldots
\Delta{\bf w}_{\tau_{l_k}}^{(i_k)}
-\right.
$$

\vspace{3mm}
$$
-\left.
\hbox{\vtop{\offinterlineskip\halign{
\hfil#\hfil\cr
{\rm l.i.m.}\cr
$\stackrel{}{{}_{N\to \infty}}$\cr
}} }\sum_{(l_1,\ldots,l_k)\in {\rm G}_k}
\phi_{j_{1}}(\tau_{l_1})
\Delta{\bf w}_{\tau_{l_1}}^{(i_1)}\ldots
\phi_{j_{k}}(\tau_{l_k})
\Delta{\bf w}_{\tau_{l_k}}^{(i_k)}\right)
+
$$

\vspace{3mm}
$$
+R_{T,t,s}^{p_1,\ldots,p_k}=
$$

\vspace{9mm}
$$
=\sum_{j_1=0}^{p_1}\ldots\sum_{j_k=0}^{p_k}
C_{j_k\ldots j_1}(s)\times
$$

\vspace{3mm}
$$
\times
\left(
\prod_{l=1}^k\zeta_{j_l}^{(i_l)}-
\hbox{\vtop{\offinterlineskip\halign{
\hfil#\hfil\cr
{\rm l.i.m.}\cr
$\stackrel{}{{}_{N\to \infty}}$\cr
}} }\sum_{(l_1,\ldots,l_k)\in {\rm G}_k}
\phi_{j_{1}}(\tau_{l_1})
\Delta{\bf w}_{\tau_{l_1}}^{(i_1)}\ldots
\phi_{j_{k}}(\tau_{l_k})
\Delta{\bf w}_{\tau_{l_k}}^{(i_k)}\right)+
$$

\vspace{4mm}
$$
+R_{T,t,s}^{p_1,\ldots,p_k}\ \ \ \hbox{w.\ p.\ 1},
$$

\vspace{9mm}
\noindent
where

$$
R_{T,t,s}^{p_1,\ldots,p_k}
=
$$

\vspace{3mm}
$$
=
\sum_{(t_1,\ldots,t_k)}
\int\limits_{t}^{T}
\ldots
\int\limits_{t}^{t_2}
\left({\bf 1}_{\{t_k<s\}} K(t_1,\ldots,t_k)-
\sum_{j_1=0}^{p_1}\ldots
\sum_{j_k=0}^{p_k}
C_{j_k\ldots j_1}(s)
\prod_{l=1}^k\phi_{j_l}(t_l)\right)\times
$$

\vspace{3mm}
$$
\times
d{\bf w}_{t_1}^{(i_1)}
\ldots
d{\bf w}_{t_k}^{(i_k)}=
$$

\vspace{3mm}

\begin{equation}
\label{hhhsx1}
=
\sum_{(t_1,\ldots,t_k)}
\int\limits_{t}^{T}
\ldots
\int\limits_{t}^{t_2}
K(t_1,\ldots,t_k){\bf 1}_{\{t_k<s\}}d{\bf w}_{t_1}^{(i_1)}
\ldots
d{\bf w}_{t_k}^{(i_k)}-             
\end{equation}

\vspace{3mm}
\begin{equation}
\label{y007xxx}
-
\sum_{(t_1,\ldots,t_k)}
\int\limits_{t}^{T}
\ldots
\int\limits_{t}^{t_2}
\sum_{j_1=0}^{p_1}\ldots
\sum_{j_k=0}^{p_k}
C_{j_k\ldots j_1}(s)
\prod_{l=1}^k\phi_{j_l}(t_l)
d{\bf w}_{t_1}^{(i_1)}
\ldots
d{\bf w}_{t_k}^{(i_k)}
\end{equation}

\vspace{4mm}
\noindent
w.~p.~1, where permutations $(t_1,\ldots,t_k)$ when summing 
in (\ref{hhhsx1})
are performed only 
in the values ${\bf 1}_{\{t_k<s\}}d{\bf w}_{t_1}^{(i_1)}
\ldots $
$d{\bf w}_{t_k}^{(i_k)}$. 
At the same time
permutations $(t_1,\ldots,t_k)$ when summing 
in (\ref{y007xxx})
are performed only 
in the values
$d{\bf w}_{t_1}^{(i_1)}
\ldots
d{\bf w}_{t_k}^{(i_k)}.$
Moreover, the indices near 
upper limits of integration in the iterated stochastic integrals 
in (\ref{hhhsx1}), (\ref{y007xxx})
are changed correspondently and if $t_r$ swapped with $t_q$ in the  
permutation $(t_1,\ldots,t_k)$, then $i_r$ swapped with $i_q$ in the 
permutation $(i_1,\ldots,i_k)$.

Let us estimate the remainder
$R_{T,t,s}^{p_1,\ldots,p_k}$ of the series.

According to Lemma 2 and (\ref{riemann}), we have

$$
{\sf M}\left\{\left(R_{T,t,s}^{p_1,\ldots,p_k}\right)^2\right\}
\le 
$$

\vspace{2mm}
$$
\le C_k
\sum_{(t_1,\ldots,t_k)}
\int\limits_{t}^{T}
\ldots
\int\limits_{t}^{t_2}
\left(K(t_1,\ldots,t_k){\bf 1}_{\{t_k<s\}}-
\sum_{j_1=0}^{p_1}\ldots
\sum_{j_k=0}^{p_k}
C_{j_k\ldots j_1}(s)
\prod_{l=1}^k\phi_{j_l}(t_l)\right)^2\times
$$

\begin{equation}
\label{yes999}
\times
dt_1
\ldots
dt_k,
\end{equation}

\vspace{4mm}
\noindent
where constant $C_k$ 
depends only
on the multiplicity $k$ of the iterated Ito stochastic integral
$J[\psi^{(k)}]_{s,t}$ and permutations $(t_1,\ldots,t_k)$ when summing 
in (\ref{yes999})
are performed only 
in the values ${\bf 1}_{\{t_k<s\}}$ and $dt_1\ldots dt_k.$
At the same time the indices near 
upper limits of integration in the iterated integrals 
in (\ref{yes999})
are changed correspondently.

Since 
$K(t_1,\ldots,t_k)\equiv 0$
if the condition $t_1<\ldots <t_k$ is not fulfilled, then

$$
{\sf M}\left\{\left(R_{T,t,s}^{p_1,\ldots,p_k}\right)^2\right\}
\le 
$$

\vspace{2mm}
$$
\le C_k
\sum_{(t_1,\ldots,t_k)}
\int\limits_{t}^{T}
\ldots
\int\limits_{t}^{t_2}
\left(K(t_1,\ldots,t_k){\bf 1}_{\{t_k<s\}}-
\sum_{j_1=0}^{p_1}\ldots
\sum_{j_k=0}^{p_k}
C_{j_k\ldots j_1}(s)
\prod_{l=1}^k\phi_{j_l}(t_l)\right)^2\times
$$

\begin{equation}
\label{yes9991}
\times
dt_1
\ldots
dt_k,
\end{equation} 

\vspace{4mm}
\noindent
where permutations $(t_1,\ldots,t_k)$ when summing 
in (\ref{yes9991})
are performed only 
in the values $dt_1\ldots dt_k.$
At the same time the indices near 
upper limits of integration in the iterated integrals 
in (\ref{yes9991})
are changed correspondently.

Then from (\ref{riemann}), (\ref{sos1zxyz}), and (\ref{yes9991}) we obtain

$$
{\sf M}\left\{\left(R_{T,t,s}^{p_1,\ldots,p_k}\right)^2\right\}
\le 
$$

\vspace{2mm}
$$
\le C_k
\sum_{(t_1,\ldots,t_k)}
\int\limits_{t}^{T}
\ldots
\int\limits_{t}^{t_2}
\left(K(t_1,\ldots,t_k){\bf 1}_{\{t_k<s\}}-
\sum_{j_1=0}^{p_1}\ldots
\sum_{j_k=0}^{p_k}
C_{j_k\ldots j_1}(s)
\prod_{l=1}^k\phi_{j_l}(t_l)\right)^2\times
$$

\vspace{1.8mm}
$$
\times
dt_1
\ldots
dt_k=
$$

\vspace{3mm}
$$
=C_k\int\limits_{[t,T]^k}
\left(\bar K(t_1,\ldots,t_k,s)-
\sum_{j_1=0}^{p_1}\ldots
\sum_{j_k=0}^{p_k}
C_{j_k\ldots j_1}(s)
\prod_{l=1}^k\phi_{j_l}(t_l)\right)^2
dt_1
\ldots
dt_k\to 0
$$
         
\vspace{4mm}
\noindent
if $p_1,\ldots,p_k\to\infty,$ where constant $C_k$ 
depends only
on the multiplicity $k$ of the iterated Ito stochastic integral
$J[\psi^{(k)}]_{s,t}$. 
Theorem 7 is proved.

\vspace{2mm}

{\bf Remark 7.} {\it Obviously from Theorem {\rm 7} for the case $s=T$
we obtain the variant of Theorem {\rm 1}.}

\vspace{2mm}

It is not difficult to see that for the case of pairwise different numbers
$i_1,\ldots,i_k=1,\ldots,m$ from Theorem 7 we obtain

\vspace{2mm}
$$
J[\psi^{(k)}]_{s,t}=
\hbox{\vtop{\offinterlineskip\halign{
\hfil#\hfil\cr
{\rm l.i.m.}\cr
$\stackrel{}{{}_{p_1,\ldots,p_k\to \infty}}$\cr
}} }\sum_{j_1=0}^{p_1}\ldots\sum_{j_k=0}^{p_k}
C_{j_k\ldots j_1}(s)\zeta_{j_1}^{(i_1)}\ldots \zeta_{j_k}^{(i_k)}.
$$

\vspace{7mm}

Consider particular cases of Theorem 7 for 
$k=1,\ldots,5$

\vspace{2mm}
$$
\label{a1zzz}
J[\psi^{(1)}]_{s,t}
=\hbox{\vtop{\offinterlineskip\halign{
\hfil#\hfil\cr
{\rm l.i.m.}\cr
$\stackrel{}{{}_{p_1\to \infty}}$\cr
}} }\sum_{j_1=0}^{p_1}
C_{j_1}(s)\zeta_{j_1}^{(i_1)},
$$

\vspace{5mm}

$$
J[\psi^{(2)}]_{s,t}
=\hbox{\vtop{\offinterlineskip\halign{
\hfil#\hfil\cr
{\rm l.i.m.}\cr
$\stackrel{}{{}_{p_1,p_2\to \infty}}$\cr
}} }\sum_{j_1=0}^{p_1}\sum_{j_2=0}^{p_2}
C_{j_2j_1}(s)\Biggl(\zeta_{j_1}^{(i_1)}\zeta_{j_2}^{(i_2)}
-{\bf 1}_{\{i_1=i_2\ne 0\}}
{\bf 1}_{\{j_1=j_2\}}\Biggr),
$$

\vspace{8mm}
$$
J[\psi^{(3)}]_{s,t}=
\hbox{\vtop{\offinterlineskip\halign{
\hfil#\hfil\cr
{\rm l.i.m.}\cr
$\stackrel{}{{}_{p_1,\ldots,p_3\to \infty}}$\cr
}} }\sum_{j_1=0}^{p_1}\sum_{j_2=0}^{p_2}\sum_{j_3=0}^{p_3}
C_{j_3j_2j_1}(s)\Biggl(
\zeta_{j_1}^{(i_1)}\zeta_{j_2}^{(i_2)}\zeta_{j_3}^{(i_3)}
-\Biggr.
$$

\vspace{2mm}
$$
-\Biggl.
{\bf 1}_{\{i_1=i_2\ne 0\}}
{\bf 1}_{\{j_1=j_2\}}
\zeta_{j_3}^{(i_3)}
-{\bf 1}_{\{i_2=i_3\ne 0\}}
{\bf 1}_{\{j_2=j_3\}}
\zeta_{j_1}^{(i_1)}-
{\bf 1}_{\{i_1=i_3\ne 0\}}
{\bf 1}_{\{j_1=j_3\}}
\zeta_{j_2}^{(i_2)}\Biggr),
$$

\vspace{7mm}

$$
J[\psi^{(4)}]_{s,t}
=
\hbox{\vtop{\offinterlineskip\halign{
\hfil#\hfil\cr
{\rm l.i.m.}\cr
$\stackrel{}{{}_{p_1,\ldots,p_4\to \infty}}$\cr
}} }\sum_{j_1=0}^{p_1}\ldots\sum_{j_4=0}^{p_4}
C_{j_4\ldots j_1}(s)\Biggl(
\prod_{l=1}^4\zeta_{j_l}^{(i_l)}
\Biggr.
-
$$

$$
-
{\bf 1}_{\{i_1=i_2\ne 0\}}
{\bf 1}_{\{j_1=j_2\}}
\zeta_{j_3}^{(i_3)}
\zeta_{j_4}^{(i_4)}
-
{\bf 1}_{\{i_1=i_3\ne 0\}}
{\bf 1}_{\{j_1=j_3\}}
\zeta_{j_2}^{(i_2)}
\zeta_{j_4}^{(i_4)}-
$$

$$
-
{\bf 1}_{\{i_1=i_4\ne 0\}}
{\bf 1}_{\{j_1=j_4\}}
\zeta_{j_2}^{(i_2)}
\zeta_{j_3}^{(i_3)}
-
{\bf 1}_{\{i_2=i_3\ne 0\}}
{\bf 1}_{\{j_2=j_3\}}
\zeta_{j_1}^{(i_1)}
\zeta_{j_4}^{(i_4)}-
$$

$$
-
{\bf 1}_{\{i_2=i_4\ne 0\}}
{\bf 1}_{\{j_2=j_4\}}
\zeta_{j_1}^{(i_1)}
\zeta_{j_3}^{(i_3)}
-
{\bf 1}_{\{i_3=i_4\ne 0\}}
{\bf 1}_{\{j_3=j_4\}}
\zeta_{j_1}^{(i_1)}
\zeta_{j_2}^{(i_2)}+
$$

$$
+
{\bf 1}_{\{i_1=i_2\ne 0\}}
{\bf 1}_{\{j_1=j_2\}}
{\bf 1}_{\{i_3=i_4\ne 0\}}
{\bf 1}_{\{j_3=j_4\}}
+
{\bf 1}_{\{i_1=i_3\ne 0\}}
{\bf 1}_{\{j_1=j_3\}}
{\bf 1}_{\{i_2=i_4\ne 0\}}
{\bf 1}_{\{j_2=j_4\}}+
$$
$$
+\Biggl.
{\bf 1}_{\{i_1=i_4\ne 0\}}
{\bf 1}_{\{j_1=j_4\}}
{\bf 1}_{\{i_2=i_3\ne 0\}}
{\bf 1}_{\{j_2=j_3\}}\Biggr),
$$

\vspace{7mm}
$$
J[\psi^{(5)}]_{s,t}
=\hbox{\vtop{\offinterlineskip\halign{
\hfil#\hfil\cr
{\rm l.i.m.}\cr
$\stackrel{}{{}_{p_1,\ldots,p_5\to \infty}}$\cr
}} }\sum_{j_1=0}^{p_1}\ldots\sum_{j_5=0}^{p_5}
C_{j_5\ldots j_1}(s)\Biggl(
\prod_{l=1}^5\zeta_{j_l}^{(i_l)}
-\Biggr.
$$

\vspace{-1mm}
$$
-
{\bf 1}_{\{i_1=i_2\ne 0\}}
{\bf 1}_{\{j_1=j_2\}}
\zeta_{j_3}^{(i_3)}
\zeta_{j_4}^{(i_4)}
\zeta_{j_5}^{(i_5)}-
{\bf 1}_{\{i_1=i_3\ne 0\}}
{\bf 1}_{\{j_1=j_3\}}
\zeta_{j_2}^{(i_2)}
\zeta_{j_4}^{(i_4)}
\zeta_{j_5}^{(i_5)}-
$$

\vspace{-1mm}
$$
-
{\bf 1}_{\{i_1=i_4\ne 0\}}
{\bf 1}_{\{j_1=j_4\}}
\zeta_{j_2}^{(i_2)}
\zeta_{j_3}^{(i_3)}
\zeta_{j_5}^{(i_5)}-
{\bf 1}_{\{i_1=i_5\ne 0\}}
{\bf 1}_{\{j_1=j_5\}}
\zeta_{j_2}^{(i_2)}
\zeta_{j_3}^{(i_3)}
\zeta_{j_4}^{(i_4)}-
$$

\vspace{-1mm}
$$
-
{\bf 1}_{\{i_2=i_3\ne 0\}}
{\bf 1}_{\{j_2=j_3\}}
\zeta_{j_1}^{(i_1)}
\zeta_{j_4}^{(i_4)}
\zeta_{j_5}^{(i_5)}-
{\bf 1}_{\{i_2=i_4\ne 0\}}
{\bf 1}_{\{j_2=j_4\}}
\zeta_{j_1}^{(i_1)}
\zeta_{j_3}^{(i_3)}
\zeta_{j_5}^{(i_5)}-
$$

\vspace{-1mm}
$$
-
{\bf 1}_{\{i_2=i_5\ne 0\}}
{\bf 1}_{\{j_2=j_5\}}
\zeta_{j_1}^{(i_1)}
\zeta_{j_3}^{(i_3)}
\zeta_{j_4}^{(i_4)}
-{\bf 1}_{\{i_3=i_4\ne 0\}}
{\bf 1}_{\{j_3=j_4\}}
\zeta_{j_1}^{(i_1)}
\zeta_{j_2}^{(i_2)}
\zeta_{j_5}^{(i_5)}-
$$

\vspace{-1mm}
$$
-
{\bf 1}_{\{i_3=i_5\ne 0\}}
{\bf 1}_{\{j_3=j_5\}}
\zeta_{j_1}^{(i_1)}
\zeta_{j_2}^{(i_2)}
\zeta_{j_4}^{(i_4)}
-{\bf 1}_{\{i_4=i_5\ne 0\}}
{\bf 1}_{\{j_4=j_5\}}
\zeta_{j_1}^{(i_1)}
\zeta_{j_2}^{(i_2)}
\zeta_{j_3}^{(i_3)}+
$$

\vspace{-2mm}
$$
+
{\bf 1}_{\{i_1=i_2\ne 0\}}
{\bf 1}_{\{j_1=j_2\}}
{\bf 1}_{\{i_3=i_4\ne 0\}}
{\bf 1}_{\{j_3=j_4\}}\zeta_{j_5}^{(i_5)}+
{\bf 1}_{\{i_1=i_2\ne 0\}}
{\bf 1}_{\{j_1=j_2\}}
{\bf 1}_{\{i_3=i_5\ne 0\}}
{\bf 1}_{\{j_3=j_5\}}\zeta_{j_4}^{(i_4)}+
$$

\vspace{-2mm}
$$
+
{\bf 1}_{\{i_1=i_2\ne 0\}}
{\bf 1}_{\{j_1=j_2\}}
{\bf 1}_{\{i_4=i_5\ne 0\}}
{\bf 1}_{\{j_4=j_5\}}\zeta_{j_3}^{(i_3)}+
{\bf 1}_{\{i_1=i_3\ne 0\}}
{\bf 1}_{\{j_1=j_3\}}
{\bf 1}_{\{i_2=i_4\ne 0\}}
{\bf 1}_{\{j_2=j_4\}}\zeta_{j_5}^{(i_5)}+
$$

\vspace{-2mm}
$$
+
{\bf 1}_{\{i_1=i_3\ne 0\}}
{\bf 1}_{\{j_1=j_3\}}
{\bf 1}_{\{i_2=i_5\ne 0\}}
{\bf 1}_{\{j_2=j_5\}}\zeta_{j_4}^{(i_4)}+
{\bf 1}_{\{i_1=i_3\ne 0\}}
{\bf 1}_{\{j_1=j_3\}}
{\bf 1}_{\{i_4=i_5\ne 0\}}
{\bf 1}_{\{j_4=j_5\}}\zeta_{j_2}^{(i_2)}+
$$

\vspace{-2mm}
$$
+
{\bf 1}_{\{i_1=i_4\ne 0\}}
{\bf 1}_{\{j_1=j_4\}}
{\bf 1}_{\{i_2=i_3\ne 0\}}
{\bf 1}_{\{j_2=j_3\}}\zeta_{j_5}^{(i_5)}+
{\bf 1}_{\{i_1=i_4\ne 0\}}
{\bf 1}_{\{j_1=j_4\}}
{\bf 1}_{\{i_2=i_5\ne 0\}}
{\bf 1}_{\{j_2=j_5\}}\zeta_{j_3}^{(i_3)}+
$$

\vspace{-2mm}
$$
+
{\bf 1}_{\{i_1=i_4\ne 0\}}
{\bf 1}_{\{j_1=j_4\}}
{\bf 1}_{\{i_3=i_5\ne 0\}}
{\bf 1}_{\{j_3=j_5\}}\zeta_{j_2}^{(i_2)}+
{\bf 1}_{\{i_1=i_5\ne 0\}}
{\bf 1}_{\{j_1=j_5\}}
{\bf 1}_{\{i_2=i_3\ne 0\}}
{\bf 1}_{\{j_2=j_3\}}\zeta_{j_4}^{(i_4)}+
$$

\vspace{-2mm}
$$
+
{\bf 1}_{\{i_1=i_5\ne 0\}}
{\bf 1}_{\{j_1=j_5\}}
{\bf 1}_{\{i_2=i_4\ne 0\}}
{\bf 1}_{\{j_2=j_4\}}\zeta_{j_3}^{(i_3)}+
{\bf 1}_{\{i_1=i_5\ne 0\}}
{\bf 1}_{\{j_1=j_5\}}
{\bf 1}_{\{i_3=i_4\ne 0\}}
{\bf 1}_{\{j_3=j_4\}}\zeta_{j_2}^{(i_2)}+
$$

\vspace{-2mm}
$$
+
{\bf 1}_{\{i_2=i_3\ne 0\}}
{\bf 1}_{\{j_2=j_3\}}
{\bf 1}_{\{i_4=i_5\ne 0\}}
{\bf 1}_{\{j_4=j_5\}}\zeta_{j_1}^{(i_1)}+
{\bf 1}_{\{i_2=i_4\ne 0\}}
{\bf 1}_{\{j_2=j_4\}}
{\bf 1}_{\{i_3=i_5\ne 0\}}
{\bf 1}_{\{j_3=j_5\}}\zeta_{j_1}^{(i_1)}+
$$
$$
+\Biggl.
{\bf 1}_{\{i_2=i_5\ne 0\}}
{\bf 1}_{\{j_2=j_5\}}
{\bf 1}_{\{i_3=i_4\ne 0\}}
{\bf 1}_{\{j_3=j_4\}}\zeta_{j_1}^{(i_1)}\Biggr),
$$

\vspace{6mm}
\noindent
where ${\bf 1}_A$ is the indicator of the set $A,$
$C_{j_k\ldots j_1}(s)$ $(k=1,\ldots,5)$ has the form (\ref{ppppaxyz}),
$s\in (t, T]$ ($s$ is fixed).

\vspace{2mm}

{\bf Remark 8.}\ {\it Note that
by analogy with the proof of estimate {\rm (\ref{ddd1})} we obtain the following
inequality

\vspace{1mm}
\begin{equation}
\label{road900}
\int\limits_{[t,T]^k}\bar K^2(t_1,\ldots,t_k,s)dt_1\ldots dt_k-
\sum_{j_1=0}^{p}\ldots \sum_{j_k=0}^{p}
C^2_{j_k\ldots j_1}(s)\le \frac{G_k(s)}{p},
\end{equation}

\vspace{4mm}
\noindent
where $\bar K(t_1,\ldots,t_k,s)$ and  $C_{j_k\ldots j_1}(s)$
are defined by the equalities {\rm (\ref{pppxyz})} and {\rm (\ref{ppppaxyz})}, respectively{\rm ;}  
constant $G_k(s)$ depends on $k$ and $s-t$ $(s\in (t, T]$, $s$ is fixed{\rm ).}}

\vspace{2mm}

The following obvious modification of Theorem {\rm 3} takes place.

\vspace{2mm}

{\bf Theorem 8} \cite{20a}, \cite{zzzzzs1}, \cite{2023xxx1}.\  {\it Suppose that
every $\psi_l(\tau)$ $(l=1,\ldots, k)$ is a continuous nonrandom function on 
$[t, T]$ and
$\{\phi_j(x)\}_{j=0}^{\infty}$ is a complete orthonormal system  
of functions in the space $L_2([t,T]),$ each function $\phi_j(x)$ of which 
for finite $j$ satisfies the condition 
$(\star)$ {\rm (}see Sect.~{\rm 4)}.
Then the estimate

\vspace{1mm}
$$
{\sf M}\left\{\left(
J[\psi^{(k)}]_{s,t}-J[\psi^{(k)}]_{s,t}^{p_1,\ldots,p_k}
\right)^2\right\}
\le 
$$

\vspace{2mm}
\begin{equation}
\label{road888}
\le k!\left(\int\limits_{[t,T]^k}
\bar K^2(t_1,\ldots,t_k,s)
dt_1\ldots dt_k -\sum_{j_1=0}^{p_1}\ldots
\sum_{j_k=0}^{p_k}C^2_{j_k\ldots j_1}(s)\right)
\end{equation}

\vspace{5mm}
\noindent
is valid for the following cases{\rm :}

\vspace{2mm}

{\rm 1.}\ $i_1,\ldots,i_k=1,\ldots,m$\ \ and\ \ $0<T-t<\infty,$

\vspace{1mm}

{\rm 2.}\ $i_1,\ldots,i_k=0, 1,\ldots,m,$\ \ $i_1^2+\ldots+i_k^2>0,$\ \
and\ \ $0<T-t<1,$

\vspace{3mm}
\noindent
where $J[\psi^{(k)}]_{s,t}$ is the stochastic integral {\rm (\ref{opr22}),}
$J[\psi^{(k)}]_{s,t}^{p_1,\ldots,p_k}$ is the 
expression on the right-hand side of {\rm (\ref{form1})} before
passing to the limit 
$\hbox{\vtop{\offinterlineskip\halign{
\hfil#\hfil\cr
{\rm l.i.m.}\cr
$\stackrel{}{{}_{p_1,\ldots,p_k\to \infty}}$\cr
}} },$ $\bar K(t_1,\ldots,t_k,s)$ and  $C_{j_k\ldots j_1}(s)$
are defined by the equalities {\rm (\ref{pppxyz})} 
and {\rm (\ref{ppppaxyz})}, respectively{\rm ;}   $s\in (t, T]$ {\rm (}$s$ is fixed{\rm );} another 
notations are the same as in Theorem {\rm 1.11}.
}

\vspace{2mm}

{\bf Remark 9.}\ {\it Combining the estimates {\rm (\ref{road900})} and {\rm (\ref{road888}),}
we obtain

\begin{equation}
\label{road1888}
{\sf M}\left\{\left(
J[\psi^{(k)}]_{s,t}-J[\psi^{(k)}]_{s,t}^{p,\ldots,p}
\right)^2\right\}
\le \frac{k! P_k(s-t)^k}{p},
\end{equation}

\vspace{4mm}
\noindent
where $i_1,\ldots,i_k=1,\ldots,m,$ constant $P_k$ depends only on $k;$
another notations are the same as in {\rm (\ref{road900})} and {\rm (\ref{road888})}.}

\vspace{2mm}
                                                            
{\bf Remark 10.}\ {\it The analogue of the estimate {\rm (\ref{2026ch1001s11})} for
the iterated Ito stochastic integral {\rm (\ref{opr22})} has the
following form

\vspace{1mm}
$$
{\sf M}\left\{\left(J[\psi^{(k)}]_{s,t}-
J[\psi^{(k)}]_{s,t}^{p_1,\ldots,p_k}\right)^{2n}\right\}\le
$$

\vspace{2mm}
$$
\le
(k!)^{n} (2n-1)^{nk}\ \times
$$
\begin{equation}
\label{agent01000}
\times\ 
\left(
\int\limits_{[t,T]^k}
\bar K^2(t_1,\ldots,t_k,s)
dt_1\ldots dt_k -\sum_{j_1=0}^{p_1}\ldots
\sum_{j_k=0}^{p_k}C^2_{j_k\ldots j_1}(s)
\right)^n,
\end{equation}

\vspace{4mm}
\noindent
where $J[\psi^{(k)}]_{s,t}^{p_1,\ldots,p_k}$ is the 
expression on the right-hand side of {\rm (\ref{form1})} before
passing to the limit 
$\hbox{\vtop{\offinterlineskip\halign{
\hfil#\hfil\cr
{\rm l.i.m.}\cr
$\stackrel{}{{}_{p_1,\ldots,p_k\to \infty}}$\cr
}} },$ $\bar K(t_1,\ldots,t_k,s)$ and  $C_{j_k\ldots j_1}(s)$
are defined by the equalities {\rm (\ref{pppxyz})} and 
{\rm (\ref{ppppaxyz})}, respectively{\rm ;}   $s\in (t, T]$ {\rm (}$s$ is fixed{\rm );}
$i_1,\ldots,i_k=1,\ldots,m.$}

\vspace{2mm}

{\bf Remark 11.} {\it The estimates
{\rm (\ref{road900})} and 
{\rm (\ref{agent01000})} imply the following 
inequality

\vspace{2mm}
$$
{\sf M}\left\{\left(
J[\psi^{(k)}]_{s,t}-J[\psi^{(k)}]_{s,t}^{p,\ldots,p}
\right)^{2n}\right\}\le 
$$

\vspace{3mm}
$$
\le (k!)^{n} (2n-1)^{nk}\
\frac{\left(P_k\right)^n (s-t)^{nk}}{p^n},
$$

\vspace{7mm}
\noindent
where $i_1,\ldots,i_k=1,\ldots,m,$ $n\in \mathbb{N},$ and 
constant $P_k$ depends only on $k$.}

\vspace{5mm}

\section{Expansion of Multiple Wiener Stochastic Integral Based on 
Generalized Multiple Fourier Series}

\vspace{5mm}

Let us consider the 
multiple stochastic integral (\ref{mult11})

\begin{equation}
\label{mult11www}
\hbox{\vtop{\offinterlineskip\halign{
\hfil#\hfil\cr
{\rm l.i.m.}\cr
$\stackrel{}{{}_{N\to \infty}}$\cr
}} }
\sum\limits_{\stackrel{j_1,\ldots,j_k=0}{{}_{j_q\ne j_r;\ q\ne r;\ 
q, r=1,\ldots, k}}}^{N-1}
\Phi\left(\tau_{j_1},\ldots,\tau_{j_k}\right)
\prod\limits_{l=1}^k
\Delta{\bf w}_{\tau_{j_l}}^{(i_l)}
\stackrel{\rm def}{=}J'[\Phi]_{T,t}^{(k)},
\end{equation}

\vspace{3mm}
\noindent
where for simplicity we assume that
$\Phi(t_1,\ldots,t_k):\ [t, T]^k\to\mathbb{R}$ is a 
continuous nonrandom
function on $[t, T]^k.$ Moreover, 
$\left\{\tau_{j}\right\}_{j=0}^{N}$ is a partition of
$[t,T],$ which satisfies the condition {\rm (\ref{1111})}.

The stochastic integral with respect to the scalar standard Wiener process
($i_1=\ldots=i_k\ne 0$)
and similar to (\ref{mult11www}) was considered in \cite{ito1951} (1951)
and is called the multiple Wiener stochastic integral \cite{ito1951}.
The case $\Phi(t_1,\ldots,t_k)\in L_2([t, T]^k)$ \cite{ito1951} 
will be considered in Sect.~15.

Consider the following theorem on expansion of the multiple
Wiener stochastic integral (\ref{mult11www}) based on generalized multiple Fourier series.

\vspace{2mm} 

{\bf Theorem 9} \cite{20a}, \cite{zzzzzs1}, \cite{2023xxx1}.\
{\it Suppose that $\Phi(t_1,\ldots,t_k):\ [t, T]^k\to\mathbb{R}$ is a 
continuous nonrandom
function on $[t, T]^k$ and
$\{\phi_j(x)\}_{j=0}^{\infty}$ is a complete orthonormal system  
of functions in the space $L_2([t,T]),$ each function $\phi_j(x)$ of which 
for finite $j$ satisfies the condition 
$(\star)$ {\rm (}see Sect.~{\rm 4)}.
Then the following expansions

$$
J'[\Phi]_{T,t}^{(k)} =
\hbox{\vtop{\offinterlineskip\halign{
\hfil#\hfil\cr
{\rm l.i.m.}\cr
$\stackrel{}{{}_{p_1,\ldots,p_k\to \infty}}$\cr
}} }\sum_{j_1=0}^{p_1}\ldots\sum_{j_k=0}^{p_k}
C_{j_k\ldots j_1}\Biggl(
\prod_{l=1}^k\zeta_{j_l}^{(i_l)} -
\Biggr.
$$

\vspace{2mm}
\begin{equation}
\label{quq50}
-\Biggl.
\hbox{\vtop{\offinterlineskip\halign{
\hfil#\hfil\cr
{\rm l.i.m.}\cr
$\stackrel{}{{}_{N\to \infty}}$\cr
}} }\sum_{(l_1,\ldots,l_k)\in {\rm G}_k}
\phi_{j_{1}}(\tau_{l_1})
\Delta{\bf w}_{\tau_{l_1}}^{(i_1)}\ldots
\phi_{j_{k}}(\tau_{l_k})
\Delta{\bf w}_{\tau_{l_k}}^{(i_k)}\Biggr),
\end{equation}

\vspace{6mm}
$$
J'[\Phi]_{T,t}^{(k)}=
\hbox{\vtop{\offinterlineskip\halign{
\hfil#\hfil\cr
{\rm l.i.m.}\cr
$\stackrel{}{{}_{p_1,\ldots,p_k\to \infty}}$\cr
}} }
\sum\limits_{j_1=0}^{p_1}\ldots
\sum\limits_{j_k=0}^{p_k}
C_{j_k\ldots j_1}\Biggl(
\prod_{l=1}^k\zeta_{j_l}^{(i_l)}+\sum\limits_{r=1}^{[k/2]}
(-1)^r \times
\Biggr.
$$

\vspace{2mm}
\begin{equation}
\label{quq11}
\times
\sum_{\stackrel{(\{\{g_1, g_2\}, \ldots, 
\{g_{2r-1}, g_{2r}\}\}, \{q_1, \ldots, q_{k-2r}\})}
{{}_{\{g_1, g_2, \ldots, 
g_{2r-1}, g_{2r}, q_1, \ldots, q_{k-2r}\}=\{1, 2, \ldots, k\}}}}
\prod\limits_{s=1}^r
{\bf 1}_{\{i_{g_{{}_{2s-1}}}=~i_{g_{{}_{2s}}}\ne 0\}}
\Biggl.{\bf 1}_{\{j_{g_{{}_{2s-1}}}=~j_{g_{{}_{2s}}}\}}
\prod_{l=1}^{k-2r}\zeta_{j_{q_l}}^{(i_{q_l})}\Biggr)
\end{equation}

\vspace{6mm}
\noindent
con\-verg\-ing in the mean-square sense are valid, 
where

$$
{\rm G}_k={\rm H}_k\backslash{\rm L}_k,\ \ \
{\rm H}_k=\bigl\{(l_1,\ldots,l_k):\ l_1,\ldots,l_k=0,\ 1,\ldots,N-1\bigr\},
$$

$$
{\rm L}_k=\bigl\{(l_1,\ldots,l_k):\ l_1,\ldots,l_k=0,\ 1,\ldots,N-1;\
l_g\ne l_r\ (g\ne r);\ g, r=1,\ldots,k\bigr\},
$$

\vspace{3mm}
\noindent
${\rm l.i.m.}$ is a limit in the mean-square sense,
$i_1,\ldots,i_k=0,1,\ldots,m,$ 

$$
\zeta_{j}^{(i)}=
\int\limits_t^T \phi_{j}(s) d{\bf w}_s^{(i)}
$$

\vspace{3mm}
\noindent
are independent standard Gaussian random variables
for various
$i$ or $j$ {\rm(}in the case when $i\ne 0${\rm),}

\begin{equation}
\label{quq12}
C_{j_k\ldots j_1}=\int\limits_{[t,T]^k}
\Phi(t_1,\ldots,t_k)\prod_{l=1}^{k}\phi_{j_l}(t_l)dt_1\ldots dt_k
\end{equation}

\vspace{3mm}
\noindent
is the Fourier coefficient,
$\Delta{\bf w}_{\tau_{j}}^{(i)}=
{\bf w}_{\tau_{j+1}}^{(i)}-{\bf w}_{\tau_{j}}^{(i)}$
$(i=0,\ 1,\ldots,m),$\
$\left\{\tau_{j}\right\}_{j=0}^{N}$ is a partition of
$[t,T],$ which satisfies the condition {\rm (\ref{1111});}
$[x]$ is an integer part of a real number $x;$
another notations are the same as in Theorem {\rm 2}.}

\vspace{2mm}

{\bf Proof.}\ Using Lemma 3 and (\ref{pobeda}), 
(\ref{pobedaxyz}),
we get the following representation

\vspace{2mm}
$$
J'[\Phi]_{T,t}^{(k)}=
$$

\vspace{2mm}
$$
=
\sum_{j_1=0}^{p_1}\ldots
\sum_{j_k=0}^{p_k}
C_{j_k\ldots j_1}
\int\limits_{t}^{T}
\ldots
\int\limits_{t}^{t_2}
\sum_{(t_1,\ldots,t_k)}\left(
\phi_{j_1}(t_1)
\ldots
\phi_{j_k}(t_k)
d{\bf w}_{t_1}^{(i_1)}
\ldots
d{\bf w}_{t_k}^{(i_k)}\right)
+
$$

\vspace{4mm}
$$
+R_{T,t}^{p_1,\ldots,p_k}=
$$

\vspace{5mm}
$$
=\sum_{j_1=0}^{p_1}\ldots
\sum_{j_k=0}^{p_k}
C_{j_k\ldots j_1}\             
\hbox{\vtop{\offinterlineskip\halign{
\hfil#\hfil\cr
{\rm l.i.m.}\cr
$\stackrel{}{{}_{N\to \infty}}$\cr
}} }
\sum\limits_{\stackrel{l_1,\ldots,l_k=0}{{}_{l_q\ne l_r;\ 
q\ne r;\ q, r=1,\ldots, k}}}^{N-1}
\phi_{j_1}(\tau_{l_1})\ldots
\phi_{j_k}(\tau_{l_k})
\Delta{\bf w}_{\tau_{l_1}}^{(i_1)}
\ldots
\Delta{\bf w}_{\tau_{l_k}}^{(i_k)}+
$$

\vspace{5mm}
$$
+R_{T,t}^{p_1,\ldots,p_k}=
$$

\vspace{5mm}
$$
=\sum_{j_1=0}^{p_1}\ldots
\sum_{j_k=0}^{p_k}
C_{j_k\ldots j_1}\left(
\hbox{\vtop{\offinterlineskip\halign{
\hfil#\hfil\cr
{\rm l.i.m.}\cr
$\stackrel{}{{}_{N\to \infty}}$\cr
}} }\sum_{l_1,\ldots,l_k=0}^{N-1}
\phi_{j_1}(\tau_{l_1})
\ldots
\phi_{j_k}(\tau_{l_k})
\Delta{\bf w}_{\tau_{l_1}}^{(i_1)}
\ldots
\Delta{\bf w}_{\tau_{l_k}}^{(i_k)}
-\right.
$$

\vspace{2mm}
$$
-\left.
\hbox{\vtop{\offinterlineskip\halign{
\hfil#\hfil\cr
{\rm l.i.m.}\cr
$\stackrel{}{{}_{N\to \infty}}$\cr
}} }\sum_{(l_1,\ldots,l_k)\in {\rm G}_k}
\phi_{j_{1}}(\tau_{l_1})
\Delta{\bf w}_{\tau_{l_1}}^{(i_1)}\ldots
\phi_{j_{k}}(\tau_{l_k})
\Delta{\bf w}_{\tau_{l_k}}^{(i_k)}\right)
+
$$

\vspace{4mm}
$$
+R_{T,t}^{p_1,\ldots,p_k}=
$$

\vspace{5mm}
$$
=\sum_{j_1=0}^{p_1}\ldots\sum_{j_k=0}^{p_k}
C_{j_k\ldots j_1}\times
$$

\vspace{2mm}
$$
\times
\left(
\prod_{l=1}^k\zeta_{j_l}^{(i_l)}-
\hbox{\vtop{\offinterlineskip\halign{
\hfil#\hfil\cr
{\rm l.i.m.}\cr
$\stackrel{}{{}_{N\to \infty}}$\cr
}} }\sum_{(l_1,\ldots,l_k)\in {\rm G}_k}
\phi_{j_{1}}(\tau_{l_1})
\Delta{\bf w}_{\tau_{l_1}}^{(i_1)}\ldots
\phi_{j_{k}}(\tau_{l_k})
\Delta{\bf w}_{\tau_{l_k}}^{(i_k)}\right)+
$$

\vspace{5mm}
$$
+R_{T,t}^{p_1,\ldots,p_k}\ \ \ \hbox{w.\ p.\ 1},
$$

\vspace{7mm}
\noindent
where

\vspace{2mm}
$$
R_{T,t}^{p_1,\ldots,p_k}
=\sum_{(t_1,\ldots,t_k)}
\int\limits_{t}^{T}
\ldots
\int\limits_{t}^{t_2}
\left(\Phi(t_1,\ldots,t_k)-
\sum_{j_1=0}^{p_1}\ldots
\sum_{j_k=0}^{p_k}
C_{j_k\ldots j_1}
\prod_{l=1}^k\phi_{j_l}(t_l)\right)\times
$$

\vspace{3mm}
$$
\times
d{\bf w}_{t_1}^{(i_1)}
\ldots
d{\bf w}_{t_k}^{(i_k)},
$$

\vspace{4mm}
\noindent
where permutations $(t_1,\ldots,t_k)$ when summing are performed only 
in the values $d{\bf w}_{t_1}^{(i_1)}
\ldots $
$d{\bf w}_{t_k}^{(i_k)}$. At the same time the indices near 
upper limits of integration in the iterated stochastic integrals 
are changed correspondently and if $t_r$ swapped with $t_q$ in the  
permutation $(t_1,\ldots,t_k)$, then $i_r$ swapped with $i_q$ in the 
permutation $(i_1,\ldots,i_k)$.

Let us estimate the remainder
$R_{T,t}^{p_1,\ldots,p_k}$ of the series using Lemma 2 and (\ref{riemann}). We have

\vspace{2mm}
$$
{\sf M}\left\{\left(R_{T,t}^{p_1,\ldots,p_k}\right)^2\right\}
\le 
$$

\vspace{2mm}
$$
\le C_k
\sum_{(t_1,\ldots,t_k)}
\int\limits_{t}^{T}
\ldots
\int\limits_{t}^{t_2}
\left(\Phi(t_1,\ldots,t_k)-
\sum_{j_1=0}^{p_1}\ldots
\sum_{j_k=0}^{p_k}
C_{j_k\ldots j_1}
\prod_{l=1}^k\phi_{j_l}(t_l)\right)^2\times
$$

\vspace{2mm}
$$
\times
dt_1
\ldots
dt_k=
$$

\vspace{3mm}
$$
=C_k\int\limits_{[t,T]^k}
\left(\Phi(t_1,\ldots,t_k)-
\sum_{j_1=0}^{p_1}\ldots
\sum_{j_k=0}^{p_k}
C_{j_k\ldots j_1}
\prod_{l=1}^k\phi_{j_l}(t_l)\right)^2
dt_1
\ldots
dt_k\to 0
$$

\vspace{5mm}
\noindent
if $p_1,\ldots,p_k\to\infty,$ where constant $C_k$ 
depends only
on the multiplicity $k$ of the multiple Wiener stochastic integral
$J'[\Phi]_{T,t}^{(k)}$.
The expansion (\ref{quq50}) is proved. Using (\ref{quq50}) and Remark 2,
we get the expansion (\ref{quq11}) (see Theorem 2).
Theorem 9 is proved.

Note that particular cases of the expansion (\ref{quq11})
are determined by the equalities (\ref{a1})--(\ref{a7}), in which the Fourier
coefficient $C_{j_k\ldots j_1}$ $(k=1,\ldots,7)$
has the form (\ref{quq12}).

\vspace{5mm}

\section{Reformulation of Theorems 1, 2, and 9 Using Hermite Polynomials}

\vspace{5mm}

In \cite{Rybakov1000} it was noted that
Theorem 3.1 (\cite{ito1951}, p.~162) can be applied
to the case of multiple Wiener stochastic integral
with respect to components of the multidimensional
Wiener process. As a result, Theorems 1, 2, and 9
can be reformulated 
using Hermite polynomials. Consider this approach 
using our notations.
Note that we derive the formula (\ref{ziko20}) (see below)
in two different ways. One of
them is not based on Theorem 3.1 \cite{ito1951} (see the proof of Theorem~20 below for details).

\vspace{2mm}

{\it We will say that the condition {\rm ($\star\star$)} is fulfilled
for the multi-index $(i_1\ldots i_k)$ $(i_1,\ldots,i_k=0, 1,\ldots, m)$ if
$m_1,\ldots,m_k$ are multiplicities of the elements $i_1,\ldots,i_k,$ respectively$,$ i.e.
$$
\{i_1,\ldots, i_k\}\hspace{-0.4mm}=\hspace{-0.4mm}\{\overbrace{{i_1, \ldots, i_1}}^{m_1},
\overbrace{{i_2, \ldots, i_2}}^{m_2},
\ldots, \overbrace{{i_r, \ldots, i_r}}^{m_r}\}\ \ \ (m_{r+1}=\ldots=m_k=0),
$$

\vspace{2mm}
\noindent
where $r=1,\ldots, k,$ braces   
mean an unordered 
set, and pa\-ren\-the\-ses mean an ordered set. At that, 
$m_1+\ldots+m_k=k,$\ $m_1,\ldots, m_k=0,1,\ldots,k,$\ 
and all elements with nonzero multiplicities are pairwise different.}

In this section, we consider the case $i_1,\ldots,i_k=0, 1,\ldots, m.$
Let the condition {\rm ($\star\star$)} is fulfilled
for the mul\-ti-\-in\-dex $(i_1 \ldots i_k).$ Then

\vspace{2mm}
$$
J'\left[\phi_{j_1}\ldots \phi_{j_k}\right]_{T,t}^{(i_1\ldots i_k)}=
J'\biggl[\underbrace{\phi_{j_{g_1}}
\ldots \phi_{j_{g_{{}_{m_1}}}}}_{m_1}
\underbrace{\phi_{j_{g_{m_1+1}}}
\ldots \phi_{j_{g_{m_1+m_2}}}}_{m_2}\ldots \biggr.
$$

\begin{equation}
\label{newe11}
\biggl.\ldots
\underbrace{\phi_{j_{g_{m_1+m_2+\ldots+m_{k-1}+1}}}\ldots
\phi_{j_{g_{m_1+m_2+\ldots+m_k}}}}_{m_k}\biggr]_{T,t}^
{(\overbrace{{}_{i_1 \ldots i_1}}^{m_1}
\overbrace{{}_{i_2 \ldots i_2}}^{m_2}
\ldots \overbrace{{}_{i_k \ldots i_k}}^{m_k})}
\end{equation}

\vspace{4mm}
\noindent
w.~p.~1, where 
$J'\left[\phi_{j_1}\ldots \phi_{j_k}\right]_{T,t}^{(i_1\ldots i_k)}$
is defined by (\ref{mult11}) (also see (\ref{mult11www})),
$\Phi(t_1,\ldots,t_k)=\phi_{j_1}(t_1)\ldots $ $\phi_{j_k}(t_k),$
$\{\phi_j(x)\}_{j=0}^{\infty}$ is a complete orthonormal system  
of functions in the space $L_2([t,T]),$ each function $\phi_j(x)$ of which 
for finite $j$ satisfies the condition 
$(\star)$ (see Sect.~4),
$\{j_{g_1},\ldots,j_{g_{m_1+m_2+\ldots+m_k}}\}=\{j_1,\ldots,j_k\}$.

From (\ref{newe11}) we have 

\vspace{-2mm}
$$
J'\left[\phi_{j_{1}}\ldots \phi_{j_{k}}\right]_{T,t}^{(i_1\ldots i_k)}=
J'\left[\phi_{j_{g_1}}
\ldots \phi_{j_{g_{{}_{m_1}}}}\right]_{T,t}^{
(\hspace{0.5mm}\overbrace{{}_{i_1 \ldots i_1}}^{m_1}\hspace{0.5mm})}
\cdot 
J'\left[\phi_{j_{g_{m_1+1}}}
\ldots \phi_{j_{g_{m_1+m_2}}}\right]_{T,t}^{
(\hspace{0.5mm}\overbrace{{}_{i_2 \ldots i_2}}^{m_2}\hspace{0.5mm})}
\cdot \ldots 
$$

\begin{equation}
\label{ziko30}
\ldots \cdot 
J'\left[\phi_{j_{g_{m_1+m_2+\ldots+m_{k-1}+1}}}\ldots
\phi_{j_{g_{m_1+m_2+\ldots+m_k}}}\right]_{T,t}^{
(\hspace{0.5mm}\overbrace{{}_{i_k \ldots i_k}}^{m_k}\hspace{0.5mm})}
\end{equation}

\vspace{4mm}
\noindent
w.~p.~1, where
we suppose that 
the condition {\rm ($\star\star$)} is fulfilled
for the multi-index $(i_1 \ldots i_k)$ and

\vspace{-1mm}
\begin{equation}
\label{ziko10}
J'\left[\phi_{j_{g_{m_1+m_2+\ldots+m_{l-1}+1}}}\ldots
\phi_{j_{g_{m_1+m_2+\ldots+m_l}}}\right]_{T,t}^{
(\hspace{0.5mm}\overbrace{{}_{i_l \ldots i_l}}^{m_l}\hspace{0.5mm})}
\stackrel{\sf def}{=}1\ \ \ \hbox{for}\ \ \ m_l=0,
\end{equation}

\vspace{4mm}
\noindent
braces   
mean an unordered 
set, and pa\-ren\-the\-ses mean an ordered set.
The detailed proof of the equality (\ref{ziko30}) will be given 
in Sect.~18 (see the proof of Theorem~20).

Let us consider the following multiple Wiener stochastic integral 

\vspace{-1mm}
$$
J'\left[\phi_{j_{m_1+m_2+\ldots+m_{l-1}+1}}\ldots
\phi_{j_{m_1+m_2+\ldots+m_l}}\right]_{T,t}^{
(\hspace{0.5mm}\overbrace{{}_{i_l \ldots i_l}}^{m_l}\hspace{0.5mm})}\ \ \ (m_l>0),
$$

\vspace{4mm}
\noindent
where we suppose that 

\vspace{-2mm}
\begin{equation}
\label{ziko999}
\left\{j_{g_{m_1+m_2+\ldots+m_{l-1}+1}}, \ldots, j_{g_{m_1+m_2+\ldots+m_{l}}}
\right\}=
\biggl\{\underbrace{j_{h_{1,l}}, \ldots, j_{h_{1,l}}}_{n_{1,l}}\ \hspace{-1mm},
\underbrace{j_{h_{2,l}}, \ldots, j_{h_{2,l}}}_{n_{2,l}}\ \hspace{-1mm}, \ldots,
\underbrace{j_{h_{d_l,l}}, \ldots, j_{h_{d_l,l}}}_{n_{d_l,l}}\biggr\},
\end{equation}

\vspace{3mm}
\noindent
where
$n_{1,l}+n_{2,l}+\ldots+n_{d_l,l}=m_l;$\ \ $n_{1,l}, n_{2,l}, \ldots, n_{d_l,l}=1,\ldots, m_l;$\ \ 
$d_l=1,\ldots,m_l;$\ \ $l=1,\ldots,k.$ Note that the numbers $m_1,\ldots,m_k,$\ $g_1,\ldots,g_k$
depend on $(i_1,\ldots ,i_k)$ and the numbers
$n_{1,l},\ldots,n_{d_l,l},$\ $h_{1,l},\ldots,h_{d_l,l},$\ $d_l$
depend on $\{j_1,\ldots,j_k\}$. Moreover, 
$\left\{j_{g_1},\ldots,j_{g_k}\right\}
=\{j_1,\ldots,j_k\}.$

Using Theorem 3.1 \cite{ito1951}, we get w.~p.~1

\vspace{1mm}
$$
J'\left[\phi_{j_{g_{m_1+m_2+\ldots+m_{l-1}+1}}}\ldots
\phi_{j_{g_{m_1+m_2+\ldots+m_l}}}
\right]_{T,t}^{(\hspace{0.5mm}\overbrace{{}_{i_l \ldots i_l}}^{m_l}\hspace{0.5mm})}
=
$$

\vspace{4mm}
\begin{equation}
\label{ziko20}
=\left\{
\begin{matrix}
H_{n_{1,l}}\left(\zeta_{j_{h_{1,l}}}^{(i_l)}\right)\ldots 
H_{n_{d_l,l}}\left(\zeta_{j_{h_{d_l,l}}}^{(i_l)}\right),\ 
&\hbox{\rm if}\ \ \ 
i_l\ne 0\cr\cr
\left(\zeta_{j_{h_{1,l}}}^{(0)}\right)^{n_{1,l}}\ldots
\left(\zeta_{j_{h_{d_l,l}}}^{(0)}\right)^{n_{d_l,l}},\  &\hbox{\rm if}\ \ \ 
i_l=0
\end{matrix}\right.\ \ \ (m_l>0),
\end{equation}

\vspace{7mm}
\noindent
where $H_n(x)$ is the Hermite polynomial of degree $n$

\vspace{1mm}
$$
H_n(x)=(-1)^n e^{x^2/2} \frac{d^n}{dx^n}\left(e^{-x^2/2}\right)
$$

\vspace{2mm}
\noindent
or

\vspace{-5mm}
\begin{equation}
\label{ziko500}
H_n(x)=n!\sum\limits_{m=0}^{[n/2]}\frac{(-1)^m x^{n-2m}}{m!(n-2m)! 2^m}\ \ \ (n\in\mathbb{N}),
\end{equation}

\vspace{4mm}
\noindent
and $\zeta_j^{(i)}$ $(i=0,1,\ldots,m,\ j=0,1,\ldots)$ is defined by (\ref{rr23}). 

For example,
$$
H_0(x)=1,\ \ \ H_1(x)=x,\ \ \ H_2(x)=x^2-1,
$$
$$
H_3(x)=x^3-3x,\ \ \ 
H_4(x)=x^4-6x^2 + 3,
$$
$$
H_5(x)=x^5-10x^3 + 15x.
$$

\vspace{5mm}

From (\ref{ziko10}) and (\ref{ziko20}) we obtain w.~p.~1

\vspace{2mm}
$$
J'\left[\phi_{j_{g_{m_1+m_2+\ldots+m_{l-1}+1}}}\ldots
\phi_{j_{g_{m_1+m_2+\ldots+m_l}}}
\right]_{T,t}^{(\hspace{0.5mm}\overbrace{{}_{i_l \ldots i_l}}^{m_l}\hspace{0.5mm})}
=
$$

\vspace{5mm}
\begin{equation}
\label{ziko40}
={\bf 1}_{\{m_l=0\}}+{\bf 1}_{\{m_l>0\}}\left\{
\begin{matrix}
H_{n_{1,l}}\left(\zeta_{j_{h_{1,l}}}^{(i_l)}\right)\ldots 
H_{n_{d_l,l}}\left(\zeta_{j_{h_{d_l,l}}}^{(i_l)}\right),\ 
&\hbox{\rm if}\ \ \ 
i_l\ne 0\cr\cr
\left(\zeta_{j_{h_{1,l}}}^{(0)}\right)^{n_{1,l}}\ldots
\left(\zeta_{j_{h_{d_l,l}}}^{(0)}\right)^{n_{d_l,l}},\  &\hbox{\rm if}\ \ \ 
i_l=0
\end{matrix}\right.,
\end{equation}

\vspace{8mm}
\noindent
where ${\bf 1}_A$ denotes the indicator of the set $A$.

Using (\ref{ziko30}) and (\ref{ziko40}), we get w.~p.~1

\vspace{1mm}
$$
J'\left[\phi_{j_1}\ldots \phi_{j_k}\right]_{T,t}^{(i_1\ldots i_k)}=
$$

\vspace{4mm}
\begin{equation}
\label{ziko50}
=\prod_{l=1}^k\left({\bf 1}_{\{m_l=0\}}+{\bf 1}_{\{m_l>0\}}\left\{
\begin{matrix}
H_{n_{1,l}}\left(\zeta_{j_{h_{1,l}}}^{(i_l)}\right)\ldots 
H_{n_{d_l,l}}\left(\zeta_{j_{h_{d_l,l}}}^{(i_l)}\right),\ 
&\hbox{\rm if}\ \ \ 
i_l\ne 0\cr\cr
\left(\zeta_{j_{h_{1,l}}}^{(0)}\right)^{n_{1,l}}\ldots
\left(\zeta_{j_{h_{d_l,l}}}^{(0)}\right)^{n_{d_l,l}},\  &\hbox{\rm if}\ \ \ 
i_l=0
\end{matrix}\right.\ \right),
\end{equation}

\vspace{7mm}
\noindent
where notations are the same as in (\ref{ziko999}) and (\ref{ziko20}).

The equality (\ref{ziko50}) allows us to reformulate Theorems 1, 2, and 9
using the Hermite polynomials.
                               
\vspace{2mm}

{\bf Theorem 10} \cite{20a}, \cite{zzzzzs1}, \cite{2023xxx1}\ 
(reformulation of Theorems 1 and 2).\ {\it Suppose that
the condition {\rm ($\star\star$)} is fulfilled
for the multi-index $(i_1 \ldots i_k)$ and the condition {\rm (\ref{ziko999})} is also 
fulfilled.
Furthermore, let 
every $\psi_l(\tau)$ $(l=$ $1,\ldots, k)$ is a continuous 
nonrandom function on 
$[t, T]$ and
$\{\phi_j(x)\}_{j=0}^{\infty}$ is a complete orthonormal system  
of functions in the space $L_2([t,T]),$ each function $\phi_j(x)$ of which 
for finite $j$ satisfies the condition 
$(\star)$ {\rm (}see Sect.~{\rm 4)}.
Then the following expansion

\vspace{1mm}
$$
J[\psi^{(k)}]_{T,t}^{(i_1\ldots i_k)}=
\hbox{\vtop{\offinterlineskip\halign{
\hfil#\hfil\cr
{\rm l.i.m.}\cr
$\stackrel{}{{}_{p_1,\ldots,p_k\to \infty}}$\cr
}} }
\sum\limits_{j_1=0}^{p_1}\ldots
\sum\limits_{j_k=0}^{p_k}
C_{j_k\ldots j_1}\times
$$

\vspace{4mm}
\begin{equation}
\label{ziko800}
\times
\prod_{l=1}^k\left({\bf 1}_{\{m_l=0\}}+{\bf 1}_{\{m_l>0\}}\left\{
\begin{matrix}
H_{n_{1,l}}\left(\zeta_{j_{h_{1,l}}}^{(i_l)}\right)\ldots 
H_{n_{d_l,l}}\left(\zeta_{j_{h_{d_l,l}}}^{(i_l)}\right),\ 
&\hbox{\rm if}\ \ \ 
i_l\ne 0\cr\cr
\left(\zeta_{j_{h_{1,l}}}^{(0)}\right)^{n_{1,l}}\ldots
\left(\zeta_{j_{h_{d_l,l}}}^{(0)}\right)^{n_{d_l,l}},\  &\hbox{\rm if}\ \ \ 
i_l=0
\end{matrix}\right.\ \right)
\end{equation}

\vspace{7mm}
\noindent
con\-verg\-ing in the mean-square sense is valid,
where $J[\psi^{(k)}]_{T,t}^{(i_1\ldots i_k)}$ is the 
iterated Ito stochastic integral {\rm (\ref{sodom20});} 
$n_{1,l}+n_{2,l}+\ldots+n_{d_l,l}=m_l;$\ \ $n_{1,l}, n_{2,l}, \ldots, n_{d_l,l}=1,\ldots, m_l;$\ \ 
$d_l=1,\ldots,m_l;$\ \ $l=1,\ldots,k;$\ \ 
$m_1+\ldots+m_k=k;$\ \
the numbers $m_1,\ldots,m_k,$\ $g_1,\ldots,g_k$
depend on $(i_1,\ldots,i_k)$ and 
the numbers $n_{1,l},\ldots,n_{d_l,l},$\ $h_{1,l},\ldots,h_{d_l,l},$\ $d_l$
depend on $\{j_1,\ldots,j_k\};$ moreover$,$ $\left\{j_{g_1},\ldots,j_{g_k}\right\}
=\{j_1,\ldots,j_k\};$
$H_n(x)$ is the Hermite polynomial {\rm (\ref{ziko500});}
another
notations are the same as in Theorem {\rm 1}.}

\vspace{2mm} 

{\bf Theorem 11} \cite{20a}, \cite{zzzzzs1}, 
\cite{2023xxx1}\  (reformulation of Theorems 9).\ {\it Suppose that
the condition {\rm ($\star\star$)} is fulfilled
for the multi-index $(i_1 \ldots i_k)$ 
and the condition {\rm (\ref{ziko999})} is also 
fulfilled.
Furthermore, let
$\Phi(t_1,\ldots,t_k):\ [t, T]^k\to\mathbb{R}$ is a 
continuous nonrandom
function on $[t, T]^k$
and
$\{\phi_j(x)\}_{j=0}^{\infty}$ is a complete orthonormal system  
of functions in the space $L_2([t,T]),$ each function $\phi_j(x)$ of which 
for finite $j$ satisfies the condition 
$(\star)$ {\rm (}see Sect.~{\rm 4)}.
Then the following expansion

\vspace{2mm}
$$
J'[\Phi]_{T,t}^{(i_1\ldots i_k)}=
\hbox{\vtop{\offinterlineskip\halign{
\hfil#\hfil\cr
{\rm l.i.m.}\cr
$\stackrel{}{{}_{p_1,\ldots,p_k\to \infty}}$\cr
}} }
\sum\limits_{j_1=0}^{p_1}\ldots
\sum\limits_{j_k=0}^{p_k}
C_{j_k\ldots j_1}\times
$$

\vspace{5mm}
$$
\times
\prod_{l=1}^k\left({\bf 1}_{\{m_l=0\}}+{\bf 1}_{\{m_l>0\}}\left\{
\begin{matrix}
H_{n_{1,l}}\left(\zeta_{j_{h_{1,l}}}^{(i_l)}\right)\ldots 
H_{n_{d_l,l}}\left(\zeta_{j_{h_{d_l,l}}}^{(i_l)}\right),\ 
&\hbox{\rm if}\ \ \ 
i_l\ne 0\cr\cr
\left(\zeta_{j_{h_{1,l}}}^{(0)}\right)^{n_{1,l}}\ldots
\left(\zeta_{j_{h_{d_l,l}}}^{(0)}\right)^{n_{d_l,l}},\  &\hbox{\rm if}\ \ \ 
i_l=0
\end{matrix}\right.\ \right)
$$

\vspace{6mm}
\noindent
con\-verg\-ing in the mean-square sense is valid,
where we denote the multiple Wiener stochastic integral 
{\rm (\ref{mult11www})}
as
$J'[\Phi]_{T,t}^{(i_1\ldots i_k)};$ 
$n_{1,l}+n_{2,l}+\ldots+n_{d_l,l}=m_l;$\ \ $n_{1,l}, n_{2,l}, \ldots, n_{d_l,l}=1,\ldots, m_l;$\ \ 
$d_l=1,\ldots,m_l;$\ \ $l=1,\ldots,k;$\ \
$m_1+\ldots+m_k=k;$\ \ 
the numbers $m_1,\ldots,m_k,$\ $g_1,\ldots,g_k$
depend on $(i_1,\ldots,i_k)$ and 
the numbers $n_{1,l},\ldots,n_{d_l,l},$\ $h_{1,l},\ldots,h_{d_l,l},$\ $d_l$
depend on $\{j_1,\ldots,j_k\};$ moreover$,$ $\left\{j_{g_1},\ldots,j_{g_k}\right\}
=\{j_1,\ldots,j_k\};$
$H_n(x)$ is the Hermite polynomial {\rm (\ref{ziko500});}
another
notations are the same as in Theorem {\rm 9}.}

\vspace{2mm}

From (\ref{ziko40}) we have w.~p.~1

\vspace{1mm}
\begin{equation}
\label{ziko80}
J'[\hspace{0.5mm}\underbrace{\phi_{j_1}\ldots
\phi_{j_1}}_{k}]_{T,t}^{(\hspace{0.5mm}\overbrace{{}_{i_1 \ldots i_1}}^{k}\hspace{0.5mm})}
=
\left\{
\begin{matrix}
H_{k}\left(\zeta_{j_{1}}^{(i_1)}\right),\ 
&\hbox{\rm if}\ \ \ 
i_1\ne 0\cr\cr
\left(\zeta_{j_{1}}^{(0)}\right)^{k},\  &\hbox{\rm if}\ \ \ 
i_1=0
\end{matrix}\right.\ \ \ (k>0).
\end{equation}

\vspace{4mm}

Let us show how the relation (\ref{ziko80}) can be obtained from Theorem 2.
To prove (\ref{ziko80}) using Theorem 2 
we choose $i_1=\ldots=i_k$ and $j_1=\ldots=j_k$ $(i_1=0, 1,\ldots,m)$
in the following formula 
(this formula follows from a 
comparison of (\ref{2023abc11}) and (\ref{leto6000}) or can be obtained
using the recurrence relation (\ref{recur1}))

\vspace{-1mm}
$$
J'[\phi_{j_1}\ldots \phi_{j_k}]_{T,t}^{(i_1\ldots i_k)}=\prod_{l=1}^k\zeta_{j_l}^{(i_l)}+
\sum\limits_{r=1}^{[k/2]}
(-1)^r \times
$$

\vspace{1mm}
\begin{equation}
\label{rezo7}
\times\sum_{\stackrel{(\{\{g_1, g_2\}, \ldots, 
\{g_{2r-1}, g_{2r}\}\}, \{q_1, \ldots, q_{k-2r}\})}
{{}_{\{g_1, g_2, \ldots, 
g_{2r-1}, g_{2r}, q_1, \ldots, q_{k-2r}\}=\{1, 2, \ldots, k\}}}}
\prod\limits_{s=1}^r
{\bf 1}_{\{i_{g_{{}_{2s-1}}}=~i_{g_{{}_{2s}}}\ne 0\}}
\Biggl.{\bf 1}_{\{j_{g_{{}_{2s-1}}}=~j_{g_{{}_{2s}}}\}}
\prod_{l=1}^{k-2r}\zeta_{j_{q_l}}^{(i_{q_l})}
\end{equation}

\vspace{4mm}
\noindent
w.~p.~1, where notations are the same as in Theorem 2.

The case $i_1=0$ of (\ref{ziko80}) is obvious.
Simple combinatorial reasoning shows that 

\vspace{2mm}
$$
\sum_{\stackrel{(\{\{g_1, g_2\}, \ldots, 
\{g_{2r-1}, g_{2r}\}\}, \{q_1, \ldots, q_{k-2r}\})}
{{}_{\{g_1, g_2, \ldots, 
g_{2r-1}, g_{2r}, q_1, \ldots, q_{k-2r}\}=\{1, 2, \ldots, k\}}}}
\prod\limits_{s=1}^r
{\bf 1}_{\{i_{g_{{}_{2s-1}}}=~i_{g_{{}_{2s}}}\ne 0\}}
\Biggl.{\bf 1}_{\{j_{g_{{}_{2s-1}}}=~j_{g_{{}_{2s}}}\}}
\prod_{l=1}^{k-2r}\zeta_{j_{q_l}}^{(i_{q_l})}=
$$

\vspace{2mm}
\begin{equation}
\label{rezo15}
=\frac{C_k^2 \cdot C_{k-2}^2 \cdot \ldots \cdot 
C_{k-(r-1)2}^2}{r!} \left(\zeta_{j_{1}}^{(i_{1})}\right)^{k-2r},
\end{equation}

\vspace{5mm}
\noindent
where $i_1=\ldots=i_k,$\ \ $j_1=\ldots=j_k$\ $(i_1=1,\ldots,m),$ and

\vspace{1mm}
$$
C_k^l=\frac{k!}{l!(k-l)!}
$$

\vspace{3mm}
\noindent
is the binomial coefficient.

We have

\vspace{-2mm}
\begin{equation}
\label{rezo16}
\frac{C_k^2 \cdot C_{k-2}^2 \cdot \ldots \cdot 
C_{k-(r-1)2}^2}{r!}=\frac{k!}{r!(k-2r)!2^r}.
\end{equation}

\vspace{3mm}

Combining (\ref{rezo7}), (\ref{rezo15}), and (\ref{rezo16}), we get w.~p.~1

\vspace{2mm}
$$
J'[\underbrace{\phi_{j_1}
\ldots \phi_{j_1}}_{k}]_{T,t}^{(\hspace{0.5mm}\overbrace{{}_{i_1 \ldots i_1}}^{k}\hspace{0.5mm})}
=\left(\zeta_{j_1}^{(i_1)}\right)^k 
+ k!\sum\limits_{r=1}^{[k/2]}
\frac{(-1)^r}{r!(k-2r)!2^r}
\left(\zeta_{j_{1}}^{(i_{1})}\right)^{k-2r}=
$$

\vspace{2mm}
$$
=k!\sum\limits_{r=0}^{[k/2]}
\frac{(-1)^r}{r!(k-2r)!2^r}
\left(\zeta_{j_{1}}^{(i_{1})}\right)^{k-2r}=H_k\left(\zeta_{j_{1}}^{(i_{1})}\right).
$$

\vspace{5mm}
\noindent
The relation (\ref{ziko80}) is proved using (\ref{rezo7}).

From (\ref{ziko50}) and (\ref{rezo7}) we obtain the following equality
for multiple Wiener stochastic integral

$$
J'[\phi_{j_1}\ldots \phi_{j_k}]_{T,t}^{(i_1\ldots i_k)}=\prod_{l=1}^k\zeta_{j_l}^{(i_l)}
+\sum\limits_{r=1}^{[k/2]}
(-1)^r \times
$$

\vspace{3mm}
$$
\times\sum_{\stackrel{(\{\{g_1, g_2\}, \ldots, 
\{g_{2r-1}, g_{2r}\}\}, \{q_1, \ldots, q_{k-2r}\})}
{{}_{\{g_1, g_2, \ldots, 
g_{2r-1}, g_{2r}, q_1, \ldots, q_{k-2r}\}=\{1, 2, \ldots, k\}}}}
\prod\limits_{s=1}^r
{\bf 1}_{\{i_{g_{{}_{2s-1}}}=~i_{g_{{}_{2s}}}\ne 0\}}
\Biggl.{\bf 1}_{\{j_{g_{{}_{2s-1}}}=~j_{g_{{}_{2s}}}\}}
\prod_{l=1}^{k-2r}\zeta_{j_{q_l}}^{(i_{q_l})}=
$$

\vspace{4mm}
\begin{equation}
\label{ziko100}
=\prod_{l=1}^k\left({\bf 1}_{\{m_l=0\}}+{\bf 1}_{\{m_l>0\}}\left\{
\begin{matrix}
H_{n_{1,l}}\left(\zeta_{j_{h_{1,l}}}^{(i_l)}\right)\ldots 
H_{n_{d_l,l}}\left(\zeta_{j_{h_{d_l,l}}}^{(i_l)}\right),\ 
&\hbox{\rm if}\ \ \ 
i_l\ne 0\cr\cr
\left(\zeta_{j_{h_{1,l}}}^{(0)}\right)^{n_{1,l}}\ldots
\left(\zeta_{j_{h_{d_l,l}}}^{(0)}\right)^{n_{d_l,l}},\  &\hbox{\rm if}\ \ \ 
i_l=0
\end{matrix}\right.\ \right)
\end{equation}

\vspace{7mm}
\noindent
w.~p.~1, where notations are the same as in Theorem 2 and (\ref{ziko999}), (\ref{ziko20}).

Let us make a remark about how it is possible to obtain the formula
(\ref{ziko20}) without using Theorem 3.1 \cite{ito1951}.
Consider the set of 
polynomials 
$H_n(x,y),$ $n=0, 1,\ldots$ defined by \cite{Ch}

\vspace{1mm}
\begin{equation}
\label{new1090}
\Biggl.H_n(x,y)=\left(\frac{d^n}{d\alpha^n} 
e^{\alpha x-\alpha^2 y/2}\right)
\Biggr|_{\alpha=0}\ \ \ (H_0(x,y)\stackrel{\sf def}{=}1).
\end{equation}

\vspace{3mm}

It is well known that polynomials $H_n(x,y)$ are connected with 
the Hermite polynomials (\ref{ziko500}) by the formula \cite{Ch}

\vspace{-1mm}
\begin{equation}
\label{ziko1000}
H_n(x,y)=y^{n/2}
H_n\left(\frac{x}{\sqrt{y}}\right)=
n!\sum\limits_{i=0}^{[n/2]}\frac{(-1)^i x^{n-2i} y^i}{i!(n-2i)! 2^i}.
\end{equation}

\vspace{3mm}

For example,
$$
H_1(x,y)
=x,
$$
$$
H_2(x,y)
=x^2-y,
$$
$$
H_3(x,y)
=x^3-3xy,
$$
$$
H_4(x,y)
=x^4-6x^2 y
+3y^2,
$$
$$
H_5(x,y)=x^5-10x^3 y+15xy^2.
$$

\vspace{3mm}

From (\ref{ziko500}) and (\ref{ziko1000}) we get

\vspace{-2mm}
\begin{equation}
\label{ziko1001}
H_n(x,1)=H_n(x).
\end{equation}

\vspace{4mm}

Obviously, without loss of generality, we can write

\vspace{-1mm}
\begin{equation}
\label{ziko1002}
\left(j_1\ldots j_k\right)=
\bigl(\underbrace{j_1 \ldots j_1}_{m_1}\
\underbrace{j_2 \ldots j_2}_{m_2}\ \ldots\ 
\underbrace{j_r \ldots j_r}_{m_r}\bigr),
\end{equation}

\vspace{3mm}
\noindent
where $m_1+\ldots+m_r=k,$\ \ $m_1,\ldots, m_r=1,\ldots,k,$\ \
$r=1,\ldots,k,$\ \ $k>0,$ and $j_1,\ldots,j_r$ are pairwise different.

Analyzing the proof of Theorem 1 and using (\ref{new1010}), (\ref{new1200})
(see the proof of Theorem~20 below), we can notice that
(we suppose that the condition (\ref{ziko1002}) is fulfilled)

\vspace{1mm}
$$
J'[\phi_{j_1}\ldots \phi_{j_k}]_{T,t}^{(i_1\ldots i_1)}=
$$

\vspace{3mm}
$$
=
\hbox{\vtop{\offinterlineskip\halign{
\hfil#\hfil\cr
{\rm l.i.m.}\cr
$\stackrel{}{{}_{N\to \infty}}$\cr
}} }
\sum\limits_{\stackrel{l_1,\ldots,l_k=0}{{}_{l_q\ne l_g;\ 
q\ne g;\ q, g=1,\ldots, k}}}^{N-1}
\phi_{j_1}(\tau_{l_1})\ldots
\phi_{j_k}(\tau_{l_k})
\Delta{\bf w}_{\tau_{l_1}}^{(i_1)}
\ldots
\Delta{\bf w}_{\tau_{l_k}}^{(i_1)}=
$$

\vspace{3mm}
$$
=
\hbox{\vtop{\offinterlineskip\halign{
\hfil#\hfil\cr
{\rm l.i.m.}\cr
$\stackrel{}{{}_{N\to \infty}}$\cr
}} }
\sum\limits_{\stackrel{l_1,\ldots,l_{m_1}=0}{{}_{l_q\ne l_g;\ 
q\ne g;\ q, g=1,\ldots, m_1}}}^{N-1}
\phi_{j_1}(\tau_{l_1})\ldots
\phi_{j_1}(\tau_{l_{m_1}})
\Delta{\bf w}_{\tau_{l_1}}^{(i_1)}
\ldots
\Delta{\bf w}_{\tau_{l_{m_1}}}^{(i_1)} \times
$$

\vspace{3mm}
$$
\times
\hbox{\vtop{\offinterlineskip\halign{
\hfil#\hfil\cr
{\rm l.i.m.}\cr
$\stackrel{}{{}_{N\to \infty}}$\cr
}} }
\sum\limits_{\stackrel{l_{m_1+1},\ldots,l_{m_1+m_2}=0}{{}_{l_q\ne l_g;\ 
q\ne g;\ q, g= m_1+1,\ldots, m_1+m_2}}}^{N-1}
\phi_{j_2}(\tau_{l_{m_1+1}})\ldots
\phi_{j_2}(\tau_{l_{m_1+m_2}})
\Delta{\bf w}_{\tau_{l_{m_1+1}}}^{(i_1)}
\ldots
\Delta{\bf w}_{\tau_{l_{m_1+m_2}}}^{(i_1)} \times 
$$

\vspace{2mm}
$$
\ldots
$$

\vspace{2mm}
$$
\times
\hbox{\vtop{\offinterlineskip\halign{
\hfil#\hfil\cr
{\rm l.i.m.}\cr
$\stackrel{}{{}_{N\to \infty}}$\cr
}} }
\sum\limits_{\stackrel{l_{k-m_r+1},\ldots,l_k=0}{{}_{l_q\ne l_g;\ 
q\ne g;\ q, g=k-m_r+1,\ldots, k}}}^{N-1}
\phi_{j_r}(\tau_{l_{k-m_r+1}})\ldots
\phi_{j_r}(\tau_{l_k})
\Delta{\bf w}_{\tau_{l_{k-m_r+1}}}^{(i_1)}
\ldots
\Delta{\bf w}_{\tau_{l_k}}^{(i_1)}=
$$

\vspace{5mm}
$$
=
\hbox{\vtop{\offinterlineskip\halign{
\hfil#\hfil\cr
{\rm l.i.m.}\cr
$\stackrel{}{{}_{N\to \infty}}$\cr
}} }
\left(
\sum\limits_{l_1=0}^{N-1}
\phi_{j_1}(\tau_{l_1})\Delta{\bf w}_{\tau_{l_1}}^{(i_1)}
\ldots
\sum\limits_{l_{m_1}=0}^{N-1}
\phi_{j_1}(\tau_{l_{m_1}})
\Delta{\bf w}_{\tau_{l_{m_1}}}^{(i_1)}
-\right.
$$

\vspace{2.5mm}
$$
\left.-
\sum_{(l_1,\ldots,l_{m_1})\in {\rm G}_{1,m_1}'}
\phi_{j_1}(\tau_{l_1})\Delta{\bf w}_{\tau_{l_1}}^{(i_1)}
\ldots
\phi_{j_1}(\tau_{l_{m_1}})
\Delta{\bf w}_{\tau_{l_{m_1}}}^{(i_1)}\right)\times
$$

\vspace{5mm}

$$
\times\hbox{\vtop{\offinterlineskip\halign{
\hfil#\hfil\cr
{\rm l.i.m.}\cr
$\stackrel{}{{}_{N\to \infty}}$\cr
}} }
\left(
\sum\limits_{l_{m_1+1}=0}^{N-1}
\phi_{j_2}(\tau_{l_{m_1}+1})\Delta{\bf w}_{\tau_{l_{m_1+1}}}^{(i_1)}
\ldots
\sum\limits_{l_{m_1+m_2}=0}^{N-1}
\phi_{j_2}(\tau_{l_{m_1+m_2}})
\Delta{\bf w}_{\tau_{l_{m_1+m_2}}}^{(i_1)}
-\right.
$$

\vspace{2.5mm}
$$
\left.-
\sum_{(l_{m_1+1},\ldots,l_{m_1+m_2})\in {\rm G}_{m_1+1,m_1+m_2}'}
\phi_{j_2}(\tau_{l_{m_1}+1})\Delta{\bf w}_{\tau_{l_{m_1+1}}}^{(i_1)}
\ldots
\phi_{j_2}(\tau_{l_{m_1+m_2}})
\Delta{\bf w}_{\tau_{l_{m_1+m_2}}}^{(i_1)}
\right)\times
$$

\vspace{1.5mm}

$$
\ldots 
$$

\vspace{1.5mm}
$$
\times\hbox{\vtop{\offinterlineskip\halign{
\hfil#\hfil\cr
{\rm l.i.m.}\cr
$\stackrel{}{{}_{N\to \infty}}$\cr
}} }
\left(
\sum\limits_{l_{k-m_r+1}=0}^{N-1}
\phi_{j_r}(\tau_{l_{k-m_r+1}})
\Delta{\bf w}_{\tau_{l_{k-m_r+1}}}^{(i_1)}
\ldots
\sum\limits_{l_{k}=0}^{N-1}
\phi_{j_r}(\tau_{l_{k}})
\Delta{\bf w}_{\tau_{l_{k}}}^{(i_1)}
-\right.
$$

\vspace{2.5mm}
$$
\left.-
\sum_{(l_{k-m_r+1},\ldots,l_{k})\in {\rm G}_{k-m_r+1,k}'}
\phi_{j_r}(\tau_{l_{k-m_r+1}})\Delta{\bf w}_{\tau_{k-m_r+1}}^{(i_1)}
\ldots
\phi_{j_r}(\tau_{l_{k}})
\Delta{\bf w}_{\tau_{l_{k}}}^{(i_1)}
\right),
$$

\vspace{6mm}
\noindent
where the set 
${\rm G}_{m,n}'$ is defined according to the same rule 
as the set ${\rm G}_{k}$ in (\ref{tyyy}). 
However, the elements of the set ${\rm G}_{m,n}'$ are the numbers 
$l_m,\ldots, l_n$ $(n>m)$, while the elements of the set ${\rm G}_{k}$ are the numbers 
$l_1,\ldots, l_k$.

We have (see the proof of Theorem 1) w.~p.~1 $(i_1\ne 0)$

\vspace{2mm}

$$
\hbox{\vtop{\offinterlineskip\halign{
\hfil#\hfil\cr
{\rm l.i.m.}\cr
$\stackrel{}{{}_{N\to \infty}}$\cr
}} }
\left(\sum\limits_{l_1=0}^{N-1}
\phi_{j_1}(\tau_{l_1})\Delta{\bf w}_{\tau_{l_1}}^{(i_1)}
\ldots
\sum\limits_{l_{m_1}=0}^{N-1}
\phi_{j_1}(\tau_{l_{m_1}})
\Delta{\bf w}_{\tau_{l_{m_1}}}^{(i_1)}\right.
-
$$

\vspace{2.5mm}
$$
\left.-
\sum_{(l_1,\ldots,l_{m_1})\in {\rm G}_{1,m_1}'}
\phi_{j_1}(\tau_{l_1})\Delta{\bf w}_{\tau_{l_1}}^{(i_1)}
\ldots
\phi_{j_1}(\tau_{l_{m_1}})
\Delta{\bf w}_{\tau_{l_{m_1}}}^{(i_1)}\right)=
$$

\vspace{5mm}
$$
=\hbox{\vtop{\offinterlineskip\halign{
\hfil#\hfil\cr
{\rm l.i.m.}\cr
$\stackrel{}{{}_{N\to \infty}}$\cr
}} }
\left(\left(\sum\limits_{l_1=0}^{N-1}\phi_{j_1}(\tau_{l_1})
\Delta {\bf w}_{\tau_{l_1}}^{(i_1)}\right)^{m_1}
+\sum\limits_{r=1}^{[m_1/2]}
(-1)^r \times\right.
$$

\vspace{2.5mm}
$$
\times\sum_{\stackrel{(\{\{g_1, g_2\}, \ldots, 
\{g_{2r-1}, g_{2r}\}\}, \{q_1, \ldots, q_{m_1-2r}\})}
{{}_{\{g_1, g_2, \ldots, 
g_{2r-1}, g_{2r}, q_1, \ldots, q_{m_1-2r}\}=\{1, 2, \ldots, m_1\}}}}
\left(\sum\limits_{l_1=0}^{N-1}\phi_{j_1}^2(\tau_{l_1})
\left(\Delta {\bf w}_{\tau_{l_1}}^{(i_1)}\right)^2\right)^r\times
$$

\vspace{2.5mm}
$$
\left.\times
\left(\sum\limits_{l_1=0}^{N-1}\phi_{j_1}(\tau_{l_1})
\Delta {\bf w}_{\tau_{l_1}}^{(i_1)}\right)^{m_1-2r}\right)=
$$

\vspace{4mm}
$$
=\hbox{\vtop{\offinterlineskip\halign{
\hfil#\hfil\cr
{\rm l.i.m.}\cr
$\stackrel{}{{}_{N\to \infty}}$\cr
}} }
\left(\left(\sum\limits_{l_1=0}^{N-1}\phi_{j_1}(\tau_{l_1})
\Delta {\bf w}_{\tau_{l_1}}^{(i_1)}\right)^{m_1}
+\sum\limits_{r=1}^{[m_1/2]}
\frac{(-1)^r m_1!}{r!(m_1-2r)!2^r}
\left(\sum\limits_{l_1=0}^{N-1}\phi_{j_1}^2(\tau_{l_1})
\left(\Delta {\bf w}_{\tau_{l_1}}^{(i_1)}\right)^2\right)^r\times\right.
$$

\vspace{2.5mm}
$$
\left.\times
\left(\sum\limits_{l_1=0}^{N-1}\phi_{j_1}(\tau_{l_1})
\Delta {\bf w}_{\tau_{l_1}}^{(i_1)}\right)^{m_1-2r}\right)=
$$

\vspace{3.5mm}
$$
=\hbox{\vtop{\offinterlineskip\halign{
\hfil#\hfil\cr
{\rm l.i.m.}\cr
$\stackrel{}{{}_{N\to \infty}}$\cr
}} }
\left(\sum\limits_{r=0}^{[m_1/2]}
\frac{(-1)^r m_1!}{r!(m_1-2r)!2^r}
\left(\sum\limits_{l_1=0}^{N-1}\phi_{j_1}^2(\tau_{l_1})
\left(\Delta {\bf w}_{\tau_{l_1}}^{(i_1)}\right)^2\right)^r
\left(\sum\limits_{l_1=0}^{N-1}\phi_{j_1}(\tau_{l_1})
\Delta {\bf w}_{\tau_{l_1}}^{(i_1)}\right)^{m_1-2r}\right)=
$$

\vspace{3.5mm}
$$
=\hbox{\vtop{\offinterlineskip\halign{
\hfil#\hfil\cr
{\rm l.i.m.}\cr
$\stackrel{}{{}_{N\to \infty}}$\cr
}} }
H_{m_1}\left(\sum\limits_{l_1=0}^{N-1}\phi_{j_1}(\tau_{l_1})
\Delta {\bf w}_{\tau_{l_1}}^{(i_1)},\
\sum\limits_{l_1=0}^{N-1}\phi_{j_1}^2(\tau_{l_1})
\left(\Delta {\bf w}_{\tau_{l_1}}^{(i_1)}\right)^2\right),
$$

\vspace{6mm}
\noindent
where notations are the same as in Theorems~1, 2.

Similarly we get w.~p.~1

\vspace{1mm}
$$
\hbox{\vtop{\offinterlineskip\halign{
\hfil#\hfil\cr
{\rm l.i.m.}\cr
$\stackrel{}{{}_{N\to \infty}}$\cr
}} }
\left(\sum\limits_{l_{m_1+1}=0}^{N-1}
\phi_{j_2}(\tau_{l_{m_1}+1})\Delta{\bf w}_{\tau_{l_{m_1+1}}}^{(i_1)}
\ldots
\sum\limits_{l_{m_1+m_2}=0}^{N-1}
\phi_{j_2}(\tau_{l_{m_1+m_2}})
\Delta{\bf w}_{\tau_{l_{m_1+m_2}}}^{(i_1)}
-\right.
$$

\vspace{2.5mm}
$$
\left.-
\sum_{(l_{m_1+1},\ldots,l_{m_1+m_2})\in {\rm G}_{m_1+1,m_1+m_2}'}
\phi_{j_2}(\tau_{l_{m_1}+1})\Delta{\bf w}_{\tau_{l_{m_1+1}}}^{(i_1)}
\ldots
\phi_{j_2}(\tau_{l_{m_1+m_2}})
\Delta{\bf w}_{\tau_{l_{m_1+m_2}}}^{(i_1)}\right)=
$$

\vspace{4mm}
$$
=\hbox{\vtop{\offinterlineskip\halign{
\hfil#\hfil\cr
{\rm l.i.m.}\cr
$\stackrel{}{{}_{N\to \infty}}$\cr
}} }
H_{m_2}\left(\sum\limits_{l_1=0}^{N-1}\phi_{j_2}(\tau_{l_1})
\Delta {\bf w}_{\tau_{l_1}}^{(i_1)},\
\sum\limits_{l_1=0}^{N-1}\phi_{j_2}^2(\tau_{l_1})
\left(\Delta {\bf w}_{\tau_{l_1}}^{(i_1)}\right)^2\right),
$$

\vspace{2mm}

$$
\ldots
$$

\vspace{2mm}
$$
\hbox{\vtop{\offinterlineskip\halign{
\hfil#\hfil\cr
{\rm l.i.m.}\cr
$\stackrel{}{{}_{N\to \infty}}$\cr
}} }
\left(\sum\limits_{l_{k-m_r+1}=0}^{N-1}
\phi_{j_r}(\tau_{l_{k-m_r+1}})
\Delta{\bf w}_{\tau_{l_{k-m_r+1}}}^{(i_1)}
\ldots
\sum\limits_{l_{k}=0}^{N-1}
\phi_{j_r}(\tau_{l_{k}})
\Delta{\bf w}_{\tau_{l_{k}}}^{(i_1)}
-\right.
$$

\vspace{2.5mm}
$$
\left.-
\sum_{(l_{k-m_r+1},\ldots,l_{k})\in {\rm G}_{k-m_r+1,k}'}
\phi_{j_r}(\tau_{l_{k-m_r+1}})\Delta{\bf w}_{\tau_{k-m_r+1}}^{(i_1)}
\ldots
\phi_{j_r}(\tau_{l_{k}})
\Delta{\bf w}_{\tau_{l_{k}}}^{(i_1)}\right)
=
$$

\vspace{4mm}
$$
=\hbox{\vtop{\offinterlineskip\halign{
\hfil#\hfil\cr
{\rm l.i.m.}\cr
$\stackrel{}{{}_{N\to \infty}}$\cr
}} }
H_{m_r}\left(\sum\limits_{l_1=0}^{N-1}\phi_{j_r}(\tau_{l_1})
\Delta {\bf w}_{\tau_{l_1}}^{(i_1)},\ 
\sum\limits_{l_1=0}^{N-1}\phi_{j_r}^2(\tau_{l_1})
\left(\Delta {\bf w}_{\tau_{l_1}}^{(i_1)}\right)^2\right).
$$

\vspace{5mm}

Then

\vspace{-1mm}
$$
J'[\phi_{j_1}\ldots \phi_{j_k}]_{T,t}^{(i_1\ldots i_1)}=
$$

\vspace{2mm}
$$
=
\hbox{\vtop{\offinterlineskip\halign{
\hfil#\hfil\cr
{\rm l.i.m.}\cr
$\stackrel{}{{}_{N\to \infty}}$\cr
}} }
H_{m_1}\left(\sum\limits_{l_1=0}^{N-1}\phi_{j_1}(\tau_{l_1})
\Delta {\bf w}_{\tau_{l_1}}^{(i_1)},\
\sum\limits_{l_1=0}^{N-1}\phi_{j_1}^2(\tau_{l_1})
\left(\Delta {\bf w}_{\tau_{l_1}}^{(i_1)}\right)^2\right)\times 
$$

\vspace{2mm}
$$
~~~~~~\times \hbox{\vtop{\offinterlineskip\halign{
\hfil#\hfil\cr
{\rm l.i.m.}\cr
$\stackrel{}{{}_{N\to \infty}}$\cr
}} }
H_{m_2}\left(\sum\limits_{l_1=0}^{N-1}\phi_{j_2}(\tau_{l_1})
\Delta {\bf w}_{\tau_{l_1}}^{(i_1)},\
\sum\limits_{l_1=0}^{N-1}\phi_{j_2}^2(\tau_{l_1})
\left(\Delta {\bf w}_{\tau_{l_1}}^{(i_1)}\right)^2\right)\times \ldots
$$

\vspace{2mm}
\begin{equation}
\label{ziko2000}
\hspace{-6mm}
\ldots \times 
\hbox{\vtop{\offinterlineskip\halign{
\hfil#\hfil\cr
{\rm l.i.m.}\cr
$\stackrel{}{{}_{N\to \infty}}$\cr
}} }
H_{m_r}\left(\sum\limits_{l_1=0}^{N-1}\phi_{j_r}(\tau_{l_1})
\Delta {\bf w}_{\tau_{l_1}}^{(i_1)},\ 
\sum\limits_{l_1=0}^{N-1}\phi_{j_r}^2(\tau_{l_1})
\left(\Delta {\bf w}_{\tau_{l_1}}^{(i_1)}\right)^2\right)
\end{equation}

\vspace{5mm}
\noindent
w.~p.~1 for $i_1\ne 0$ and

\vspace{-1mm}
$$
J'[\phi_{j_1}\ldots \phi_{j_k}]_{T,t}^{(0\ldots 0)}=
$$

\vspace{2mm}

$$
=
\lim\limits_{N\to\infty}
\left(\sum\limits_{l_1=0}^{N-1}\phi_{j_1}(\tau_{l_1})
\Delta \tau_{l_1}\right)^{m_1}
\ldots 
\left(\sum\limits_{l_r=0}^{N-1}\phi_{j_r}(\tau_{l_r})
\Delta \tau_{l_r}\right)^{m_r}=
$$

\vspace{3mm}
$$
=\left(\int\limits_t^T\phi_{j_1}(s)
ds\right)^{m_1}
\ldots 
\left(\int\limits_t^T\phi_{j_r}(s)
ds\right)^{m_r}=
$$

\vspace{3mm}
\begin{equation}
\label{ziko2001}
=\left(\zeta_{j_1}^{(0)}\right)^{m_1}
\ldots 
\left(\zeta_{j_r}^{(0)}\right)^{m_r}
\end{equation}

\vspace{5mm}
\noindent
for $i_1=0,$ where we suppose that the condition 
(\ref{ziko1002}) is fulfilled; also we use in (\ref{ziko2000}) and (\ref{ziko2001}) 
the same notations as in the proof of Theorem 1.

Applying (\ref{ziko1000}), (\ref{ziko1001}), Lemma~3, and Remark~2 to the right-hand side of
(\ref{ziko2000}), we finally obtain w.~p.~1

\vspace{-2mm}
$$
J'[\phi_{j_1}\ldots \phi_{j_k}]_{T,t}^{(i_1\ldots i_1)}=
H_{m_1}\left(\int\limits_t^T\phi_{j_1}(s)
d{\bf w}_s^{(i_1)},
\int\limits_t^T\phi_{j_1}^2(s)ds
\right)\times
$$

\vspace{2mm}
$$
\times
H_{m_2}\left(\int\limits_t^T\phi_{j_2}(s)
d{\bf w}_s^{(i_1)},
\int\limits_t^T\phi_{j_2}^2(s)ds
\right)
\ldots 
H_{m_r}\left(\int\limits_t^T\phi_{j_r}(s)
d{\bf w}_s^{(i_1)},
\int\limits_t^T\phi_{j_r}^2(s)ds\right)=
$$

\vspace{3mm}
$$
=
H_{m_1}\left(\zeta_{j_1}^{(i_1)},
1\right)H_{m_2}\left(\zeta_{j_2}^{(i_1)},
1\right)
\ldots 
H_{m_r}\left(\zeta_{j_r}^{(i_1)},1\right)=
$$

\vspace{3mm}
$$
=
H_{m_1}\left(\zeta_{j_1}^{(i_1)}\right)
H_{m_2}\left(\zeta_{j_2}^{(i_1)}\right)
\ldots 
H_{m_r}\left(\zeta_{j_r}^{(i_1)}\right)
$$

\vspace{6mm}
\noindent
for $i_1\ne 0,$ 
where we suppose that the condition 
(\ref{ziko1002}) is fulfilled.
An equality similar to (\ref{ziko20}) was proved without using
Theorem 3.1 \cite{ito1951}.

Consider particular cases of the equality (\ref{ziko100}) for
$k=1,\ldots,4$ and $i_1,\ldots,i_4=1,\ldots,m$ (see (\ref{a1})--(\ref{a4})). We have
w.~p.~1

\vspace{1mm}
$$
J'[\phi_{j_1}]_{T,t}^{(i_1)}=
\zeta_{j_1}^{(i_1)}=H_1\left(\zeta_{j_1}^{(i_1)}\right);
$$

\vspace{5mm}
$$
J'[\phi_{j_1}\phi_{j_2}]_{T,t}^{(i_1i_2)}=\zeta_{j_1}^{(i_1)}\zeta_{j_2}^{(i_2)}
-{\bf 1}_{\{i_1=i_2\}}
{\bf 1}_{\{j_1=j_2\}}=
$$

\vspace{2mm}
\begin{equation}
\label{ziko301}
=\left\{
\begin{matrix}
H_2\left(\zeta_{j_1}^{(i_1)}\right)H_0\left(\zeta_{j_2}^{(i_2)}\right),\ 
&\hbox{\rm if}\ \ \ 
i_1=i_2,\ j_1=j_2\cr\cr
H_1\left(\zeta_{j_1}^{(i_1)}\right)H_1\left(\zeta_{j_2}^{(i_2)}\right),\  &\hbox{\rm otherwise}
\end{matrix}\right.;
\end{equation}

\vspace{10mm}
$$
J'[\phi_{j_1}\phi_{j_2}\phi_{j_3}]_{T,t}^{(i_1i_1i_1)}=
\zeta_{j_1}^{(i_1)}\zeta_{j_2}^{(i_1)}\zeta_{j_3}^{(i_1)}
-
{\bf 1}_{\{j_1=j_2\}}
\zeta_{j_3}^{(i_1)}
-
{\bf 1}_{\{j_2=j_3\}}
\zeta_{j_1}^{(i_1)}
-
$$

$$
-
{\bf 1}_{\{j_1=j_3\}}
\zeta_{j_2}^{(i_1)}=
$$

\vspace{2mm}
\begin{equation}
\label{ziko302}
=\left\{
\begin{matrix}
H_3\left(\zeta_{j_1}^{(i_1)}\right)H_0\left(\zeta_{j_2}^{(i_1)}\right)
H_0\left(\zeta_{j_3}^{(i_1)}\right),\ 
&\hbox{\rm if}\ \ \ 
j_1=j_2=j_3
\cr\cr\cr
H_2\left(\zeta_{j_1}^{(i_1)}\right)H_0\left(\zeta_{j_2}^{(i_1)}\right)
H_1\left(\zeta_{j_3}^{(i_1)}\right),\ 
&\hbox{\rm if}\ \ \ 
j_1=j_2\ne j_3
\cr\cr\cr
H_1\left(\zeta_{j_1}^{(i_1)}\right)H_2\left(\zeta_{j_2}^{(i_1)}\right)
H_0\left(\zeta_{j_3}^{(i_1)}\right),\ 
&\hbox{\rm if}\ \ \ 
j_2=j_3\ne j_1
\cr\cr\cr
H_0\left(\zeta_{j_1}^{(i_1)}\right)H_1\left(\zeta_{j_2}^{(i_1)}\right)
H_2\left(\zeta_{j_3}^{(i_1)}\right),\ 
&\hbox{\rm if}\ \ \ 
j_1=j_3\ne j_2
\cr\cr\cr
H_1\left(\zeta_{j_1}^{(i_1)}\right)H_1\left(\zeta_{j_2}^{(i_1)}\right)
H_1\left(\zeta_{j_3}^{(i_1)}\right),\ 
&\hbox{\rm if}\ \ \ 
j_1\ne j_2,\ 
j_2\ne j_3, j_1\ne j_3
\end{matrix}\right.;
\end{equation}

\vspace{10mm}
$$
J'[\phi_{j_1}\phi_{j_2}\phi_{j_3}]_{T,t}^{(i_1i_2i_3)}=
\zeta_{j_1}^{(i_1)}\zeta_{j_2}^{(i_2)}\zeta_{j_3}^{(i_2)}=
$$

\vspace{1mm}
$$
=
H_1\left(\zeta_{j_1}^{(i_1)}\right)H_1\left(\zeta_{j_2}^{(i_2)}\right)
H_1\left(\zeta_{j_3}^{(i_3)}\right),
$$

\vspace{5mm}
\noindent
where $i_1,i_2,i_3$ are pairwise different;

\vspace{4mm}
$$
J'[\phi_{j_1}\phi_{j_2}\phi_{j_3}]_{T,t}^{(i_1i_1i_3)}=
\zeta_{j_1}^{(i_1)}\zeta_{j_2}^{(i_1)}\zeta_{j_3}^{(i_3)}
-
{\bf 1}_{\{j_1=j_2\}}
\zeta_{j_3}^{(i_3)}=
$$

\vspace{2mm}
$$
=\left(\zeta_{j_1}^{(i_1)}\zeta_{j_2}^{(i_1)}
-
{\bf 1}_{\{j_1=j_2\}}\right)
\zeta_{j_3}^{(i_3)}=J'[\phi_{j_1}\phi_{j_2}]_{T,t}^{(i_1i_1)}J'[\phi_{j_3}]_{T,t}^{(i_3)}=
$$

\vspace{3mm}
$$
=\left\{
\begin{matrix}
H_2\left(\zeta_{j_1}^{(i_1)}\right)H_0\left(\zeta_{j_2}^{(i_1)}\right)
H_1\left(\zeta_{j_3}^{(i_3)}\right),\ 
&\hbox{\rm if}\ \ \ 
j_1=j_2\cr\cr\cr
H_1\left(\zeta_{j_1}^{(i_1)}\right)H_1\left(\zeta_{j_2}^{(i_1)}\right)
H_1\left(\zeta_{j_3}^{(i_3)}\right),\  &\hbox{\rm if}\ \ \ 
j_1\ne j_2
\end{matrix}\right.,
$$

\vspace{3mm}
\noindent
where $i_1=i_2\ne i_3$;

\vspace{4mm}
$$
J'[\phi_{j_1}\phi_{j_2}\phi_{j_3}]_{T,t}^{(i_1i_2i_2)}=
\zeta_{j_1}^{(i_1)}\zeta_{j_2}^{(i_2)}\zeta_{j_3}^{(i_2)}
-
{\bf 1}_{\{j_2=j_3\}}
\zeta_{j_1}^{(i_1)}=
$$

\vspace{2mm}
$$
=\zeta_{j_1}^{(i_1)}\left(\zeta_{j_2}^{(i_2)}\zeta_{j_3}^{(i_2)}
-
{\bf 1}_{\{j_2=j_3\}}\right)
=J'[\phi_{j_1}]_{T,t}^{(i_1)}J'[\phi_{j_2}\phi_{j_3}]_{T,t}^{(i_2i_2)}=
$$

\vspace{3mm}
$$
=\left\{
\begin{matrix}
H_1\left(\zeta_{j_1}^{(i_1)}\right)H_2\left(\zeta_{j_2}^{(i_2)}\right)
H_0\left(\zeta_{j_3}^{(i_2)}\right),\ 
&\hbox{\rm if}\ \ \ 
j_2=j_3\cr\cr\cr
H_1\left(\zeta_{j_1}^{(i_1)}\right)H_1\left(\zeta_{j_2}^{(i_2)}\right)
H_1\left(\zeta_{j_3}^{(i_2)}\right),\  &\hbox{\rm if}\ \ \ 
j_1\ne j_2
\end{matrix}\right.,
$$

\vspace{3mm}
\noindent
where $i_1\ne i_2=i_3$;

\vspace{4mm}
$$
J'[\phi_{j_1}\phi_{j_2}\phi_{j_3}]_{T,t}^{(i_1i_2i_1)}=
\zeta_{j_1}^{(i_1)}\zeta_{j_2}^{(i_2)}\zeta_{j_3}^{(i_1)}
-
{\bf 1}_{\{j_1=j_3\}}
\zeta_{j_2}^{(i_2)}=
$$

\vspace{2mm}
$$
=\zeta_{j_2}^{(i_2)}\left(\zeta_{j_1}^{(i_1)}\zeta_{j_3}^{(i_1)}
-
{\bf 1}_{\{j_1=j_3\}}\right)
=J'[\phi_{j_2}]_{T,t}^{(i_2)}J'[\phi_{j_1}\phi_{j_3}]_{T,t}^{(i_1i_1)}=
$$

\vspace{3mm}
$$
=\left\{
\begin{matrix}
H_2\left(\zeta_{j_1}^{(i_1)}\right)H_1\left(\zeta_{j_2}^{(i_2)}\right)
H_0\left(\zeta_{j_3}^{(i_1)}\right),\ 
&\hbox{\rm if}\ \ \ 
j_1=j_3\cr\cr\cr
H_1\left(\zeta_{j_1}^{(i_1)}\right)H_1\left(\zeta_{j_2}^{(i_2)}\right)
H_1\left(\zeta_{j_3}^{(i_1)}\right),\  &\hbox{\rm if}\ \ \ 
j_1\ne j_3
\end{matrix}\right.,
$$

\vspace{3mm}
\noindent
where $i_1=i_3\ne i_2$;

\vspace{2mm}

$$
J'[\phi_{j_1}\phi_{j_2}\phi_{j_3}\phi_{j_4}]_{T,t}^{(i_1i_1i_1i_1)}=
\prod\limits_{l=1}^4 \zeta_{j_l}^{(i_1)}-
$$
$$
-
{\bf 1}_{\{j_1=j_2\}}
\zeta_{j_3}^{(i_1)}
\zeta_{j_4}^{(i_1)}
-
{\bf 1}_{\{j_1=j_3\}}
\zeta_{j_2}^{(i_1)}
\zeta_{j_4}^{(i_1)}-
$$
$$
-
{\bf 1}_{\{j_1=j_4\}}
\zeta_{j_2}^{(i_1)}
\zeta_{j_3}^{(i_1)}
-
{\bf 1}_{\{j_2=j_3\}}
\zeta_{j_1}^{(i_1)}
\zeta_{j_4}^{(i_1)}-
$$
$$
-
{\bf 1}_{\{j_2=j_4\}}
\zeta_{j_1}^{(i_1)}
\zeta_{j_3}^{(i_1)}
-
{\bf 1}_{\{j_3=j_4\}}
\zeta_{j_1}^{(i_1)}
\zeta_{j_2}^{(i_1)}+
$$
$$
+
{\bf 1}_{\{j_1=j_2\}}
{\bf 1}_{\{j_3=j_4\}}
+
{\bf 1}_{\{j_1=j_3\}}
{\bf 1}_{\{j_2=j_4\}}+
$$
$$
+\Biggl.
{\bf 1}_{\{j_1=j_4\}}
{\bf 1}_{\{j_2=j_3\}}=
$$
$$
=\left\{
\begin{matrix}
H_4\left(\zeta_{j_1}^{(i_1)}\right)H_0\left(\zeta_{j_2}^{(i_1)}\right)
H_0\left(\zeta_{j_3}^{(i_1)}\right)H_0\left(\zeta_{j_4}^{(i_1)}\right),\ 
&\hbox{\rm if}\ \ \ {\rm (I)}
\cr\cr
H_1\left(\zeta_{j_1}^{(i_1)}\right)H_1\left(\zeta_{j_2}^{(i_1)}\right)
H_1\left(\zeta_{j_3}^{(i_1)}\right)H_1\left(\zeta_{j_4}^{(i_1)}\right),\ 
&\hbox{\rm if}\ \ \ {\rm (II)}
\cr\cr
H_2\left(\zeta_{j_1}^{(i_1)}\right)H_0\left(\zeta_{j_2}^{(i_1)}\right)
H_1\left(\zeta_{j_3}^{(i_1)}\right)H_1\left(\zeta_{j_4}^{(i_1)}\right),\ 
&\hbox{\rm if}\ \ \ {\rm (III)}
\cr\cr
H_0\left(\zeta_{j_1}^{(i_1)}\right)H_1\left(\zeta_{j_2}^{(i_1)}\right)
H_2\left(\zeta_{j_3}^{(i_1)}\right)H_1\left(\zeta_{j_4}^{(i_1)}\right),\ 
&\hbox{\rm if}\ \ \ {\rm (IV)}
\cr\cr
H_0\left(\zeta_{j_1}^{(i_1)}\right)H_1\left(\zeta_{j_2}^{(i_1)}\right)
H_1\left(\zeta_{j_3}^{(i_1)}\right)H_2\left(\zeta_{j_4}^{(i_1)}\right),\ 
&\hbox{\rm if}\ \ \ {\rm (V)}
\cr\cr
H_1\left(\zeta_{j_1}^{(i_1)}\right)H_0\left(\zeta_{j_2}^{(i_1)}\right)
H_2\left(\zeta_{j_3}^{(i_1)}\right)H_1\left(\zeta_{j_4}^{(i_1)}\right),\ 
&\hbox{\rm if}\ \ \ {\rm (VI)}
\cr\cr
H_1\left(\zeta_{j_1}^{(i_1)}\right)H_0\left(\zeta_{j_2}^{(i_1)}\right)
H_1\left(\zeta_{j_3}^{(i_1)}\right)H_2\left(\zeta_{j_4}^{(i_1)}\right),\ 
&\hbox{\rm if}\ \ \ {\rm (VII)}
\cr\cr
H_1\left(\zeta_{j_1}^{(i_1)}\right)H_1\left(\zeta_{j_2}^{(i_1)}\right)
H_0\left(\zeta_{j_3}^{(i_1)}\right)H_2\left(\zeta_{j_4}^{(i_1)}\right),\ 
&\hbox{\rm if}\ \ \ {\rm (VIII)}
\cr\cr
H_3\left(\zeta_{j_1}^{(i_1)}\right)H_0\left(\zeta_{j_2}^{(i_1)}\right)
H_0\left(\zeta_{j_3}^{(i_1)}\right)H_1\left(\zeta_{j_4}^{(i_1)}\right),\ 
&\hbox{\rm if}\ \ \ {\rm (IX)}
\cr\cr
H_1\left(\zeta_{j_1}^{(i_1)}\right)H_3\left(\zeta_{j_2}^{(i_1)}\right)
H_0\left(\zeta_{j_3}^{(i_1)}\right)H_0\left(\zeta_{j_4}^{(i_1)}\right),\ 
&\hbox{\rm if}\ \ \ {\rm (X)}
\cr\cr
H_0\left(\zeta_{j_1}^{(i_1)}\right)H_0\left(\zeta_{j_2}^{(i_1)}\right)
H_1\left(\zeta_{j_3}^{(i_1)}\right)H_3\left(\zeta_{j_4}^{(i_1)}\right),\ 
&\hbox{\rm if}\ \ \ {\rm (XI)}
\cr\cr
H_0\left(\zeta_{j_1}^{(i_1)}\right)H_1\left(\zeta_{j_2}^{(i_1)}\right)
H_0\left(\zeta_{j_3}^{(i_1)}\right)H_3\left(\zeta_{j_4}^{(i_1)}\right),\ 
&\hbox{\rm if}\ \ \ {\rm (XII)}
\cr\cr
H_2\left(\zeta_{j_1}^{(i_1)}\right)H_0\left(\zeta_{j_2}^{(i_1)}\right)
H_0\left(\zeta_{j_3}^{(i_1)}\right)H_2\left(\zeta_{j_4}^{(i_1)}\right),\ 
&\hbox{\rm if}\ \ \ {\rm (XIII)}
\cr\cr
H_2\left(\zeta_{j_1}^{(i_1)}\right)H_2\left(\zeta_{j_2}^{(i_1)}\right)
H_0\left(\zeta_{j_3}^{(i_1)}\right)H_0\left(\zeta_{j_4}^{(i_1)}\right),\ 
&\hbox{\rm if}\ \ \ {\rm (XIV)}
\cr\cr
H_2\left(\zeta_{j_1}^{(i_1)}\right)H_0\left(\zeta_{j_2}^{(i_1)}\right)
H_2\left(\zeta_{j_3}^{(i_1)}\right)H_0\left(\zeta_{j_4}^{(i_1)}\right),\ 
&\hbox{\rm if}\ \ \ {\rm (XV)}
\end{matrix}\right.,
$$

\vspace{8mm}
\noindent
where $H_n(x)$ is the Hermite polynomial (\ref{ziko500}) of degree $n$  
and (I)--(XV) are the following conditions

\vspace{3mm}

(I).\ $j_1=j_2=j_3=j_4,$

\vspace{1mm}

(II).\ $j_1, j_2, j_3, j_4$ are pairwise different,

\vspace{1mm}

(III).\ $j_1=j_2\ne j_3, j_4;\ j_3\ne j_4,$
 
\vspace{1mm}

(IV).\ $j_1=j_3\ne j_2, j_4;\ j_2\ne j_4,$

\vspace{1mm}

(V).\ $j_1=j_4\ne j_2, j_3;\ j_2\ne j_3,$

\vspace{1mm}

(VI).\ $j_2=j_3\ne j_1, j_4;\ j_1\ne j_4,$

\vspace{1mm}

(VII).\ $j_2=j_4\ne j_1, j_3;\ j_1\ne j_3,$

\vspace{1mm}

(VIII).\ $j_3=j_4\ne j_1, j_2;\ j_1\ne j_2,$

\vspace{1mm}

(IX).\ $j_1=j_2=j_3\ne j_4,$

\vspace{1mm}

(X).\ $j_2=j_3=j_4\ne j_1,$

\vspace{1mm}

(XI).\ $j_1=j_2=j_4\ne j_3,$

\vspace{1mm}

(XII).\ $j_1=j_3=j_4\ne j_2,$

\vspace{1mm}

(XIII).\ $j_1=j_2\ne j_3=j_4,$

\vspace{1mm}

(XIV).\ $j_1=j_3\ne j_2=j_4,$

\vspace{1mm}

(XV).\ $j_1=j_4\ne j_2=j_3.$

\vspace{3mm}

\noindent
Moreover, from (\ref{ziko30}) we have w.~p.~1

\vspace{1mm}
$$
J'[\phi_{j_1}\phi_{j_2}\phi_{j_3}\phi_{j_4}]_{T,t}^{(i_1i_2i_3i_4)}=
H_1\left(\zeta_{j_1}^{(i_1)}\right)H_1\left(\zeta_{j_2}^{(i_2)}\right)
H_1\left(\zeta_{j_3}^{(i_3)}\right)H_1\left(\zeta_{j_4}^{(i_4)}\right),
$$

\vspace{5mm}
\noindent
where $i_1,\ldots,i_4$ are pairwise different;

\vspace{3mm}
\begin{equation}
\label{ziko201}
J'[\phi_{j_1}\phi_{j_2}\phi_{j_3}\phi_{j_4}]_{T,t}^{(i_1i_1i_3i_4)}=
J'[\phi_{j_1}\phi_{j_2}]_{T,t}^{(i_1i_1)}
H_1\left(\zeta_{j_3}^{(i_3)}\right)H_1\left(\zeta_{j_4}^{(i_4)}\right)\ \ \ 
(i_1=i_2\ne i_3, i_4;\ i_3\ne i_4);
\end{equation}

\vspace{3mm}
\begin{equation}
\label{ziko202}
J'[\phi_{j_1}\phi_{j_2}\phi_{j_3}\phi_{j_4}]_{T,t}^{(i_1i_2i_1i_4)}=
J'[\phi_{j_1}\phi_{j_3}]_{T,t}^{(i_1i_1)}
H_1\left(\zeta_{j_2}^{(i_2)}\right)H_1\left(\zeta_{j_4}^{(i_4)}\right)\ \ \ 
(i_1=i_3\ne i_2, i_4;\ i_2\ne i_4);
\end{equation}

\vspace{3mm}
\begin{equation}
\label{ziko203}
J'[\phi_{j_1}\phi_{j_2}\phi_{j_3}\phi_{j_4}]_{T,t}^{(i_1i_2i_3i_1)}=
J'[\phi_{j_1}\phi_{j_4}]_{T,t}^{(i_1i_1)}
H_1\left(\zeta_{j_2}^{(i_2)}\right)H_1\left(\zeta_{j_3}^{(i_3)}\right)\ \ \ 
(i_1=i_4\ne i_2, i_3;\ i_2\ne i_3);
\end{equation}

\vspace{3mm}
\begin{equation}
\label{ziko204}
J'[\phi_{j_1}\phi_{j_2}\phi_{j_3}\phi_{j_4}]_{T,t}^{(i_1i_2i_2i_4)}=
J'[\phi_{j_2}\phi_{j_3}]_{T,t}^{(i_2i_2)}
H_1\left(\zeta_{j_1}^{(i_1)}\right)H_1\left(\zeta_{j_4}^{(i_4)}\right)\ \ \ 
(i_2=i_3\ne i_1, i_4;\ i_1\ne i_4);
\end{equation}

\vspace{3mm}
\begin{equation}
\label{ziko205}
J'[\phi_{j_1}\phi_{j_2}\phi_{j_3}\phi_{j_4}]_{T,t}^{(i_1i_2i_3i_2)}=
J'[\phi_{j_2}\phi_{j_4}]_{T,t}^{(i_2i_2)}
H_1\left(\zeta_{j_1}^{(i_1)}\right)H_1\left(\zeta_{j_3}^{(i_3)}\right)\ \ \ 
(i_2=i_4\ne i_1, i_3;\ i_1\ne i_3);
\end{equation}

\vspace{3mm}
\begin{equation}
\label{ziko206}
J'[\phi_{j_1}\phi_{j_2}\phi_{j_3}\phi_{j_4}]_{T,t}^{(i_1i_2i_3i_3)}=
J'[\phi_{j_3}\phi_{j_4}]_{T,t}^{(i_3i_3)}
H_1\left(\zeta_{j_1}^{(i_1)}\right)H_1\left(\zeta_{j_2}^{(i_2)}\right)\ \ \
(i_3=i_4\ne i_1, i_2;\ i_1\ne i_2);
\end{equation}

\vspace{3mm}
\begin{equation}
\label{ziko207}
J'[\phi_{j_1}\phi_{j_2}\phi_{j_3}\phi_{j_4}]_{T,t}^{(i_1i_1i_1i_4)}=
J'[\phi_{j_1}\phi_{j_2}\phi_{j_3}]_{T,t}^{(i_1i_1i_1)}
H_1\left(\zeta_{j_4}^{(i_4)}\right)\ \ \
(i_1=i_2=i_3\ne i_4);
\end{equation}

\vspace{3mm}
\begin{equation}
\label{ziko208}
J'[\phi_{j_1}\phi_{j_2}\phi_{j_3}\phi_{j_4}]_{T,t}^{(i_1i_2i_2i_2)}=
J'[\phi_{j_2}\phi_{j_3}\phi_{j_4}]_{T,t}^{(i_2i_2i_2)}
H_1\left(\zeta_{j_1}^{(i_1)}\right)\ \ \ 
(i_2=i_3=i_4\ne i_1);
\end{equation}

\vspace{3mm}
\begin{equation}
\label{ziko209}
J'[\phi_{j_1}\phi_{j_2}\phi_{j_3}\phi_{j_4}]_{T,t}^{(i_1i_1i_3i_1)}=
J'[\phi_{j_1}\phi_{j_2}\phi_{j_4}]_{T,t}^{(i_1i_1i_1)}
H_1\left(\zeta_{j_3}^{(i_3)}\right)\ \ \ 
(i_1=i_2=i_4\ne i_3);
\end{equation}

\vspace{3mm}
\begin{equation}
\label{ziko210}
J'[\phi_{j_1}\phi_{j_2}\phi_{j_3}\phi_{j_4}]_{T,t}^{(i_1i_2i_1i_1)}=
J'[\phi_{j_1}\phi_{j_3}\phi_{j_4}]_{T,t}^{(i_1i_1i_1)}
H_1\left(\zeta_{j_2}^{(i_2)}\right)\ \ \ 
(i_1=i_3=i_4\ne i_2);
\end{equation}

\vspace{3mm}
\begin{equation}
\label{ziko211}
J'[\phi_{j_1}\phi_{j_2}\phi_{j_3}\phi_{j_4}]_{T,t}^{(i_1i_1i_3i_3)}=
J'[\phi_{j_1}\phi_{j_2}]_{T,t}^{(i_1i_1)}J'[\phi_{j_3}\phi_{j_4}]_{T,t}^{(i_3i_3)}\ \ \
(i_1=i_2\ne i_3=i_4);
\end{equation}

\vspace{3mm}
\begin{equation}
\label{ziko212}
J'[\phi_{j_1}\phi_{j_2}\phi_{j_3}\phi_{j_4}]_{T,t}^{(i_1i_2i_1i_2)}=
J'[\phi_{j_1}\phi_{j_3}]_{T,t}^{(i_1i_1)}J'[\phi_{j_2}\phi_{j_4}]_{T,t}^{(i_2i_2)}\ \ \
(i_1=i_3\ne i_2=i_4);
\end{equation}

\vspace{3mm}
\begin{equation}
\label{ziko213}
J'[\phi_{j_1}\phi_{j_2}\phi_{j_3}\phi_{j_4}]_{T,t}^{(i_1i_2i_2i_1)}=
J'[\phi_{j_1}\phi_{j_4}]_{T,t}^{(i_1i_1)}J'[\phi_{j_2}\phi_{j_3}]_{T,t}^{(i_2i_2)}\ \ \
(i_1=i_4\ne i_2=i_3).
\end{equation}

\vspace{6mm}

Note that the right-hand sides of (\ref{ziko201})--(\ref{ziko213})
contain multiple Wiener stochastic integrals of multiplicities 2 and 3.
These integrals are considered in detail in (\ref{ziko301}), (\ref{ziko302}).

It should be noted that the formulas (\ref{leto6000}) (Theorem~2) and 
(\ref{ziko800}) (Theorem~10) are interesting from various points of view.
The formulas (\ref{a1})--(\ref{a6}) (these formulas are particular cases of 
(\ref{leto6000}) for $k=1,\ldots,6$)
are convenient for numerical modeling of iterated Ito stochastic integrals 
of multiplicities 1 to 6.
For example, in \cite{Kuz-Kuz} and \cite{Mikh-1}, approximations 
of iterated Ito stochastic integrals of multiplicities 1 to 6 
in the Python programming 
language were successfully implemented using (\ref{a1})--(\ref{a6}) and 
Legendre polynomials.

On the other hand, the equality (\ref{ziko800}) is interesting 
by a number of reasons.
Firstly, this equality connects Ito's results on multiple Wiener stochastic 
integrals (\cite{ito1951}, Theorem~3.1) with the theory of 
mean-square approximation of iterated Ito stochastic integrals 
presented in this paper and in the book \cite{20a}.
Secondly, the equality (\ref{ziko800}) is based on the 
Hermite polynomials, which have the 
orthogonality property on $\mathbb{R}$ with a Gaussian weight. 
This feature opens up new possibilities in the study 
of iterated Ito stochastic integrals.

\vspace{5mm}

\section{A Generalization of Theorems 1, 2, 10, and 11 to the Case of an Arbitrary 
Complete Ortho\-nor\-mal System of Functions in the Space $L_2([t, T])$
and $\psi_1(\tau),$ $\ldots,\psi_k(\tau)\in L_2([t, T]),$
$\Phi(t_1,\ldots,t_k)\in L_2([t, T]^k)$}

\vspace{5mm}

In this section, we will use the definition of the multiple Wiener 
stochastic integral from \cite{ito1951}, \cite{Kuo} to generalize Theorems 
1, 2, 10, and 11 to the case of an arbitrary 
complete orthonormal system of functions in the space $L_2([t, T])$
and $\psi_1(\tau),$ $\ldots,\psi_k(\tau)\in L_2([t, T]),$
$\Phi(t_1,\ldots,t_k)\in L_2([t, T]^k)$.

Consider the following step function on the hypercube $[t, T]^k$

\vspace{-2mm}
\begin{equation}
\label{chain3}
\Phi_N(t_1,\ldots,t_k)=\sum\limits_{l_1,\ldots,l_k=0}^{N-1}
a_{l_1 \ldots l_k} {\bf 1}_{[\tau_{l_1},\tau_{l_1+1})}(t_1) \ldots
{\bf 1}_{[\tau_{l_k},\tau_{l_k+1})}(t_k),
\end{equation}

\vspace{3mm}
\noindent
where $a_{l_1 \ldots l_k}\in\mathbb{R}$ and such that 
$a_{l_1 \ldots l_k}=0$ if $l_p=l_q$ for some $p\ne q,$

\vspace{1mm}
$$
{\bf 1}_A (\tau)=\left\{
\begin{matrix}
1\ &{\rm if}\ \tau\in A \cr\cr
0\ &\hbox{\rm otherwise}
\end{matrix}\right.,
$$

\vspace{4mm}
\noindent
$N\in\mathbb{N},$ $\left\{\tau_{j}\right\}_{j=0}^{N}$ is a partition of
$[t,T],$ which satisfies the condition (\ref{1111}):

\vspace{1mm}
\begin{equation}
\label{1111xxx1}
t=\tau_0<\ldots <\tau_N=T,\ \ \
\Delta_N=
\hbox{\vtop{\offinterlineskip\halign{
\hfil#\hfil\cr
{\rm max}\cr
$\stackrel{}{{}_{0\le j\le N-1}}$\cr
}} }\Delta\tau_j\to 0\ \ \hbox{if}\ \ N\to \infty,\ \ \ 
\Delta\tau_j=\tau_{j+1}-\tau_j.
\end{equation}

\vspace{4mm}

Let us define the multiple Wiener stochastic integral for $\Phi_N(t_1,\ldots,t_k)$ 
\cite{ito1951}, \cite{Kuo}

\vspace{1mm}
\begin{equation}
\label{chain9}
J'[\Phi_N]_{T,t}^{(i_1\ldots i_k)}\stackrel{\sf def}{=}
\sum\limits_{l_1,\ldots,l_k=0}^{N-1}
a_{l_1 \ldots l_k}
\Delta{\bf w}_{\tau_{l_1}}^{(i_1)}\ldots \Delta{\bf w}_{\tau_{l_k}}^{(i_k)},
\end{equation}

\vspace{4mm}
\noindent
where $\Delta{\bf w}_{\tau_{j}}^{(i)}=
{\bf w}_{\tau_{j+1}}^{(i)}-{\bf w}_{\tau_{j}}^{(i)},$\
$i=0, 1,\ldots,m,$\ ${\bf w}_{\tau}^{(0)}=\tau.$

It is known (see \cite{Kuo}, Lemma~9.6.4)
that for any $\Phi(t_1,\ldots,t_k)\in L_2([t, T]^k)$ 
there exists a sequence of step functions $\Phi_N(t_1,\ldots,t_k)$ of the form (\ref{chain3})
such that

\vspace{1mm}
\begin{equation}
\label{chain15}
\lim\limits_{N\to\infty} \int\limits_{[t,T]^k}
\left(\Phi(t_1,\ldots,t_k)-\Phi_N(t_1,\ldots,t_k)\right)^2 dt_1\ldots dt_k=0.
\end{equation}

\vspace{4mm}

We have

\vspace{-1mm}
$$
\Phi_N(t_1,\ldots,t_k)=\sum\limits_{l_1,\ldots,l_k=0}^{N-1}
a_{l_1 \ldots l_k} {\bf 1}_{[\tau_{l_1},\tau_{l_1+1})}(t_1) \ldots
{\bf 1}_{[\tau_{l_k},\tau_{l_k+1})}(t_k)=
$$

\vspace{1mm}
\begin{equation}
\label{chain5}
=\sum\limits_{(l_1,\ldots,l_k)}
\sum_{\stackrel{l_1,\ldots,l_k=0}{{}_{l_1<l_2<\ldots < l_k}}}^{N-1}
a_{l_1 \ldots l_k} {\bf 1}_{[\tau_{l_1},\tau_{l_1+1})}(t_1) \ldots
{\bf 1}_{[\tau_{l_k},\tau_{l_k+1})}(t_k),
\end{equation}

\vspace{4mm}
\noindent
where permutations $(l_1,\ldots,l_k)$ when summing are 
performed only in the expression $l_1<l_2<\ldots < l_k$
(recall that $a_{l_1 \ldots l_k}=0$ if $l_p=l_q$ for some $p\ne q$).

Using (\ref{chain5}), we get

\vspace{-1mm}
\begin{equation}
\label{chain30}
\sum_{(t_1,\ldots,t_k)}
\int\limits_{t}^{T}
\ldots
\int\limits_{t}^{t_2}
\Phi_N(t_1,\ldots,t_k)d{\bf w}_{t_1}^{(i_1)}
\ldots
d{\bf w}_{t_k}^{(i_k)}=
\end{equation}

\vspace{1mm}
$$
=\sum\limits_{(l_1,\ldots,l_k)}
\sum_{\stackrel{l_1,\ldots,l_k=0}{{}_{l_1<l_2<\ldots < l_k}}}^{N-1}
a_{l_1 \ldots l_k} 
\Delta{\bf w}_{\tau_{l_1}}^{(i_1)} \ldots \Delta{\bf w}_{\tau_{l_k}}^{(i_k)}=
$$

\vspace{1mm}
$$
=\sum\limits_{\stackrel{l_1,\ldots,l_k=0}{{}_{l_q\ne l_r;\ q\ne r;\ 
q, r=1,\ldots, k}}}^{N-1}
a_{l_1 \ldots l_k} 
\Delta{\bf w}_{\tau_{l_1}}^{(i_1)} \ldots \Delta{\bf w}_{\tau_{l_k}}^{(i_k)}=
$$

\vspace{1mm}
\begin{equation}
\label{chain10}
=J'[\Phi_N]_{T,t}^{(i_1\ldots i_k)}\ \ \ \hbox{w.\ p.\ 1},
\end{equation}

\vspace{3mm}
\noindent
where permutations $(t_1,\ldots,t_k)$ when summing are 
performed only in the values
$d{\bf w}_{t_1}^{(i_1)}
\ldots $
$d{\bf w}_{t_k}^{(i_k)}$ 
and permutations $(l_1,\ldots,l_k)$ when summing are 
performed only in the expression $l_1<l_2<\ldots < l_k.$
At the same time the indices near 
upper 
limits of integration in the iterated stochastic integrals in (\ref{chain30}) are changed 
correspondently and if $t_r$ swapped with $t_q$ in the  
permutation $(t_1,\ldots,t_k)$, then $i_r$ swapped with $i_q$ in 
the permutation $(i_1,\ldots,i_k)$ (see (\ref{chain30})).
In addition, the multiple Wiener stochastic integral 
$J'[\Phi_N]_{T,t}^{(i_1\ldots i_k)}$ is defined by (\ref{chain9})
and 

\vspace{1mm}
$$
\int\limits_{t}^{T}
\ldots
\int\limits_{t}^{t_2}
\Phi_N(t_1,\ldots,t_k)d{\bf w}_{t_1}^{(i_1)}
\ldots
d{\bf w}_{t_k}^{(i_k)}
$$

\vspace{4mm}
\noindent
is the iterated Ito stochastic integral.

Using (\ref{chain15}), (\ref{chain10}), Lemma 2 for $\Phi(t_1,\ldots,t_k)\in L_2(D_k)$, 
and (\ref{riemann}) for Lebesgue integrals, 
we have

\vspace{1mm}
$$
{\sf M}\left\{\left(J'[\Phi_N]_{T,t}^{(i_1\ldots i_k)}-
J'[\Phi_M]_{T,t}^{(i_1\ldots i_k)}\right)^2\right\}\le
$$

\vspace{2mm}
$$
\le C_k 
\sum_{(t_1,\ldots,t_k)}
\int\limits_{t}^{T}
\ldots
\int\limits_{t}^{t_2}
\left(\Phi_N(t_1,\ldots,t_k)-\Phi_M(t_1,\ldots,t_k)\right)^2 dt_1
\ldots dt_k=
$$

\vspace{2mm}
$$
=C_k 
\int\limits_{[t,T]^k}
\left(\Phi_N(t_1,\ldots,t_k)-\Phi_M(t_1,\ldots,t_k)\right)^2 dt_1
\ldots dt_k=
$$

\vspace{2mm}
$$
=C_k\left\Vert \Phi_N-\Phi_M\right\Vert_{L_2([t, T]^k)}^2\le
$$

\vspace{2mm}
$$
\le 2 C_k \left(\left\Vert \Phi_N-\Phi\right\Vert_{L_2([t, T]^k)}^2+
\left\Vert \Phi-\Phi_M\right\Vert_{L_2([t, T]^k)}^2\right)^2\ \to 0
$$

\vspace{5mm}
\noindent
if $N,M\to\infty,$ 
where constant $C_k$ 
depends only
on the multiplicity $k$ of the multiple Wiener stochastic integral.

Thus, there exists the limit 

$$
\hbox{\vtop{\offinterlineskip\halign{
\hfil#\hfil\cr
{\rm l.i.m.}\cr
$\stackrel{}{{}_{N\to \infty}}$\cr
}} }J'[\Phi_N]_{T,t}^{(i_1\ldots i_k)}.
$$

\vspace{4mm}

We will define the multiple Wiener stochastic integral for $\Phi(t_1,\ldots,t_k)\in L_2([t, T]^k)$ 
by the formula \cite{ito1951}, \cite{Kuo}
\begin{equation}
\label{WiI}
J'[\Phi]_{T,t}^{(i_1\ldots i_k)}\stackrel{\sf def}{=}
\hbox{\vtop{\offinterlineskip\halign{
\hfil#\hfil\cr
{\rm l.i.m.}\cr
$\stackrel{}{{}_{N\to \infty}}$\cr
}} }J'[\Phi_N]_{T,t}^{(i_1\ldots i_k)}=
\hbox{\vtop{\offinterlineskip\halign{
\hfil#\hfil\cr
{\rm l.i.m.}\cr
$\stackrel{}{{}_{N\to \infty}}$\cr
}} }
\sum\limits_{l_1,\ldots,l_k=0}^{N-1}
a_{l_1 \ldots l_k}
\Delta{\bf w}_{\tau_{l_1}}^{(i_1)}\ldots \Delta{\bf w}_{\tau_{l_k}}^{(i_k)},
\end{equation}

\vspace{4mm}
\noindent
where $\Phi_N(t_1,\ldots,t_k)$ is defined by 
(\ref{chain3}),
$\Delta{\bf w}_{\tau_{j}}^{(i)}=
{\bf w}_{\tau_{j+1}}^{(i)}-{\bf w}_{\tau_{j}}^{(i)},$\
$i=0, 1,\ldots,m,$\ ${\bf w}_{\tau}^{(0)}=\tau.$

It is easy to see that the above definition coincides with 
(\ref{mult11}) if the function 
$\Phi(t_1,\ldots,t_k):\ [t, T]^k\to\mathbb{R}$ is continuous in the hypercube
$[t, T]^k$.

Let us prove the following equality 

\begin{equation}
\label{Wi110}
J'[\Phi]_{T,t}^{(i_1\ldots i_k)}=\sum_{(t_1,\ldots,t_k)}
\int\limits_{t}^{T}
\ldots
\int\limits_{t}^{t_2}
\Phi(t_1,\ldots,t_k)d{\bf w}_{t_1}^{(i_1)}
\ldots
d{\bf w}_{t_k}^{(i_k)}\ \ \ \hbox{w.\ p.\ 1},
\end{equation}

\vspace{3mm}
\noindent
where permutations $(t_1,\ldots,t_k)$ when summing are 
performed only in the values
$d{\bf w}_{t_1}^{(i_1)}
\ldots $
$d{\bf w}_{t_k}^{(i_k)}.$ At the same time the indices near 
upper 
limits of integration in the iterated stochastic integrals are changed 
correspondently and if $t_r$ swapped with $t_q$ in the  
permutation $(t_1,\ldots,t_k)$, then $i_r$ swapped with $i_q$ in 
the permutation $(i_1,\ldots,i_k).$ 
In addition, the multiple Wiener stochastic integral 
$J'[\Phi]_{T,t}^{(i_1\ldots i_k)}$ is defined by (\ref{WiI})
and 

$$
\int\limits_{t}^{T}
\ldots
\int\limits_{t}^{t_2}
\Phi(t_1,\ldots,t_k)d{\bf w}_{t_1}^{(i_1)}
\ldots
d{\bf w}_{t_k}^{(i_k)}
$$

\vspace{3mm}
\noindent
is the iterated Ito stochastic integral.

The equality (\ref{Wi110}) has already been proved for the case 
$\Phi(t_1,\ldots,t_k)=\Phi_N(t_1,\ldots,t_k)$ (see (\ref{chain10})).
From (\ref{chain10}) we have
$$
J'[\Phi_N]_{T,t}^{(i_1\ldots i_k)}=
\sum_{(t_1,\ldots,t_k)}
\int\limits_{t}^{T}
\ldots
\int\limits_{t}^{t_2}
\Phi_N(t_1,\ldots,t_k)d{\bf w}_{t_1}^{(i_1)}
\ldots
d{\bf w}_{t_k}^{(i_k)}=
$$

\vspace{2mm}
$$
=\sum_{(t_1,\ldots,t_k)}
\int\limits_{t}^{T}
\ldots
\int\limits_{t}^{t_2}
\Phi(t_1,\ldots,t_k)d{\bf w}_{t_1}^{(i_1)}
\ldots
d{\bf w}_{t_k}^{(i_k)}+
$$

\vspace{2mm}
\begin{equation}
\label{chain11}
+\sum_{(t_1,\ldots,t_k)}
\int\limits_{t}^{T}
\ldots
\int\limits_{t}^{t_2}
\left(\Phi_N(t_1,\ldots,t_k)-\Phi(t_1,\ldots,t_k)\right)d{\bf w}_{t_1}^{(i_1)}
\ldots
d{\bf w}_{t_k}^{(i_k)}\ \ \ \hbox{w.~p.~1.}
\end{equation}

\vspace{5mm}

Passing to the limit $\hbox{\vtop{\offinterlineskip\halign{
\hfil#\hfil\cr
{\rm l.i.m.}\cr
$\stackrel{}{{}_{N\to \infty}}$\cr
}} }$ in the equality (\ref{chain11}), we obtain

\vspace{1mm}
$$
J'[\Phi]_{T,t}^{(i_1\ldots i_k)}=
\sum_{(t_1,\ldots,t_k)}
\int\limits_{t}^{T}
\ldots
\int\limits_{t}^{t_2}
\Phi(t_1,\ldots,t_k)d{\bf w}_{t_1}^{(i_1)}
\ldots
d{\bf w}_{t_k}^{(i_k)}+
$$

\vspace{2mm}
\begin{equation}
\label{chain12}
+\hbox{\vtop{\offinterlineskip\halign{
\hfil#\hfil\cr
{\rm l.i.m.}\cr
$\stackrel{}{{}_{N\to \infty}}$\cr
}} }\sum_{(t_1,\ldots,t_k)}
\int\limits_{t}^{T}
\ldots
\int\limits_{t}^{t_2}
\left(\Phi_N(t_1,\ldots,t_k)-\Phi(t_1,\ldots,t_k)\right)d{\bf w}_{t_1}^{(i_1)}
\ldots
d{\bf w}_{t_k}^{(i_k)}\ \ \ \hbox{w.~p.~1.}
\end{equation}

\vspace{5mm}

Using Lemma 2 for $\Phi(t_1,\ldots,t_k)\in L_2(D_k)$, 
(\ref{riemann}) for Lebesgue integrals, and (\ref{chain15}), we get

\vspace{1mm}
$$
{\sf M}\left\{\left(
\sum_{(t_1,\ldots,t_k)}
\int\limits_{t}^{T}
\ldots
\int\limits_{t}^{t_2}
\left(\Phi_N(t_1,\ldots,t_k)-\Phi(t_1,\ldots,t_k)\right)d{\bf w}_{t_1}^{(i_1)}
\ldots
d{\bf w}_{t_k}^{(i_k)}\right)^2\right\}\le
$$

\vspace{2mm}
$$
\le C_k 
\sum_{(t_1,\ldots,t_k)}
\int\limits_{t}^{T}
\ldots
\int\limits_{t}^{t_2}
\left(\Phi_N(t_1,\ldots,t_k)-\Phi(t_1,\ldots,t_k)\right)^2 dt_1
\ldots dt_k=
$$

\vspace{2mm}
\begin{equation}
\label{chain20}
=C_k 
\int\limits_{[t,T]^k}
\left(\Phi_N(t_1,\ldots,t_k)-\Phi(t_1,\ldots,t_k)\right)^2 dt_1
\ldots dt_k\ \to 0
\end{equation}

\vspace{5mm}
\noindent
if $N\to\infty,$ 
where constant $C_k$ 
depends only
on the multiplicity $k$ of the multiple Wiener stochastic integral.

The relations (\ref{chain12}) and (\ref{chain20}) prove the equality 
(\ref{Wi110}).
Using (\ref{Wi110}) and the isometry property of the Ito stochastic integral, we have

\begin{equation}
\label{wi1001}
J[\psi^{(k)}]_{T,t}^{(i_1\ldots i_k)}=\int\limits_t^T\psi_k(t_k) \ldots \int\limits_t^{t_{2}}
\psi_1(t_1) d{\bf w}_{t_1}^{(i_1)}\ldots
d{\bf w}_{t_k}^{(i_k)}=J'[K]_{T,t}^{(i_1\ldots i_k)}\ \ \ \hbox{w.\ p.\ 1},
\end{equation}

\vspace{3mm}
\noindent
where 
$K=K(t_1,\ldots,t_k)$ is defined by (\ref{ppp}), i.e.

\begin{equation}
\label{chain200}
K(t_1,\ldots,t_k)=
\left\{\begin{matrix}
\psi_1(t_1)\ldots \psi_k(t_k),\ &t_1<\ldots<t_k\cr\cr\cr
0,\ &\hbox{\rm otherwise}
\end{matrix}
\right.,
\end{equation}

\vspace{4mm}
\noindent
where $\psi_1(\tau),\ldots,\psi_k(\tau)\in L_2([t,T]),$\ $t_1,\ldots,t_k\in [t, T]$ $(k\ge 2)$ and 
$K(t_1)\equiv\psi_1(t_1)$ for $t_1\in[t, T].$

Applying (\ref{wi1001}) and the linearity property of the Ito stochastic integral, we obtain

\vspace{1mm}
$$
J[\psi^{(k)}]_{T,t}^{(i_1\ldots i_k)}=J'[K]_{T,t}^{(i_1\ldots i_k)}=
$$

\vspace{2mm}
\begin{equation}
\label{chain102}
=\sum_{j_1=0}^{p_1}\ldots
\sum_{j_k=0}^{p_k}
C_{j_k\ldots j_1}
J'[\phi_{j_1}\ldots \phi_{j_k}]_{T,t}^{(i_1\ldots i_k)}+
J'[R_{p_1\ldots p_k}]_{T,t}^{(i_1\ldots i_k)}\ \ \ \hbox{w.~p.~1,}
\end{equation}

\vspace{4mm}
\noindent
where
\begin{equation}
\label{chain30001}
R_{p_1\ldots p_k}(t_1,\ldots,t_k)\stackrel{{\rm def}}{=}
K(t_1,\ldots,t_k)-
\sum_{j_1=0}^{p_1}\ldots
\sum_{j_k=0}^{p_k}
C_{j_k\ldots j_1}
\prod_{l=1}^k\phi_{j_l}(t_l)
\end{equation}

\vspace{3mm}
\noindent
and
\begin{equation}
\label{chain300}
C_{j_k\ldots j_1}=\int\limits_{[t,T]^k}
K(t_1,\ldots,t_k)\prod_{l=1}^{k}\phi_{j_l}(t_l)dt_1\ldots dt_k
\end{equation}

\vspace{4mm}
\noindent
is the Fourier coefficient corresponding to $K(t_1,\ldots,t_k).$

Using the Ito formula, we have

\vspace{2mm}
$$
\sum\limits_{(j_1,\ldots,j_q)}\int\limits_t^T \phi_{j_q}(t_q)\ldots
\int\limits_t^{t_2}\phi_{j_1}(t_1)
d{\bf w}_{t_1}^{(i_1)}\ldots {\bf w}_{t_q}^{(i_q)}\times
$$

\vspace{4mm}
$$
\times \sum\limits_{(j_1',\ldots,j_n')}
\int\limits_t^T \phi_{j_n'}(t_n')\ldots 
\int\limits_t^{t_2'}\phi_{j_1'}(t_1')d{\bf w}_{t_1'}^{(g)}\ldots {\bf w}_{t_n'}^{(g)}=
$$

\vspace{4mm}
$$
=\sum\limits_{(j_1,\ldots,j_q,j_1',\ldots,j_n')}
\int\limits_t^T \phi_{j_q}(t_q)\ldots \int\limits_t^{t_2}\phi_{j_1}(t_1)
\int\limits_t^{t_1}\phi_{j_n'}(t_n')\ldots \int\limits_t^{t_2'}
\phi_{j_1'}(t_1')\times
$$

\vspace{5mm}
\begin{equation}
\label{new1000}
\times d{\bf w}_{t_1'}^{(g)}\ldots d{\bf w}_{t_n'}^{(g)}d{\bf w}_{t_1}^{(i_1)}\ldots d{\bf w}_{t_q}^{(i_q)}
\end{equation}

\vspace{6mm}
\noindent
w.~p.~1, where $g=0$ or $g=1,$\  $n, q\in \mathbb{N},$\ $i_1,\ldots,i_q\ne 0,\ 1,$

$$
\sum\limits_{(j_1,\ldots,j_k)}
$$

\vspace{3mm}
\noindent
means the sum with respect to all possible permutations $(j_1,\ldots,j_k)$.
At the same time if $j_r$ swapped with $j_d$ in the permutation $(j_1,\ldots,j_k)$, then
$i_r$ swapped with $i_d$ in the permutation $(i_1,\ldots,i_k).$

The detailed proof of (\ref{new1000}) will be given in Sect.~18 (see the proof of Theorem~20).
The equality (\ref{new1000}) means that (see (\ref{Wi110}))

\vspace{2mm}
$$
J'[\phi_{j_1}\ldots\phi_{j_q}]_{T,t}^{(i_1\ldots i_q)}
\cdot J'[\phi_{j_1'}\ldots\phi_{j_n'}]^{(g\ldots g)}_{T,t}=
$$

\vspace{3mm}
\begin{equation}
\label{2023abc111}
=J'[\phi_{j_1}\ldots\phi_{j_q}\phi_{j_1'}\ldots\phi_{j_n'}]_{T,t}^{(i_1\ldots i_q g\ldots g)}
\end{equation}

\vspace{5mm}
\noindent
w.~p.~1, where $g=0$ or $g=1$,\ $n, q\in \{0\}\cup \mathbb{N},$\ $i_1,\ldots,i_q\ne 0,\ 1,$ 
and $J'[\phi_{j_1}\ldots\phi_{j_q}]_{T,t}^{(i_1\ldots i_q)}\stackrel{\sf def}{=}1$ for $q=0.$

Using the equality (\ref{2023abc111}), we get (\ref{ziko30}) for the case
of an arbitrary complete orthonormal system $\left\{\phi_j(x)\right\}_{j=0}^{\infty}$
of functions in $L_2([t,T])$.

Using Theorem~9.6.9 \cite{Kuo} (also see 
\cite{ito1951}, Theorem~3.1) and (\ref{ziko100}) (also see Theorem~21 below), we get

$$
J'[\phi_{j_1}\ldots \phi_{j_k}]_{T,t}^{(i_1\ldots i_k)}=
$$

\vspace{5mm}
$$
=\prod_{l=1}^k\left({\bf 1}_{\{m_l=0\}}+{\bf 1}_{\{m_l>0\}}\left\{
\begin{matrix}
H_{n_{1,l}}\left(\zeta_{j_{h_{1,l}}}^{(i_l)}\right)\ldots 
H_{n_{d_l,l}}\left(\zeta_{j_{h_{d_l,l}}}^{(i_l)}\right),\ 
&\hbox{\rm if}\ \ \ 
i_l\ne 0\cr\cr
\left(\zeta_{j_{h_{1,l}}}^{(0)}\right)^{n_{1,l}}\ldots
\left(\zeta_{j_{h_{d_l,l}}}^{(0)}\right)^{n_{d_l,l}},\  &\hbox{\rm if}\ \ \ 
i_l=0
\end{matrix}\right.\ \right)=
$$

\vspace{8mm}
$$
=\prod_{l=1}^k\zeta_{j_l}^{(i_l)}
+\sum\limits_{r=1}^{[k/2]}
(-1)^r \times
$$

\vspace{2mm}
\begin{equation}
\label{chain401}
\times\sum_{\stackrel{(\{\{g_1, g_2\}, \ldots, 
\{g_{2r-1}, g_{2r}\}\}, \{q_1, \ldots, q_{k-2r}\})}
{{}_{\{g_1, g_2, \ldots, 
g_{2r-1}, g_{2r}, q_1, \ldots, q_{k-2r}\}=\{1, 2, \ldots, k\}}}}
\prod\limits_{s=1}^r
{\bf 1}_{\{i_{g_{{}_{2s-1}}}=~i_{g_{{}_{2s}}}\ne 0\}}
\Biggl.{\bf 1}_{\{j_{g_{{}_{2s-1}}}=~j_{g_{{}_{2s}}}\}}
\prod_{l=1}^{k-2r}\zeta_{j_{q_l}}^{(i_{q_l})}
\end{equation}

\vspace{6mm}
\noindent
w.~p.~1, where notations are the same as in Theorems 2 and 10;
the multiple Wiener stochastic integral 
$J'[\phi_{j_1}\ldots \phi_{j_k}]_{T,t}^{(i_1\ldots i_k)}$ is defined by 
(\ref{WiI}).

Again applying (\ref{Wi110}), we have

\vspace{2mm}
$$
J'[R_{p_1\ldots p_k}]_{T,t}^{(i_1\ldots i_k)}
=
\sum_{(t_1,\ldots,t_k)}
\int\limits_{t}^{T}
\ldots
\int\limits_{t}^{t_2}
\Biggl(K(t_1,\ldots,t_k)-\Biggr.
$$

\vspace{4mm}
\begin{equation}
\label{wi2005}
-\Biggl.
\sum_{j_1=0}^{p_1}\ldots
\sum_{j_k=0}^{p_k}
C_{j_k\ldots j_1}
\prod_{l=1}^k\phi_{j_l}(t_l)\Biggr)
d{\bf w}_{t_1}^{(i_1)}
\ldots
d{\bf w}_{t_k}^{(i_k)},
\end{equation}

\vspace{5mm}
\noindent
where permutations $(t_1,\ldots,t_k)$ when summing are performed only 
in the values $d{\bf w}_{t_1}^{(i_1)}
\ldots $
$d{\bf w}_{t_k}^{(i_k)}$. At the same time the indices near 
upper limits of integration in the iterated stochastic integrals 
are changed correspondently and if $t_r$ swapped with $t_q$ in the  
permutation $(t_1,\ldots,t_k)$, then $i_r$ swapped with $i_q$ in the 
permutation $(i_1,\ldots,i_k).$
In addition, the multiple Wiener stochastic integral
$J'[R_{p_1\ldots p_k}]_{T,t}^{(i_1\ldots i_k)}$ is defined by 
(\ref{WiI}).

According to Lemma 2 for $\Phi(t_1,\ldots,t_k)\in L_2(D_k)$, 
(\ref{sos1z}), and (\ref{riemann}) for Lebesgue integrals, we have

\vspace{2mm}
$$
{\sf M}\left\{\left(J'[R_{p_1\ldots p_k}]_{T,t}^{(i_1\ldots i_k)}\right)^2\right\}
\le 
$$

\vspace{3mm}
$$
\le C_k
\sum_{(t_1,\ldots,t_k)}
\int\limits_{t}^{T}
\ldots
\int\limits_{t}^{t_2}
\left(K(t_1,\ldots,t_k)-
\sum_{j_1=0}^{p_1}\ldots
\sum_{j_k=0}^{p_k}
C_{j_k\ldots j_1}
\prod_{l=1}^k\phi_{j_l}(t_l)\right)^2
dt_1
\ldots
dt_k=
$$

\vspace{3mm}
\begin{equation}
\label{chain7771}
=C_k\int\limits_{[t,T]^k}
\left(K(t_1,\ldots,t_k)-
\sum_{j_1=0}^{p_1}\ldots
\sum_{j_k=0}^{p_k}
C_{j_k\ldots j_1}
\prod_{l=1}^k\phi_{j_l}(t_l)\right)^2
dt_1
\ldots
dt_k\to 0
\end{equation}

\vspace{5mm}
\noindent
if $p_1,\ldots,p_k\to\infty,$ where constant $C_k$ 
depends only
on the multiplicity $k$ of the 
iterated  Ito stochastic integral
$J[\psi^{(k)}]_{T,t}^{(i_1\ldots i_k)}$.

Thus, the following theorem is proved.

\vspace{2mm}
         
{\bf Theorem 12}\ \cite{20a}, \cite{2023xxx1}\ (generalization of Theorems 1, 2, and 10).\ 
{\it Suppose that
the condition {\rm ($\star\star$)} is fulfilled
for the multi-index $(i_1 \ldots i_k)$ {\rm (}see Sect.~{\rm 14)} 
and the condition {\rm (\ref{ziko999})} is also 
fulfilled.
Furthermore$,$ let 
$\psi_l(\tau)\in L_2([t, T])$ $(l=$ $1,\ldots, k)$ and
$\{\phi_j(x)\}_{j=0}^{\infty}$ is an arbitrary complete orthonormal system  
of functions in the space $L_2([t,T]).$
Then the following expansions

\vspace{1mm}
\begin{equation}
\label{new9999}
J[\psi^{(k)}]_{T,t}^{(i_1\ldots i_k)}=
\hbox{\vtop{\offinterlineskip\halign{
\hfil#\hfil\cr
{\rm l.i.m.}\cr
$\stackrel{}{{}_{p_1,\ldots,p_k\to \infty}}$\cr
}} }
\sum\limits_{j_1=0}^{p_1}\ldots
\sum\limits_{j_k=0}^{p_k}
C_{j_k\ldots j_1}\times
$$

\vspace{4mm}
$$
\times
\prod_{l=1}^k\left({\bf 1}_{\{m_l=0\}}+{\bf 1}_{\{m_l>0\}}\left\{
\begin{matrix}
H_{n_{1,l}}\left(\zeta_{j_{h_{1,l}}}^{(i_l)}\right)\ldots 
H_{n_{d_l,l}}\left(\zeta_{j_{h_{d_l,l}}}^{(i_l)}\right),\ 
&\hbox{\rm if}\ \ \ 
i_l\ne 0\cr\cr
\left(\zeta_{j_{h_{1,l}}}^{(0)}\right)^{n_{1,l}}\ldots
\left(\zeta_{j_{h_{d_l,l}}}^{(0)}\right)^{n_{d_l,l}},\  &\hbox{\rm if}\ \ \ 
i_l=0
\end{matrix}\right.\ \right),
\end{equation}

\vspace{7mm}
$$
J[\psi^{(k)}]_{T,t}^{(i_1\ldots i_k)}=
\hbox{\vtop{\offinterlineskip\halign{
\hfil#\hfil\cr
{\rm l.i.m.}\cr
$\stackrel{}{{}_{p_1,\ldots,p_k\to \infty}}$\cr
}} }
\sum\limits_{j_1=0}^{p_1}\ldots
\sum\limits_{j_k=0}^{p_k}
C_{j_k\ldots j_1}\Biggl(
\prod_{l=1}^k\zeta_{j_l}^{(i_l)}+\sum\limits_{r=1}^{[k/2]}
(-1)^r \times
\Biggr.
$$

\vspace{3mm}
\begin{equation}
\label{razzar1}
\times
\sum_{\stackrel{(\{\{g_1, g_2\}, \ldots, 
\{g_{2r-1}, g_{2r}\}\}, \{q_1, \ldots, q_{k-2r}\})}
{{}_{\{g_1, g_2, \ldots, 
g_{2r-1}, g_{2r}, q_1, \ldots, q_{k-2r}\}=\{1, 2, \ldots, k\}}}}
\prod\limits_{s=1}^r
{\bf 1}_{\{i_{g_{{}_{2s-1}}}=~i_{g_{{}_{2s}}}\ne 0\}}
\Biggl.{\bf 1}_{\{j_{g_{{}_{2s-1}}}=~j_{g_{{}_{2s}}}\}}
\prod_{l=1}^{k-2r}\zeta_{j_{q_l}}^{(i_{q_l})}\Biggr)
\end{equation}

\vspace{5mm}
\noindent
con\-verg\-ing in the mean-square sense are valid$,$
where $[x]$ is an integer part of a real number $x;$\ \
$n_{1,l}+n_{2,l}+\ldots+n_{d_l,l}=m_l;$\ \ $n_{1,l}, n_{2,l}, \ldots, n_{d_l,l}=1,\ldots, m_l;$\ \ 
$d_l=1,\ldots,m_l;$\ \ $l=1,\ldots,k;$\ \
$m_1+\ldots+m_k=k;$\ \ 
the numbers $m_1,\ldots,m_k,$\ $g_1,\ldots,g_k$
depend on $(i_1,\ldots,i_k)$ and 
the numbers $n_{1,l},\ldots,n_{d_l,l},$\ $h_{1,l},\ldots,h_{d_l,l},$\ $d_l$
depend on $\{j_1,\ldots,j_k\};$ moreover$,$ $\left\{j_{g_1},\ldots,j_{g_k}\right\}
=\{j_1,\ldots,j_k\};$
$H_n(x)$ is the Hermite polynomial {\rm (\ref{ziko500});}
another
notations are the same as in Theorems {\rm 1, 2, and 10}.}

\vspace{2mm}

Replacing the function $K(t_1,\ldots,t_k)$ by 
$\Phi(t_1,\ldots,t_k)$ we get the following theorem.

\vspace{2mm}

{\bf Theorem 13}\ \cite{20a}, \cite{2023xxx1}\ (generalization of Theorem 11).\ {\it Suppose that
the condition {\rm ($\star\star$)} is fulfilled
for the multi-index $(i_1 \ldots i_k)$ {\rm (}see Sect.~{\rm 14)} 
and the condition {\rm (\ref{ziko999})} is also 
fulfilled.
Furthermore$,$ let 
$\Phi(t_1,\ldots,t_k)\in L_2([t, T]^k)$ and
$\{\phi_j(x)\}_{j=0}^{\infty}$ is an arbitrary complete orthonormal system  
of functions in the space $L_2([t,T]).$
Then the following expansions

\vspace{2mm}
$$
J'[\Phi]_{T,t}^{(i_1\ldots i_k)}=
\hbox{\vtop{\offinterlineskip\halign{
\hfil#\hfil\cr
{\rm l.i.m.}\cr
$\stackrel{}{{}_{p_1,\ldots,p_k\to \infty}}$\cr
}} }
\sum\limits_{j_1=0}^{p_1}\ldots
\sum\limits_{j_k=0}^{p_k}
C_{j_k\ldots j_1}\times
$$

\vspace{4mm}
\begin{equation}
\label{lllnew1}
\times
\prod_{l=1}^k\left({\bf 1}_{\{m_l=0\}}+{\bf 1}_{\{m_l>0\}}\left\{
\begin{matrix}
H_{n_{1,l}}\left(\zeta_{j_{h_{1,l}}}^{(i_l)}\right)\ldots 
H_{n_{d_l,l}}\left(\zeta_{j_{h_{d_l,l}}}^{(i_l)}\right),\ 
&\hbox{\rm if}\ \ \ 
i_l\ne 0\cr\cr
\left(\zeta_{j_{h_{1,l}}}^{(0)}\right)^{n_{1,l}}\ldots
\left(\zeta_{j_{h_{d_l,l}}}^{(0)}\right)^{n_{d_l,l}},\  &\hbox{\rm if}\ \ \ 
i_l=0
\end{matrix}\right.\ \right),
\end{equation}

\vspace{7mm}
$$
J'[\Phi]_{T,t}^{(i_1\ldots i_k)}=
\hbox{\vtop{\offinterlineskip\halign{
\hfil#\hfil\cr
{\rm l.i.m.}\cr
$\stackrel{}{{}_{p_1,\ldots,p_k\to \infty}}$\cr
}} }
\sum\limits_{j_1=0}^{p_1}\ldots
\sum\limits_{j_k=0}^{p_k}
C_{j_k\ldots j_1}\Biggl(
\prod_{l=1}^k\zeta_{j_l}^{(i_l)}+\sum\limits_{r=1}^{[k/2]}
(-1)^r \times
\Biggr.
$$

\vspace{5mm}
$$
\times
\sum_{\stackrel{(\{\{g_1, g_2\}, \ldots, 
\{g_{2r-1}, g_{2r}\}\}, \{q_1, \ldots, q_{k-2r}\})}
{{}_{\{g_1, g_2, \ldots, 
g_{2r-1}, g_{2r}, q_1, \ldots, q_{k-2r}\}=\{1, 2, \ldots, k\}}}}
\prod\limits_{s=1}^r
{\bf 1}_{\{i_{g_{{}_{2s-1}}}=~i_{g_{{}_{2s}}}\ne 0\}}
\Biggl.{\bf 1}_{\{j_{g_{{}_{2s-1}}}=~j_{g_{{}_{2s}}}\}}
\prod_{l=1}^{k-2r}\zeta_{j_{q_l}}^{(i_{q_l})}\Biggr)
$$

\vspace{5mm}
\noindent
con\-verg\-ing in the mean-square sense are valid$,$
where
$n_{1,l}+n_{2,l}+\ldots+n_{d_l,l}=m_l;$\ \ $n_{1,l}, n_{2,l}, \ldots, n_{d_l,l}=1,\ldots, m_l;$\ \ 
$d_l=1,\ldots,m_l;$\ \ $l=1,\ldots,k;$\ \
$m_1+\ldots+m_k=k;$\ \ 
the numbers $m_1,\ldots,m_k,$\ $g_1,\ldots,g_k$
depend on $(i_1,\ldots,i_k)$ and 
the numbers $n_{1,l},\ldots,n_{d_l,l},$\ $h_{1,l},\ldots,h_{d_l,l},$\ $d_l$
depend on $\{j_1,\ldots,j_k\};$ moreover$,$ $\left\{j_{g_1},\ldots,j_{g_k}\right\}
=\{j_1,\ldots,j_k\};$
the multiple Wiener stochastic integral
$J'[\Phi]_{T,t}^{(i_1\ldots i_k)}$ is defined by 
{\rm (\ref{WiI});} $H_n(x)$ is the Hermite polynomial {\rm (\ref{ziko500});}
another
notations are the same as in Theorem {\rm 9,  11}.}

\vspace{2mm}

It should be noted that an analogue of the expansion (\ref{lllnew1}) was obtained
in \cite{Rybakov1000} for the case $i_1,\ldots,i_k=1,\ldots,m.$  The proof in \cite{Rybakov1000} is 
different from the proof given in this section and Sect.~18.
Note that the results of work \cite{Rybakov1000}, as well as 
the results of this article, are based on our idea 
\cite{7} (2006) on the expansion of the kernel (\ref{ppp}) (or $\Phi(t_1,\ldots,t_k)\in L_2([t,T]^k)$)
into a generalized multiple Fourier series 
(see \cite{7}, Chapter~5, Theorem~5.1, pp.~235-245 
or \cite{20a}, Sect.~1.1.3 for details).

\vspace{5mm}

\section{Exact Calculation of the Mean-Square Error in Theorems~1, 2, and 12}

\vspace{5mm}

In this section, we will use the multiple Wiener 
stochastic integral 
with respect to the components of a multidimensional Wiener
process
to generalize theorem on the 
exact calculation of the mean-square error in Theorems~1, 2.
More precisely, we will generalize the following theorem.

\vspace{2mm}

{\bf Theorem 14}\ \cite{20a}, \cite{2023xxx1}, \cite{26}.\
{\it Suppose that
every $\psi_l(\tau)$ $(l=1,\ldots, k)$ is a continuous nonrandom function on 
$[t, T]$ and
$\{\phi_j(x)\}_{j=0}^{\infty}$ is a complete orthonormal system  
of functions in the space $L_2([t,T]),$ 
each function $\phi_j(x)$ of which 
for finite $j$ satisfies the condition 
$(\star)$ {\rm (}see Sect.~{\rm 4)}.
Then

\vspace{1mm}
$$
{\sf M}\left\{\left(J[\psi^{(k)}]_{T,t}-
J[\psi^{(k)}]_{T,t}^p\right)^2\right\}
= \int\limits_{[t,T]^k} K^2(t_1,\ldots,t_k)
dt_1\ldots dt_k - 
$$

\vspace{1mm}
\begin{equation}
\label{tttr11}
- \sum_{j_1=0}^{p}\ldots\sum_{j_k=0}^{p}
C_{j_k\ldots j_1}
{\sf M}\left\{J[\psi^{(k)}]_{T,t}
\sum\limits_{(j_1,\ldots,j_k)}
\int\limits_t^T \phi_{j_k}(t_k)
\ldots
\int\limits_t^{t_{2}}\phi_{j_{1}}(t_{1})
d{\bf f}_{t_1}^{(i_1)}\ldots
d{\bf f}_{t_k}^{(i_k)}\right\},
\end{equation}

\vspace{4mm}
\noindent
where
$$
J[\psi^{(k)}]_{T,t}=\int\limits_t^T\psi_k(t_k) \ldots \int\limits_t^{t_{2}}
\psi_1(t_1) d{\bf f}_{t_1}^{(i_1)}\ldots
d{\bf f}_{t_k}^{(i_k)},
$$

\vspace{2mm}
\begin{equation}
\label{yeee2}
J[\psi^{(k)}]_{T,t}^p=
\sum_{j_1=0}^{p}\ldots\sum_{j_k=0}^{p}
C_{j_k\ldots j_1}\left(
\prod_{l=1}^k\zeta_{j_l}^{(i_l)}-S_{j_1,\ldots,j_k}^{(i_1\ldots i_k)}
\right),
\end{equation}

\vspace{3mm}

\begin{equation}
\label{ppp1}
S_{j_1,\ldots,j_k}^{(i_1\ldots i_k)}=
\hbox{\vtop{\offinterlineskip\halign{
\hfil#\hfil\cr
{\rm l.i.m.}\cr
$\stackrel{}{{}_{N\to \infty}}$\cr
}} }\sum_{(l_1,\ldots,l_k)\in {\rm G}_k}
\phi_{j_{1}}(\tau_{l_1})
\Delta{\bf f}_{\tau_{l_1}}^{(i_1)}\ldots
\phi_{j_{k}}(\tau_{l_k})
\Delta{\bf f}_{\tau_{l_k}}^{(i_k)},
\end{equation}

\vspace{6mm}
\noindent
the Fourier coefficient $C_{j_k\ldots j_1}$ has the form {\rm (\ref{ppppa})},

\vspace{-1mm}
\begin{equation}
\label{rr232}
\zeta_{j}^{(i)}=
\int\limits_t^T \phi_{j}(s) d{\bf f}_s^{(i)}
\end{equation}

\vspace{3mm}
\noindent
are independent standard Gaussian random variables
for various
$i$ or $j$ $(i=1,\ldots,m),$

\vspace{-1mm}
$$
\sum\limits_{(j_1,\ldots,j_k)}
$$ 

\vspace{3mm}
\noindent
means the sum with respect to all
possible permutations 
$(j_1,\ldots,j_k).$ At the same time if 
$j_r$ swapped with $j_q$ in the permutation $(j_1,\ldots,j_k)$,
then $i_r$ swapped with $i_q$ in the permutation
$(i_1,\ldots,i_k)$ {\rm (}see {\rm (\ref{tttr11}));}
another notations are the same as in Theorem {\rm 1.}}

\vspace{2mm}

Let us generalize Theorem~14 
to the case of an arbitrary 
complete orthonormal system of functions in the space $L_2([t, T])$
and $\psi_1(\tau),$ $\ldots,\psi_k(\tau)\in L_2([t, T]).$

\vspace{2mm}

{\bf Theorem 15}\ \cite{20a}, \cite{2023xxx1}.\
{\it Suppose that
$\psi_1(\tau),\ldots,\psi_k(\tau)\in L_2([t, T])$ 
and
$\{\phi_j(x)\}_{j=0}^{\infty}$ is an arbitrary complete orthonormal system  
of functions in the space $L_2([t,T]).$ 
Then

\vspace{3mm}
$$
{\sf M}\left\{\left(J[\psi^{(k)}]_{T,t}-
J[\psi^{(k)}]_{T,t}^p\right)^2\right\}
= \int\limits_{[t,T]^k} K^2(t_1,\ldots,t_k)
dt_1\ldots dt_k - 
$$

\vspace{1mm}
\begin{equation}
\label{chain100}
- \sum_{j_1=0}^{p}\ldots\sum_{j_k=0}^{p}
C_{j_k\ldots j_1}
{\sf M}\left\{J[\psi^{(k)}]_{T,t}
\sum\limits_{(j_1,\ldots,j_k)}
\int\limits_t^T \phi_{j_k}(t_k)
\ldots
\int\limits_t^{t_{2}}\phi_{j_{1}}(t_{1})
d{\bf f}_{t_1}^{(i_1)}\ldots
d{\bf f}_{t_k}^{(i_k)}\right\},
\end{equation}

\vspace{6mm}
\noindent
where
$$
J[\psi^{(k)}]_{T,t}=\int\limits_t^T\psi_k(t_k) \ldots \int\limits_t^{t_{2}}
\psi_1(t_1) d{\bf f}_{t_1}^{(i_1)}\ldots
d{\bf f}_{t_k}^{(i_k)},
$$

\vspace{2mm}
\begin{equation}
\label{chain101}
J[\psi^{(k)}]_{T,t}^p=
\sum_{j_1=0}^{p}\ldots\sum_{j_k=0}^{p}
C_{j_k\ldots j_1} J'[\phi_{j_1}\ldots \phi_{j_k}]_{T,t}^{(i_1\ldots i_k)},
\end{equation}

\vspace{5mm}

\noindent
$J'[\phi_{j_1}\ldots \phi_{j_k}]_{T,t}^{(i_1\ldots i_k)}$
is the multiple Wiener stochastic integral 
defined by {\rm (\ref{WiI}),}
the Fourier coefficient $C_{j_k\ldots j_1}$ has the form {\rm (\ref{chain300}),}
$K(t_1,\ldots,t_k)$ is defined by {\rm (\ref{chain200}),}

\vspace{-1mm}
$$
\zeta_{j}^{(i)}=
\int\limits_t^T \phi_{j}(s) d{\bf f}_s^{(i)}
$$

\vspace{2mm}
\noindent
are independent standard Gaussian random variables
for various
$i$ or $j$ $(i=1,\ldots,m),$

\vspace{-1mm}
$$
\sum\limits_{(j_1,\ldots,j_k)}
$$ 

\vspace{2mm}
\noindent
means the sum with respect to all
possible permutations 
$(j_1,\ldots,j_k).$ At the same time if 
$j_r$ swapped with $j_q$ in the permutation $(j_1,\ldots,j_k)$,
then $i_r$ swapped with $i_q$ in the permutation
$(i_1,\ldots,i_k)$ {\rm (}see {\rm (\ref{chain100})).}}

\vspace{2mm}

{\bf Proof.}\ First, note that the formula (\ref{chain101}) 
appears due to the equality (\ref{chain102}).
Using the equality (\ref{Wi110}), we get

\begin{equation}
\label{chain103}
J'[\phi_{j_1}\ldots \phi_{j_k}]_{T,t}^{(i_1\ldots i_k)}=\sum_{(t_1,\ldots,t_k)}
\int\limits_t^T \phi_{j_k}(t_k)
\ldots
\int\limits_t^{t_{2}}\phi_{j_{1}}(t_{1})
d{\bf f}_{t_1}^{(i_1)}\ldots
d{\bf f}_{t_k}^{(i_k)}\ \ \ \hbox{w.\ p.\ 1},
\end{equation}

\vspace{3mm}
\noindent
where permutations $(t_1,\ldots,t_k)$ when summing are 
performed only in the values
$d{\bf f}_{t_1}^{(i_1)}
\ldots $
$d{\bf f}_{t_k}^{(i_k)}.$ At the same time the indices near 
upper 
limits of integration in the iterated stochastic integrals are changed 
correspondently and if $t_r$ swapped with $t_q$ in the  
permutation $(t_1,\ldots,t_k)$, then $i_r$ swapped with $i_q$ in 
the permutation $(i_1,\ldots,i_k).$

It is easy to see that the equality (\ref{chain103}) can be written in the form

\vspace{-1mm}
\begin{equation}
\label{chain104}
J'[\phi_{j_1}\ldots \phi_{j_k}]_{T,t}^{(i_1\ldots i_k)}=
\sum\limits_{(j_1,\ldots,j_k)}
\int\limits_t^T \phi_{j_k}(t_k)
\ldots
\int\limits_t^{t_{2}}\phi_{j_{1}}(t_{1})
d{\bf f}_{t_1}^{(i_1)}\ldots
d{\bf f}_{t_k}^{(i_k)}\ \ \ \hbox{w.\ p.\ 1},
\end{equation}

\vspace{4mm}
\noindent
where 
$$
\sum\limits_{(j_1,\ldots,j_k)}
$$ 

\vspace{3mm}
\noindent
means the sum with respect to all
possible permutations 
$(j_1,\ldots,j_k).$ At the same time if 
$j_r$ swapped with $j_q$ in the permutation $(j_1,\ldots,j_k)$,
then $i_r$ swapped with $i_q$ in the permutation $(i_1,\ldots,i_k)$.

Further proof of Theorem~15 is based on the equality 
(\ref{chain104}) and is similar to the proof of Theorem 14 in \cite{20a}, \cite{26}.
Theorem~15 is proved.

The equalities (\ref{chain401}) and (\ref{chain104}) 
allow us to formulate the following modification of Theorem~12.

\vspace{2mm}

{\bf Theorem 16.}\ 
{\it Suppose that
$\psi_1(\tau),\ldots,\psi_k(\tau)\in L_2([t, T])$ and
$\{\phi_j(x)\}_{j=0}^{\infty}$ is an arbitrary complete orthonormal system  
of functions in the space $L_2([t,T]).$
Then the following expansion

\vspace{2mm}
$$
J[\psi^{(k)}]_{T,t}^{(i_1\ldots i_k)}=
\hbox{\vtop{\offinterlineskip\halign{
\hfil#\hfil\cr
{\rm l.i.m.}\cr
$\stackrel{}{{}_{p_1,\ldots,p_k\to \infty}}$\cr
}} }
\sum\limits_{j_1=0}^{p_1}\ldots
\sum\limits_{j_k=0}^{p_k}
C_{j_k\ldots j_1}J'[\phi_{j_1}\ldots \phi_{j_k}]_{T,t}^{(i_1\ldots i_k)}=
$$

\vspace{3mm}
$$
=\hbox{\vtop{\offinterlineskip\halign{
\hfil#\hfil\cr
{\rm l.i.m.}\cr
$\stackrel{}{{}_{p_1,\ldots,p_k\to \infty}}$\cr
}} }
\sum\limits_{j_1=0}^{p_1}\ldots
\sum\limits_{j_k=0}^{p_k}
C_{j_k\ldots j_1}
\sum\limits_{(j_1,\ldots,j_k)}
\int\limits_t^T \phi_{j_k}(t_k)
\ldots
\int\limits_t^{t_{2}}\phi_{j_{1}}(t_{1})
d{\bf f}_{t_1}^{(i_1)}\ldots
d{\bf f}_{t_k}^{(i_k)}
$$

\vspace{7mm}
\noindent
con\-verg\-ing in the mean-square sense is valid$,$
where $i_1,\ldots,i_k=1,\ldots,m;$
another
notations are the same as in Theorems {\rm 1, 12}.}

\vspace{2mm}

Consider the following obvious generalization of Theorem 3.

\vspace{2mm}

{\bf Theorem 17}\ \cite{20a}, \cite{2023xxx1}.\
{\it Suppose that
$\psi_1(\tau),\ldots,\psi_k(\tau)\in L_2([t, T])$ 
and
$\{\phi_j(x)\}_{j=0}^{\infty}$ is an arbitrary complete orthonormal system  
of functions in the space $L_2([t,T]).$ 
Then the estimate

\vspace{1mm}
$$
{\sf M}\left\{\left(
J[\psi^{(k)}]_{T,t}-J[\psi^{(k)}]_{T,t}^{p_1,\ldots,p_k}
\right)^2\right\}
\le 
$$

\vspace{1mm}
\begin{equation}
\label{z1new}
\le k!\left(~\int\limits_{[t,T]^k}
K^2(t_1,\ldots,t_k)
dt_1\ldots dt_k -\sum_{j_1=0}^{p_1}\ldots
\sum_{j_k=0}^{p_k}C^2_{j_k\ldots j_1}\right)
\end{equation}

\vspace{5mm}
\noindent
is valid for the following cases{\rm :}

\vspace{2mm}

{\rm 1.}\ $i_1,\ldots,i_k=1,\ldots,m$\ \ and\ \ $0<T-t<\infty,$

\vspace{1mm}

{\rm 2.}\ $i_1,\ldots,i_k=0, 1,\ldots,m,$\ \ $i_1^2+\ldots+i_k^2>0,$\ \
and\ \ $0<T-t<1,$

\vspace{3mm}
\noindent
where $J[\psi^{(k)}]_{T,t}$ is the stochastic integral {\rm (\ref{sodom20}),}
$J[\psi^{(k)}]_{T,t}^{p_1,\ldots,p_k}$ is the 
expression on the right-hand side of {\rm (\ref{razzar1})} before
passing to the limit 
$\hbox{\vtop{\offinterlineskip\halign{
\hfil#\hfil\cr
{\rm l.i.m.}\cr
$\stackrel{}{{}_{p_1,\ldots,p_k\to \infty}}$\cr
}} };$ another 
notations are the same as in Theorem {\rm 1, 2, 12}.
}

\vspace{2mm}

In addition, under the conditions of Theorem 17 we have the estimate
(also see (\ref{2026ch1001s11}))

\vspace{1mm}
$$
{\sf M}\left\{\left(J[\psi^{(k)}]_{T,t}-
J[\psi^{(k)}]_{T,t}^{p_1,\ldots,p_k}\right)^{2n}\right\}\le
$$

\vspace{3mm}
$$
\le
(k!)^{n} (2n-1)^{nk}\  \times
$$

\vspace{1mm}
$$
\times\ 
\left(~
\int\limits_{[t,T]^k}
K^2(t_1,\ldots,t_k)
dt_1\ldots dt_k -\sum_{j_1=0}^{p_1}\ldots
\sum_{j_k=0}^{p_k}C^2_{j_k\ldots j_1}
\right)^n.
$$

\vspace{5mm}

\section{Generalization of Theorems 4, 5 to the Case of an Arbitrary 
Complete Ortho\-nor\-mal with Weight $r(x)\ge 0$ System of Functions in the Space $L_2([t, T])$
and $\psi_1(x)\sqrt{r(x)},$ $\ldots,$ $\psi_k(x)\sqrt{r(x)}$ $\in $ $L_2([t, T])$}

\vspace{5mm}

In this section, we will use the multiple Wiener 
stochastic integral 
with respect 
to the components of a multidimensional Wiener
process
to generalize Theorems 
4, 5 to the case of an arbitrary 
complete ortho\-nor\-mal with weight $r(x)\ge 0$ system of functions in the 
space $L_2([t, T])$
and $\psi_1(x)\sqrt{r(x)},$ $\ldots,$ $\psi_k(x)\sqrt{r(x)}$ $\in $ $L_2([t, T])$.
From the results of Sect.~8, 15 we obtain the following two theorems.

\vspace{2mm}

{\bf Theorem 18}\ \cite{20a}, \cite{2023xxx1}.\
{\it Suppose that $\psi_1(x)\sqrt{r(x)},\ldots,
\psi_k(x)\sqrt{r(x)}\in L_2([t, T]),$ 
where $r(x)\ge 0.$
Moreover$,$ let 
$$
\left\{\Psi_j(x)\sqrt{r(x)}\right\}_{j=0}^{\infty}
$$

\vspace{2mm}
\noindent 
is an arbitrary complete orthonormal
system of functions in the space $L_2([t,T]).$
Then$,$ for the iterated Ito stochastic integral

\begin{equation}
\label{fifi1}
{\tilde J}[\psi^{(k)}]_{T,t}=\int\limits_t^T\psi_k(t_k)\sqrt{r(t_k)} 
\ldots \int\limits_t^{t_{2}}
\psi_1(t_1)\sqrt{r(t_1)} d{\bf w}_{t_1}^{(i_1)}\ldots
d{\bf w}_{t_k}^{(i_k)}
\end{equation}

\vspace{4mm}
\noindent
the following expansion 

$$
{\tilde J}[\psi^{(k)}]_{T,t}=
\hbox{\vtop{\offinterlineskip\halign{
\hfil#\hfil\cr
{\rm l.i.m.}\cr
$\stackrel{}{{}_{p_1,\ldots,p_k\to \infty}}$\cr
}} }
\sum\limits_{j_1=0}^{p_1}\ldots
\sum\limits_{j_k=0}^{p_k}
\tilde C_{j_k\ldots j_1}\Biggl(
\prod_{l=1}^k {\tilde \zeta}_{j_l}^{(i_l)}+\sum\limits_{r=1}^{[k/2]}
(-1)^r \times
\Biggr.
$$

\vspace{3mm}
\begin{equation}
\label{fifi2}
\times
\sum_{\stackrel{(\{\{g_1, g_2\}, \ldots, 
\{g_{2r-1}, g_{2r}\}\}, \{q_1, \ldots, q_{k-2r}\})}
{{}_{\{g_1, g_2, \ldots, 
g_{2r-1}, g_{2r}, q_1, \ldots, q_{k-2r}\}=\{1, 2, \ldots, k\}}}}
\prod\limits_{s=1}^r
{\bf 1}_{\{i_{g_{{}_{2s-1}}}=~i_{g_{{}_{2s}}}\ne 0\}}
\Biggl.{\bf 1}_{\{j_{g_{{}_{2s-1}}}=~j_{g_{{}_{2s}}}\}}
\prod_{l=1}^{k-2r} {\tilde \zeta}_{j_{q_l}}^{(i_{q_l})}\Biggr)
\end{equation}

\vspace{6mm}
\noindent
that converges in the mean-square
sense   
is valid, where 
$i_1,\ldots,i_k=0,1,\ldots,m,$ 
$$
{\tilde \zeta}_{j}^{(i)}=
\int\limits_t^T \Psi_{j}(s)\sqrt{r(s)}d{\bf w}_s^{(i)}
$$

\vspace{2mm}
\noindent
are independent standard Gaussian random variables
for various
$i$ or $j$ {\rm(}in the case when $i\ne 0${\rm),}

$$
{\tilde C}_{j_k\ldots j_1}=\int\limits_{[t,T]^k}
K(t_1,\ldots,t_k)
\prod_{l=1}^{k}\biggl(\Psi_{j_l}(t_l)r(t_l)\biggr)dt_1\ldots dt_k
$$

\vspace{4mm}
\noindent
is the Fourier coefficient$,$
$K(t_1,\ldots,t_k)$ is defined by {\rm (\ref{ppp});}
another notations are the same as in Theorems {\rm 1, 2, 4.}
}

\vspace{2mm}

{\bf Theorem 19}\ \cite{20a}, \cite{2023xxx1}.\
{\it Under the conditions of Theorem~{\rm 18}
the following estimate

\vspace{1mm}
$$
{\sf M}\left\{\left(
{\tilde J}[\psi^{(k)}]_{T,t}-{\tilde J}[\psi^{(k)}]_{T,t}^{p_1,\ldots,p_k}
\right)^2\right\}
\le 
$$

\vspace{3mm}
$$
~ \le k!\left(~\int\limits_{[t,T]^k}
K^2(t_1,\ldots,t_k)\left(\prod_{l=1}^k r(t_l)\right)
dt_1\ldots dt_k -\sum_{j_1=0}^{p_1}\ldots
\sum_{j_k=0}^{p_k}{\tilde C}^2_{j_k\ldots j_1}\right)
$$

\vspace{6mm}
\noindent
is valid for the following cases{\rm :}

\vspace{2mm}

{\rm 1.}\ $i_1,\ldots,i_k=1,\ldots,m$\ \ and\ \ $0<T-t<\infty,$

\vspace{1mm}

{\rm 2.}\ $i_1,\ldots,i_k=0, 1,\ldots,m,$\ \ $i_1^2+\ldots+i_k^2>0,$\ \
and\ \ $0<T-t<1,$

\vspace{3mm}
\noindent
where ${\tilde J}[\psi^{(k)}]_{T,t}$ is the 
stochastic integral {\rm (\ref{fifi1}),}
${\tilde J}[\psi^{(k)}]_{T,t}^{p_1,\ldots,p_k}$ is the 
expression on the right-hand side of {\rm (\ref{fifi2})} before
passing to the limit 
$\hbox{\vtop{\offinterlineskip\halign{
\hfil#\hfil\cr
{\rm l.i.m.}\cr
$\stackrel{}{{}_{p_1,\ldots,p_k\to \infty}}$\cr
}} };$ another 
notations are the same as in Theorems {\rm 4, 5, 18}.
}

\vspace{5mm}

\section{Proof of Theorems~12 and 13 Based on the Ito Formula and 
Without Explicit Use of the Multiple Wiener Stochastic Integral}

\vspace{5mm}

Note that Theorems~12 and 13 can also be proved without explicit use of
the multiple Wiener stochastic integral.
To do this, we introduce the following sum
of iterated Ito stochastic integrals

\vspace{1mm}
\begin{equation}
\label{chain10100}
J''[\Phi]_{T,t}^{(i_1\ldots i_k)}\stackrel{\sf def}{=}\sum_{(t_1,\ldots,t_k)}
\int\limits_{t}^{T}
\ldots
\int\limits_{t}^{t_2}
\Phi(t_1,\ldots,t_k)d{\bf w}_{t_1}^{(i_1)}
\ldots
d{\bf w}_{t_k}^{(i_k)},
\end{equation}

\vspace{3mm}
\noindent
where $\Phi(t_1,\ldots,t_k)\in L_2([t, T]^k),$ $i_1,\ldots,i_k=0,1,\ldots,m,$\ $d{\bf w}_{\tau}^{(0)}
\stackrel{\sf def}{=}d\tau;$
another notations are the same as in (\ref{Wi110}).

Further, using the isometry property of the Ito stochastic integral as well as
the linearity property of this integral, we have

\vspace{1mm}
$$
J[\psi^{(k)}]_{T,t}^{(i_1\ldots i_k)}=J''[K]_{T,t}^{(i_1\ldots i_k)}=
$$

\vspace{2mm}
\begin{equation}
\label{chain10200}
=\sum_{j_1=0}^{p_1}\ldots
\sum_{j_k=0}^{p_k}
C_{j_k\ldots j_1}
J''[\phi_{j_1}\ldots \phi_{j_k}]_{T,t}^{(i_1\ldots i_k)}+
J''[R_{p_1\ldots p_k}]_{T,t}^{(i_1\ldots i_k)}\ \ \ \hbox{w.~p.~1,}
\end{equation}

\vspace{4mm}
\noindent
where $K(t_1,\ldots,t_k)$ and $R_{p_1\ldots p_k}(t_1,\ldots,t_k)$
are defined by (\ref{chain200}) and (\ref{chain30001})
correspondingly. Moreover, 
$J''[\phi_{j_1}\ldots \phi_{j_k}]_{T,t}^{(i_1\ldots i_k)}$
and $J''[R_{p_1\ldots p_k}]_{T,t}^{(i_1\ldots i_k)}$
are defined by (\ref{chain10100}). Obviously, we can consider
an analogue of (\ref{chain10200}) for $\Phi(t_1,\ldots,t_k)$
instead of $K(t_1,\ldots,t_k)$.

Passing to the limit $\hbox{\vtop{\offinterlineskip\halign{
\hfil#\hfil\cr
{\rm l.i.m.}\cr
$\stackrel{}{{}_{p_1,\ldots,p_k\to \infty}}$\cr
}} }$ in (\ref{chain10200}) and using (\ref{wi2005}), (\ref{chain7771}), (\ref{chain10100}), we obtain

\vspace{2mm}
$$
J[\psi^{(k)}]_{T,t}^{(i_1\ldots i_k)}=
\hbox{\vtop{\offinterlineskip\halign{
\hfil#\hfil\cr
{\rm l.i.m.}\cr
$\stackrel{}{{}_{p_1,\ldots,p_k\to \infty}}$\cr
}} }
\sum\limits_{j_1=0}^{p_1}\ldots
\sum\limits_{j_k=0}^{p_k}
C_{j_k\ldots j_1}J''[\phi_{j_1}\ldots \phi_{j_k}]_{T,t}^{(i_1\ldots i_k)}=
$$

\vspace{3mm}
\begin{equation}
\label{chain4002}
=\hbox{\vtop{\offinterlineskip\halign{
\hfil#\hfil\cr
{\rm l.i.m.}\cr
$\stackrel{}{{}_{p_1,\ldots,p_k\to \infty}}$\cr
}} }
\sum\limits_{j_1=0}^{p_1}\ldots
\sum\limits_{j_k=0}^{p_k}
C_{j_k\ldots j_1}
\sum\limits_{(t_1,\ldots,t_k)}
\int\limits_t^T \phi_{j_k}(t_k)
\ldots
\int\limits_t^{t_{2}}\phi_{j_{1}}(t_{1})
d{\bf w}_{t_1}^{(i_1)}\ldots
d{\bf w}_{t_k}^{(i_k)},
\end{equation}

\vspace{4mm}
\noindent
where permutations $(t_1,\ldots,t_k)$ when summing are 
performed only in the values
$d{\bf w}_{t_1}^{(i_1)}
\ldots $
$d{\bf w}_{t_k}^{(i_k)}.$ At the same time the indices near 
upper 
limits of integration in the iterated stochastic integrals are changed 
correspondently and if $t_r$ swapped with $t_q$ in the  
permutation $(t_1,\ldots,t_k)$, then $i_r$ swapped with $i_q$ in 
the permutation $(i_1,\ldots,i_k).$

It is easy to see that the equality (\ref{chain4002}) can be written as

\vspace{1mm}
$$
J[\psi^{(k)}]_{T,t}^{(i_1\ldots i_k)}=
\hbox{\vtop{\offinterlineskip\halign{
\hfil#\hfil\cr
{\rm l.i.m.}\cr
$\stackrel{}{{}_{p_1,\ldots,p_k\to \infty}}$\cr
}} }
\sum\limits_{j_1=0}^{p_1}\ldots
\sum\limits_{j_k=0}^{p_k}
C_{j_k\ldots j_1}\times
$$

\vspace{2mm}
\begin{equation}
\label{chain7878}
\times \sum\limits_{(j_1,\ldots,j_k)}
\int\limits_t^T \phi_{j_k}(t_k)
\ldots
\int\limits_t^{t_{2}}\phi_{j_{1}}(t_{1})
d{\bf w}_{t_1}^{(i_1)}\ldots
d{\bf w}_{t_k}^{(i_k)},
\end{equation}

\vspace{3mm}
\noindent
where 
$$
\sum\limits_{(j_1,\ldots,j_k)}
$$ 

\vspace{2mm}
\noindent
means the sum with respect to all
possible permutations 
$(j_1,\ldots,j_k).$ At the same time if 
$j_r$ swapped with $j_q$ in the permutation $(j_1,\ldots,j_k)$,
then $i_r$ swapped with $i_q$ in the permutation $(i_1,\ldots,i_k)$.

Further, using the Ito formula, we can prove the following equality

\vspace{1mm}
$$
\sum\limits_{(j_1,\ldots,j_k)}
\int\limits_t^T \phi_{j_k}(t_k)
\ldots
\int\limits_t^{t_{2}}\phi_{j_{1}}(t_{1})
d{\bf w}_{t_1}^{(i_1)}\ldots
d{\bf w}_{t_k}^{(i_k)}
=\prod_{l=1}^k\zeta_{j_l}^{(i_l)}
+\sum\limits_{r=1}^{[k/2]}
(-1)^r \times
$$

\begin{equation}
\label{chain405xx}
\times\sum_{\stackrel{(\{\{g_1, g_2\}, \ldots, 
\{g_{2r-1}, g_{2r}\}\}, \{q_1, \ldots, q_{k-2r}\})}
{{}_{\{g_1, g_2, \ldots, 
g_{2r-1}, g_{2r}, q_1, \ldots, q_{k-2r}\}=\{1, 2, \ldots, k\}}}}
\prod\limits_{s=1}^r
{\bf 1}_{\{i_{g_{{}_{2s-1}}}=~i_{g_{{}_{2s}}}\ne 0\}}
\Biggl.{\bf 1}_{\{j_{g_{{}_{2s-1}}}=~j_{g_{{}_{2s}}}\}}
\prod_{l=1}^{k-2r}\zeta_{j_{q_l}}^{(i_{q_l})}
\end{equation}

\vspace{4mm}
\noindent
w.~p.~1, where notations are the same as in Theorem 2 and (\ref{chain7878}).

The main difficulty in proving (\ref{chain405xx}) 
using the Ito formula is related to the need to 
take into account various combinations of indices $i_1,\ldots,i_k=0,1,\ldots,m.$
To avoid this difficulty, consider another approach, also based
on the Ito formula.

First, we prove the following modification and generalization 
of Theorem~3.1 from \cite{ito1951} (1951) for the case $i_1,\ldots,i_k=0, 1,\ldots,m$
using the Ito formula and without explicit use of the multiple Wiener 
stochastic integral.

\vspace{2mm}                   

{\bf Theorem~20}\ \cite{20a}.\ {\it Suppose that
the condition {\rm ($\star\star$)} is fulfilled
for the multi-index $(i_1 \ldots i_k)$ {\rm (}see Sect.~{\rm 14)} 
and the condition {\rm (\ref{ziko999})} is also 
fulfilled.
Furthermore$,$ let 
$\{\phi_j(x)\}_{j=0}^{\infty}$ is an arbitrary complete orthonormal system  
of functions in the space $L_2([t,T]).$
Then

$$
J''[\phi_{j_1}\ldots \phi_{j_k}]_{T,t}^{(i_1\ldots i_k)}=
$$

\vspace{2mm}
\begin{equation}
\label{new6000}
=\prod_{l=1}^k\left({\bf 1}_{\{m_l=0\}}+{\bf 1}_{\{m_l>0\}}\left\{
\begin{matrix}
H_{n_{1,l}}\left(\zeta_{j_{h_{1,l}}}^{(i_l)}\right)\ldots 
H_{n_{d_l,l}}\left(\zeta_{j_{h_{d_l,l}}}^{(i_l)}\right),\ 
&\hbox{\rm if}\ \ \ 
i_l\ne 0\cr\cr
\left(\zeta_{j_{h_{1,l}}}^{(0)}\right)^{n_{1,l}}\ldots
\left(\zeta_{j_{h_{d_l,l}}}^{(0)}\right)^{n_{d_l,l}},\  &\hbox{\rm if}\ \ \ 
i_l=0
\end{matrix}\right.\ \right)
\end{equation}

\vspace{5mm}
\noindent
w.~p.~{\rm 1,} where $i_1,\ldots,i_k=0, 1,\ldots,m;$\ \
$n_{1,l}+n_{2,l}+\ldots+n_{d_l,l}=m_l;$\ \ $n_{1,l}, n_{2,l}, \ldots, n_{d_l,l}=1,\ldots, m_l;$\ \ 
$d_l=1,\ldots,m_l;$\ \ $l=1,\ldots,k;$\ \ $m_1+\ldots+m_k=k;$\ \  
the numbers $m_1,\ldots,m_k,$\ $g_1,\ldots,g_k$
depend on $(i_1,\ldots,i_k)$ and 
the numbers $n_{1,l},\ldots,n_{d_l,l},$\ $h_{1,l},\ldots,h_{d_l,l},$\ $d_l$
depend on $\{j_1,\ldots,j_k\};$ moreover$,$ $\left\{j_{g_1},\ldots,j_{g_k}\right\}
=\{j_1,\ldots,j_k\};$
$H_n(x)$ is the Hermite polynomial {\rm (\ref{ziko500});}
another
notations are the same as in Theorem {\rm 10}.}

\vspace{2mm}

{\bf Proof.}\ First, 
consider the case $i_1=\ldots=i_k=1,\ldots, m$ and $j_1,\ldots,j_k\in \{0\}\cup \mathbb{N}$.
By induction, we prove the following equality

\vspace{1mm}
$$
p! \int\limits_t^T \phi_l(t_p)\ldots \int\limits_t^{t_2}
\phi_l(t_1)d{\bf w}_{t_1}^{(1)}\ldots d{\bf w}_{t_p}^{(1)}\times
$$

\vspace{2mm}
$$
\times \sum\limits_{(j_1,\ldots,j_q)}
\int\limits_t^T \phi_{j_q}(t_q)\ldots \int\limits_t^{t_2}
\phi_{j_1}(t_1)d{\bf w}_{t_1}^{(1)}\ldots d{\bf w}_{t_q}^{(1)}=
$$

\vspace{2mm}
$$
=\sum\limits_{(j_1,\ldots,j_q, \underbrace{{}_{l, \ldots ,l}}_{p})}
\int\limits_t^T \phi_{j_q}(t_q)\ldots \int\limits_t^{t_2}
\phi_{j_1}(t_1)
\int\limits_t^{t_1} \phi_{l}(t_p')\ldots \int\limits_t^{t_2'}
\phi_{l}(t_1')\times
$$

\begin{equation}
\label{new1010}
\times
d{\bf w}_{t_1'}^{(1)}\ldots d{\bf w}_{t_p'}^{(1)}
d{\bf w}_{t_1}^{(1)}\ldots d{\bf w}_{t_q}^{(1)}
\end{equation}

\vspace{7mm}
\noindent
w.~p.~1, where $p\in\mathbb{N},$\ $l\ne j_1,\ldots,j_q,$ and
$$
\sum\limits_{(q_1,\ldots, q_n)}
$$

\vspace{2mm}
\noindent
means the sum with respect to all possible permutations
$(q_1,\ldots, q_n)$.

Consider the case $p=1.$ Using the Ito formula, we get w.~p.~1 for
$s\in[t, T]$

\vspace{1mm}
$$
\int\limits_t^s \phi_l(\tau)
d{\bf w}_{\tau}^{(1)}
\int\limits_t^s \phi_{j_q}(t_q)\ldots \int\limits_t^{t_2}
\phi_{j_1}(t_1)d{\bf w}_{t_1}^{(1)}\ldots d{\bf w}_{t_q}^{(1)}=
$$

\vspace{2mm}
$$
=\int\limits_t^s \phi_l(\tau)\phi_{j_q}(\tau)
\int\limits_t^{\tau} \phi_{j_{q-1}}(t_{q-1})\ldots \int\limits_t^{t_2}
\phi_{j_1}(t_1)d{\bf w}_{t_1}^{(1)}\ldots d{\bf w}_{t_{q-1}}^{(1)}d\tau+
$$

\vspace{2mm}
$$
+\int\limits_t^s \phi_l(\tau)\
\int\limits_t^{\tau} \phi_{j_q}(t_q)\ldots \int\limits_t^{t_2}
\phi_{j_1}(t_1)d{\bf w}_{t_1}^{(1)}\ldots d{\bf w}_{t_q}^{(1)}d{\bf w}_{\tau}^{(1)}+
$$

\vspace{2mm}
\begin{equation}
\label{new1011}
+\int\limits_t^s \phi_{j_q}(\tau)\hspace{-0.5mm}
\left(\int\limits_t^{\tau} \phi_{l}(\theta)
d{\bf w}_{\theta}^{(1)}
\int\limits_t^{\tau} \phi_{j_{q-1}}(t_{q-1})\ldots \int\limits_t^{t_2}
\phi_{j_1}(t_1)d{\bf w}_{t_1}^{(1)}\ldots d{\bf w}_{t_{q-1}}^{(1)}\right)
\hspace{-0.5mm}d{\bf w}_{\tau}^{(1)}\hspace{-0.2mm}.
\end{equation}

\vspace{5mm}

Hereinafter in this section always $s\in [t, T].$
Differentiating by the Ito formula the expression in parentheses
on the right-hand side of equality (\ref{new1011}) and combining the result 
of differentiation with (\ref{new1011}), we obtain w.~p.~1

\vspace{1mm}
$$
J_{(l)s,t} J_{(j_q\ldots j_1)s,t}=
$$

\vspace{2mm}
$$
=\int\limits_t^s \phi_l(\tau)\phi_{j_q}(\tau)
\int\limits_t^{\tau} \phi_{j_{q-1}}(t_{q-1})\ldots \int\limits_t^{t_2}
\phi_{j_1}(t_1)d{\bf w}_{t_1}^{(1)}\ldots d{\bf w}_{t_{q-1}}^{(1)}d\tau+
$$

\vspace{2mm}
$$
+
J_{(l j_q\ldots j_1)s,t}+
$$

\vspace{2mm}
$$
+
\int\limits_t^s \hspace{-0.2mm}\phi_{j_q}(\tau)
\hspace{-0.2mm}
\int\limits_t^{\tau}\hspace{-0.2mm}
\phi_{l}(\theta)\phi_{j_{q-1}}(\theta)\hspace{-0.1mm}\int\limits_t^{\theta}
\hspace{-0.2mm}\phi_{j_{q-2}}(t_{q-2})
\ldots \int\limits_t^{t_2}\hspace{-0.2mm}
\phi_{j_1}(t_1)d{\bf w}_{t_1}^{(1)}\ldots d{\bf w}_{t_{q-2}}^{(1)}d\theta
d{\bf w}_{\tau}^{(1)}+
$$

\vspace{2mm}
$$
+
J_{(j_q l j_{q-1}\ldots j_1)s,t}+
$$

\vspace{2mm}
$$
+\int\limits_t^s \phi_{j_q}(\tau)
\int\limits_t^{\tau} \phi_{j_{q-1}}(\theta)\times
$$

\vspace{2mm}
\begin{equation}
\label{new1012}
\times\left(\int\limits_t^{\theta} \phi_l(u)\
d{\bf w}_{u}^{(1)}\int\limits_t^{\theta} 
\phi_{j_{q-2}}(t_{q-2})\ldots \int\limits_t^{t_2}
\phi_{j_1}(t_1)d{\bf w}_{t_1}^{(1)}\ldots d{\bf w}_{t_{q-2}}^{(1)}
\right)
d{\bf w}_{\theta}^{(1)}d{\bf w}_{\tau}^{(1)},
\end{equation}

\vspace{5mm}
\noindent
where
$$
\int\limits_t^s \phi_{j_q}(t_q)\ldots \int\limits_t^{t_2}
\phi_{j_1}(t_1)d{\bf w}_{t_1}^{(1)}\ldots d{\bf w}_{t_q}^{(1)}
\stackrel{\sf def}{=}J_{(j_q \ldots j_1)s,t}.
$$

\vspace{4mm}

Continuing the process of iterative application of the Ito formula, we have w.~p.~1

\vspace{1mm}
$$
J_{(l)s,t} J_{(j_q\ldots j_1)s,t}=
$$

\vspace{2mm}
$$
=J_{(l j_q\ldots j_1)s,t}+ J_{(j_q l j_{q-1}\ldots j_1)s,t}+\ldots + J_{(j_q\ldots j_1 l)s,t}+
$$

\vspace{2mm}
$$
+\int\limits_t^s \phi_l(\tau)\phi_{j_q}(\tau)
\int\limits_t^{\tau} \phi_{j_{q-1}}(t_{q-1})\ldots \int\limits_t^{t_2}
\phi_{j_1}(t_1)d{\bf w}_{t_1}^{(1)}\ldots d{\bf w}_{t_{q-1}}^{(1)}d\tau+\ldots
$$

\vspace{2mm}
\begin{equation}
\label{new1040ss}
\ldots +
\int\limits_t^{s} \phi_{j_{q}}(t_{q})\ldots \int\limits_t^{t_3}
\phi_{j_2}(t_2)\int\limits_t^{t_2} \phi_{l}(\tau)\phi_{j_1}(\tau)      
d\tau d{\bf w}_{t_2}^{(1)}\ldots d{\bf w}_{t_{q}}^{(1)}.
\end{equation}

\vspace{5mm}

Summing the equality (\ref{new1040ss}) over permutations $(j_1,\ldots, j_q)$, we get

\begin{equation}
\label{new1025}
\sum\limits_{(j_1,\ldots, j_q)}J_{(l)s,t} J_{(j_q\ldots j_1)s,t}=
\sum\limits_{(j_1,\ldots, j_q,l)}J_{(l j_q\ldots j_1)s,t}+ S(s)
\end{equation}

\vspace{3mm}
\noindent
w.~p.~1, where
$$
S(s)=
$$

\vspace{2mm}
$$
=\sum\limits_{(j_1,\ldots, j_q)}\left(\int\limits_t^s \phi_l(\tau)\phi_{j_q}(\tau)
\int\limits_t^{\tau} \phi_{j_{q-1}}(t_{q-1})\ldots \int\limits_t^{t_2}
\phi_{j_1}(t_1)d{\bf w}_{t_1}^{(1)}\ldots d{\bf w}_{t_{q-1}}^{(1)}d\tau+\ldots\right.
$$

\vspace{2mm}
\begin{equation}
\label{new1040}
\left.\ldots +
\int\limits_t^{s} \phi_{j_{q}}(t_{q})\ldots \int\limits_t^{t_3}
\phi_{j_2}(t_2)\int\limits_t^{t_2} \phi_{l}(\tau)\phi_{j_1}(\tau)      
d\tau d{\bf w}_{t_2}^{(1)}\ldots d{\bf w}_{t_{q}}^{(1)}\right).
\end{equation}

\vspace{5mm}

Consider 
$$
\int\limits_t^s \phi_l(\tau)\phi_{j_q}(\tau)d\tau
\int\limits_t^{s} \phi_{j_{q-1}}(t_{q-1})\ldots \int\limits_t^{t_2}
\phi_{j_1}(t_1)d{\bf w}_{t_1}^{(1)}\ldots d{\bf w}_{t_{q-1}}^{(1)}.
$$

\vspace{4mm}

Applying the Ito formula, we get w.~p.~1

$$
\int\limits_t^s \phi_l(\tau)\phi_{j_q}(\tau)d\tau
\int\limits_t^{s} \phi_{j_{q-1}}(t_{q-1})\ldots \int\limits_t^{t_2}
\phi_{j_1}(t_1)d{\bf w}_{t_1}^{(1)}\ldots d{\bf w}_{t_{q-1}}^{(1)}=
$$

\vspace{2mm}
$$
=\int\limits_t^s \phi_l(\tau)\phi_{j_q}(\tau)
\int\limits_t^{\tau} \phi_{j_{q-1}}(t_{q-1})\ldots \int\limits_t^{t_2}
\phi_{j_1}(t_1)d{\bf w}_{t_1}^{(1)}\ldots d{\bf w}_{t_{q-1}}^{(1)}d\tau+
$$

\vspace{2mm}
$$
+\int\limits_t^s \phi_{j_{q-1}}(t_{q-1})\times
$$
$$
\times\left(\int\limits_t^{t_{q-1}}\phi_l(\tau)\phi_{j_q}(\tau)d\tau
\int\limits_t^{t_{q-1}} \phi_{j_{q-2}}(t_{q-2})\ldots \int\limits_t^{t_2}
\phi_{j_1}(t_1)d{\bf w}_{t_1}^{(1)}\ldots d{\bf w}_{t_{q-2}}^{(1)}\right)d{\bf w}_{t_{q-1}}^{(1)}.
$$

\vspace{5mm}

By iterative application of the Ito formula (as above), we obtain w.~p.~1

\vspace{1mm}
$$
\int\limits_t^s \phi_l(\tau)\phi_{j_q}(\tau)d\tau
\int\limits_t^{s} \phi_{j_{q-1}}(t_{q-1})\ldots \int\limits_t^{t_2}
\phi_{j_1}(t_1)d{\bf w}_{t_1}^{(1)}\ldots d{\bf w}_{t_{q-1}}^{(1)}=
$$

\vspace{2mm}
$$
=\int\limits_t^s \phi_l(\tau)\phi_{j_q}(\tau)
\int\limits_t^{\tau} \phi_{j_{q-1}}(t_{q-1})\ldots \int\limits_t^{t_2}
\phi_{j_1}(t_1)d{\bf w}_{t_1}^{(1)}\ldots d{\bf w}_{t_{q-1}}^{(1)}d\tau+\ldots
$$

\vspace{2mm}
\begin{equation}
\label{newx1020ss}
\ldots +
\int\limits_t^{s} \phi_{j_{q-1}}(t_{q-1})\ldots \int\limits_t^{t_2}
\phi_{j_1}(t_1)\int\limits_t^{t_1} \phi_{l}(\tau)\phi_{j_q}(\tau)      
d\tau d{\bf w}_{t_1}^{(1)}\ldots d{\bf w}_{t_{q-1}}^{(1)}.
\end{equation}

\vspace{5mm}

Summing the equality (\ref{newx1020ss}) over permutations $(j_1,\ldots, j_q)$, we get

\begin{equation}
\label{new1020}
\sum\limits_{(j_1,\ldots, j_q)}\int\limits_t^s \phi_l(\tau)\phi_{j_q}(\tau)d\tau
\int\limits_t^{s} \phi_{j_{q-1}}(t_{q-1})\ldots \int\limits_t^{t_2}
\phi_{j_1}(t_1)d{\bf w}_{t_1}^{(1)}\ldots d{\bf w}_{t_{q-1}}^{(1)}=S_1(s)
\end{equation}

\vspace{4mm}
\noindent
w.~p.~1, where
$$
S_1(s)=
$$

\vspace{2mm}
$$
=\sum\limits_{(j_1,\ldots, j_q)}\left(\int\limits_t^s \phi_l(\tau)\phi_{j_q}(\tau)
\int\limits_t^{\tau} \phi_{j_{q-1}}(t_{q-1})\ldots \int\limits_t^{t_2}
\phi_{j_1}(t_1)d{\bf w}_{t_1}^{(1)}\ldots d{\bf w}_{t_{q-1}}^{(1)}d\tau+\ldots\right.
$$

\vspace{2mm}
\begin{equation}
\label{newx1020}
\left.\ldots +
\int\limits_t^{s} \phi_{j_{q-1}}(t_{q-1})\ldots \int\limits_t^{t_2}
\phi_{j_1}(t_1)\int\limits_t^{t_1} \phi_{l}(\tau)\phi_{j_q}(\tau)      
d\tau d{\bf w}_{t_1}^{(1)}\ldots d{\bf w}_{t_{q-1}}^{(1)}\right).
\end{equation}

\vspace{5mm}

It is not difficult to see that
         
\vspace{-3mm}
\begin{equation}
\label{new1021}
S(s)=S_1(s)\ \ \ \hbox{w.~p.~1.}
\end{equation}

\vspace{4mm}

Moreover, due to the orthogonality of $\{\phi_j(x)\}_{j=0}^{\infty}$
and (\ref{new1020}), (\ref{new1021}), we have

\vspace{1mm}
\begin{equation}
\label{new1026}
S(T)=S_1(T)=0\ \ \ \hbox{w.~p.~1.}
\end{equation}

\vspace{4mm}

Thus (see (\ref{new1025}), (\ref{new1026})), the equality (\ref{new1010}) is proved for the case
$p=1.$
Let us assume that the equality (\ref{new1010}) is true for $p=2, 3, \ldots, k-1$, and prove
its validity for $p=k.$

From (\ref{new1025}) for the case $q=k-1,$ $j_1=\ldots=j_{k-1}=l$ we obtain

\begin{equation}
\label{new1042}
\left(J_1\right)_{s,t} (k-1)! \left(J_{k-1}\right)_{s,t}=k! \left(J_{k}\right)_{s,t} + S_2(s)
\end{equation}

\vspace{2mm}
\noindent
w.~p.~1, where 

\vspace{2mm}
$$
S_2(s)=S(s)\biggl|_{j_1=\ldots=j_q=l,\ q=k-1}\biggr.\ \ (k\ge 2)\ \ \ \hbox{and}\ \ \
S_2(s)\stackrel{\sf def}{=}0\ \ (q=k-1,\ k=1),
$$

\vspace{5mm}
$$
\int\limits_t^s \phi_{l}(t_r)\ldots \int\limits_t^{t_2}
\phi_{l}(t_1)d{\bf w}_{t_1}^{(1)}\ldots d{\bf w}_{t_r}^{(1)}
\stackrel{\sf def}{=}\left(J_{r}\right)_{s,t}\ \ (r\in \mathbb{N})\ \ \ \hbox{and}\ \ \
\left(J_{0}\right)_{s,t}\stackrel{\sf def}{=}1.
$$

\vspace{5mm}

Taking into account (\ref{new1040}), (\ref{new1020})--(\ref{new1021})
and the orthonormality of $\{\phi_j(x)\}_{j=0}^{\infty}$, we have

\vspace{0.5mm}
\begin{equation}
\label{new1043}
S_2(T)=(k-1)!\left(J_{k-2}\right)_{T,t}.
\end{equation}

\vspace{4mm}

Combining (\ref{new1042}) and (\ref{new1043}), we obtain the following recurrence 
relation

\begin{equation}
\label{new1044}
~k! \left(J_{k}\right)_{T,t}=\left(J_1\right)_{T,t} (k-1)! \left(J_{k-1}\right)_{T,t}-
(k-1)!\left(J_{k-2}\right)_{T,t}
\end{equation}

\vspace{4mm}
\noindent
w.~p.~1.

Using (\ref{new1044}) and the induction hypothesis, we get w.~p.~1

$$
k! \int\limits_t^T \phi_l(t_k)\ldots \int\limits_t^{t_2}
\phi_l(t_1)d{\bf w}_{t_1}^{(1)}\ldots d{\bf w}_{t_k}^{(1)}\times
$$

\vspace{2mm}
$$
\times \sum\limits_{(j_1,\ldots,j_q)}
\int\limits_t^T \phi_{j_q}(t_q)\ldots \int\limits_t^{t_2}
\phi_{j_1}(t_1)d{\bf w}_{t_1}^{(1)}\ldots d{\bf w}_{t_q}^{(1)}=
$$

\vspace{2mm}
$$
=\int\limits_t^T \phi_l(\tau)\
d{\bf w}_{\tau}^{(1)}\Biggl(
(k-1)!\int\limits_t^T \phi_l(t_{k-1})\ldots \int\limits_t^{t_2}
\phi_l(t_1)d{\bf w}_{t_1}^{(1)}\ldots d{\bf w}_{t_{k-1}}^{(1)}\times\Biggr.
$$

\vspace{2mm}
$$
\Biggl.\times \sum\limits_{(j_1,\ldots,j_q)}
\int\limits_t^T \phi_{j_q}(t_q)\ldots \int\limits_t^{t_2}
\phi_{j_1}(t_1)d{\bf w}_{t_1}^{(1)}\ldots d{\bf w}_{t_q}^{(1)}\Biggr)-
$$

\vspace{2mm}
$$
-(k-1)!\int\limits_t^T \phi_l(t_{k-2})\ldots \int\limits_t^{t_2}
\phi_l(t_1)d{\bf w}_{t_1}^{(1)}\ldots d{\bf w}_{t_{k-2}}^{(1)}\times
$$

\vspace{2mm}
$$
\times \sum\limits_{(j_1,\ldots,j_q)}
\int\limits_t^T \phi_{j_q}(t_q)\ldots \int\limits_t^{t_2}
\phi_{j_1}(t_1)d{\bf w}_{t_1}^{(1)}\ldots d{\bf w}_{t_q}^{(1)}=
$$

\vspace{2mm}
$$
=\int\limits_t^T \phi_l(\tau)\
d{\bf w}_{\tau}^{(1)}
\sum\limits_{(j_1,\ldots,j_q, \underbrace{{}_{l, \ldots ,l}}_{k-1})}
\int\limits_t^T \phi_{j_q}(t_q)\ldots \int\limits_t^{t_2}
\phi_{j_1}(t_1)
\int\limits_t^{t_1} \phi_{l}(t_{k-1}')\ldots \int\limits_t^{t_2'}
\phi_{l}(t_1')\times
$$

\vspace{2mm}
$$
\times
d{\bf w}_{t_1'}^{(1)}\ldots d{\bf w}_{t_{k-1}'}^{(1)}
d{\bf w}_{t_1}^{(1)}\ldots d{\bf w}_{t_{q}}^{(1)}-
$$

\vspace{2mm}
$$
-(k-1)\sum\limits_{(j_1,\ldots,j_q, \underbrace{{}_{l, \ldots ,l}}_{k-2})}
\int\limits_t^T \phi_{j_q}(t_q)\ldots \int\limits_t^{t_2}
\phi_{j_1}(t_1)
\int\limits_t^{t_1} \phi_{l}(t_{k-2}')\ldots \int\limits_t^{t_2'}
\phi_{l}(t_1')\times
$$

\begin{equation}
\label{new1050a}
\times
d{\bf w}_{t_1'}^{(1)}\ldots d{\bf w}_{t_{k-2}'}^{(1)}
d{\bf w}_{t_1}^{(1)}\ldots d{\bf w}_{t_{q}}^{(1)}.
\end{equation}

\vspace{5mm}

Let $\fbox{\it l}$ be the symbol $l$ which does not participate
in the following sum with respect to permutations

$$
\sum\limits_{(j_1,\ldots,j_q, \underbrace{{}_{l, \ldots ,l}}_{k-1})}.
$$

\vspace{2mm}

Applying (\ref{new1025}), we have w.~p.~1

$$
\int\limits_t^s \phi_{l}(\tau)\
d{\bf w}_{\tau}^{(1)}
\sum\limits_{(j_1,\ldots,j_q, \underbrace{{}_{l, \ldots ,l}}_{k-1})}
\int\limits_t^s \phi_{j_q}(t_q)\ldots \int\limits_t^{t_2}
\phi_{j_1}(t_1)
\int\limits_t^{t_1} \phi_{l}(t_{k-1}')\ldots \int\limits_t^{t_2'}
\phi_{l}(t_1')\times
$$

\vspace{2mm}
$$
\times
d{\bf w}_{t_1'}^{(1)}\ldots d{\bf w}_{t_{k-1}'}^{(1)}
d{\bf w}_{t_1}^{(1)}\ldots d{\bf w}_{t_{q}}^{(1)}=
$$

\vspace{2mm}
$$
=
\sum\limits_{(j_1,\ldots,j_q, \underbrace{{}_{l, \ldots ,l}}_{k-1})}
\int\limits_t^s \phi_{\small{\fbox{\it l}}}(\tau)\
d{\bf w}_{\tau}^{(1)}\int\limits_t^s \phi_{j_q}(t_q)\ldots \int\limits_t^{t_2}
\phi_{j_1}(t_1)
\int\limits_t^{t_1} \phi_{l}(t_{k-1}')\ldots \int\limits_t^{t_2'}
\phi_{l}(t_1')\times
$$

$$
\times
d{\bf w}_{t_1'}^{(1)}\ldots d{\bf w}_{t_{k-1}'}^{(1)}
d{\bf w}_{t_1}^{(1)}\ldots d{\bf w}_{t_{q}}^{(1)}=
$$

\vspace{2mm}
$$
=\sum\limits_{(j_1,\ldots,j_q, \underbrace{{}_{l, \ldots ,l}}_{k-1})}
\left(J_{(\small{\fbox{\it l}}j_q\ldots j_1 \underbrace{l \ldots l}_{k-1})s,t}+
J_{(\small{j_q\fbox{\it l}}j_{q-1}\ldots j_1 \underbrace{l \ldots l}_{k-1})s,t}+\ldots\right.
$$

\vspace{2mm}
$$
\left.\ldots +J_{(j_q\ldots j_1 \small{\fbox{\it l}}\underbrace{l \ldots l}_{k-1})s,t}+
J_{(j_q\ldots j_1 l\small{\fbox{\it l}}\underbrace{l \ldots l}_{k-2})s,t}+
\ldots + J_{(j_q\ldots j_1 \small{\underbrace{l \ldots l}_{k-1}}\small{\fbox{\it l}})s,t}\right)+
S_3(s)=
$$

\vspace{5mm}
\begin{equation}
\label{new1050b}
=\sum\limits_{(j_1,\ldots,j_q, \underbrace{{}_{l, \ldots ,l}}_{k})}
J_{(j_q\ldots j_1 \small{\underbrace{l \ldots l}_{k}})s,t}
+S_3(s),
\end{equation}

\vspace{3mm}
\noindent
where

\vspace{-2mm}
$$
S_3(s)=
$$

\vspace{1mm}
$$
=
\sum\limits_{(j_1,\ldots,j_q, \underbrace{{}_{l, \ldots ,l}}_{k-1})}
\Biggl(\int\limits_t^s \phi_{\small{\fbox{\it l}}}(\tau)\phi_{j_q}(\tau)
\int\limits_t^{\tau} \phi_{j_{q-1}}(t_{q-1})\ldots \int\limits_t^{t_2}
\phi_{j_1}(t_1)
\times\Biggr.
$$

$$
\times\int\limits_t^{t_1} \phi_{l}(t_{k-1}')\ldots \int\limits_t^{t_2'}
\phi_{l}(t_1')
d{\bf w}_{t_1'}^{(1)}\ldots d{\bf w}_{t_{k-1}'}^{(1)}
d{\bf w}_{t_1}^{(1)}\ldots d{\bf w}_{t_{q-1}}^{(1)}d\tau +  \ldots
$$

\vspace{3mm}
$$
+\ldots 
\int\limits_t^{s} \phi_{j_{q}}(t_{q})\ldots \int\limits_t^{t_3}
\phi_{j_2}(t_2)\int\limits_t^{t_2} \phi_{\small{\fbox{\it l}}}(\tau)\phi_{j_1}(\tau)
\times\Biggr.
$$

\vspace{2mm}
$$
\times\int\limits_t^{\tau} \phi_{l}(t_{k-1}')\ldots \int\limits_t^{t_2'}
\phi_{l}(t_1')
d{\bf w}_{t_1'}^{(1)}\ldots d{\bf w}_{t_{k-1}'}^{(1)}
d\tau d{\bf w}_{t_2}^{(1)}\ldots d{\bf w}_{t_{q}}^{(1)}+
$$

\vspace{2mm}
$$
+
\int\limits_t^{s} \phi_{j_{q}}(t_{q})\ldots \int\limits_t^{t_2}
\phi_{j_1}(t_1)\int\limits_t^{t_1} \phi_{\small{\fbox{\it l}}}(\tau)\phi_{l}(\tau)
\times\Biggr.
$$

\vspace{2mm}
$$
\times\int\limits_t^{\tau} \phi_{l}(t_{k-2}')\ldots \int\limits_t^{t_2'}
\phi_{l}(t_1')
d{\bf w}_{t_1'}^{(1)}\ldots d{\bf w}_{t_{k-2}'}^{(1)}d\tau
d{\bf w}_{t_1}^{(1)}\ldots d{\bf w}_{t_{q}}^{(1)}+ \ldots
$$

\vspace{2mm}
$$
\ldots +
\int\limits_t^{s} \phi_{j_{q}}(t_{q})\ldots \int\limits_t^{t_2}
\phi_{j_1}(t_1)
\times\Biggr.
$$

\vspace{2mm}
$$
\Biggl.\times\int\limits_t^{t_1} \phi_{l}(t_{k-1}')\ldots \int\limits_t^{t_3'}
\phi_{l}(t_2')\int\limits_t^{t_2'} \phi_{\small{\fbox{\it l}}}(\tau)\phi_{l}(\tau)
d\tau d{\bf w}_{t_2'}^{(1)}\ldots d{\bf w}_{t_{k-1}'}^{(1)}
d{\bf w}_{t_1}^{(1)}\ldots d{\bf w}_{t_{q}}^{(1)}\Biggr).
$$

\vspace{6mm}

Using (\ref{new1040}), (\ref{new1020})--(\ref{new1021}), we get w.~p.~1

$$
S_3(s)=
$$

\vspace{2mm}
$$
=\sum\limits_{(j_1,\ldots,j_q, \underbrace{{}_{l, \ldots ,l}}_{k-1})}
\int\limits_t^{s} \phi_{\small{\fbox{\it l}}}(\tau)\phi_{l}(\tau)d\tau
\int\limits_t^{s} \phi_{j_{q}}(t_{q})\ldots \int\limits_t^{t_2}
\phi_{j_1}(t_1)\times
$$

\vspace{2mm}
$$
\times
\int\limits_t^{t_1}\phi_{l}(t_{k-2}')\ldots \int\limits_t^{t_2'}
\phi_{l}(t_1')
d{\bf w}_{t_1'}^{(1)}\ldots d{\bf w}_{t_{k-2}'}^{(1)}
d{\bf w}_{t_1}^{(1)}\ldots d{\bf w}_{t_{q}}^{(1)} =
$$

\vspace{2mm}
$$
=(k-1)\sum\limits_{(j_1,\ldots,j_q, \underbrace{{}_{l, \ldots ,l}}_{k-2})}
\int\limits_t^{s} \phi_{\small{\fbox{\it l}}}(\tau)\phi_{l}(\tau)d\tau
\int\limits_t^{s} \phi_{j_{q}}(t_{q})\ldots \int\limits_t^{t_2}
\phi_{j_1}(t_1)\times
$$

\vspace{2mm}
$$
\times
\int\limits_t^{t_1}\phi_{l}(t_{k-2}')\ldots \int\limits_t^{t_2'}
\phi_{l}(t_1')
d{\bf w}_{t_1'}^{(1)}\ldots d{\bf w}_{t_{k-2}'}^{(1)}
d{\bf w}_{t_1}^{(1)}\ldots d{\bf w}_{t_{q}}^{(1)} +
$$

\vspace{2mm}
$$
+\sum\limits_{(j_1,\ldots,j_{q-1}, \underbrace{{}_{l, \ldots ,l}}_{k-1})}
\int\limits_t^{s} \phi_{\small{\fbox{\it l}}}(\tau)\phi_{j_q}(\tau)d\tau
\int\limits_t^{s} \phi_{j_{q-1}}(t_{q-1})\ldots \int\limits_t^{t_2}
\phi_{j_1}(t_1)\times
$$

\vspace{2mm}
$$
\times
\int\limits_t^{t_1}\phi_{l}(t_{k-1}')\ldots \int\limits_t^{t_2'}
\phi_{l}(t_1')
d{\bf w}_{t_1'}^{(1)}\ldots d{\bf w}_{t_{k-1}'}^{(1)}
d{\bf w}_{t_1}^{(1)}\ldots d{\bf w}_{t_{q-1}}^{(1)} + 
$$

\vspace{2mm}
$$
+\sum\limits_{(j_1,\ldots,j_{q-2}, j_q \underbrace{{}_{l, \ldots ,l}}_{k-1})}
\int\limits_t^{s} \phi_{\small{\fbox{\it l}}}(\tau)\phi_{j_{q-1}}(\tau)d\tau
\int\limits_t^{s} \phi_{j_{q}}(t_{q}) \int\limits_t^{t_q} \phi_{j_{q-2}}(t_{q-2})\ldots \int\limits_t^{t_2}
\phi_{j_1}(t_1)\times
$$

\vspace{2mm}
$$
\times
\int\limits_t^{t_1}\phi_{l}(t_{k-1}')\ldots \int\limits_t^{t_2'}
\phi_{l}(t_1')
d{\bf w}_{t_1'}^{(1)}\ldots d{\bf w}_{t_{k-1}'}^{(1)}
d{\bf w}_{t_1}^{(1)}\ldots d{\bf w}_{t_{q-2}}^{(1)}d{\bf w}_{t_{q}}^{(1)} + \ldots
$$

\vspace{2mm}
$$
\ldots
$$

\vspace{2mm}
$$
\ldots +\sum\limits_{(j_2,\ldots,j_q \underbrace{{}_{l, \ldots ,l}}_{k-1})}
\int\limits_t^{s} \phi_{\small{\fbox{\it l}}}(\tau)\phi_{j_{1}}(\tau)d\tau
\int\limits_t^{s} \phi_{j_{q}}(t_{q}) \ldots \int\limits_t^{t_3}
\phi_{j_2}(t_2)\times
$$

\vspace{2mm}
\begin{equation}
\label{new1047}
\times
\int\limits_t^{t_2}\phi_{l}(t_{k-1}')\ldots \int\limits_t^{t_2'}
\phi_{l}(t_1')
d{\bf w}_{t_1'}^{(1)}\ldots d{\bf w}_{t_{k-1}'}^{(1)}
d{\bf w}_{t_2}^{(1)}\ldots d{\bf w}_{t_{q}}^{(1)}.
\end{equation}

\vspace{6mm}

Applying (\ref{new1047}) and the orthonormality of $\{\phi_j(x)\}_{j=0}^{\infty}$, we finally have

$$
S_3(T)=(k-1)\sum\limits_{(j_1,\ldots,j_q, \underbrace{{}_{l, \ldots ,l}}_{k-2})}
\int\limits_t^{T} \phi_{j_{q}}(t_{q})\ldots \int\limits_t^{t_2}
\phi_{j_1}(t_1)
\times
$$

\vspace{2mm}
\begin{equation}
\label{new1048}
\times
\int\limits_t^{t_1}\phi_{l}(t_{k-2}')\ldots \int\limits_t^{t_2'}
\phi_{l}(t_1')
d{\bf w}_{t_1'}^{(1)}\ldots d{\bf w}_{t_{k-2}'}^{(1)}
d{\bf w}_{t_1}^{(1)}\ldots d{\bf w}_{t_{q}}^{(1)}.
\end{equation}

\vspace{6mm}

Combining (\ref{new1050a}), (\ref{new1050b}), (\ref{new1048}), we obtain w.~p.~1

$$
k! \int\limits_t^T \phi_l(t_k)\ldots \int\limits_t^{t_2}
\phi_l(t_1)d{\bf w}_{t_1}^{(1)}\ldots d{\bf w}_{t_k}^{(1)}\times
$$

\vspace{2mm}
$$
\times \sum\limits_{(j_1,\ldots,j_q)}
\int\limits_t^T \phi_{j_q}(t_q)\ldots \int\limits_t^{t_2}
\phi_{j_1}(t_1)d{\bf w}_{t_1}^{(1)}\ldots d{\bf w}_{t_q}^{(1)}=
$$

\vspace{2mm}
$$
=\sum\limits_{(\underbrace{{}_{l, \ldots ,l}}_{k})}
\int\limits_t^T \phi_l(t_k)\ldots \int\limits_t^{t_2}
\phi_l(t_1)d{\bf w}_{t_1}^{(1)}\ldots d{\bf w}_{t_k}^{(1)}\times
$$

\vspace{2mm}
$$
\times \sum\limits_{(j_1,\ldots,j_q)}
\int\limits_t^T \phi_{j_q}(t_q)\ldots \int\limits_t^{t_2}
\phi_{j_1}(t_1)d{\bf w}_{t_1}^{(1)}\ldots d{\bf w}_{t_q}^{(1)}=
$$

\vspace{2mm}
$$
=\sum\limits_{(j_1,\ldots,j_q, \underbrace{{}_{l, \ldots ,l}}_{k})}
\int\limits_t^T \phi_{j_q}(t_q)\ldots \int\limits_t^{t_2}
\phi_{j_1}(t_1)
\int\limits_t^{t_1} \phi_{l}(t_k')\ldots \int\limits_t^{t_2'}
\phi_{l}(t_1')\times
$$

\vspace{2mm}
\begin{equation}
\label{new1060}
\times
d{\bf w}_{t_1'}^{(1)}\ldots d{\bf w}_{t_k'}^{(1)}
d{\bf w}_{t_1}^{(1)}\ldots d{\bf w}_{t_q}^{(1)},
\end{equation}

\vspace{5mm}
\noindent
where $l\ne j_1,\ldots, j_q.$

The equality (\ref{new1010}) is proved. From the other hand, (\ref{new1060}) means that

\begin{equation}
\label{new1061}
J''[\phi_{j_1}\ldots \phi_{j_q}\underbrace{\phi_{l}\ldots \phi_{l}}_{n}]^
{(\hspace{0.5mm}\small{\overbrace{1 \ldots 1}^{q+n}}\hspace{0.5mm})}_{T,t}=
J''[\underbrace{\phi_{l}\ldots \phi_{l}}_{n}]^
{(\hspace{0.5mm}\small{\overbrace{1 \ldots 1}^{n}}\hspace{0.5mm})}_{T,t}
\cdot J''[\phi_{j_1}\ldots \phi_{j_q}]^
{(\hspace{0.5mm}\small{\overbrace{1 \ldots 1}^{q}}\hspace{0.5mm})}_{T,t}
\end{equation}

\vspace{4mm}
\noindent
w.~p.~1, where $n, q=0,1,2\ldots;$\ $l\ne j_1,\ldots, j_q$ and

$$
J''[\phi_{j_1}\ldots \phi_{j_q}]^
{(\hspace{0.5mm}\small{\overbrace{1 \ldots 1}^{q}}\hspace{0.5mm})}_{T,t}\stackrel{\sf def}{=}1
$$

\vspace{3mm}
\noindent
for $q=0.$

Note that \cite{Ch} 

\vspace{-1mm}
$$
\int\limits_t^T \phi_l(t_n)\ldots \int\limits_t^{t_2}
\phi_l(t_1)d{\bf w}_{t_1}^{(1)}\ldots d{\bf w}_{t_n}^{(1)}=
$$

\vspace{2mm}
$$
=
\frac{1}{n!} H_n\left(\int\limits_t^T \phi_l(\tau)d{\bf w}_{\tau}^{(1)},
\int\limits_t^T \phi_l^2(\tau)d\tau\right)=
$$

\vspace{2mm}
\begin{equation}
\label{new1100}
=
\frac{1}{n!} H_n\left(\int\limits_t^T \phi_l(\tau)d{\bf w}_{\tau}^{(1)},1\right)=
\frac{1}{n!} H_n\left(\int\limits_t^T \phi_l(\tau)d{\bf w}_{\tau}^{(1)}\right)
\end{equation}

\vspace{5mm}
\noindent
w.~p.~1, where $n\in \mathbb{N},$ $H_n(x, y)$ is defined by (\ref{new1090})
(also see (\ref{ziko1000})), and 
$H_n(x)$ is the Hermite polynomial (\ref{ziko500}).

From (\ref{new1100}) we have w.~p.~1

\vspace{1mm}
$$
J''[\underbrace{\phi_{l}\ldots \phi_{l}}_{n}]^
{(\hspace{0.5mm}\small{\overbrace{1 \ldots 1}^{n}}\hspace{0.5mm})}_{T,t}
=n! \int\limits_t^T \phi_l(t_n)\ldots \int\limits_t^{t_2}
\phi_l(t_1)d{\bf w}_{t_1}^{(1)}\ldots d{\bf w}_{t_n}^{(1)}=
$$

\vspace{2mm}
\begin{equation}
\label{new1101}
=n! \frac{1}{n!} H_n\left(\int\limits_t^T \phi_l(\tau)d{\bf w}_{\tau}^{(1)}\right)=
H_n\left(\int\limits_t^T \phi_l(\tau)d{\bf w}_{\tau}^{(1)}\right),
\end{equation}

\vspace{4mm}
\noindent
where $n\in\mathbb{N}.$

Combining (\ref{new1061}) and (\ref{new1101}), we obtain

\begin{equation}
\label{new1102}
J''[\phi_{j_1}\ldots \phi_{j_q}\underbrace{\phi_{l}\ldots \phi_{l}}_{n}]^
{(\hspace{0.5mm}\small{\overbrace{1 \ldots 1}^{q+n}}\hspace{0.5mm})}_{T,t}=
H_n\left(\int\limits_t^T \phi_l(\tau)d{\bf w}_{\tau}^{(1)}\right)
\cdot J''[\phi_{j_1}\ldots \phi_{j_q}]^
{(\hspace{0.5mm}\small{\overbrace{1 \ldots 1}^{q}}\hspace{0.5mm})}_{T,t}
\end{equation}

\vspace{4mm}
\noindent
w.~p.~1, where $n, q=0,1,2\ldots;$\ $l\ne j_1,\ldots, j_q.$ 

The iterated application of the formula (\ref{new1102})
completes the proof of Theorem~20 for the case 
$i_1=\ldots=i_k=1,\ldots, m$ and $j_1,\ldots,j_k\in \{0\}\cup \mathbb{N}$.

To prove Theorem~20 for the case $i_1=\ldots=i_k=0, 1,\ldots, m$ and 
$j_1,\ldots,j_k\in \{0\}\cup \mathbb{N}$, we need to prove
the following formula in addition to the previous proof

\vspace{1mm}
$$
p! \int\limits_t^T \phi_l(t_p)\ldots \int\limits_t^{t_2}
\phi_l(t_1)dt_1\ldots dt_p
\sum\limits_{(j_1,\ldots,j_q)}
\int\limits_t^T \phi_{j_q}(t_q)\ldots \int\limits_t^{t_2}
\phi_{j_1}(t_1)dt_1\ldots dt_q=
$$

\vspace{2mm}
\begin{equation}
\label{new1200}
=\sum\limits_{(j_1,\ldots,j_q, \underbrace{{}_{l, \ldots ,l}}_{p})}
\int\limits_t^T \phi_{j_q}(t_q)\ldots \int\limits_t^{t_2}
\phi_{j_1}(t_1)
\int\limits_t^{t_1} \phi_{l}(t_p')\ldots \int\limits_t^{t_2'}
\phi_{l}(t_1')
dt_1'\ldots dt_p'
dt_1\ldots dt_q,
\end{equation}

\vspace{5mm}
\noindent
where $p\in\mathbb{N}$,
$$
\sum\limits_{(j_1,\ldots, j_{d})}
$$

\vspace{3mm}
\noindent
means the sum with respect to all possible permutations $(j_1,\ldots,j_{d}).$ 

First, consider the case $p=1.$ We have 

$$
d\left(\int\limits_t^s \phi_l(\theta)d\theta
\int\limits_t^s \phi_{j_q}(t_q)\ldots \int\limits_t^{t_2}
\phi_{j_1}(t_1)dt_1\ldots dt_q\right)=
$$

\vspace{2mm}
$$
=\phi_l(s)
\int\limits_t^s \phi_{j_q}(t_q)\ldots \int\limits_t^{t_2}
\phi_{j_1}(t_1)dt_1\ldots dt_q ds+
$$

\vspace{2mm}
$$
+\phi_{j_q}(s)\left(\int\limits_t^{s}\phi_{j_{q-1}}(t_{q-1})\ldots \int\limits_t^{t_2}
\phi_{j_1}(t_1)dt_{1}\ldots dt_{q-1}
\cdot
\int\limits_t^s \phi_l(\theta)d\theta\right)ds.
$$

\vspace{5mm}

Then
$$
\int\limits_t^s \phi_l(\theta)d\theta
\int\limits_t^s \phi_{j_q}(t_q)\ldots \int\limits_t^{t_2}
\phi_{j_1}(t_1)dt_1\ldots dt_q=
$$

\vspace{2mm}
$$
=I_{(l j_q \ldots j_1)s,t}+
$$

\vspace{2mm}
$$
+\int\limits_t^s
\phi_{j_q}(\tau)\left(\int\limits_t^{\tau}\phi_{j_{q-1}}(t_{q-1})\ldots \int\limits_t^{t_2}
\phi_{j_1}(t_1)dt_{1}\ldots dt_{q-1}
\cdot
\int\limits_t^{\tau} \phi_l(\theta)d\theta\right)d\tau,
$$

\vspace{5mm}
\noindent
where
\begin{equation}
\label{new1701}
\int\limits_t^s \phi_{j_r}(t_r)\ldots \int\limits_t^{t_2}
\phi_{j_1}(t_1)dt_1\ldots dt_r
\stackrel{\sf def}{=}I_{(j_r \ldots j_1)s,t}.
\end{equation}

\vspace{5mm}

Continuing this process, we get

\begin{equation}
\label{new1301}
\int\limits_t^s \phi_l(\theta)d\theta\sum\limits_{(j_1,\ldots,j_q)}
I_{(j_q \ldots j_1)s,t}=
\sum\limits_{(j_1,\ldots,j_q, l)}
I_{(l j_q \ldots j_1)s,t},
\end{equation}

\vspace{4mm}
\noindent
where
$$
\sum\limits_{(j_1,\ldots,j_d)}
$$

\vspace{2mm}
\noindent
means the sum with respect to all possible permutations $(j_1,\ldots,j_d)$.

The equality (\ref{new1200}) is proved for the case $p=1.$
Let us assume that the equality (\ref{new1200}) is true for $p=2, 3, \ldots, k-1$, and prove
its validity for $p=k.$

From (\ref{new1301}) for $j_1=\ldots=j_q=l,$ $q=k-1$ we have

\begin{equation}
\label{new1400}
\left(I_1\right)_{s,t} (k-1)! \left(I_{k-1}\right)_{s,t}=k!\left(I_{k}\right)_{s,t},
\end{equation}

\vspace{4mm}
\noindent
where $k\in \mathbb{N}$ and

$$
\int\limits_t^s \phi_{l}(t_k)\ldots \int\limits_t^{t_2}
\phi_{l}(t_1)dt_1\ldots dt_k
\stackrel{\sf def}{=}\left(I_k\right)_{s,t},\ \ \ \left(I_0\right)_{s,t}\stackrel{\sf def}{=}1.
$$

\vspace{5mm}

Using (\ref{new1400}) and the induction hypothesis, we obtain 

\vspace{2mm}
$$
k! \left(I_k\right)_{s,t}
\sum\limits_{(j_1,\ldots,j_q)}
I_{(j_q\ldots j_1)s,t}=
\left(I_1\right)_{s,t} (k-1)! \left(I_{k-1}\right)_{s,t}
\sum\limits_{(j_1,\ldots,j_q)}I_{(j_q\ldots j_1)s,t}=
$$

\vspace{2mm}
\begin{equation}
\label{new1499}
=I_{(l)s,t} 
\sum\limits_{(j_1,\ldots,j_q, \underbrace{{}_{l, \ldots ,l}}_{k-1})}
I_{(j_q \ldots j_1 \underbrace{{}_{l, \ldots ,l}}_{k-1})s,t}=
\sum\limits_{(j_1,\ldots,j_q, \underbrace{{}_{l, \ldots ,l}}_{k-1})}
I_{(\small{\fbox{\it l}})s,t} 
I_{(j_q \ldots j_1 \underbrace{{}_{l, \ldots ,l}}_{k-1})s,t},
\end{equation}

\vspace{5mm}
\noindent
where $I_{(j_r \ldots j_1)s,t}$ is defined by (\ref{new1701})
and $\fbox{\it l}$ is the symbol $l$ which does not participate
in the following sum with respect to permutations

$$
\sum\limits_{(j_1,\ldots,j_q, \underbrace{{}_{l, \ldots ,l}}_{k-1})}.
$$

\vspace{3mm}

By analogy with (\ref{new1301}) we have 

$$
\sum\limits_{(j_1,\ldots,j_q, \underbrace{{}_{l, \ldots ,l}}_{k-1})}
I_{(\small{\fbox{\it l}})s,t}
I_{(j_q \ldots j_1 \underbrace{{}_{l, \ldots ,l}}_{k-1})s,t}=
\sum\limits_{(j_1,\ldots,j_q, \underbrace{{}_{l, \ldots ,l}}_{k-1})}
\left(I_{(\small{\fbox{\it l}}j_q\ldots j_1 \underbrace{l \ldots l}_{k-1})s,t}+
I_{(\small{j_q\fbox{\it l}}j_{q-1}\ldots j_1 \underbrace{l \ldots l}_{k-1})s,t}
+\ldots\right.
$$

\vspace{1mm}
$$
\left.\ldots +I_{(j_q\ldots j_1 \small{\fbox{\it l}}\underbrace{l \ldots l}_{k-1})s,t}
+
I_{(j_q\ldots j_1 l\small{\fbox{\it l}}\underbrace{l \ldots l}_{k-2})s,t}
+
\ldots + I_{(j_q\ldots j_1 \small{\underbrace{l \ldots l}_{k-1}}\small{\fbox{\it l}})s,t}
\right)=
$$

\vspace{3mm}
\begin{equation}
\label{new1500}
=\sum\limits_{(j_1,\ldots,j_q, \underbrace{{}_{l, \ldots ,l}}_{k})}
I_{(j_q\ldots j_1 \small{\underbrace{l \ldots l}_{k}})s,t}.
\end{equation}
                                                          
\vspace{5mm}

Substituting $s=T$ into (\ref{new1499}), (\ref{new1500})
and combining (\ref{new1499}), (\ref{new1500}), we conlude that the equality (\ref{new1200}) 
is proved for $p=k.$ The equality (\ref{new1200}) 
is proved.

Note that

\begin{equation}
\label{new1505}
n! \int\limits_t^T \phi_l(t_n)\ldots \int\limits_t^{t_2}
\phi_l(t_1)dt_1\ldots dt_n = n! \frac{1}{n!}
\left(\int\limits_t^T \phi_l(\tau)d\tau\right)^n=
\left(\int\limits_t^T \phi_l(\tau)d\tau\right)^n,
\end{equation}

\vspace{5mm}
\noindent
where $n\in\mathbb{N}.$

After substituting (\ref{new1505}) into (\ref{new1200}), we have for $p=n$

\begin{equation}
\label{new1506}
\left(\int\limits_t^T \phi_l(\tau)d\tau\right)^n
\sum\limits_{(j_1,\ldots,j_q)}
J_{(j_q\ldots j_1)T,t}
=\sum\limits_{(j_1,\ldots,j_q, \underbrace{{}_{l, \ldots ,l}}_{n})}
J_{(j_q\ldots j_1 \small{\underbrace{l \ldots l}_{n}})T,t}.
\end{equation}

\vspace{4mm}

The equality (\ref{new1506}) means that

\begin{equation}
\label{new1507b}
J''[\phi_{j_1}\ldots \phi_{j_q}\underbrace{\phi_{l}\ldots \phi_{l}}_{n}]^
{(\hspace{0.5mm}\small{\overbrace{0 \ldots 0}^{q+n}}\hspace{0.5mm})}_{T,t}=
\left(\int\limits_t^T \phi_l(\tau)d\tau\right)^n
\cdot J''[\phi_{j_1}\ldots \phi_{j_q}]^
{(\hspace{0.5mm}\small{\overbrace{0 \ldots 0}^{q}}\hspace{0.5mm})}_{T,t},
\end{equation}

\vspace{4mm}
\noindent
where $n, q=0,1,2\ldots $ 
and $J''[\phi_{j_1}\ldots \phi_{j_q}]^
{(0 \ldots 0)}_{T,t}\stackrel{\sf def}{=}1$
for $q=0.$

The relations (\ref{new1102}) and (\ref{new1507b}) prove 
Theorem~20 for the case $i_1=\ldots=i_k=0, 1,\ldots, m$ and $j_1,\ldots,j_k\in \{0\}\cup \mathbb{N}$.

\vspace{2mm}

{\bf Remark~12.}\ {\it Note that the equality
{\rm (\ref{new1200})} can be obtained in another way.
Let 

$$
D_q=\left\{(t_1,\ldots,t_q)\in [t, T]^q:\
\exists\ i\ne j\ \hbox{such that}\ t_i=t_j\right\}
$$ 

\vspace{4mm}
\noindent
be the "diagonal set" of $[t,T]^q$
$(q=2,3,\ldots)$
{\rm \cite{Kuo}}. Since the Lebesgue meashure of the set $D_q$ is equal to zero {\rm \cite{Kuo}}$,$ 
then {\rm(}see {\rm (\ref{chain10100}))}

\vspace{1mm}
\begin{equation}
\label{new9000a}
J''[\phi_{j_1}\ldots \phi_{j_q}]^
{(\hspace{0.5mm}\small{\overbrace{0 \ldots 0}^{q}}\hspace{0.5mm})}_{T,t}=
\int\limits_{[t,T]^q}\phi_{j_1}(t_1)\ldots \phi_{j_q}(t_q)dt_1\ldots dt_q.
\end{equation}

\vspace{4mm}

From {\rm (\ref{new9000a})} we have

$$
J''[\phi_{l}\ldots \phi_{l}]^
{(\hspace{0.5mm}\small{\overbrace{0 \ldots 0}^{p}}\hspace{0.5mm})}_{T,t}\cdot
J''[\phi_{j_1}\ldots \phi_{j_q}]^
{(\hspace{0.5mm}\small{\overbrace{0 \ldots 0}^{q}}\hspace{0.5mm})}_{T,t}=
$$

\vspace{3mm}
$$
=\int\limits_{[t,T]^q}\phi_{j_1}(t_1)\ldots \phi_{j_q}(t_q)dt_1\ldots dt_q
\int\limits_{[t,T]^p}\phi_{l}(t_1)\ldots \phi_{l}(t_p)dt_1\ldots dt_p=
$$

\vspace{3mm}
$$
=\int\limits_{[t,T]^{p+q}}\phi_{j_1}(t_1)\ldots \phi_{j_q}(t_q)
\phi_l(t_1')\ldots \phi_l(t_p')dt_1'\ldots dt_p'dt_1\ldots dt_q=
$$

\begin{equation}
\label{new100000}
=J''[\phi_{j_1}\ldots \phi_{j_q} \phi_{l}\ldots \phi_{l}]^
{(\hspace{0.5mm}\small{\overbrace{0 \ldots 0}^{p+q}}\hspace{0.5mm})}_{T,t}.
\end{equation}

\vspace{5mm}

It is not difficult to see that the equality {\rm (\ref{new100000})} is nothing but
the equality {\rm (\ref{new1200})} written in another form.}

\vspace{2mm}

To complete the proof of Theorem~20, we need to consider the case
$i_1,\ldots,i_k=0, 1,\ldots, m$ and $j_1,\ldots,j_k\in \{0\}\cup \mathbb{N}$.

Obviously, the proof of Theorem~20 will be completed if we prove the following equalities

$$
\sum\limits_{(j_1,\ldots,j_q)}\int\limits_t^T \phi_{j_q}(t_q)\ldots
\int\limits_t^{t_2}\phi_{j_1}(t_1)
d{\bf w}_{t_1}^{(i_1)}\ldots d{\bf w}_{t_q}^{(i_q)}\times
$$

\vspace{2mm}
$$
\times \sum\limits_{(j_1',\ldots,j_n')}
\int\limits_t^T 
\phi_{j_n'}(t_n')\ldots \int\limits_t^{t_2'}\phi_{j_1'}(t_1')d{\bf w}_{t_1'}^{(1)}\ldots 
d{\bf w}_{t_n'}^{(1)}=
$$

\vspace{2mm}
$$
=\sum\limits_{(j_1,\ldots,j_q,j_1',\ldots,j_n')}
\int\limits_t^T \phi_{j_q}(t_q)\ldots \int\limits_t^{t_2}\phi_{j_1}(t_1)
\int\limits_t^{t_1}\phi_{j_n'}(t_n')\ldots \int\limits_t^{t_2'}
\phi_{j_1'}(t_1')\times
$$

\vspace{2mm}
\begin{equation}
\label{new1600}
\times d{\bf w}_{t_1'}^{(1)}\ldots d{\bf w}_{t_n'}^{(1)}d{\bf w}_{t_1}^{(i_1)}\ldots d{\bf w}_{t_q}^{(i_q)},
\end{equation}

\vspace{4mm}
$$
\sum\limits_{(j_1,\ldots,j_q)}\int\limits_t^T \phi_{j_q}(t_q)\ldots
\int\limits_t^{t_2}\phi_{j_1}(t_1)
d{\bf w}_{t_1}^{(i_1)}\ldots d{\bf w}_{t_q}^{(i_q)}\times
$$

\vspace{2mm}
$$
\times \sum\limits_{(j_1',\ldots,j_n')}
\int\limits_t^T \phi_{j_n'}(t_n')\ldots \int\limits_t^{t_2'}\phi_{j_1'}(t_1')
d{\bf w}_{t_1'}^{(0)}\ldots d{\bf w}_{t_n'}^{(0)}=
$$

\vspace{2mm}
$$
=\sum\limits_{(j_1,\ldots,j_q,j_1',\ldots,j_n')}
\int\limits_t^T \phi_{j_q}(t_q)\ldots \int\limits_t^{t_2}\phi_{j_1}(t_1)
\int\limits_t^{t_1}\phi_{j_n'}(t_n')\ldots \int\limits_t^{t_2'}
\phi_{j_1'}(t_1')\times
$$

\vspace{3mm}
\begin{equation}
\label{new1600a}
\times d{\bf w}_{t_1'}^{(0)}\ldots d{\bf w}_{t_n'}^{(0)}d{\bf w}_{t_1}^{(i_1)}\ldots d{\bf w}_{t_q}^{(i_q)}
\end{equation}

\vspace{5mm}
\noindent
w.~p.~1, where $n, q\in \mathbb{N},$\ $d{\bf w}_{\tau}^{(0)}
\stackrel{\sf def}{=}d\tau,$\ $i_1,\ldots,i_q\ne 1$ in (\ref{new1600})
and $i_1,\ldots,i_q\ne 0$ in (\ref{new1600a}),

\vspace{-1mm}
$$
\sum\limits_{(j_1,\ldots,j_g)}
$$

\vspace{2mm}
\noindent
means the sum with respect to all possible permutations $(j_1,\ldots,j_g)$.
At the same time if $j_r$ swapped with $j_d$ in the permutation $(j_1,\ldots,j_g)$, then
$i_r$ swapped with $i_d$ in the permutation $(i_1,\ldots,i_g).$

The equalities (\ref{new1600}) and (\ref{new1600a}) mean that

\vspace{1mm}
\begin{equation}
\label{new1601}
J''[\phi_{j_1}\ldots\phi_{j_q}\phi_{j_1'}\ldots\phi_{j_n'}]_{T,t}^{(i_1\ldots i_q 1\ldots 1)}=
J''[\phi_{j_1}\ldots\phi_{j_q}]_{T,t}^{(i_1\ldots i_q)}
\cdot J''[\phi_{j_1'}\ldots\phi_{j_n'}]^{(1\ldots 1)}_{T,t},
\end{equation}

\vspace{2mm}
\begin{equation}
\label{new1601a}
J''[\phi_{j_1}\ldots\phi_{j_q}\phi_{j_1'}\ldots\phi_{j_n'}]_{T,t}^{(i_1\ldots i_q 0\ldots 0)}=
J''[\phi_{j_1}\ldots\phi_{j_q}]_{T,t}^{(i_1\ldots i_q)}
\cdot J''[\phi_{j_1'}\ldots\phi_{j_n'}]^{(0\ldots 0)}_{T,t}
\end{equation}

\vspace{5mm}
\noindent
w.~p.~1, where $i_1,\ldots,i_q\ne 1$ in (\ref{new1601}) and
$i_1,\ldots,i_q\ne 0$ in (\ref{new1601a}).

First, we prove the equality (\ref{new1600}). Consider the case $n=1.$
Using the Ito formula, we get w.~p.~1

\vspace{1mm}
$$
\int\limits_t^{s}\phi_{j_1'}(\theta)d{\bf w}_{\theta}^{(1)}
\int\limits_t^s \phi_{j_q}(t_q)\ldots
\int\limits_t^{t_2}\phi_{j_1}(t_1)
d{\bf w}_{t_1}^{(i_1)}\ldots d{\bf w}_{t_q}^{(i_q)}=
$$

\vspace{3mm}
$$
=J_{(j_1' j_q\ldots j_1)s,t}^{(1 i_q\ldots i_1)}+
$$

\vspace{2mm}
$$
+
\int\limits_t^s \phi_{j_q}(\tau)\left(\int\limits_t^{\tau}
\phi_{j_{q-1}}(t_{q-1})\ldots
\int\limits_t^{t_2}\phi_{j_1}(t_1)
d{\bf w}_{t_1}^{(i_1)}\ldots d{\bf w}_{t_{q-1}}^{(i_{q-1})}
\int\limits_t^{\tau}\phi_{j_1'}(\theta)d{\bf w}_{\theta}^{(1)}\right)
d{\bf w}_{\tau}^{(i_{q})}=
$$

\vspace{2mm}
$$
=\ldots =
$$

\vspace{1mm}
\begin{equation}
\label{new1700}
=J_{(j_1' j_q\ldots j_1)s,t}^{(1 i_q\ldots i_1)}+J_{(j_q j_1' j_{q-1}\ldots j_1)s,t}
^{(i_q 1 i_{q-1}\ldots i_1)}+\ldots +J_{(j_q\ldots j_1 j_1')s,t}^{(i_q\ldots i_1 1)},
\end{equation}

\vspace{5mm}
\noindent
where
\begin{equation}
\label{new100001}
\int\limits_t^s \phi_{j_r}(t_r)\ldots \int\limits_t^{t_2}
\phi_{j_1}(t_1)d{\bf w}_{t_1}^{(i_1)}\ldots d{\bf w}_{t_r}^{(i_r)}
\stackrel{\sf def}{=}J_{(j_r \ldots j_1)s,t}^{(i_r\ldots i_1)},
\end{equation}

\vspace{4mm}
\noindent
$i_1,\ldots,i_r=0,1,\ldots,m.$

From (\ref{new1700}) we obtain

$$
\int\limits_t^{s}\phi_{j_1'}(\theta)d{\bf w}_{\theta}^{(1)}
\sum\limits_{(j_1,\ldots,j_q)}\int\limits_t^s \phi_{j_q}(t_q)\ldots
\int\limits_t^{t_2}\phi_{j_1}(t_1)
d{\bf w}_{t_1}^{(i_1)}\ldots d{\bf w}_{t_q}^{(i_q)}=
$$

\vspace{3mm}
$$
=\sum\limits_{(j_1,\ldots,j_q)}\int\limits_t^{s}\phi_{j_1'}(\theta)d{\bf w}_{\theta}^{(1)}
\int\limits_t^s \phi_{j_q}(t_q)\ldots
\int\limits_t^{t_2}\phi_{j_1}(t_1)
d{\bf w}_{t_1}^{(i_1)}\ldots d{\bf w}_{t_q}^{(i_q)}=
$$

\vspace{5mm}
$$
=
\sum\limits_{(j_1,\ldots,j_q)}\left(
J_{(j_1' j_q\ldots j_1)s,t}^{(1 i_q\ldots i_1)}+J_{(j_q j_1' j_{q-1}\ldots j_1)s,t}
^{(i_q 1 i_{q-1}\ldots i_1)}+\ldots +J_{(j_q\ldots j_1 j_1')s,t}^{(i_q\ldots i_1 1)}\right)=
$$

\vspace{3mm}
\begin{equation}
\label{new1700a}
=
\sum\limits_{(j_1,\ldots,j_q, j_1')}
J_{(j_q\ldots j_1 j_1')s,t}^{(i_q\ldots i_1 1)}
\end{equation}

\vspace{5mm}
\noindent
w.~p.~1, where $J_{(j_r\ldots j_1)s,t}^{(i_r\ldots i_1)}$
is defined by (\ref{new100001}).
The equality (\ref{new1600}) is proved for the case $n=1.$

Let us assume that the equality (\ref{new1600}) is true for $n=2, 3, \ldots, k-1$, and prove
its validity for $n=k.$

Applying (\ref{new1025}), (\ref{new1040}), (\ref{new1020})--(\ref{new1021}), we obtain w.~p.~1

$$
\sum\limits_{(j_1',\ldots,j_k')}
\int\limits_t^s 
\phi_{j_k'}(t_k')\ldots \int\limits_t^{t_2'}\phi_{j_1'}(t_1')d{\bf w}_{t_1'}^{(1)}\ldots 
d{\bf w}_{t_k'}^{(1)}=
$$

\vspace{2mm}
$$
=\int\limits_t^{s}\phi_{j_k'}(\theta)d{\bf w}_{\theta}^{(1)}
\sum\limits_{(j_1',\ldots,j_{k-1}')}
\int\limits_t^s 
\phi_{j_{k-1}'}(t_{k-1})\ldots \int\limits_t^{t_2}\phi_{j_1'}(t_1)d{\bf w}_{t_1}^{(1)}\ldots 
d{\bf w}_{t_{k-1}}^{(1)}-
$$

\vspace{2mm}
\begin{equation}
\label{new1800}
-\sum\limits_{(j_1',\ldots,j_{k-1}')}
\int\limits_t^s \phi_{j_k'}(\theta)\phi_{j_{k-1}'}(\theta)d\theta
\int\limits_t^s 
\phi_{j_{k-2}'}(t_{k-2})\ldots \int\limits_t^{t_2}\phi_{j_1'}(t_1)d{\bf w}_{t_1}^{(1)}\ldots 
d{\bf w}_{t_{k-2}}^{(1)}.
\end{equation}

\vspace{5mm}

After substituting $s=T$ in (\ref{new1800})
and applying the orthonormality of 
$\{\phi_j(x)\}_{j=0}^{\infty}$, we get w.~p.~1

$$
\sum\limits_{(j_1',\ldots,j_k')}
\int\limits_t^T
\phi_{j_k'}(t_k')\ldots \int\limits_t^{t_2'}\phi_{j_1'}(t_1')d{\bf w}_{t_1'}^{(1)}\ldots 
d{\bf w}_{t_k'}^{(1)}=
$$

\vspace{2mm}
$$
=\int\limits_t^{T}\phi_{j_k'}(\theta)d{\bf w}_{\theta}^{(1)}
\sum\limits_{(j_1',\ldots,j_{k-1}')}
\int\limits_t^T 
\phi_{j_{k-1}'}(t_{k-1})\ldots \int\limits_t^{t_2}\phi_{j_1'}(t_1)d{\bf w}_{t_1}^{(1)}\ldots 
d{\bf w}_{t_{k-1}}^{(1)}-
$$

\vspace{2mm}
\begin{equation}
\label{new1801}
-\sum\limits_{(j_1',\ldots,j_{k-1}')}
{\bf 1}_{\{j_k'=j_{k-1}'\}}
\int\limits_t^T
\phi_{j_{k-2}'}(t_{k-2})\ldots \int\limits_t^{t_2}\phi_{j_1'}(t_1)d{\bf w}_{t_1}^{(1)}\ldots 
d{\bf w}_{t_{k-2}}^{(1)},
\end{equation}

\vspace{5mm}
\noindent
where ${\bf 1}_{A}$ is the indicator of the set $A$.

Using (\ref{new1801}) and the induction hypothesis, we obtain w.~p.~1

$$
\sum\limits_{(j_1',\ldots,j_k')}
\int\limits_t^T 
\phi_{j_k'}(t_k)\ldots \int\limits_t^{t_2}\phi_{j_1'}(t_1)d{\bf w}_{t_1}^{(1)}\ldots 
d{\bf w}_{t_k}^{(1)}\times
$$

\vspace{2mm}
$$
\times\sum\limits_{(j_1,\ldots,j_q)}\int\limits_t^T \phi_{j_q}(t_q)\ldots
\int\limits_t^{t_2}\phi_{j_1}(t_1)
d{\bf w}_{t_1}^{(i_1)}\ldots d{\bf w}_{t_q}^{(i_q)}=
$$

\vspace{2mm}
$$
=\int\limits_t^{T}\phi_{j_k'}(\theta)d{\bf w}_{\theta}^{(1)}
\sum\limits_{(j_1',\ldots,j_{k-1}')}
\int\limits_t^T 
\phi_{j_{k-1}'}(t_{k-1})\ldots \int\limits_t^{t_2}\phi_{j_1'}(t_1)d{\bf w}_{t_1}^{(1)}\ldots 
d{\bf w}_{t_{k-1}}^{(1)}\times
$$

\vspace{2mm}
$$
\times\sum\limits_{(j_1,\ldots,j_q)}\int\limits_t^T \phi_{j_q}(t_q)\ldots
\int\limits_t^{t_2}\phi_{j_1}(t_1)
d{\bf w}_{t_1}^{(i_1)}\ldots d{\bf w}_{t_q}^{(i_q)}-
$$

\vspace{2mm}
$$
-\sum\limits_{(j_1',\ldots,j_{k-1}')}
{\bf 1}_{\{j_k'=j_{k-1}'\}}
\int\limits_t^T
\phi_{j_{k-2}'}(t_{k-2})\ldots \int\limits_t^{t_2}\phi_{j_1'}(t_1)d{\bf w}_{t_1}^{(1)}\ldots 
d{\bf w}_{t_{k-2}}^{(1)}\times
$$

\vspace{2mm}
$$
\times\sum\limits_{(j_1,\ldots,j_q)}\int\limits_t^T \phi_{j_q}(t_q)\ldots
\int\limits_t^{t_2}\phi_{j_1}(t_1)
d{\bf w}_{t_1}^{(i_1)}\ldots d{\bf w}_{t_q}^{(i_q)}=
$$

\vspace{2mm}
$$
=\int\limits_t^{T}\phi_{j_k'}(\theta)d{\bf w}_{\theta}^{(1)}\times
$$
$$
\times
\sum\limits_{(j_1,\ldots,j_q,j_1',\ldots,j_{k-1}')}
\int\limits_t^T \phi_{j_q}(t_q)\ldots \int\limits_t^{t_2}\phi_{j_1}(t_1)
\int\limits_t^{t_1}\phi_{j_{k-1}'}(t_{k-1}')\ldots \int\limits_t^{t_2'}
\phi_{j_1'}(t_1')\times
$$

\vspace{4mm}
$$
\times d{\bf w}_{t_1'}^{(1)}\ldots 
d{\bf w}_{t_{k-1}'}^{(1)}d{\bf w}_{t_1}^{(i_1)}\ldots d{\bf w}_{t_q}^{(i_q)}-
$$

\vspace{2mm}
$$
-\sum\limits_{(j_1',\ldots,j_{k-1}')}
{\bf 1}_{\{j_k'=j_{k-1}'\}}
\int\limits_t^T
\phi_{j_{k-2}'}(t_{k-2})\ldots \int\limits_t^{t_2}\phi_{j_1'}(t_1)d{\bf w}_{t_1}^{(1)}\ldots 
d{\bf w}_{t_{k-2}}^{(1)}\times
$$

\vspace{2mm}
\begin{equation}
\label{newi9000}
\times\sum\limits_{(j_1,\ldots,j_q)}\int\limits_t^T \phi_{j_q}(t_q)\ldots
\int\limits_t^{t_2}\phi_{j_1}(t_1)
d{\bf w}_{t_1}^{(i_1)}\ldots d{\bf w}_{t_q}^{(i_q)}.
\end{equation}

\vspace{5mm}

Further, applying the induction hypothesis, we have w.~p.~1

$$
\sum\limits_{(j_1',\ldots,j_{k-1}')}
{\bf 1}_{\{j_k'=j_{k-1}'\}}
\int\limits_t^T
\phi_{j_{k-2}'}(t_{k-2})\ldots \int\limits_t^{t_2}\phi_{j_1'}(t_1)d{\bf w}_{t_1}^{(1)}\ldots 
d{\bf w}_{t_{k-2}}^{(1)}\times
$$

\vspace{2mm}
$$
\times\sum\limits_{(j_1,\ldots,j_q)}\int\limits_t^T \phi_{j_q}(t_q)\ldots
\int\limits_t^{t_2}\phi_{j_1}(t_1)
d{\bf w}_{t_1}^{(i_1)}\ldots d{\bf w}_{t_q}^{(i_q)}=
$$

\vspace{2mm}
$$
=\Biggl(
\sum\limits_{(j_1',\ldots,j_{k-2}')}
{\bf 1}_{\{j_k'=j_{k-1}'\}}
\int\limits_t^T
\phi_{j_{k-2}'}(t_{k-2})\ldots \int\limits_t^{t_2}\phi_{j_1'}(t_1)d{\bf w}_{t_1}^{(1)}\ldots 
d{\bf w}_{t_{k-2}}^{(1)}+
\Biggr.
$$

\vspace{2mm}
$$
+\sum\limits_{(j_1',\ldots,j_{k-3}', j_{k-1}')}
{\bf 1}_{\{j_k'=j_{k-2}'\}}
\int\limits_t^T
\phi_{j_{k-1}'}(t_{k-2})
\int\limits_t^{t_{k-2}}
\phi_{j_{k-3}'}(t_{k-3})\ldots \int\limits_t^{t_2}\phi_{j_1'}(t_1)\times
$$

\vspace{4mm}
$$
\times d{\bf w}_{t_1}^{(1)}\ldots 
d{\bf w}_{t_{k-3}}^{(1)}d{\bf w}_{t_{k-2}}^{(1)}+\ldots
$$

\vspace{2mm}
$$
\ldots+\sum\limits_{(j_2',\ldots,j_{k-1}')}
{\bf 1}_{\{j_k'=j_{1}'\}}
\int\limits_t^T
\phi_{j_{k-2}'}(t_{k-2})\ldots \int\limits_t^{t_3}\phi_{j_2'}(t_2)
\int\limits_t^{t_2}
\phi_{j_{k-1}'}(t_{1})\times
$$

\vspace{2mm}
$$
\Biggl.\times
d{\bf w}_{t_1}^{(1)}d{\bf w}_{t_2}^{(1)}\ldots 
d{\bf w}_{t_{k-2}}^{(1)}\Biggr)\times
$$

\vspace{2mm}
$$
\times\sum\limits_{(j_1,\ldots,j_q)}\int\limits_t^T \phi_{j_q}(t_q)\ldots
\int\limits_t^{t_2}\phi_{j_1}(t_1)
d{\bf w}_{t_1}^{(i_1)}\ldots d{\bf w}_{t_q}^{(i_q)}=
$$

\vspace{2mm}
$$
=\Biggl(
{\bf 1}_{\{j_k'=j_{k-1}'\}}\sum\limits_{(j_1',\ldots,j_{k-2}')}
\int\limits_t^T
\phi_{j_{k-2}'}(t_{k-2})\ldots \int\limits_t^{t_2}\phi_{j_1'}(t_1)d{\bf w}_{t_1}^{(1)}\ldots 
d{\bf w}_{t_{k-2}}^{(1)}+
\Biggr.
$$

\vspace{2mm}
$$
+{\bf 1}_{\{j_k'=j_{k-2}'\}}\sum\limits_{(j_1',\ldots,j_{k-3}', j_{k-1}')}
\int\limits_t^T
\phi_{j_{k-1}'}(t_{k-2})
\int\limits_t^{t_{k-2}}
\phi_{j_{k-3}'}(t_{k-3})\ldots \int\limits_t^{t_2}\phi_{j_1'}(t_1)\times
$$

\vspace{4mm}
$$
\times d{\bf w}_{t_1}^{(1)}\ldots 
d{\bf w}_{t_{k-3}}^{(1)}d{\bf w}_{t_{k-2}}^{(1)}+\ldots
$$

\vspace{2mm}
$$
\ldots +{\bf 1}_{\{j_k'=j_{1}'\}}\sum\limits_{(j_2',\ldots,j_{k-1}')}
\int\limits_t^T
\phi_{j_{k-2}'}(t_{k-2})\ldots \int\limits_t^{t_3}\phi_{j_2'}(t_2)
\int\limits_t^{t_2}
\phi_{j_{k-1}'}(t_{1})\times
$$

\vspace{2mm}
$$
\Biggl.\times
d{\bf w}_{t_1}^{(1)}d{\bf w}_{t_2}^{(1)}\ldots 
d{\bf w}_{t_{k-2}}^{(1)}\Biggr)\times
$$

\vspace{2mm}
$$
\times\sum\limits_{(j_1,\ldots,j_q)}\int\limits_t^T \phi_{j_q}(t_q)\ldots
\int\limits_t^{t_2}\phi_{j_1}(t_1)
d{\bf w}_{t_1}^{(i_1)}\ldots d{\bf w}_{t_q}^{(i_q)}=
$$

\vspace{2mm}
$$
={\bf 1}_{\{j_k'=j_{k-1}'\}}\sum\limits_{(j_1,\ldots,j_q,j_1',\ldots,j_{k-2}')}
\int\limits_t^T \phi_{j_q}(t_q)\ldots \int\limits_t^{t_2}\phi_{j_1}(t_1)
\int\limits_t^{t_1}\phi_{j_{k-2}'}(t_{k-2}')\ldots \int\limits_t^{t_2'}
\phi_{j_1'}(t_1')\times
$$

\vspace{4mm}
$$
\times d{\bf w}_{t_1'}^{(1)}\ldots d{\bf w}_{t_{k-2}'}^{(1)}
d{\bf w}_{t_1}^{(i_1)}\ldots d{\bf w}_{t_q}^{(i_q)}+
$$

\vspace{2mm}
$$
+{\bf 1}_{\{j_k'=j_{k-2}'\}}\sum\limits_{(j_1,\ldots,j_q,j_1',\ldots,j_{k-3}',j_{k-1}')}
\int\limits_t^T \phi_{j_q}(t_q)\ldots \int\limits_t^{t_2}\phi_{j_1}(t_1)
\int\limits_t^{t_1}\phi_{j_{k-1}'}(t_{k-2}')\times
$$

\vspace{4mm}
$$
\times\int\limits_t^{t_{k-2}'}
\phi_{j_{k-3}'}(t_{k-3}')
\ldots \int\limits_t^{t_2'}
\phi_{j_1'}(t_1')d{\bf w}_{t_1'}^{(1)}\ldots d{\bf w}_{t_{k-3}'}^{(1)}d{\bf w}_{t_{k-2}'}^{(1)}
d{\bf w}_{t_1}^{(i_1)}\ldots d{\bf w}_{t_q}^{(i_q)}+\ldots
$$

\vspace{1mm}
$$
\ldots
$$

\vspace{1mm}
$$
\ldots +{\bf 1}_{\{j_k'=j_{1}'\}}\sum\limits_{(j_1,\ldots,j_q,j_2',\ldots,j_{k-1}')}
\int\limits_t^T \phi_{j_q}(t_q)\ldots \int\limits_t^{t_2}\phi_{j_1}(t_1)\times
$$

\vspace{4mm}
$$
\times
\int\limits_t^{t_1}\phi_{j_{k-2}'}(t_{k-2}')\ldots \int\limits_t^{t_{3}'}
\phi_{j_{2}'}(t_{2}')
\int\limits_t^{t_2'}
\phi_{j_{k-1}'}(t_1')d{\bf w}_{t_1'}^{(1)}d{\bf w}_{t_2'}^{(1)}\ldots d{\bf w}_{t_{k-2}'}^{(1)}
d{\bf w}_{t_1}^{(i_1)}\ldots d{\bf w}_{t_q}^{(i_q)}\stackrel{\sf def}{=}
$$

\vspace{2mm}
\begin{equation}
\label{new1900}
\stackrel{\sf def}{=}S_4(T).
\end{equation}

\vspace{6mm}

By analogy with (\ref{newx1020ss}) we obtain w.~p.~1

\vspace{2mm}
$$
\int\limits_t^T \phi_l(\tau)\phi_{j_r}(\tau)d\tau
\int\limits_t^{T} \phi_{j_{r-1}}(t_{r-1})\ldots \int\limits_t^{t_2}
\phi_{j_1}(t_1)d{\bf w}_{t_1}^{(i_1)}\ldots d{\bf w}_{t_{r-1}}^{(i_{r-1})}=
$$

\vspace{2mm}
$$
=\int\limits_t^T \phi_l(\tau)\phi_{j_r}(\tau)
\int\limits_t^{\tau} \phi_{j_{r-1}}(t_{r-1})\ldots \int\limits_t^{t_2}
\phi_{j_1}(t_1)d{\bf w}_{t_1}^{(i_1)}\ldots d{\bf w}_{t_{r-1}}^{(i_{r-1})}d\tau+\ldots
$$

\vspace{2mm}
\begin{equation}
\label{new5001}
\ldots +
\int\limits_t^{T} \phi_{j_{r-1}}(t_{r-1})\ldots \int\limits_t^{t_2}
\phi_{j_1}(t_1)\int\limits_t^{t_1} \phi_{l}(\tau)\phi_{j_r}(\tau)      
d\tau d{\bf w}_{t_1}^{(i_1)}\ldots d{\bf w}_{t_{r-1}}^{(i_{r-1})},
\end{equation}

\vspace{5mm}
\noindent
where $i_1,\ldots,i_{r-1}=0,1,\ldots,m.$

Using iteratively the Ito formula, as well as (\ref{new5001})
and combinatorial reasoning, we obtain w.~p.~1
(see Remark~13 below for details)

$$
\int\limits_t^{T}\phi_{j_k'}(\theta)d{\bf w}_{\theta}^{(1)}\times
$$

\vspace{2mm}
$$
\times
\sum\limits_{(j_1,\ldots,j_q,j_1',\ldots,j_{k-1}')}
\int\limits_t^T \phi_{j_q}(t_q)\ldots \int\limits_t^{t_2}\phi_{j_1}(t_1)
\int\limits_t^{t_1}\phi_{j_{k-1}'}(t_{k-1}')\ldots \int\limits_t^{t_2'}
\phi_{j_1'}(t_1')\times
$$

\vspace{4mm}
$$
\times d{\bf w}_{t_1'}^{(1)}\ldots 
d{\bf w}_{t_{k-1}'}^{(1)}d{\bf w}_{t_1}^{(i_1)}\ldots d{\bf w}_{t_q}^{(i_q)}=
$$

\vspace{2mm}
$$
=\sum\limits_{(j_1,\ldots,j_q,j_1',\ldots,j_{k}')}
\int\limits_t^T \phi_{j_q}(t_q)\ldots \int\limits_t^{t_2}\phi_{j_1}(t_1)
\int\limits_t^{t_1}\phi_{j_{k}'}(t_{k}')\ldots \int\limits_t^{t_2'}
\phi_{j_1'}(t_1')\times
$$

\vspace{4mm}
$$
\times d{\bf w}_{t_1'}^{(1)}\ldots 
d{\bf w}_{t_{k}'}^{(1)}d{\bf w}_{t_1}^{(i_1)}\ldots d{\bf w}_{t_q}^{(i_q)}+
$$

\vspace{2mm}
$$
+\sum\limits_{(j_1,\ldots,j_q,j_1',\ldots,j_{k-1}')}
\Biggl(\int\limits_t^T \phi_{j_q}(t_q)\ldots \int\limits_t^{t_2}\phi_{j_1}(t_1)
\int\limits_t^{t_1}\phi_{j_{k}'}(\theta)\phi_{j_{k-1}'}(\theta)\int\limits_t^{\theta}
\phi_{j_{k-2}'}(t_{k-2}')\ldots\Biggr.
$$

\vspace{2mm}
$$
\Biggl.\ldots \int\limits_t^{t_2'}\phi_{j_1'}(t_1')
d{\bf w}_{t_1'}^{(1)}\ldots 
d{\bf w}_{t_{k-2}'}^{(1)}d{\bf w}_{\theta}^{(0)}d{\bf w}_{t_1}^{(i_1)}\ldots d{\bf w}_{t_q}^{(i_q)}+
$$

\vspace{2mm}
$$
+
\int\limits_t^T \phi_{j_q}(t_q)\ldots \int\limits_t^{t_2}
\phi_{j_1}(t_1)\int\limits_t^{t_1}\phi_{j_{k-1}'}(t_{k-1}')
\int\limits_t^{t_{k-1}'}\phi_{j_{k}'}(\theta)\phi_{j_{k-2}'}(\theta)\int\limits_t^{\theta}
\phi_{j_{k-3}'}(t_{k-3}')\ldots
$$

\vspace{2mm}
$$
\ldots
\int\limits_t^{t_2'}
\phi_{j_{1}'}(t_{1}')d{\bf w}_{t_1'}^{(1)}\ldots 
d{\bf w}_{t_{k-3}'}^{(1)}d{\bf w}_{\theta}^{(0)}d{\bf w}_{t_{k-1}'}^{(1)}
d{\bf w}_{t_1}^{(i_1)}\ldots d{\bf w}_{t_q}^{(i_q)}+\ldots
$$

\vspace{2mm}
$$
\ldots +
\int\limits_t^T \phi_{j_q}(t_q)\ldots \int\limits_t^{t_2}
\phi_{j_1}(t_1)\int\limits_t^{t_1}\phi_{j_{k-1}'}(t_{k-1}')\ldots
\int\limits_t^{t_{3}'}\phi_{j_{2}'}(t_2')\int\limits_t^{t_{2}'}
\phi_{j_{k}'}(\theta)\phi_{j_{1}'}(\theta)d{\bf w}_{\theta}^{(0)}\times
$$

\vspace{2mm}
$$
\Biggl.\times
d{\bf w}_{t_2'}^{(1)}\ldots 
d{\bf w}_{t_{k-1}'}^{(1)}
d{\bf w}_{t_1}^{(i_1)}\ldots d{\bf w}_{t_q}^{(i_q)}\Biggr)=
$$

\vspace{2mm}
$$
=\sum\limits_{(j_1,\ldots,j_q,j_1',\ldots,j_{k}')}
\int\limits_t^T \phi_{j_q}(t_q)\ldots \int\limits_t^{t_2}\phi_{j_1}(t_1)
\int\limits_t^{t_1}\phi_{j_{k}'}(t_{k}')\ldots \int\limits_t^{t_2'}
\phi_{j_1'}(t_1')\times
$$

\vspace{4mm}
$$
\times d{\bf w}_{t_1'}^{(1)}\ldots 
d{\bf w}_{t_{k}'}^{(1)}d{\bf w}_{t_1}^{(i_1)}\ldots d{\bf w}_{t_q}^{(i_q)}+
$$

\vspace{2mm}
$$
+\sum\limits_{(j_1,\ldots,j_q,j_1',\ldots,j_{k-2}')}
\Biggl\{\int\limits_t^T 
\phi_{j_{k}'}(\theta)\phi_{j_{k-1}'}(\theta)\int\limits_t^{\theta}
\phi_{j_q}(t_q)\ldots \int\limits_t^{t_2}\phi_{j_1}(t_1)
\int\limits_t^{t_1}
\phi_{j_{k-2}'}(t_{k-2}')\ldots\Biggr.
$$

\vspace{2mm}
$$
\Biggl.\ldots \int\limits_t^{t_2'}\phi_{j_1'}(t_1')
d{\bf w}_{t_1'}^{(1)}\ldots 
d{\bf w}_{t_{k-2}'}^{(1)}d{\bf w}_{t_1}^{(i_1)}\ldots d{\bf w}_{t_q}^{(i_q)}d{\bf w}_{\theta}^{(0)}+\ldots
$$

\vspace{2mm}
$$
\ldots+\int\limits_t^T 
\phi_{j_q}(t_q)\ldots \int\limits_t^{t_2}\phi_{j_1}(t_1)
\int\limits_t^{t_1}
\phi_{j_{k-2}'}(t_{k-2}')\ldots \int\limits_t^{t_2'}\phi_{j_1'}(t_1')
\int\limits_t^{t_1'} 
\phi_{j_{k}'}(\theta)\phi_{j_{k-1}'}(\theta)d{\bf w}_{\theta}^{(0)}\times
$$

\vspace{4mm}
$$
\Biggl.
\times d{\bf w}_{t_1'}^{(1)}\ldots 
d{\bf w}_{t_{k-2}'}^{(1)}d{\bf w}_{t_1}^{(i_1)}\ldots d{\bf w}_{t_q}^{(i_q)}\Biggr\}+
$$

\vspace{4mm}
$$
+\sum\limits_{(j_1,\ldots,j_q,j_1',\ldots,j_{k-3}',j_{k-1}')}
\Biggl\{\int\limits_t^T 
\phi_{j_{k}'}(\theta)\phi_{j_{k-2}'}(\theta)\int\limits_t^{\theta}
\phi_{j_q}(t_q)\ldots \int\limits_t^{t_2}\phi_{j_1}(t_1)
\int\limits_t^{t_1}
\phi_{j_{k-1}'}(t_{k-1}')\times\Biggr.
$$

\vspace{2mm}
$$
\times\int\limits_t^{t_{k-1}'}
\phi_{j_{k-3}'}(t_{k-3}')
\ldots \int\limits_t^{t_2'}\phi_{j_1'}(t_1')
d{\bf w}_{t_1'}^{(1)}\ldots 
d{\bf w}_{t_{k-3}'}^{(1)}d{\bf w}_{t_{k-1}'}^{(1)}
d{\bf w}_{t_1}^{(i_1)}\ldots d{\bf w}_{t_q}^{(i_q)}
d{\bf w}_{\theta}^{(0)}+\ldots
$$

\vspace{2mm}
$$
\ldots+\int\limits_t^T 
\phi_{j_q}(t_q)\ldots \int\limits_t^{t_2}\phi_{j_1}(t_1)
\int\limits_t^{t_1}
\phi_{j_{k-1}'}(t_{k-1}')\int\limits_t^{t_{k-1}'}
\phi_{j_{k-3}'}(t_{k-3}')\ldots
\int\limits_t^{t_2'}\phi_{j_1'}(t_1')
\times
$$

\vspace{4mm}
$$
\Biggl.\times
\int\limits_t^{t_1'} 
\phi_{j_{k}'}(\theta)\phi_{j_{k-2}'}(\theta)d{\bf w}_{\theta}^{(0)}
d{\bf w}_{t_1'}^{(1)}\ldots 
d{\bf w}_{t_{k-3}'}^{(1)}d{\bf w}_{t_{k-1}'}^{(1)}
d{\bf w}_{t_1}^{(i_1)}\ldots d{\bf w}_{t_q}^{(i_q)}\Biggr\}+\ldots
$$

\vspace{4mm}
$$
\ldots +\sum\limits_{(j_1,\ldots,j_q,j_2',\ldots,j_{k-1}')}
\Biggl\{\int\limits_t^T 
\phi_{j_{k}'}(\theta)\phi_{j_{1}'}(\theta)\int\limits_t^{\theta}
\phi_{j_q}(t_q)\ldots \int\limits_t^{t_2}\phi_{j_1}(t_1)
\int\limits_t^{t_1}
\phi_{j_{k-1}'}(t_{k-1}')\ldots\Biggr.
$$

\vspace{2mm}
$$
\ldots\int\limits_t^{t_{3}'}
\phi_{j_2'}(t_2')
d{\bf w}_{t_2'}^{(1)}\ldots 
d{\bf w}_{t_{k-1}'}^{(1)}
d{\bf w}_{t_1}^{(i_1)}\ldots d{\bf w}_{t_q}^{(i_q)}
d{\bf w}_{\theta}^{(0)}+\ldots
$$

\vspace{2mm}
$$
\ldots+\int\limits_t^T 
\phi_{j_q}(t_q)\ldots \int\limits_t^{t_2}\phi_{j_1}(t_1)
\int\limits_t^{t_1}
\phi_{j_{k-1}'}(t_{k-1}')\ldots
\int\limits_t^{t_3'}\phi_{j_2'}(t_2')\int\limits_t^{t_2'}
\phi_{j_{k}'}(\theta)\phi_{j_{1}'}(\theta)d{\bf w}_{\theta}^{(0)}\times
$$

\vspace{4mm}
$$
\Biggl.\times
d{\bf w}_{t_2'}^{(1)}\ldots 
d{\bf w}_{t_{k-1}'}^{(1)}
d{\bf w}_{t_1}^{(i_1)}\ldots d{\bf w}_{t_q}^{(i_q)}\Biggr\}=
$$

\vspace{4mm}
$$
=\sum\limits_{(j_1,\ldots,j_q,j_1',\ldots,j_{k}')}
\int\limits_t^T \phi_{j_q}(t_q)\ldots \int\limits_t^{t_2}\phi_{j_1}(t_1)
\int\limits_t^{t_1}\phi_{j_{k}'}(t_{k}')\ldots \int\limits_t^{t_2'}
\phi_{j_1'}(t_1')\times
$$

\vspace{4mm}
$$
\times d{\bf w}_{t_1'}^{(1)}\ldots 
d{\bf w}_{t_{k}'}^{(1)}d{\bf w}_{t_1}^{(i_1)}\ldots d{\bf w}_{t_q}^{(i_q)}+
$$

\vspace{4mm}
$$
+\int\limits_t^T 
\phi_{j_{k}'}(\theta)\phi_{j_{k-1}'}(\theta)d\theta\sum\limits_{(j_1,\ldots,j_q,j_1',\ldots,j_{k-2}')}
\int\limits_t^{T}
\phi_{j_q}(t_q)\ldots \int\limits_t^{t_2}\phi_{j_1}(t_1)
\int\limits_t^{t_1}
\phi_{j_{k-2}'}(t_{k-2}')\ldots
$$

\vspace{4mm}
$$
\ldots \int\limits_t^{t_2'}\phi_{j_1'}(t_1')
d{\bf w}_{t_1'}^{(1)}\ldots 
d{\bf w}_{t_{k-2}'}^{(1)}d{\bf w}_{t_1}^{(i_1)}\ldots d{\bf w}_{t_q}^{(i_q)}+
$$

\vspace{4mm}
$$
+\int\limits_t^T 
\phi_{j_{k}'}(\theta)\phi_{j_{k-2}'}(\theta)d\theta
\sum\limits_{(j_1,\ldots,j_q,j_1',\ldots,j_{k-3}',j_{k-1}')}
\int\limits_t^{T}
\phi_{j_q}(t_q)\ldots \int\limits_t^{t_2}\phi_{j_1}(t_1)
\int\limits_t^{t_1}
\phi_{j_{k-1}'}(t_{k-1}')\times\Biggr.
$$

\vspace{4mm}
$$
\times\int\limits_t^{t_{k-1}'}
\phi_{j_{k-3}'}(t_{k-3}')
\ldots \int\limits_t^{t_2'}\phi_{j_1'}(t_1')
d{\bf w}_{t_1'}^{(1)}\ldots 
d{\bf w}_{t_{k-3}'}^{(1)}d{\bf w}_{t_{k-1}'}^{(1)}
d{\bf w}_{t_1}^{(i_1)}\ldots d{\bf w}_{t_q}^{(i_q)}
+\ldots
$$

\vspace{4mm}
$$
\ldots +\int\limits_t^T 
\phi_{j_{k}'}(\theta)\phi_{j_{1}'}(\theta)d\theta\sum\limits_{(j_1,\ldots,j_q,j_2',\ldots,j_{k-1}')}
\int\limits_t^{T}
\phi_{j_q}(t_q)\ldots \int\limits_t^{t_2}\phi_{j_1}(t_1)
\int\limits_t^{t_1}
\phi_{j_{k-1}'}(t_{k-1}')\ldots
$$

\vspace{4mm}
$$
\ldots\int\limits_t^{t_{3}'}
\phi_{j_2'}(t_2')
d{\bf w}_{t_2'}^{(1)}\ldots 
d{\bf w}_{t_{k-1}'}^{(1)}
d{\bf w}_{t_1}^{(i_1)}\ldots d{\bf w}_{t_q}^{(i_q)}=
$$

\vspace{2mm}
$$
=\sum\limits_{(j_1,\ldots,j_q,j_1',\ldots,j_{k}')}
\int\limits_t^T \phi_{j_q}(t_q)\ldots \int\limits_t^{t_2}\phi_{j_1}(t_1)
\int\limits_t^{t_1}\phi_{j_{k}'}(t_{k}')\ldots \int\limits_t^{t_2'}
\phi_{j_1'}(t_1')\times
$$

\vspace{2mm}
\begin{equation}
\label{new1901}
\times d{\bf w}_{t_1'}^{(1)}\ldots 
d{\bf w}_{t_{k}'}^{(1)}d{\bf w}_{t_1}^{(i_1)}\ldots d{\bf w}_{t_q}^{(i_q)}+S_4(T).
\end{equation}

\vspace{5mm}

From (\ref{newi9000}), (\ref{new1900}), and (\ref{new1901}) we conclude that
the equality (\ref{new1600}) is proved for $n=k.$
The equality (\ref{new1600}) is proved.

\vspace{2mm} 

{\bf Remark~13.}\ {\it It should be noted that the sums with respect to 
permutations 

\vspace{1mm}
$$
\sum\limits_{(j_1,\ldots,j_q,j_1',\ldots,j_{k-1}')}
$$

\vspace{4mm}
\noindent
in {\rm (\ref{new1901}),} containing the expressions 
$\phi_{j_{k}'}(\theta)\phi_{j_{k-1}'}(\theta),\ldots,
\phi_{j_{k}'}(\theta)\phi_{j_{1}'}(\theta),$
should be understood in a special way.
Let us explain this rule on the basis of the sum

\vspace{1mm}
$$
\sum\limits_{(j_1,\ldots,j_q,j_1',\ldots,j_{k-1}')}
\int\limits_t^T \phi_{j_q}(t_q)\ldots \int\limits_t^{t_2}\phi_{j_1}(t_1)
\int\limits_t^{t_1}\phi_{j_{k}'}(\theta)\phi_{j_{k-1}'}(\theta)\int\limits_t^{\theta}
\phi_{j_{k-2}'}(t_{k-2}')\ldots
$$

\vspace{2.5mm}
\begin{equation}
\label{new777100}
\Biggl.\ldots \int\limits_t^{t_2'}\phi_{j_1'}(t_1')
d{\bf w}_{t_1'}^{(1)}\ldots 
d{\bf w}_{t_{k-2}'}^{(1)}d{\bf w}_{\theta}^{(0)}d{\bf w}_{t_1}^{(i_1)}\ldots d{\bf w}_{t_q}^{(i_q)}.
\end{equation}

\vspace{4mm}

More precisely, permutations $\left(j_1,\ldots,j_q,j_1',\ldots,j_{k-1}'\right)$ 
when summing in {\rm (\ref{new777100})}
are performed in such a way that if
$j_r^{*}$ swapped with $j_d^{*}$ in the  
permutation 

\vspace{1mm}
$$
\left(j_{q+k-1}^{*},\ldots,j_1^{*}\right)=
\left(j_q,\ldots,j_1,j_{k-1}',j_{k-2}',\ldots,j_{1}'\right),
$$ 

\vspace{4mm}
\noindent
then $i_r^{*}$ swapped with $i_d^{*}$ in 
the permutation 

\vspace{1mm}
$$
\left(i_{q+k-1}^{*},\ldots,i_1^{*}\right)=
\bigl(i_q,\ldots,i_1,0,\underbrace{1, \ldots ,1}_{k-2}\bigr).
$$

\vspace{4mm}
\noindent
Moreover, 
$\bar \phi_{j_r^{*}}$ swapped with $\bar \phi_{j_d^{*}}$
in the permutation 

\vspace{1mm}
$$
\bigl(\bar \phi_{j_{q+k-1}^{*}},\ldots,\bar \phi_{j_1^{*}}\bigr)=
\bigl(\phi_{j_q},\ldots,\phi_{j_1},\hspace{1.5mm} \phi_{j_{k}'}\hspace{-0.5mm}\cdot\hspace{-0.5mm}
\phi_{j_{k-1}'},\hspace{1.5mm}
\phi_{j_{k-2}'},\ldots, \phi_{j_{1}'}\bigr).
$$

\vspace{4mm}
\noindent
A similar rule should be applied to all other sums with respect to permutations

\vspace{-1mm}
$$
\sum\limits_{(j_1,\ldots,j_q,j_1',\ldots,j_{k-1}')}
$$

\vspace{3mm}
\noindent
in {\rm (\ref{new1901})} that contain the expressions
$\phi_{j_{k}'}(\theta)\phi_{j_{k-2}'}(\theta),\ldots,
\phi_{j_{k}'}(\theta)\phi_{j_{1}'}(\theta).$}

\vspace{2mm}

Let us prove the equality (\ref{new1600a}). Consider the case  $n=1.$
By analogy with (\ref{new1700}) and (\ref{new1700a}) we obtain 

\vspace{1mm}
$$
\int\limits_t^{s}\phi_{j_1'}(\theta)d{\bf w}_{\theta}^{(0)}
\sum\limits_{(j_1,\ldots,j_q)}\int\limits_t^s \phi_{j_q}(t_q)\ldots
\int\limits_t^{t_2}\phi_{j_1}(t_1)
d{\bf w}_{t_1}^{(i_1)}\ldots {\bf w}_{t_q}^{(i_q)}=
$$

\vspace{2.5mm}
$$
=
\sum\limits_{(j_1,\ldots,j_q, j_1')}
J_{(j_q\ldots j_1 j_1')s,t}^{(i_q\ldots i_1 0)}
$$

\vspace{5mm}
\noindent
w.~p.~1, where $J_{(j_r\ldots j_1)s,t}^{(i_r\ldots i_1)}$
is defined by (\ref{new100001}).
The equality (\ref{new1600a}) is proved for the case $n=1.$

Let us assume that the equality (\ref{new1600a}) is true for $n=2, 3, \ldots, k-1$, and prove
its validity for $n=k.$

In complete analogy with (\ref{new1301}) we get

\vspace{1mm}
$$
\int\limits_t^{s}\phi_{j_k'}(\theta)d\theta
\int\limits_t^s \phi_{j_{k-1}'}(t_{k-1})\ldots
\int\limits_t^{t_2}\phi_{j_1'}(t_1)
dt_1\ldots dt_{k-1}=
$$

\vspace{2.5mm}
\begin{equation}
\label{new5000}
=
J_{(j_k'j_{k-1}'\ldots j_1')s,t}^{(0\ldots  0)}
+J_{(j_{k-1}'j_k' j_{k-2}'\ldots j_1')s,t}^{(0\ldots 0)}+
\ldots + J_{(j_{k-1}'\ldots j_1' j_k')s,t}^{(0\ldots 0)}.
\end{equation}

\vspace{6mm}

Applying (\ref{new5000}), we have

\vspace{1mm}
$$
\sum\limits_{(j_1',\ldots,j_k')}
\int\limits_t^T \phi_{j_k'}(t_k')\ldots \int\limits_t^{t_2'}\phi_{j_1'}(t_1')
d{\bf w}_{t_1'}^{(0)}\ldots d{\bf w}_{t_k'}^{(0)}=
$$

\vspace{4.5mm}
$$
=\sum\limits_{(j_1',\ldots,j_{k-1}')}
\left(J_{(j_k'j_{k-1}'\ldots j_1')s,t}^{(0\ldots  0)}
+J_{(j_{k-1}'j_k' j_{k-2}'\ldots j_1')s,t}^{(0\ldots 0)}+
\ldots + J_{(j_{k-1}'\ldots j_1' j_k')s,t}^{(0\ldots 0)}\right)=
$$

\vspace{2.5mm}
\begin{equation}
\label{new3005}
=\int\limits_t^{T}\phi_{j_k'}(\theta)d\theta\sum\limits_{(j_1',\ldots,j_{k-1}')}
\int\limits_t^T \phi_{j_{k-1}'}(t_{k-1})\ldots \int\limits_t^{t_2'}\phi_{j_1'}(t_1)
d{\bf w}_{t_1}^{(0)}\ldots d{\bf w}_{t_{k-1}}^{(0)}.
\end{equation}

\vspace{6mm}

Using (\ref{new3005}) and the induction hypothesis, we obtain w.~p.~1

\vspace{1mm}
$$
\sum\limits_{(j_1',\ldots,j_k')}
\int\limits_t^T 
\phi_{j_k'}(t_k)\ldots \int\limits_t^{t_2}\phi_{j_1'}(t_1)d{\bf w}_{t_1}^{(0)}\ldots 
d{\bf w}_{t_k}^{(0)}\times
$$

\vspace{2.5mm}
$$
\times\sum\limits_{(j_1,\ldots,j_q)}\int\limits_t^T \phi_{j_q}(t_q)\ldots
\int\limits_t^{t_2}\phi_{j_1}(t_1)
d{\bf w}_{t_1}^{(i_1)}\ldots d{\bf w}_{t_q}^{(i_q)}=
$$

\vspace{2.5mm}
$$
=\int\limits_t^{T}\phi_{j_k'}(\theta)d\theta\sum\limits_{(j_1',\ldots,j_{k-1}')}
\int\limits_t^T \phi_{j_{k-1}'}(t_{k-1}')\ldots \int\limits_t^{t_2'}\phi_{j_1'}(t_1')
d{\bf w}_{t_1'}^{(0)}\ldots d{\bf w}_{t_{k-1}'}^{(0)}\times
$$

\vspace{2.5mm}
$$
\times\sum\limits_{(j_1,\ldots,j_q)}\int\limits_t^T \phi_{j_q}(t_q)\ldots
\int\limits_t^{t_2}\phi_{j_1}(t_1)
d{\bf w}_{t_1}^{(i_1)}\ldots d{\bf w}_{t_q}^{(i_q)}=
$$

\vspace{2.5mm}
$$
=\int\limits_t^{T}\phi_{j_k'}(\theta)d\theta
\sum\limits_{(j_1,\ldots,j_q,j_1',\ldots,j_{k-1}')}
\int\limits_t^T \phi_{j_q}(t_q)\ldots \int\limits_t^{t_2}\phi_{j_1}(t_1)\times
$$

\vspace{2.5mm}
$$
\times
\int\limits_t^{t_1}\phi_{j_{k-1}'}(t_{k-1}')\ldots \int\limits_t^{t_2'}
\phi_{j_1'}(t_1')d{\bf w}_{t_1'}^{(0)}\ldots d{\bf w}_{t_{k-1}'}^{(0)}
d{\bf w}_{t_1}^{(i_1)}\ldots d{\bf w}_{t_q}^{(i_q)}=
$$

\vspace{2.5mm}
$$
=
\sum\limits_{(j_1,\ldots,j_q,j_1',\ldots,j_{k-1}')}\int\limits_t^{T}\phi_{j_k'}(\theta)d\theta
\int\limits_t^T \phi_{j_q}(t_q)\ldots \int\limits_t^{t_2}\phi_{j_1}(t_1)\times
$$

\vspace{2.5mm}
\begin{equation}
\label{new4000}
\times
\int\limits_t^{t_1}\phi_{j_{k-1}'}(t_{k-1}')\ldots \int\limits_t^{t_2'}
\phi_{j_1'}(t_1')d{\bf w}_{t_1'}^{(0)}\ldots d{\bf w}_{t_{k-1}'}^{(0)}
d{\bf w}_{t_1}^{(i_1)}\ldots d{\bf w}_{t_q}^{(i_q)}.
\end{equation}

\vspace{6mm}

An iterative application of the Ito formula leads to the following equality

\vspace{1mm}
$$
\int\limits_t^{T}\phi_{j_k'}(\theta)d\theta
\int\limits_t^T \phi_{j_q}(t_q)\ldots \int\limits_t^{t_2}\phi_{j_1}(t_1)\times
$$

\vspace{2.5mm}
$$
\times
\int\limits_t^{t_1}\phi_{j_{k-1}'}(t_{k-1}')\ldots \int\limits_t^{t_2'}
\phi_{j_1'}(t_1')d{\bf w}_{t_1'}^{(0)}\ldots d{\bf w}_{t_{k-1}'}^{(0)}
d{\bf w}_{t_1}^{(i_1)}\ldots d{\bf w}_{t_q}^{(i_q)}=
$$

\vspace{4.5mm}
$$
=J_{(j_k'j_q \ldots j_1 j_{k-1}'\ldots j_1')T,t}^{(0 i_q\ldots i_1 0\ldots 0)}+
J_{(j_q j_k'j_{q-1} \ldots j_1 j_{k-1}'\ldots j_1')T,t}^{(i_q 0 i_{q-1}\ldots i_1 0\ldots 0)}+\ldots
J_{(j_q \ldots j_1 j_k' j_{k-1}'\ldots j_1')T,t}^{(i_q\ldots i_1 0\ldots 0)}+
$$

\vspace{2.5mm}
\begin{equation}
\label{new4001}
+J_{(j_q \ldots j_1 j_{k-1}' j_k' j_{k-2}'\ldots j_1')T,t}^{(i_q\ldots i_1 0\ldots 0)}+\ldots
+J_{(j_q \ldots j_1 j_{k-1}' \ldots j_1' j_k')T,t}^{(i_q\ldots i_1 0\ldots 0)}
\end{equation}

\vspace{6mm}
\noindent
w.~p.~1.

Combining (\ref{new4000}) and (\ref{new4001}) we finally obtain w.~p.~1

\vspace{1mm}
$$
\sum\limits_{(j_1,\ldots,j_q)}\int\limits_t^T \phi_{j_q}(t_q)\ldots
\int\limits_t^{t_2}\phi_{j_1}(t_1)
d{\bf w}_{t_1}^{(i_1)}\ldots d{\bf w}_{t_q}^{(i_q)}\times
$$

\vspace{2.5mm}
$$
\times \sum\limits_{(j_1',\ldots,j_k')}
\int\limits_t^T \phi_{j_k'}(t_k')\ldots \int\limits_t^{t_2'}\phi_{j_1'}(t_1')
d{\bf w}_{t_1'}^{(0)}\ldots d{\bf w}_{t_k'}^{(0)}=
$$

\vspace{2.5mm}
$$
=\sum\limits_{(j_1,\ldots,j_q,j_1',\ldots,j_k')}
\int\limits_t^T \phi_{j_q}(t_q)\ldots \int\limits_t^{t_2}\phi_{j_1}(t_1)
\int\limits_t^{t_1}\phi_{j_k'}(t_k')\ldots \int\limits_t^{t_2'}
\phi_{j_1'}(t_1')\times
$$

\vspace{4.5mm}
$$
\times d{\bf w}_{t_1'}^{(0)}\ldots d{\bf w}_{t_k'}^{(0)}d{\bf w}_{t_1}^{(i_1)}\ldots d{\bf w}_{t_q}^{(i_q)}.
$$

\vspace{6mm}

The equality (\ref{new1600a}) is proved for $n=k.$
The equality (\ref{new1600a}) is proved. Theorem~20 is proved.

To complete the proof of Theorems~12 and 13, we prove the following theorem.

\vspace{2mm}   

{\bf Theorem 21}\ \cite{20a}.\ {\it Suppose that
$\{\phi_j(x)\}_{j=0}^{\infty}$ is an arbitrary complete orthonormal system  
of func\-ti\-ons in the space $L_2([t,T]).$
Then the following representation
        
\vspace{2mm}
$$
J''[\phi_{j_1}\ldots\phi_{j_k}]^{(i_1 \ldots i_k)}=
\prod_{l=1}^k\zeta_{j_l}^{(i_l)}+\sum\limits_{r=1}^{[k/2]}
(-1)^r \times \Biggr.
$$

\vspace{2mm}
\begin{equation}
\label{leto6000xxa}
\times
\sum_{\stackrel{(\{\{g_1, g_2\}, \ldots, 
\{g_{2r-1}, g_{2r}\}\}, \{q_1, \ldots, q_{k-2r}\})}
{{}_{\{g_1, g_2, \ldots, 
g_{2r-1}, g_{2r}, q_1, \ldots, q_{k-2r}\}=\{1,2, \ldots, k\}}}}
\prod\limits_{s=1}^r
{\bf 1}_{\{i_{g_{{}_{2s-1}}}=~i_{g_{{}_{2s}}}\ne 0\}}
\Biggl.{\bf 1}_{\{j_{g_{{}_{2s-1}}}=~j_{g_{{}_{2s}}}\}}
\prod_{l=1}^{k-2r}\zeta_{j_{q_l}}^{(i_{q_l})}
\end{equation}

\vspace{4mm}
\noindent 
is valid w.~p.~{\rm 1,} where $i_1,\ldots,i_k=0,1,\ldots,m,$
$[x]$ is an integer part of a real number $x,$
the sum in the second line of the formula {\rm (\ref{leto6000xxa})} 
is the sum with respect to all possible
partitions {\rm (\ref{leto5008}),} 
$\prod\limits_{\emptyset}\stackrel{\sf def}{=}1,$ $\sum\limits_{\emptyset}
\stackrel{\sf def}{=}0;$
another notations are the same as in Theorems~{\rm 1, 2.}}

\vspace{2mm}

{\bf Remark~14.}\ {\it It should be noted that the formulas {\rm (\ref{new1010}),}
{\rm (\ref{new100000}),} {\rm (\ref{new1601}),} {\rm (\ref{new1601a})}
follow from {\rm (\ref{leto6000xxa}).}
It is only necessary to set the values
of the corresponding indicators of the form ${\bf 1}_A$ from the formula
{\rm (\ref{leto6000xxa})} equal to $0$ or $1.$}

\vspace{2mm}

{\bf Proof.}\ The proof of Theorem~21 is carried out
by induction using the following recurrence 
relation

\vspace{2mm}
$$
J''[\phi_{j_1}\ldots\phi_{j_k}]^{(i_1 \ldots i_k)}_{T,t}=
J''[\phi_{j_k}]^{(i_k)}_{T,t}\cdot
J''[\phi_{j_1}\ldots\phi_{j_{k-1}}]^{(i_1 \ldots i_{k-1})}_{T,t}-
$$

\vspace{1mm}
\begin{equation}
\label{recur1}
-\sum\limits_{l=1}^{k-1}{\bf 1}_{\{i_l=i_k\ne 0\}}
{\bf 1}_{\{j_l=j_k\}}\cdot 
J''[\phi_{j_1}\ldots\phi_{j_{l-1}}\phi_{j_{l+1}}\ldots\phi_{j_{k-1}}]^{(i_1
\ldots  i_{l-1}i_{l+1}\ldots i_{k-1})}_{T,t}
\end{equation}

\vspace{3mm}
\noindent
w.~p.~1.

Let us prove the recurrence relation (\ref{recur1}).
Using iteratively the Ito formula, the orthonormality of $\{\phi_j(x)\}_{j=0}^{\infty}$,
as well as (\ref{new5001}) and 
combinatorial reasoning, we obtain w.~p.~1 (see Remark~15 below for details)

\vspace{1mm}
$$
J''[\phi_{j_k}]^{(i_k)}_{T,t}\cdot
J''[\phi_{j_1}\ldots\phi_{j_{k-1}}]^{(i_1 \ldots i_{k-1})}_{T,t}=
$$

\vspace{2.5mm}
$$
=\int\limits_t^{T}\phi_{j_k}(\theta)d{\bf w}_{\theta}^{(i_k)}
\sum\limits_{(j_1,\ldots,j_{k-1})}
\int\limits_t^T \phi_{j_{k-1}}(t_{k-1})\ldots \int\limits_t^{t_2}\phi_{j_1}(t_1)
d{\bf w}_{t_1}^{(i_1)}\ldots d{\bf w}_{t_{k-1}}^{(i_{k-1})}=
$$

\vspace{2.5mm}
$$
=\sum\limits_{(j_1,\ldots,j_{k-1})}\int\limits_t^{T}\phi_{j_k}(\theta)d{\bf w}_{\theta}^{(i_k)}
\int\limits_t^T \phi_{j_{k-1}}(t_{k-1})\ldots \int\limits_t^{t_2}\phi_{j_1}(t_1)
d{\bf w}_{t_1}^{(i_1)}\ldots d{\bf w}_{t_{k-1}}^{(i_{k-1})}=
$$

\vspace{2.5mm}
$$
=\sum\limits_{(j_1,\ldots,j_{k})}
\int\limits_t^T \phi_{j_{k}}(t_{k})\ldots \int\limits_t^{t_2}\phi_{j_1}(t_1)
d{\bf w}_{t_1}^{(i_1)}\ldots d{\bf w}_{t_{k}}^{(i_{k})}+
$$

\vspace{2.5mm}
$$
+\sum\limits_{(j_1,\ldots,j_{k-1})}
\Biggl({\bf 1}_{\{i_k=i_{k-1}\ne 0\}}
\int\limits_t^T \phi_{j_{k}}(\theta)\phi_{j_{k-1}}(\theta)
\int\limits_t^{\theta}
\phi_{j_{k-2}}(t_{k-2})
\ldots \int\limits_t^{t_2}\phi_{j_1}(t_1)\times\Biggr.
$$

\vspace{2.5mm}
$$
\times
d{\bf w}_{t_1}^{(i_1)}\ldots d{\bf w}_{t_{k-2}}^{(i_{k-2})}d{\bf w}_{\theta}^{(0)}+
$$

\vspace{2.5mm}
$$
+{\bf 1}_{\{i_k=i_{k-2}\ne 0\}}
\int\limits_t^T \phi_{j_{k-1}}(t_{k-1})
\int\limits_t^{t_{k-1}} \phi_{j_{k}}(\theta)\phi_{j_{k-2}}(\theta)
\int\limits_t^{\theta}
\phi_{j_{k-3}}(t_{k-3})
\ldots \int\limits_t^{t_2}\phi_{j_1}(t_1)\times\Biggr.
$$

\vspace{5mm}
$$
\times
d{\bf w}_{t_1}^{(i_1)}\ldots d{\bf w}_{t_{k-3}}^{(i_{k-3})}d{\bf w}_{\theta}^{(0)}
d{\bf w}_{t_{k-1}}^{(i_{k-1})}+ \ldots
$$

\vspace{2.5mm}
$$
\ldots +{\bf 1}_{\{i_k=i_{1}\ne 0\}}
\int\limits_t^T \phi_{j_{k-1}}(t_{k-1})\ldots 
\int\limits_t^{t_3}
\phi_{j_{2}}(t_{2})
\int\limits_t^{t_{2}} \phi_{j_{k}}(\theta)\phi_{j_{1}}(\theta)\times
$$

\vspace{2.5mm}
$$
\Biggl.\times
d{\bf w}_{\theta}^{(0)}d{\bf w}_{t_2}^{(i_2)}\ldots d{\bf w}_{t_{k-1}}^{(i_{k-1})}
\Biggr)=
$$

\vspace{2.5mm}
$$
=\sum\limits_{(j_1,\ldots,j_{k})}
\int\limits_t^T \phi_{j_{k}}(t_{k})\ldots \int\limits_t^{t_2}\phi_{j_1}(t_1)
d{\bf w}_{t_1}^{(i_1)}\ldots d{\bf w}_{t_{k}}^{(i_{k})}+
$$

\vspace{2.5mm}
$$
+\sum\limits_{(j_1,\ldots,j_{k-2})}
{\bf 1}_{\{i_k=i_{k-1}\ne 0\}}
\Biggl\{\int\limits_t^T \phi_{j_{k}}(\theta)\phi_{j_{k-1}}(\theta)
\int\limits_t^{\theta}
\phi_{j_{k-2}}(t_{k-2})
\ldots \int\limits_t^{t_2}\phi_{j_1}(t_1)\times\Biggr.
$$

\vspace{4.5mm}
$$
\times
d{\bf w}_{t_1}^{(i_1)}\ldots d{\bf w}_{t_{k-2}}^{(i_{k-2})}d{\bf w}_{\theta}^{(0)}+\ldots
$$

\vspace{2.5mm}
$$
\Biggl.\ldots +\int\limits_t^T 
\phi_{j_{k-2}}(t_{k-2})
\ldots \int\limits_t^{t_2}\phi_{j_1}(t_1)
\int\limits_t^{t_1}
\phi_{j_{k}}(\theta)\phi_{j_{k-1}}(\theta)
d{\bf w}_{\theta}^{(0)}
d{\bf w}_{t_1}^{(i_1)}\ldots d{\bf w}_{t_{k-2}}^{(i_{k-2})}\Biggr\}+
$$

\vspace{2.5mm}
$$
+\sum\limits_{(j_1,\ldots,j_{k-3},j_{k-1})}
{\bf 1}_{\{i_k=i_{k-2}\ne 0\}}
\Biggl\{\int\limits_t^T \phi_{j_{k}}(\theta)\phi_{j_{k-2}}(\theta)
\int\limits_t^{\theta}
\phi_{j_{k-1}}(t_{k-1})
\int\limits_t^{t_{k-1}}
\phi_{j_{k-3}}(t_{k-3})
\ldots \Biggr.
$$

\vspace{2.5mm}
$$
\ldots \int\limits_t^{t_2}\phi_{j_1}(t_1)
d{\bf w}_{t_1}^{(i_1)}\ldots d{\bf w}_{t_{k-3}}^{(i_{k-3})}d{\bf w}_{t_{k-1}}^{(i_{k-1})}
d{\bf w}_{\theta}^{(0)}+\ldots
$$

\vspace{2.5mm}
$$
\ldots +\int\limits_t^T 
\phi_{j_{k-1}}(t_{k-1})
\int\limits_t^{t_{k-1}}
\phi_{j_{k-3}}(t_{k-3})
\ldots \int\limits_t^{t_2}\phi_{j_1}(t_1)
\int\limits_t^{t_1}
\phi_{j_{k}}(\theta)\phi_{j_{k-2}}(\theta)\times
$$

\vspace{2.5mm}
$$
\Biggl.\times d{\bf w}_{\theta}^{(0)}
d{\bf w}_{t_1}^{(i_1)}\ldots d{\bf w}_{t_{k-3}}^{(i_{k-3})}d{\bf w}_{t_{k-1}}^{(i_{k-1})}\Biggr\}+\ldots
$$

\vspace{2.5mm}
$$
\ldots +\sum\limits_{(j_2,\ldots,j_{k-1})}
{\bf 1}_{\{i_k=i_{1}\ne 0\}}
\Biggl\{\int\limits_t^T \phi_{j_{k}}(\theta)\phi_{j_{1}}(\theta)
\int\limits_t^{\theta}
\phi_{j_{k-1}}(t_{k-1})
\ldots \int\limits_t^{t_3}\phi_{j_2}(t_2)\times\Biggr.
$$

\vspace{4.5mm}
$$
\times
d{\bf w}_{t_2}^{(i_2)}\ldots d{\bf w}_{t_{k-1}}^{(i_{k-1})}d{\bf w}_{\theta}^{(0)}+\ldots
$$

\vspace{2.5mm}
$$
\Biggl.\ldots +\int\limits_t^T 
\phi_{j_{k-1}}(t_{k-1})
\ldots \int\limits_t^{t_3}\phi_{j_2}(t_2)
\int\limits_t^{t_2}
\phi_{j_{k}}(\theta)\phi_{j_{1}}(\theta)
d{\bf w}_{\theta}^{(0)}
d{\bf w}_{t_2}^{(i_2)}\ldots d{\bf w}_{t_{k-1}}^{(i_{k-1})}\Biggr\}=
$$

\vspace{2.5mm}
$$
=\sum\limits_{(j_1,\ldots,j_{k})}
\int\limits_t^T \phi_{j_{k}}(t_{k})\ldots \int\limits_t^{t_2}\phi_{j_1}(t_1)
d{\bf w}_{t_1}^{(i_1)}\ldots d{\bf w}_{t_{k}}^{(i_{k})}+
$$

\vspace{2.5mm}
$$
+\int\limits_t^T \phi_{j_{k}}(\theta)\phi_{j_{k-1}}(\theta)d\theta
\sum\limits_{(j_1,\ldots,j_{k-2})}
{\bf 1}_{\{i_k=i_{k-1}\ne 0\}}
\int\limits_t^{T}
\phi_{j_{k-2}}(t_{k-2})
\ldots \int\limits_t^{t_2}\phi_{j_1}(t_1)\times\Biggr.
$$

\vspace{4.5mm}
$$
\times
d{\bf w}_{t_1}^{(i_1)}\ldots d{\bf w}_{t_{k-2}}^{(i_{k-2})}+
$$

\vspace{2.5mm}
$$
+\int\limits_t^T \phi_{j_{k}}(\theta)\phi_{j_{k-2}}(\theta)d\theta
\sum\limits_{(j_1,\ldots,j_{k-3},j_{k-1})}
{\bf 1}_{\{i_k=i_{k-2}\ne 0\}}
\int\limits_t^{T}
\phi_{j_{k-1}}(t_{k-1})
\int\limits_t^{t_{k-1}}
\phi_{j_{k-3}}(t_{k-3})
\ldots 
$$

\vspace{3.5mm}
$$
\ldots \int\limits_t^{t_2}\phi_{j_1}(t_1)
d{\bf w}_{t_1}^{(i_1)}\ldots d{\bf w}_{t_{k-3}}^{(i_{k-3})}d{\bf w}_{t_{k-1}}^{(i_{k-1})}
+\ldots
$$

\vspace{2.5mm}

$$
\ldots +\int\limits_t^T \phi_{j_{k}}(\theta)\phi_{j_{1}}(\theta)d\theta
\sum\limits_{(j_2,\ldots,j_{k-1})}
{\bf 1}_{\{i_k=i_{1}\ne 0\}}
\int\limits_t^{T}
\phi_{j_{k-1}}(t_{k-1})
\ldots \int\limits_t^{t_3}\phi_{j_2}(t_2)\times
$$

\vspace{4.5mm}
$$
\times
d{\bf w}_{t_2}^{(i_2)}\ldots d{\bf w}_{t_{k-1}}^{(i_{k-1})}=
$$

\vspace{3.5mm}

$$
=J''[\phi_{j_1}\ldots\phi_{j_k}]^{(i_1 \ldots i_k)}_{T,t}+
{\bf 1}_{\{i_k=i_{k-1}\ne 0\}}{\bf 1}_{\{j_k=j_{k-1}\}}
\cdot J''[\phi_{j_1}\ldots \phi_{j_{k-2}}]_{T,t}^{(i_1\ldots i_{k-2})}+
$$

\vspace{3.5mm}
$$
+
{\bf 1}_{\{i_k=i_{k-2}\ne 0\}}{\bf 1}_{\{j_k=j_{k-2}\}}
\cdot J''[\phi_{j_1}\ldots \phi_{j_{k-3}}\phi_{j_{k-1}}]_{T,t}^{(i_1\ldots i_{k-3}i_{k-1})}+\ldots
$$

\vspace{3.5mm}
$$
\ldots +
{\bf 1}_{\{i_k=i_{1}\ne 0\}}{\bf 1}_{\{j_k=j_{1}\}}
\cdot J''[\phi_{j_2}\ldots \phi_{j_{k-1}}]_{T,t}^{(i_2\ldots i_{k-1})}=
$$

\vspace{3.5mm}
$$
=J''[\phi_{j_1}\ldots\phi_{j_k}]^{(i_1 \ldots i_k)}_{T,t}+
$$

\vspace{1mm}
\begin{equation}
\label{new00002}
+\sum\limits_{l=1}^{k-1}{\bf 1}_{\{i_l=i_k\ne 0\}}
{\bf 1}_{\{j_l=j_k\}}\cdot 
J''[\phi_{j_1}\ldots\phi_{j_{l-1}}\phi_{j_{l+1}}\ldots\phi_{j_{k-1}}]^{(i_1
\ldots  i_{l-1}i_{l+1}\ldots i_{k-1})}_{T,t}.
\end{equation}

\vspace{6mm}

The equality (\ref{recur1}) is proved. Theorem~21 is proved.

\vspace{2mm}

{\bf Remark~15.}\ {\it It should be noted that the sums with respect to 
permutations 

\vspace{1mm}
$$
\sum\limits_{(j_1,\ldots,j_{k-1})}
$$

\vspace{3mm}
\noindent
in {\rm (\ref{new00002}),} containing the expressions 

\vspace{1mm}
$$
{\bf 1}_{\{i_k=i_{k-1}\ne 0\}}\phi_{j_{k}}(\theta)\phi_{j_{k-1}}(\theta),\ldots,
{\bf 1}_{\{i_k=i_{1}\ne 0\}}\phi_{j_{k}}(\theta)\phi_{j_{1}}(\theta),
$$

\vspace{5mm}
\noindent
should be understood in a special way.
Let us explain this rule on the basis of the sum

\vspace{1mm}
$$
\sum\limits_{(j_1,\ldots,j_{k-1})}
{\bf 1}_{\{i_k=i_{k-1}\ne 0\}}
\int\limits_t^T \phi_{j_{k}}(\theta)\phi_{j_{k-1}}(\theta)
\int\limits_t^{\theta}
\phi_{j_{k-2}}(t_{k-2})
\ldots \int\limits_t^{t_2}\phi_{j_1}(t_1)\times
$$

\vspace{2.5mm}
\begin{equation}
\label{new00003}
\times
d{\bf w}_{t_1}^{(i_1)}\ldots d{\bf w}_{t_{k-2}}^{(i_{k-2})}d{\bf w}_{\theta}^{(0)}.
\end{equation}

\vspace{6mm}

More precisely, permutations $(j_1,\ldots,j_{k-1})$ 
when summing in {\rm (\ref{new00003})}
are performed in such a way that if
$j_r$ swapped with $j_d$ in the  
permutation $(j_1,\ldots,j_{k-1}),$ 
then $i_r$ swapped with $i_d$ in 
the permutation $(i_1,\ldots,i_{k-2} i_{k-1})$\ {\rm (}note that $i_{k-1}=0${\rm )}.
Moreover, 
$\bar \phi_{j_r}$ swapped with $\bar \phi_{j_d}$
in the permutation 

\vspace{1mm}
$$
\bigl(\bar \phi_{j_{1}},\ldots,\bar\phi_{j_{k-1}}\bigr)=
\bigl(\phi_{j_1},\ldots,\phi_{j_{k-2}},\hspace{1.5mm} 
{\bf 1}_{\{i_k=i_{k-1}\ne 0\}}\cdot \phi_{j_k}\cdot \phi_{j_{k-1}}\bigr),
$$

\vspace{5mm}
\noindent
where $\bar\phi_{j_{k-1}}(\tau)=
{\bf 1}_{\{i_k=i_{k-1}\ne 0\}}\phi_{j_k}(\tau)\phi_{j_{k-1}}(\tau).$

A similar rule should be applied to all other sums with respect to permutations

\vspace{1mm}
$$
\sum\limits_{(j_1,\ldots,j_{k-1})}
$$

\vspace{4mm}
\noindent
in {\rm (\ref{new00002})} that contain the expressions

\vspace{1mm}
$$
{\bf 1}_{\{i_k=i_{k-2}\ne 0\}}\phi_{j_{k}}(\theta)\phi_{j_{k-2}}(\theta),\ldots,
{\bf 1}_{\{i_k=i_{1}\ne 0\}}\phi_{j_{k}}(\theta)\phi_{j_{1}}(\theta).
$$
}

\vspace{7mm}

The relations (\ref{chain4002}), (\ref{new6000}), (\ref{leto6000xxa})
prove Theorem~12. An analogue of the formula (\ref{chain4002})
for the function $\Phi(t_1,\ldots,t_k)$ instead of $K(t_1,\ldots,t_k)$ and
(\ref{new6000}), (\ref{leto6000xxa}) prove Theorem~13.

We also note a number of works \cite{ito1951}, \cite{Kuo}-\cite{major1} in which the properties
of multiple Wiener stochastic integrals were studied using
measure theory, in particular, the formulas for the product
of such integrals were obtained.

First of all, let us compare Theorem~21 with Proposition~5.1 from \cite{fox}.
An analogue of the right-hand side of (\ref{leto6000xxa})
for nonrandom $x_1,\ldots,x_k$
is constructed in \cite{fox} using diagrams (see the formula (5.1) in \cite{fox}).
This means that the application of the formula (5.1) from \cite{fox},
unlike the formula (\ref{leto6000xxa}), is difficult when
performing algebraic transformations.

Further, we note that the formula (5.1) from \cite{fox}
was applied to the representation of the multiple Wiener stochastic integral
somewhat differently than the formula (\ref{leto6000xxa}).
Namely, using Proposition~5.1 \cite{fox}.
Let us expain this difference in more detail.

Proposition~5.1 from \cite{fox} in our degree of generality 
and in our notations can be written as

\vspace{2mm}
$$
J''\left[\phi_{j_1}\ldots \phi_{j_k}\right]_{T,t}^{(i_1\ldots i_k)}=
$$

\vspace{3mm}
$$
=
J''\biggl[\underbrace{\phi_{j_1}
\ldots \phi_{j_1}}_{m_1}
\underbrace{\phi_{j_2}
\ldots \phi_{j_2}}_{m_2}\ldots 
\underbrace{\phi_{j_p}\ldots
\phi_{j_p}}_{m_p}\biggr]_{T,t}^
{(\overbrace{\hspace{0.5mm}{}_{i_1 \ldots i_{m_1}}}^{m_1}
\overbrace{\hspace{0.3mm}{}_{i_{m_1+1} \ldots i_{m_2}}}^{m_2}
\ldots \overbrace{\hspace{0.3mm}{}_{i_{m_1+\ldots +m_{p-1}+1} \ldots i_k}}^{m_p})}=
$$

\vspace{3mm}
\begin{equation}
\label{new54321}
=J''\left[\phi_{j_1}
\hspace{-0.3mm}\ldots \hspace{-0.3mm}\phi_{j_1}\right]
_{T,t}^
{(\overbrace{\hspace{0.5mm}{}_{i_1 \ldots i_{m_1}}}^{m_1})}\hspace{-0.3mm}\cdot
J''\left[\phi_{j_2}\hspace{-0.3mm}
\ldots \hspace{-0.3mm}\phi_{j_2}\right]_{T,t}^
{(\overbrace{\hspace{0.3mm}{}_{i_{m_1+1} \ldots i_{m_2}}}^{m_2})}\hspace{-0.3mm}\cdot \ldots
\cdot J''\left[\phi_{j_p}
\hspace{-0.3mm}\ldots \hspace{-0.3mm}\phi_{j_p}\right]_{T,t}^
{(\overbrace{\hspace{0.3mm}{}_{i_{m_1+\ldots +m_{p-1}+1} \ldots i_k}}^{m_p})}
\end{equation}

\vspace{6mm}
\noindent
w.~p.~1, where

\vspace{-1mm}
$$
J''\left[\phi_{j_1}
\hspace{-0.3mm}\ldots \hspace{-0.3mm}\phi_{j_1}\right]
_{T,t}^
{(\overbrace{\hspace{0.5mm}{}_{i_1 \ldots i_{m_1}}}^{m_1})}\hspace{-1.5mm},
J''\left[\phi_{j_2}\hspace{-0.3mm}
\ldots \hspace{-0.3mm}\phi_{j_2}\right]_{T,t}^
{(\overbrace{\hspace{0.3mm}{}_{i_{m_1+1} \ldots i_{m_2}}}^{m_2})}\hspace{-1mm},\ldots,
J''\left[\phi_{j_p}
\hspace{-0.3mm}\ldots \hspace{-0.3mm}\phi_{j_p}\right]_{T,t}^
{(\overbrace{\hspace{0.3mm}{}_{i_{m_1+\ldots +m_{p-1}+1} \ldots i_k}}^{m_p})}
$$

\vspace{6mm}
\noindent
are defined by the right-hand side of the formula (5.1) from \cite{fox},
$m_1+\ldots +m_p=k,$ $m_1,\ldots, m_p>0,$ $j_q\ne j_d$ $(q\ne d,\ q,d=1,\ldots,p),$
$i_1,\ldots,i_k=1,\ldots,m.$

This actually means that in \cite{fox} an analogue of the formula
(\ref{leto6000xxa}) is constructed for the special case
$j_1=\ldots=j_k$. Moreover, the specified analogue 
is based on the formula (5.1) \cite{fox} obtained using diagrams.

Comparing the formulas (\ref{leto6000xxa}) and (\ref{new54321}) (or (5.1) from \cite{fox}), it is easy
to understand that the transition from 
(\ref{leto6000xxa}) and (\ref{new54321}) is obvious.
It is only necessary to set the values
of the corresponding indicators of the form ${\bf 1}_A$ from the formula
(\ref{leto6000xxa}) equal to $0$ or $1.$
The reverse transition from the formula (\ref{new54321})
to the formula (\ref{leto6000xxa}) is not obvious.
Note that the formula 
(\ref{leto6000xxa}) (not the formula (\ref{new54321})) is convenient for the  numerical
integration of Ito stochastic differential equations (see \cite{20a}, 
Chapter~5 and \cite{Kuz-Kuz}, \cite{Mikh-1} for details).

Let us turn to the comparison of Theorem~21 with another interesting work \cite{major2} (2019).
As it turned out, a version of Theorem~21 was obtained in terms 
of Wick polynomials and for the case of vector valued random measures 
in \cite{major2} (see Theorem~7.2, p.~69).
However, much earlier the formula (\ref{leto6000xxa}) (Theorem~21) is obtained
in our monograph \cite{10} (2009) as part of the formula
(5.30) (see \cite{10}, p.~220).
Moreover, particular cases of the formula (\ref{leto6000xxa}) were obtained
even earlier in our works \cite{7} (2006) and \cite{9} (2007).
More precisely, partiular cases $k=1,\ldots,5$ of the formula (\ref{leto6000xxa})
were obtained in \cite{7} (2006) as parts of the formulas 
on the pages 243-244 and partiular cases $k=1,\ldots,7$ of the formula (\ref{leto6000xxa})
were obtained in \cite{9} (2007) as parts of the formulas 
on the pages 208-218.

We also note that we have found an explicit expression for the 
Wick polynomial of degree $k$ of the arguments $\zeta_{j_1}^{(i_1)},\ldots,\zeta_{j_k}^{(i_k)}$ 
(see the formula (\ref{leto6000xxa})),
which is very convenient for the numerical simulation of
iterated Ito stochastic integrals (\ref{sodom20})  \cite{Kuz-Kuz}, \cite{Mikh-1}.
Note that the representation of the Wick polynomial
of the arguments $\zeta_{j_1}^{(i_1)},\ldots,\zeta_{j_k}^{(i_k)}$ 
in terms of the product of Hermite polynomials
is less convenient for the numerical simulation of
iterated Ito stochastic integrals (\ref{sodom20}).
For example, the expression for $J''[\phi_{j_1}\phi_{j_2}\phi_{j_3}\phi_{j_4}]^{(i_1 i_2 i_3 i_4)}_{T,t}$
in terms of the product of Hermite polynomials,
even under the condition $i_1=i_2=i_3=i_4$, already contains
15 different expressions (see Sect.~14).
At the same time, all these 15 expressions are contained 
in one formula (\ref{leto6000xxa}) provided that $k=4$ and $i_1=i_2=i_3=i_4$.
It is very convenient, since in computer simulation
using the formula (\ref{leto6000xxa}), in addition to
modeling of random variables $\zeta_{j_1}^{(i_1)},\ldots,\zeta_{j_k}^{(i_k)}$,
it remains only to set  
the values
of the corresponding indicators of the form ${\bf 1}_A$ from the formula
(\ref{leto6000xxa}) equal to $0$ or $1.$

It should be noted that in \cite{major} (Theorem~6.1)
a diagram formula was obtained for the product
of two multiple Wiener stochastic integrals
with respect to vector valued random measures. 
The formula (\ref{new1600}) can be derived from the diagram formula \cite{major}.
Although the proof of the diagram formula \cite{major}
is much more complicated than our proof of the formula (\ref{new1600}).

To conclude this section, we say a few words about expansions 
(\ref{new9999}) and (\ref{razzar1}).
The transition from the expansion (\ref{razzar1}) to the expansion 
(\ref{new9999}) is obvious. It is only necessary to set the values
of the corresponding indicators of the form ${\bf 1}_A$ from the formula
(\ref{razzar1}) equal to $0$ or $1.$
The reverse transition from the formula (\ref{new9999})
to the formula (\ref{razzar1}) is also possible but not obvious.
However, Theorems~20 and 21 provide a transition from 
(\ref{new9999}) to (\ref{razzar1}) and vice versa.
Note that the expansion (\ref{new9999}) is interesting from the point of 
view of studying the structure of the expansion of iterated Ito
stochastic integrals. On the orther hand, 
the expansion (\ref{razzar1}) is exceptionally convenient 
for applications (see
\cite{Kuz-Kuz}, \cite{Mikh-1}).

\vspace{5mm}

\section{Generalization of Theorem~7 to the Case of an Arbitrary 
Complete Ortho\-nor\-mal System of Functions in the Space $L_2([t, T])$
and $\psi_1(\tau),$ $\ldots,\psi_k(\tau)\in L_2([t, T])$}

\vspace{5mm}

Suppose that $\psi_1(\tau),$ $\ldots,\psi_k(\tau)\in L_2([t, T])$.
Define the following function on the hypercube $[t, T]^k$

$$
\bar K(t_1,\ldots,t_k,s)={\bf 1}_{\{t_k<s\}}K(t_1,\ldots,t_k),
$$

\vspace{4mm}
\noindent
where the function $K(t_1,\ldots,t_k)$ has the form
(\ref{ppp}), $s\in (t, T]$ ($s$ is fixed), 
and ${\bf 1}_A$ is the indicator of the set $A.$

Further, we have (see (\ref{ppp}))

$$
\bar K(t_1,\ldots,t_k,s)=
{\bf 1}_{\{t_1<\ldots <t_k<s\}}\psi_1(t_1)\ldots \psi_k(t_k)=
$$

\vspace{2mm}
$$
=
\left\{\begin{matrix}
\psi_1(t_1)\ldots \psi_k(t_k),\ &t_1<\ldots<t_k<s\cr\cr\cr
0,\ &\hbox{\rm otherwise}
\end{matrix}
\right.,
$$

\vspace{5mm}
\noindent
where $\bar K(t_1,\ldots,t_k,s)\in L_2([t,T]^k),$ $k\ge 1, $ $t_1,\ldots,t_k\in [t, T],$ and 
$s\in (t, T]$.

Note that

\vspace{-1mm}
$$
J[\psi^{(k)}]_{s,t}=\int\limits_t^s\psi_k(t_k) \ldots \int\limits_t^{t_{2}}
\psi_1(t_1) d{\bf w}_{t_1}^{(i_1)}\ldots
d{\bf w}_{t_k}^{(i_k)}=
$$

\begin{equation}
\label{strange700}
=
\int\limits_t^T {\bf 1}_{\{t_k<s\}}\psi_k(t_k) \ldots \int\limits_t^{t_{2}}
\psi_1(t_1) d{\bf w}_{t_1}^{(i_1)}\ldots
d{\bf w}_{t_k}^{(i_k)}\ \ \ \hbox{w.~p.~1},
\end{equation}

\vspace{4mm}
\noindent
where $s\in (t, T]$ ($s$ is fixed), $i_1,\ldots,i_k=0,1,\ldots,m.$ 

Applying Theorem~12 to the iterated Ito stochastic integral
(\ref{strange700}), we obtain the following ge\-ne\-ra\-li\-za\-tion 
of Theorem~7 to the case of an arbitrary 
complete ortho\-nor\-mal system of functions in the space $L_2([t, T])$
and $\psi_1(\tau),$ $\ldots,\psi_k(\tau)\in L_2([t, T]).$

\vspace{2mm}

{\bf Theorem~22.}\ {\it Suppose that
$\psi_1(\tau),$ $\ldots,\psi_k(\tau)\in L_2([t, T])$ and
$\{\phi_j(x)\}_{j=0}^{\infty}$ is an arbitrary complete ortho\-nor\-mal system  
of functions in the space $L_2([t,T]).$
Then$,$ the following expansion

\vspace{1mm}
$$
J[\psi^{(k)}]_{s,t}^{(i_1\ldots i_k)}=
\hbox{\vtop{\offinterlineskip\halign{
\hfil#\hfil\cr
{\rm l.i.m.}\cr
$\stackrel{}{{}_{p_1,\ldots,p_k\to \infty}}$\cr
}} }
\sum\limits_{j_1=0}^{p_1}\ldots
\sum\limits_{j_k=0}^{p_k}
C_{j_k\ldots j_1}(s)\Biggl(
\prod_{l=1}^k\zeta_{j_l}^{(i_l)}+\sum\limits_{r=1}^{[k/2]}
(-1)^r \times
\Biggr.
$$

\vspace{2mm}
$$
\times
\sum_{\stackrel{(\{\{g_1, g_2\}, \ldots, 
\{g_{2r-1}, g_{2r}\}\}, \{q_1, \ldots, q_{k-2r}\})}
{{}_{\{g_1, g_2, \ldots, 
g_{2r-1}, g_{2r}, q_1, \ldots, q_{k-2r}\}=\{1, 2, \ldots, k\}}}}
\prod\limits_{s=1}^r
{\bf 1}_{\{i_{g_{{}_{2s-1}}}=~i_{g_{{}_{2s}}}\ne 0\}}
\Biggl.{\bf 1}_{\{j_{g_{{}_{2s-1}}}=~j_{g_{{}_{2s}}}\}}
\prod_{l=1}^{k-2r}\zeta_{j_{q_l}}^{(i_{q_l})}\Biggr)
$$

\vspace{5mm}
\noindent
con\-verg\-ing in the mean-square sense is valid$,$
where $[x]$ is an integer part of a real number $x,$

$$
C_{j_k\ldots j_1}(s)=\int\limits_{[t,T]^k}
\bar K(t_1,\ldots,t_k,s)\prod_{l=1}^{k}\phi_{j_l}(t_l)dt_1\ldots dt_k=
$$

\vspace{1mm}
$$
=\int\limits_t^s\psi_k(t_k)\phi_{j_k}(t_k)\ldots
\int\limits_t^{t_2}
\psi_1(t_1)\phi_{j_1}(t_1)
dt_1\ldots dt_k
$$

\vspace{4mm}
\noindent
is the Fourier coefficient$,$ $\prod\limits_{\emptyset}
\stackrel{\sf def}{=}1,$ $\sum\limits_{\emptyset}
\stackrel{\sf def}{=}0;$ another
notations are the same as in Theorem {\rm 2}.}

\vspace{2mm}

Note that the estimates (\ref{road888}) and (\ref{agent01000})
will also be valid under the conditions of Theorem~22.

\end{document}